\documentclass[11pt,reqno]{amsart}
\usepackage{tikz,a4wide}
\usetikzlibrary{arrows.meta,patterns.meta}
\usetikzlibrary{math,calc}
\usepackage{longtable}
\usepackage{tikz-cd}
\usepackage{tkz-euclide}
\usepackage{float}
\usepackage{amsmath}
\usepackage{xfrac}
\usepackage{stmaryrd}
\SetSymbolFont{stmry}{bold}{U}{stmry}{m}{n}
\usepackage{mathrsfs,bm,amsthm,mathtools,yfonts,amssymb,color}
\usepackage{esint}
\usepackage{xcolor}
\usepackage[pagebackref,citecolor=black,linkcolor=blue,colorlinks=true]{hyperref}
\usepackage{tikz-3dplot}
\usepackage{epstopdf}
\usepackage{thmtools}
\usepackage[capitalize,noabbrev,nosort]{cleveref}
\usepackage{csquotes}
\usepackage{comment}
\MakeOuterQuote{"} 
\usepackage{enumitem}
\usetikzlibrary{shapes.geometric, positioning, arrows.meta, calc}
\usepackage{caption}

\makeatletter 
\DeclarePairedDelimiter{\@tmpabs}{\lvert}{\rvert}
\newcommand{\@absstar}[1]{{\@tmpabs*{#1}}}
\newcommand{\@absnostar}[2][]{{\@tmpabs[#1]{#2}}}
\newcommand{\abs}{\@ifstar\@absstar\@absnostar}
\makeatother

\makeatletter 
\DeclarePairedDelimiter{\@tmpnorm}{\lVert}{\rVert}
\newcommand{\@normstar}[1]{{\@tmpnorm*{#1}}}
\newcommand{\@normnostar}[2][]{{\@tmpnorm[#1]{#2}}}
\newcommand{\norm}{\@ifstar\@normstar\@normnostar}
\makeatother

\begingroup
\newtheorem{theorem}{Theorem}[section]
\newtheorem{lemma}[theorem]{Lemma}
\newtheorem{proposition}[theorem]{Proposition}
\newtheorem{corollary}[theorem]{Corollary}
\newtheorem{lista}[theorem]{List}
\endgroup

\theoremstyle{definition}
\newtheorem{definition}[theorem]{Definition}
\newtheorem{remark}[theorem]{Remark}

\newtheorem{problem}[theorem]{Problem}

\crefname{lista}{List}{Lists}
\Crefname{lista}{List}{Lists}

\numberwithin{equation}{section}

\DeclareMathAlphabet{\mathdutchcal}{U}{dutchcal}{m}{n} 
\newcommand\scrs{\mathdutchcal{s}} 
\newcommand\scrh{\mathdutchcal{h}} 

\newcommand{\N}{\mathbb{N}} 
\newcommand{\R}{\mathbb{R}} 
\newcommand{\Sf}{\mathbb{S}^2} 
\newcommand{\Ha}{\mathcal{H}} 
\newcommand{\dist}{\mathrm{dist}} 
\newcommand{\darts}{\mathrm{Darts}} 
\newcommand{\rev}{\mathrm{rev}} 
\renewcommand{\next}{\mathrm{next}} 
\newcommand\TS{T^1\mathbb S^2} 
\DeclareMathOperator{\rot}{rot} 
\newcommand{\SO}{\operatorname{SO}}
\newcommand{\Lip}{\mathrm{Lip}}
\newcommand{\emptyparam}{\ensuremath{\,\cdot\,}}
\newcommand{\eps}{\varepsilon}

\newcommand\res{\mathop{\hbox{\vrule height 7pt width .3pt depth 0pt\vrule height .3pt width 5pt depth 0pt}}\nolimits}
\newcommand\len{\mathcal{L}}
\newcommand\inter{\mathcal{I}}

\usepackage{dsfont}
\newcommand{\id}{\ensuremath{\mathds{1}}}

\newcommand{\ssubset}{\subset\joinrel\subset}

\title{Boundary Plateau laws I}
\author[C. De Lellis]{Camillo De Lellis}
\address{School of Mathematics, Institute for Advanced Study, 1 Einstein Dr., Princeton NJ 08540, USA}
\address{Gran Sasso Science Institute, viale Francesco Crispi, 7, 67100 L’Aquila, Italy}
\email{camillo.delellis@ias.edu}
\author[F. Glaudo]{Federico Glaudo}
\email{federico.glaudo@gmail.com}

\begin{document}

\begin{abstract}
We give a complete list of the possible local shapes of soap films where they wet the whole wire. In rigorous terms, given a straight line $\ell \subset \mathbb R^3$, we study $2$-dimensional cones in $\mathbb R^3$ which contain $\ell$ and are Almgren minimizing with respect to $\ell$.
\end{abstract}

\maketitle

\tableofcontents

\section{Introduction}

\subsection{Informal explanation of the main theorem}\label{s:informal}
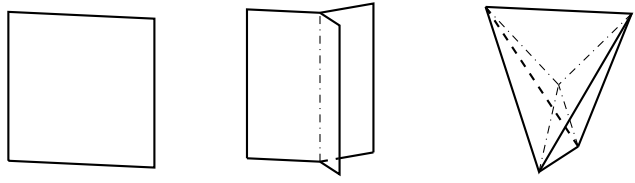
\begin{figure}
\tdplotsetmaincoords{80}{15}
\begin{tikzpicture}[tdplot_main_coords]
\draw[thick] (-7,0,-1) -- (-7,0,1) -- (-5,0,1) -- (-5,0,-1) -- (-7,0,-1);
\draw[thick] (-4,1,-1) -- (-3,1,-1);
\draw[thick] (-3,1,1) -- (-4,1,1) -- (-4,1,-1);
\draw[dashdotted] (-3,1,-1) -- (-3,1,1);
\draw[thick] (-3,1,1) -- (-2.5,{sqrt(3)/2+1},1) -- (-2.5,{sqrt(3)/2+1},-1);
\draw[thick] (-2.5,{sqrt(3)/2+1},-1) -- (-0.35*2.5 -0.65*3,{0.35*(sqrt(3)/2 + 1) + 0.65},-1);
\draw[thick, dashed] (-3,1,-1) -- (-0.35*2.5 -0.65*3,{0.35*(sqrt(3)/2 + 1) + 0.65},-1);
\draw[thick] (-3,1,1) -- (-2.5,{-sqrt(3)/2+1},1) -- (-2.5,{-sqrt(3)/2+1},-1) -- (-3,1,-1);
\draw[thick] (-1,2,1) -- (0,1,-1) -- (1,2,1) -- (-1,2,1);
\draw[thick] (0,1,-1) -- (0,3,-1);
\draw[thick, dashed] (0,3,-1) -- (-1,2,1);
\draw[thick] (0,3,-1) -- (1,2,1);
\draw[dashdotted] (-1,2,1) -- (0,2,0);
\draw[dashdotted] (0,1,-1) -- (0,2,0);
\draw[dashdotted] (1,2,1) -- (0,2,0);
\draw[dashdotted] (0,3,-1) -- (0,2,0);
\end{tikzpicture}
\caption{The three possible local geometries of a soap film away from the boundary.}\label{f:Plateau}
\end{figure}

Formulated in the nineteenth century by the Belgian physicist Joseph Plateau \cite{Plateau1873}, the classical Plateau laws describe the local geometry of a soap film at points away from the boundary. They state that the film is locally diffeomorphic to one of the following: a plane, three half-planes meeting at 120$^\circ$ angles along a line, or the cone over the 1-skeleton of a regular tetrahedron (see \cref{f:Plateau}).
In this context, the ``boundary'' corresponds to the physical frame -- usually a metal or plastic wire -- dipped into soapy water. The film attaches to this wire, which dictates the boundary conditions.

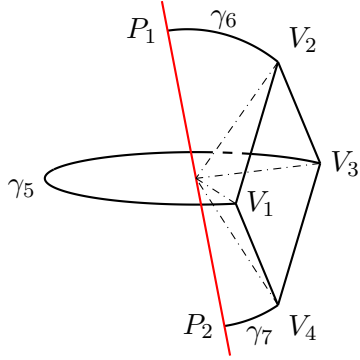
\begin{figure}[H]
\tdplotsetmaincoords{80}{20}
\begin{tikzpicture}[tdplot_main_coords, scale=2]
    \pgfmathsetmacro{\xA}{sqrt(3)/3}
    \pgfmathsetmacro{\yA}{-sqrt(6)/3}
    \pgfmathsetmacro{\zA}{0}
    \pgfmathsetmacro{\xB}{sqrt(3)/3}
    \pgfmathsetmacro{\yB}{0}
    \pgfmathsetmacro{\zB}{sqrt(6)/3}
    \pgfmathsetmacro{\xC}{sqrt(3)/3}
    \pgfmathsetmacro{\yC}{sqrt(6)/3}
    \pgfmathsetmacro{\zC}{0}
    \pgfmathsetmacro{\xD}{sqrt(3)/3}
    \pgfmathsetmacro{\yD}{0}
    \pgfmathsetmacro{\zD}{-sqrt(6)/3}
    \pgfmathsetmacro{\xE}{-1}
    \pgfmathsetmacro{\yE}{0}
    \pgfmathsetmacro{\zE}{0}
    \pgfmathsetmacro{\xP}{(-0.2)}
    \pgfmathsetmacro{\yP}{0}
    \pgfmathsetmacro{\zP}{sqrt(1-0.2*0.2)}

\pgfmathsetmacro{\angle}{acos(\xA*\xB + \yA*\yB + \zA*\zB)}        \pgfmathsetmacro{\angleBC}{acos(\xB*\xC + \yB*\yC + \zB*\zC)}    \pgfmathsetmacro{\angleCD}{acos(\xC*\xD + \yC*\yD + \zC*\zD)} \pgfmathsetmacro{\angleDA}{acos(\xD*\xA + \yD*\yA + \zD*\zA)}    \pgfmathsetmacro{\angleAE}{acos(\xA*\xE + \yA*\yE + \zA*\zE)}    \pgfmathsetmacro{\angleCE}{acos(\xC*\xE + \yC*\yE + \zC*\zE)}
\pgfmathsetmacro{\angleBP}{acos(\xP*\xB + \yP*\yB + \zP*\zB)}
\pgfmathsetmacro{\angleDQ}{acos(-\xP*\xD - \yP*\yD - \zP*\zD)}


\draw[thick] (\xA,\yA,\zA) -- (\xB,\yB,\zB);


\draw[thick] (\xB,\yB,\zB) -- (\xC,\yC,\zC);


\draw[thick] (\xC,\yC,\zC) -- (\xD,\yD,\zD);


\draw[thick] (\xD,\yD,\zD) -- (\xA,\yA,\zA);

\draw[thick, samples=100, domain=0:1] plot[variable=\t]
    (
        {(sin((1-\t)*\angleAE)/sin(\angleAE)) * \xA + (sin(\t*\angleAE)/sin(\angleAE)) * \xE},
        {(sin((1-\t)*\angleAE)/sin(\angleAE)) * \yA + (sin(\t*\angleAE)/sin(\angleAE)) * \yE},
        {(sin((1-\t)*\angleAE)/sin(\angleAE)) * \zA + (sin(\t*\angleAE)/sin(\angleAE)) * \zE}
    );

\draw[thick, samples=61, domain=0.39:1] plot[variable=\t]
    (
        {(sin((1-\t)*\angleCE)/sin(\angleCE)) * \xC + (sin(\t*\angleCE)/sin(\angleCE)) * \xE},
        {(sin((1-\t)*\angleCE)/sin(\angleCE)) * \yC + (sin(\t*\angleCE)/sin(\angleCE)) * \yE},
        {(sin((1-\t)*\angleCE)/sin(\angleCE)) * \zC + (sin(\t*\angleCE)/sin(\angleCE)) * \zE}
    );

\draw[thick, samples=24, domain=0:0.27] plot[variable=\t]
    (
        {(sin((1-\t)*\angleCE)/sin(\angleCE)) * \xC + (sin(\t*\angleCE)/sin(\angleCE)) * \xE},
        {(sin((1-\t)*\angleCE)/sin(\angleCE)) * \yC + (sin(\t*\angleCE)/sin(\angleCE)) * \yE},
        {(sin((1-\t)*\angleCE)/sin(\angleCE)) * \zC + (sin(\t*\angleCE)/sin(\angleCE)) * \zE}
    );

\draw[thick, samples=3, domain=0.27:0.3] plot[variable=\t]
    (
        {(sin((1-\t)*\angleCE)/sin(\angleCE)) * \xC + (sin(\t*\angleCE)/sin(\angleCE)) * \xE},
        {(sin((1-\t)*\angleCE)/sin(\angleCE)) * \yC + (sin(\t*\angleCE)/sin(\angleCE)) * \yE},
        {(sin((1-\t)*\angleCE)/sin(\angleCE)) * \zC + (sin(\t*\angleCE)/sin(\angleCE)) * \zE}
    );

\draw[thick, samples=3, domain=0.33:0.36] plot[variable=\t]
    (
        {(sin((1-\t)*\angleCE)/sin(\angleCE)) * \xC + (sin(\t*\angleCE)/sin(\angleCE)) * \xE},
        {(sin((1-\t)*\angleCE)/sin(\angleCE)) * \yC + (sin(\t*\angleCE)/sin(\angleCE)) * \yE},
        {(sin((1-\t)*\angleCE)/sin(\angleCE)) * \zC + (sin(\t*\angleCE)/sin(\angleCE)) * \zE}
    );

\draw[thick, samples=100, domain=0:1] plot[variable=\t]
    (
        {(sin((1-\t)*\angleBP)/sin(\angleBP)) * \xB + (sin(\t*\angleBP)/sin(\angleBP)) * \xP},
        {(sin((1-\t)*\angleBP)/sin(\angleBP)) * \yB + (sin(\t*\angleBP)/sin(\angleBP)) * \yP},
        {(sin((1-\t)*\angleBP)/sin(\angleBP)) * \zB + (sin(\t*\angleBP)/sin(\angleBP)) * \zP}
    );

\draw[thick, samples=100, domain=0:1] plot[variable=\t]
    (
        {(sin((1-\t)*\angleDQ)/sin(\angleDQ)) * \xD - (sin(\t*\angleDQ)/sin(\angleDQ)) * \xP},
        {(sin((1-\t)*\angleDQ)/sin(\angleDQ)) * \yD - (sin(\t*\angleDQ)/sin(\angleDQ)) * \yP},
        {(sin((1-\t)*\angleDQ)/sin(\angleDQ)) * \zD - (sin(\t*\angleDQ)/sin(\angleDQ)) * \zP}
    );

\draw[thick, color=red] (1.2*\xP,1.2*\yP,1.2*\zP) -- (-1.2*\xP,-1.2*\yP,-1.2*\zP);

\draw[dashdotted] (\xA,\yA,\zA) -- (0,0,0);
\draw[dashdotted] (\xB,\yB,\zB) -- (0,0,0);
\draw[dashdotted] (\xC,\yC,\zC) -- (0,0,0);
\draw[dashdotted] (\xD,\yD,\zD) -- (0,0,0);

\node[left] at (\xP,\yP,\zP) {$P_1$};
\node[left] at (-\xP,-\yP,-\zP) {$P_2$};
\node[right] at (\xA,\yA,\zA) {$V_1$};
\node[right] at (\xC,\yC,\zC) {$V_3$};
\node[above right] at (\xB,\yB,\zB) {$V_2$};
\node[below right] at (\xD,\yD,\zD) {$V_4$};
\node[below left] at (\xE,\yE,\zE) {$\gamma_5\;$};
\node[above] at (0.2,0,0.95) {$\gamma_6$};
\node[below] at (0.45,0,-0.9) {$\gamma_7$};
\end{tikzpicture}
\caption{A cone $C$ in the thin fish family. In the picture we can see the spherical arcs $\gamma_5$, $\gamma_6$, and $\gamma_7$ of Definition \ref{d:thin-fish}, while in place of the arcs $\gamma_1, \gamma_2, \gamma_3$, and $\gamma_4$ the vertices $V_1, V_2, V_3, V_4$ are joined by straight edges (which form a square). The length of the arc $\gamma_7$ is allowed to vary in a range between $0$ (when $V_4=P_2$) and a maximal length $m_0$ (when $V_2=P_1$ and the arc $\gamma_6$ vanishes).}\label{f:PlateauB4}
\end{figure}

\begin{figure}
\tdplotsetmaincoords{80}{15}
\begin{tikzpicture}[tdplot_main_coords]
\draw[thick] (-9,0,-1) -- (-8,0,-1) -- (-8,0,1) -- (-9,0,1);
\draw[thick, color=red] (-9,0,-1.2) -- (-9,0,1.2);
\draw[thick] (-6,0,-1) -- (-5,0,-1) -- (-5,0,1) -- (-6,0,1);
\draw[thick, color=red] (-6,0,-1.2) -- (-6,0,1.2);
\draw[thick] (-6,0,-1) -- ({-6-sqrt(2)*cos(70)}, {-sqrt(2)*sin(70)},-1)
-- ({-6-sqrt(2)*cos(70)}, {-sqrt(2)*sin(70)},1) -- (-6,0,1);
\draw[thick, color=red] (-3,1,-1.2) -- (-3,1,1.2);
\draw[thick] (-4,1,-1) -- (-3,1,-1);
\draw[thick] (-3,1,1) -- (-4,1,1) -- (-4,1,-1);
\draw[thick] (-3,1,1) -- (-2.5,{sqrt(3)/2+1},1) -- (-2.5,{sqrt(3)/2+1},-1);
\draw[thick] (-2.5,{sqrt(3)/2+1},-1) -- (-0.35*2.5 -0.65*3,{0.35*(sqrt(3)/2 + 1) + 0.65},-1);
\draw[thick, dashed] (-3,1,-1) -- (-0.35*2.5 -0.65*3,{0.35*(sqrt(3)/2 + 1) + 0.65},-1);
\draw[thick] (-3,1,1) -- (-2.5,{-sqrt(3)/2+1},1) -- (-2.5,{-sqrt(3)/2+1},-1) -- (-3,1,-1);
\end{tikzpicture}
\caption{The local boundary geometries (a), (b), and (c). The line $\ell$ is in red.}\label{f:PlateauB1}
\end{figure}
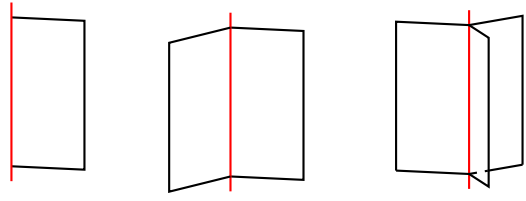

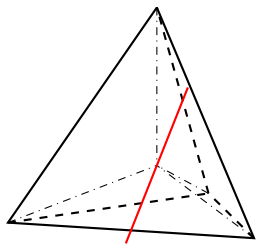
\begin{figure}
\tdplotsetmaincoords{80}{20}
\begin{tikzpicture}[tdplot_main_coords]
\draw[thick] (0,0, {2*sqrt(2)}) -- ({-sqrt(3)},-1,0) -- ({sqrt(3)}, -1, 0) -- (0,0,{2*sqrt(2)});
\draw[thick, dashed] ({sqrt(3)},-1,0) -- (0,2,0) -- (0,0,{2*sqrt(2)});
\draw[thick, dashed] ({-sqrt(3)},-1,0) -- (0,2,0);
\draw[dashdotted] ({-sqrt(3)},-1,0) -- (0,0,{sqrt(2)/2});
\draw[dashdotted] ({sqrt(3)},-1,0) -- (0,0,{sqrt(2)/2});
\draw[dashdotted] (0,2,0) -- (0,0,{sqrt(2)/2});
\draw[dashdotted] (0,0,{2*sqrt(2)}) -- (0,0,{sqrt(2)/2});
\draw[thick, red] (0,-1.2,{-0.1*sqrt(2)}) -- (0,1.2,{1.1*sqrt(2)});
\end{tikzpicture}
\caption{The local boundary geometry (d). The line $\ell$ is going through the midpoints of two opposite sides of the $1$-skeleton of the tetrahedron.}\label{f:PlateauB2}
\end{figure}

\begin{figure}
\tdplotsetmaincoords{80}{30}
\begin{tikzpicture}[tdplot_main_coords]
\tikzmath{\a=sqrt(2);}
\draw[thick] (-\a,-\a,-\a) -- (\a,-\a,-\a) -- (\a,-\a,\a) -- (-\a,-\a,\a) -- (-\a,-\a,-\a);
\draw[thick] (-\a,-\a,\a) -- (-\a,\a,\a) -- (\a,\a,\a) -- (\a,-\a,\a);
\draw[thick] (\a,-\a,-\a) -- (\a,\a,-\a) -- (\a,\a,\a);
\draw[thick, dashed] (-\a,-\a,-\a) -- (-\a,\a,-\a) -- (\a,\a,-\a);
\draw[thick, dashed] (-\a,\a,-\a) -- (-\a,\a,\a);
\draw[dashdotted] (-\a,-\a,-\a) -- (\a,\a,\a);
\draw[dashdotted] (-\a,\a,-\a) -- (\a,-\a,\a);
\draw[dashdotted] (-\a,-\a,\a) -- (\a,\a,-\a);
\draw[dashdotted] (\a,-\a,-\a) -- (-\a,\a,\a);
\draw[thick, red] (1.1*\a,1.1*\a,1.1*\a) -- (-1.1*\a,-1.1*\a,-1.1*\a);
\node[below right] at (\a,\a,\a) {$P$};
\tikzmath{\a=sqrt(2);}
\draw[thick] (5-\a,3-\a,-\a) -- (5+\a,3-\a,-\a)
-- (5+\a,3-\a,\a) -- (5-\a,3-\a,\a) -- (5-\a,3-\a,-\a);
\draw[thick] (5-\a,3-\a,\a) -- (5-\a,3+\a,\a) -- (5+\a,3+\a,\a) -- (5+\a,3-\a,\a);
\draw[thick] (5+\a,3-\a,-\a) -- (5+\a,3+\a,-\a) -- (5+\a,3+\a,\a);
\draw[thick, dashed] (5-\a,3-\a,-\a) -- (5-\a,3+\a,-\a) -- (5+\a,3+\a,-\a);
\draw[thick, dashed] (5-\a,3+\a,-\a) -- (5-\a,3+\a,\a);
\draw[dashdotted] (5-\a,3-\a,-\a) -- (5+\a,3+\a,\a);
\draw[dashdotted] (5-\a,3+\a,-\a) -- (5+\a,3-\a,\a);
\draw[dashdotted] (5-\a,3-\a,\a) -- (5+\a,3+\a,-\a);
\draw[dashdotted] (5+\a,3-\a,-\a) -- (5-\a,3+\a,\a);
\draw[thick, red] (5-1.2*0.6*\a,3-1.2*\a,-1.2*\a) -- (5+1.2*0.6*\a,3+1.2*\a,1.2*\a);
\node[above] at (5+0.6*\a,3+\a,\a) {$P$};
\end{tikzpicture}
\caption{The local boundary geometry (e). If $P$ and $-P$ are the intersections of the line $\ell$ with the skeleton of the cube, then either $P$ is a vertex or it has to be sufficiently close to a vertex.}\label{f:PlateauB3}
\end{figure}
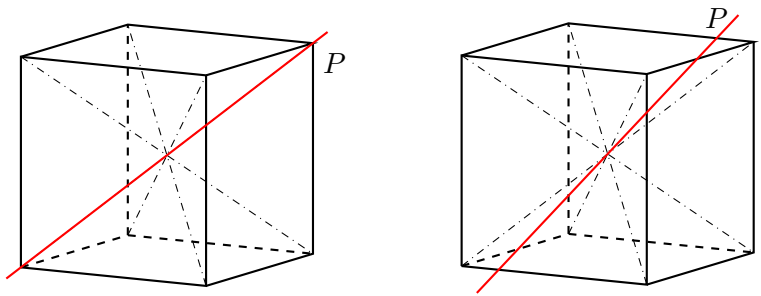

The aim of this paper is to classify the possible shapes at a point $p$ located at the boundary $\Gamma$, under the assumption that the film wets all of $\Gamma$ close to $p$. Assuming that $\Gamma$ is regular we can
move $p$ to the origin of a Cartesian system of coordinates and zoom in. This allows us to treat the boundary as an infinite straight line $\ell$ going through the origin and the film as an area-minimizing infinite cone containing $\ell$ (cf. \cref{c:blow-up} for the rigorous statement). Our list of possible shapes at $p$ comprises six families, depicted in \cref{f:PlateauB1,f:PlateauB2,f:PlateauB3,f:PlateauB4}.

\begin{lista}\label{lista}$\phantom{acapo}$
\begin{itemize}
\item[(a)] A halfplane with boundary $\ell$.
\item[(b)] Two halfplanes with boundary $\ell$ forming a convex angle $\geq \frac{2\pi}{3}$.
\item[(c)] Three halfplanes bounded by $\ell$ forming angles of $\frac{2\pi}{3}$.
\item[(d)] The cone over the skeleton of a regular tetrahedron (which contains $3$ distinct lines).
\item[(e)] The cone over the skeleton $Q$ of a cube.
\item[(f)] An element of a family which we call ``thin fish'', cf. \cref{f:PlateauB4} and \cref{d:thin-fish}.
\end{itemize}
\end{lista}

\subsection{Rigorous statements} In order to state our results in rigorous terms we first introduce the relevant objects in a more general setting.

\begin{definition}[Almgren minimizing sets]\label{d:Almgren-minimizing}
Let $k<N$ be two given positive integers, $U\subset \mathbb R^N$ an open subset, and $\Gamma$, $E$ two relatively closed subsets of $U$. $E$ is a $k$-dimensional Almgren minimizing set with respect to $\Gamma$ if $E$ has locally finite $\mathcal{H}^k$ measure in $U$ and the inequality
    \begin{equation}\label{e:Almgren-minimizing}
    \Ha^k (E\cap V) \leq \Ha^k (\Phi (E) \cap V)\,
    \end{equation}
    holds for every open set $V\ssubset U$ and every Lipschitz map $\Phi$ with $\overline{W\cup \Phi (W)} \ssubset V\setminus \Gamma$, where $W= \{x: x\neq \Phi (x)\}$.
\end{definition}

With the latter definition in hand, the focus of this paper is the classification of all $2$-dimensional cones which are Almgren minimizing with respect to a line $\ell$ contained in them.

\begin{theorem}\label{t:main}
Assume $\ell\subset \mathbb R^3$ is a line and $E\supset \ell$ a closed cone (namely invariant under the homotheties $x\mapsto \lambda x$, $\lambda > 0$) which is Almgren minimizing with respect to $\ell$. Then $E$ belongs to \cref{lista}. 
\end{theorem}

To complete the statement of \cref{t:main} we describe in \cref{d:thin-fish} the cones in \cref{lista}(f). \cref{t:main} has been achieved through a three year project, during which many steps were completed with the assistance of computers in a ``traditional sense'', namely through numerical simulations to guide our intuition and through the implementation of suitable computer algorithms to check automatically a large number of cases, after developing theoretical tools which allow a reduction to a finite case analysis. Except for helping with some visualization tools and for a few small fine-tuning of some arguments, we did not use any AI tool to complete this work. 

\begin{remark}\label{r:phys-math} 
Up until a few days before the public release of this manuscript, we were convinced that all the elements of \cref{lista} were in fact Almgren minimizing, with the exception of those elements in (e) for which the line $\ell$ is close to the midpoint of an edge of $Q$. Our opinion was motivated by both real life and numerical experiments. For instance we dipped wires shaped as the links in \cref{lista}(f) and in \cref{lista}(e) in soapy water on several occasions and we invariably saw the films forming cones (with the exceptions mentioned above among the links in \cref{lista}(e)). The same behavior has been witnessed also by Guy David, who ran independent experiments with different models of the most symmetric element of \cref{lista}(f), see \cite{Da}. Similarly, all the numerical simulations run by us using Brakke's surface evolver confirmed what we saw in the real life experiments. 

A few days before releasing this manuscript we doublechecked one of our conjectures, prompting an LLM with the case of the most symmetric element in \cref{lista}(f). The LLM suggested a relatively simple way to construct a competitor for it and, after understanding the relevant construction and the computations, we then concluded that the same idea could be applied to prove that no elements in \cref{lista}(e)\&(f) are Almgren minimizing. Since we can prove that all the others in \cref{lista} are minimizing, this in fact gives a {\em complete} classification. We also realized that a combination of this new insight and some of the ideas of this paper would give a much shorter proof which is entirely ``pen and paper'', namely with no need for a computer assisted part. We then asked LLMs to complete a second manuscript with this proof, giving a rather detailed plan of work for the overall strategy and closely steering them during the writing process. In this way, we were able to produce the companion paper \cite{DG2} in a few hours. 

Why then releasing this manuscript? First of all, several parts of it have independent interest and can be used to prove other statements. For instance we go here a long way in the direction of classifying certain type of stationary cone, and in fact some of the tools introduced here are used crucially in \cite{DG2}. Furthermore, it may have some historical interest as a record of human mathematical work compared to LLM-assisted work during this transitional period. But there is also another, perhaps more important, reason. The stronger classification in paper \cite{DG2} is at odds with what can be observed experimentally. It does seem that the boundary Plateau laws for soap films need to include shapes as in \cref{lista}(e)\&(f). How do we reconcile this with the main theorem of \cite{DG2} which excludes their Almgren minimality? We propose a way in \cref{s:phys-math}. If our suggestions are correct, then the much more complicated tools introduced in this manuscript do provide a better picture of the physical behavior of soap films than the shorter and more elegant arguments given in \cite{DG2}.
\end{remark}

\begin{definition}[The thin fish family]\label{d:thin-fish}
A cone $C\subset \mathbb R^3$ is called a {\em thin fish} if its link $C\cap \partial B_1$ is the union of geodesic arcs $\gamma_i$ of $\partial B_1$ arranged in the following way.
\begin{itemize}
\item $\gamma_1$, $\gamma_2$, $\gamma_3$, and $\gamma_4$ are congruent and form a spherical quadrilateral\footnote{Note that $F$ is the projection on the sphere $\partial B_1$ of the face of a cube inscribed in $\partial B_1$.} $F$ with angles $\frac{2\pi}{3}$.
\item $\gamma_5$ is the longest geodesic arc connecting two opposite vertices of $F$ (say $V_1$ and $V_3$).
\item $\gamma_6$ and $\gamma_7$ are two arcs which connect the other vertices, $V_2$ and $V_4$, to two antipodal points $P_1$ and $P_2$, lying outside $F$. $\gamma_6$ can be trivial (i.e. $P_1=V_2$ is allowed).
\item Any pair of arcs with a common endpoint forms an angle of $\frac{2\pi}{3}$.
\item The reference line $\ell$ is the one passing through $P_1$ and $P_2$.
\end{itemize}
\end{definition}

Note that the Almgren minimality property is invariant under addition of an $\mathcal{H}^k$-null set; thus, to be precise, the statement of \cref{t:main} must be appropriately modified. This is done by imposing on $E$ the following ``normalization'':
\begin{equation}\label{e:normalization}
\mathcal{H}^k (E\cap B_r (x)) > 0 \qquad \forall x\in E \quad\mbox{and} \quad \forall r >0 \, .
\end{equation}
Closed sets as in \cref{d:Almgren-minimizing} are actually called $(\mathbf{M}, 0, \infty)$-minimal by Almgren and in much of the subsequent literature (cf. \cite{Almgren-1976} and \cite{Taylor}). In the rest of this work we will simplify our terminology further and just use the term ``minimizing''. Furthermore it is a theorem by Almgren, cf.
\cite[Theorem II.3, Conclusion 9]{Almgren-1976}, that such sets are in fact necessarily (countably) rectifiable.
To streamline our exposition, we henceforth adopt the standing assumption that ``minimizing'' refers to a rectifiable set as in \cref{d:Almgren-minimizing} which satisfies the condition \eqref{e:normalization}.

Note that the problem of classifying cones which are minimizing with respect to a line $\ell$ can be stated more generally without assuming $\ell\subset E$. However, since $E$ is a cone, there are only two other possibilities: $\ell \cap E = \{0\}$ and $\ell \cap E$ is a halfline. The first case essentially reduces to the classical theorem of Taylor in \cite{Taylor}, cf. the discussion below, while we plan to study the remaining one some time in the future, as the methods introduced in this work can be adapted to address it.

\subsection{Almgren minimality and physical Plateau Laws}\label{s:phys-math}
Taylor's theorem \cite{Taylor} immediately implies that the cones in (c) or (d) in \cref{lista} are indeed minimizing, whereas the minimality of (a) and (b) is elementary. On the other hand, as noted in \cref{r:phys-math}, the cones in (e) and (f) are not minimizing. 

As explained in the previous section, our study is motivated by the related problems of understanding the local structure of soap films and area-minimizing sets. We refrain from discussing whether minimizing sets are a good approximation of the shape of a soap film and rather refer the reader to other works, such as \cite{AlTa} and \cite{KMS}. Concerning the abstract variational problem of minimizing the area subject to appropriate boundary conditions, we just point out that the minimizers of all variational frameworks which we know of (and aim to capture the behavior of soap films) are in fact subsets of the class introduced in \cref{d:Almgren-minimizing}, see for instance \cite{DGM}. Moreover, real-world experiments of soap films attaching to several printed models of a boundary $\Gamma$ (thanks to the assistance of the research group of Michael Niedermeier at the University of Passau, cf. \cref{s:acknowledgments}) have been of great importance in proving our theorems and have shown a remarkable agreement among experiments, theory, and numerical simulations (performed using Brakke's surface evolver, see \cite{Brakke}). The only exceptions are that, while we can disprove the minimality of cones in (e) and (f) (cf. \cite{DG2}), we do observe them in real life and numerical experiments. How can we explain such discrepancy?

One obvious possibility is that our experiments are too crude. While our wires were 3d printed in the lab of Michael Niedermeier's research group, these structures are formed by plastic tubes. In fact we used strings of nylon for the line $\ell$, upon a suggestion of David that the thickness might have an effect. Ideally our findings should be confirmed by much more accurate experiments.   

A second explanation, the one we propose here, is that when the wire has the shape of one of the links in \cref{lista}(e)\&(f), even though the actual minimum is different than the corresponding cone, it is however too close to it to be observed in an experiment (with the exception of some cases in \cref{lista}(e)).

First of all, it is possible to use the tools of this paper to show the following: for any stationary cone (cf. \cref{d:stazionarie}) which is not listed in Theorem \ref{t:main} there is a competitor which gains a ``substantial amount'' in area and coincides with the cone outside the unit ball. Without altering the structure of the manuscript, but scaling appropriately some of the needed numerical computations, this substantial amount can be quantified as being around $10^{-3}$. The thickness of a soap film in a classroom experiment ranges between 0.1 and 2 micrometers, namely between $10^{-5}$ cm and $2\cdot 10^{-4}$ cm, while the film has a typical size of a few centimeters. This substantial gain in area can also be shown when $E$ is as in \cref{lista}(e) and the line $\ell$ contains the midpoint of an edge of $Q$, cf. \cref{s:davvero-minimi}. And in this case we can see experimentally that the actual film will not form a cone, while Brakke's surface evolver does show also a good competitor. 

On the other hand, consider the case of the most symmetric element in \cref{lista}(f). The authors tried themselves and asked an LLM to optimize the ideas in the competitor construction of \cite{DG2}. The best area gain achieved in the unit ball is approximately $2\cdot 10^{-15}$. When we consider that the size of a water molecule is roughly $3\cdot 10^{-8}$ cm, such a gain would be absolutely meaningless even if Almgren minimizers were extremely accurate models of real soap films.

Let $E$ be a cone as in \cref{lista}(e)\&(f) and let $\mathscr{C} (E)$ be the family of closed sets $\Phi (E)$ for any $\Phi$ as in \cref{d:Almgren-minimizing} with $U=B_1$ and $\Gamma = \ell$. Define
\begin{align}
m (E) := \mathcal{H}^2 (E\cap B_1) - \inf \{ \mathcal{H}^2 (E' \cap B_1) : E'\in \mathscr{C} (E)\}\, .
\end{align}
We conjecture that $m (E)$ is a very small number, unless $E$ is an element of \cref{lista}(e) for which $\ell$ is close to the midpoint of an edge of $Q$. We therefore propose the following problem   

\begin{problem}
Give an upper bound for $m (E)$ for $E$ in \cref{lista}(e)\&(f).
\end{problem}
 
A related question is estimating the distance of the ``actual minimizer'' with the same boundary data from the cone $E$. In order to define things properly, consider any minimizing sequence $\{E_k\}\subset \mathscr{C} (E)$, define the measures $\mu_k := \mathcal{H}^2 \res (E_k \cap B_1)$, and let $\mu$ be the weak$^*$ limit of any subsequence. According to \cite{DGM}, the closed set $E_m := {\rm spt}\, (\mu)$ is an Almgren minimal set. We expect that the Hausdorff distance $d_H(E, E_m)$ between $E_m$ and $E\cap \overline{B}_1$ is very small in the cases under discussions.

\begin{problem}
Give an upper bound for $d_H (E, E_m)$ for $E$ in \cref{lista}(e)\&(f) and any $E_m$ as above.
\end{problem}

In particular, let $E$ be the most symmetric element of \cref{lista}(f). If it were possible to prove that $d_H (E, E_m) \leq 10^{-6}$, then clearly $E_m$ and $E$ would be indistinguishable at the physical scales which are of interest for soap films, given that the ratio between the thickness and the size of the film is in fact larger. 

\subsection{Tangent cones}
The relation between \cref{t:main} and the local structure of minimizing sets is exemplified by the seminal work \cite{Taylor}. In that paper Taylor proved that:
\begin{itemize}
\item[(i)] If we replace $\ell$ with the empty set $\emptyset$ in the statement of \cref{t:main}, then there are precisely three possibilities: $E$ is either a $2$-dimensional plane, or the union of three half-planes meeting at $120^\circ$ angles along a common line, or the cone over the $1$-skeleton of a regular tetrahedron;
\item[(ii)] if $E$ is minimizing with respect to a boundary $\Gamma$ and $x\in E\setminus \Gamma$, then there is a neighborhood $U$ of $x$ in which $E$ is diffeomorphic to one of the three cones in (i).
\end{itemize}
While we do not address the problem of showing further regularity properties (e.g. uniqueness of the tangent cone), our classification theorem, combined with Allard's classical boundary monotonicity formula and the compactness arguments of \cite{All} (or more recent and elementary ones, see e.g. \cite{DGM}), immediately implies the following corollary (cf. \cref{c:blow-up-2} and \cref{c:monot}).

\begin{corollary}\label{c:blow-up}
Assume $U\subset \mathbb R^3$ is an open set, $\Gamma \subset U$ a smooth arc containing the origin, and $E\subset U$ a minimizing set containing $\Gamma \cap U$. For every sequence $r_k\downarrow 0$ there is a subsequence, not relabeled, such that the rescaled sets
\[
E_{0,r_k} :=\{x: r_k x\in E\}
\]
converge, locally in the Hausdorff distance, to a cone $E$ which contains the line $\ell$ tangent to $\Gamma$ at $0$ and is minimizing with respect to it. Therefore $E$ belongs to \cref{lista}.
\end{corollary}

\subsection{Previous results}
The problem of finding ``boundary'' area-minimizing cones is a very natural one and is mentioned in several places in the literature, see e.g. \cite{Brakke,DavidBoundary,LM,MorganBook}. It is well known that a half-plane with boundary $\ell$ is in fact a minimizing cone with boundary $\ell$.
Allard's classical regularity theory in \cite{AllB} implies then that, if $\Gamma$ is of class $C^2$, $x\in \Gamma$, and one tangent cone to $E$ at $x$ is a half-plane (case (a) in \cref{lista}), then $E$ is diffeomorphic to it in some neighborhood of $x$. The only other result of this type is due to David in \cite{DavidBoundary}. For a notion of minimizing set which is stronger than the one of Almgren, David gives a complete description of the local behavior of $E$ at $x$ if one tangent cone is as in \cref{lista}(b) (one cannot always conclude that $E$ is diffeomorphic to the cone in a neighborhood of $x$ since it can be proved that there are some other possible behaviors, as it is explained in detail in \cite{DavidBoundary}). He also proves that if one tangent cone at $x$ is as in \cref{lista}(c) then the tangent cone is unique (i.e. the cone is the limit of the whole family $r^{-1} (E-x)$ as $r\downarrow 0$); in this case, however, he does not fully characterize the local topology of $E$. The paper \cite{LM} correctly identifies the geometries in (a), (b), (c), and (d) as minimizing, while the case of the cone over the 1-skeleton of the cube with the line $\ell$ going through two vertices is suggested in \cite{DavidBoundary} as very likely to be minimizing. The remaining cases of (e) are not mentioned anywhere, while the family in (f) appears to be a novel discovery.

\subsection{Acknowledgments}\label{s:acknowledgments}
The first author acknowledges the support of the Simons Foundation through the Simons Initative on the Geometry of Flows (Grant Award ID BD-Targeted-00017375, CDeL). Both authors are very grateful to Michael Niedermeier and D\"orte R\"uweler, who printed several models of candidate minimizing cones. The experiments involving dipping these models in soapy water have been pivotal in guiding our intuition. The first author is also grateful to Guy David, Jonas Hirsch, and Simone Steinbr\"uchel for several discussions on the present work.

\section{Plan of the paper}\label{s:plan}

The remainder of this paper is devoted to the proof of \cref{t:main}, with the exception of \cref{s:davvero-minimi}, where we discuss which cones in \cref{lista} are proven to be area-minimizing and motivate our conjectures for the remaining ones. The proof of \cref{t:main} is split into several steps, many of which rely on computer-assisted verification. In this section, we outline these steps and their logical interdependencies.

\medskip

Fix a line $\ell\subset \mathbb R^3$ and a $2$-dimensional cone $E$ containing $\ell$, which is minimizing with respect to it. Let $\mathcal{N} = E\cap \partial B_1$ denote the link of the cone. In \cref{s:reti-con-pin} we use the interior regularity theory for Almgren-minimizing sets and a classical blow-up argument to prove the following structural properties of $\mathcal{N}$ (cf. \cref{p:stationary}):
\begin{itemize}
\item $\mathcal{N}$ is the union of finitely many closed geodesic arcs $\gamma_i$.
\item Distinct edges $\gamma_i$ can only intersect at their endpoints.
\item If $p$ is an endpoint of an arc $\gamma_i$, then either $p\in \{P_1, P_2\}= \ell \cap \partial B_1$, or it is a junction where three arcs meet at equal ($120^\circ$) angles.
\item $P_1$ and $P_2$ are both endpoints of at least one arc and at most three.
\item Any pair of arcs sharing $P_1$ (or $P_2$) as an endpoint form an angle in the range $[\frac{2\pi}{3}, \pi]$.
\end{itemize}
Given their distinct role, we call the points $P_1$ and $P_2$ ``pins''. Sets $\mathcal{N}$ satisfying these conditions are called {\em stationary nets with two antipodal pins} (cf. \cref{d:reti}).

Any such $\mathcal{N}$ is naturally associated to an abstract graph embedded in $\mathbb S^2$ (cf. \cref{d:graph,d:graph-0}): the vertices correspond to the endpoints of the geodesic segments, and the edges correspond to the segments themselves. Two nets are said to have the same topology if their underlying graphs are isomorphic and their embeddings in $\mathbb S^2$ are homeomorphically equivalent (cf. \cref{d:iso-graph}).

\medskip

 Roughly speaking, our proof of \cref{t:main} relies on the following three main tasks:
\begin{itemize}
\item[(I)] ``Topological classification''. Identify meaningful restrictions on the topology of $\mathcal{N}$ and list all topologies satisfying them. These are the {\em candidate topologies}.
\item[(II)] ``Geometric realization''. For each candidate topology, classify all its possible embeddings as stationary nets with two antipodal pins. These will be called {\em realizations} of the candidate topologies.
\item[(III)] ``Area minimization''. For each realization $\mathcal{N}$ not included in \cref{lista}, show that there is a suitable competitor decreasing its area.
\end{itemize}

\subsection{Task (I)}
For any net $\mathcal{N}$, we construct a spherical embedding $\mathcal N''$ of a cubic graph (e.g., all vertices have degree three). This purely topological object is called ``twice-derived spherical graph'', see \cref{d:derived-two} (the attentive reader may observe that in this introduction we are skipping the concept of once-derived net, which is an important intermediate step between the original net and the twice-derived spherical graph, cf. \cref{d:once-derived}).

Theoretical considerations in \cref{s:derivate} imply that the number of valid twice-derived spherical graphs is finite, albeit rather large.
In order to classify them, we describe in \cref{s:primo-programma} an algorithm which produces a computer-assisted enumeration, cf. \cref{p:first-topological-list}.
Once we have classified the twice-derived spherical graphs, it is a relatively easy job to deduce a finite list of candidate topologies (because, given a twice-derived spherical graph, it is simple to enumerate the possible candidate topologies that could have generated it): this is also explained in \cref{s:primo-programma}.
Some additional geometric ``sieves'' -- simple criteria based on spherical geometry that filter some (geometrically non-realizable) topologies (see \cref{s:sieves}) -- are however used to further restrict the list of candidate topologies. These criteria are based on some important considerations in spherical geometry, which are also used in many arguments in the subsequent parts of the paper, collected in \cref{s:spher-poly}. The final list of possibilities, to which we arrive at the conclusion of this part of our paper, is given in \cref{p:second-topological-classification}.

\subsection{Task (II)}
Finding realizations of a fixed topology can be reduced to a system of algebraic equations where the unknowns are suitable trigonometric functions of the edge lengths. This system is generally overdetermined; thus, most topologies have no realization. However, the system's high degree makes a direct resolution computationally unfeasible and forces us to proceed differently.

The realizations of some exceptional topologies will be studied with pen and paper using considerations in spherical geometries, and some Gr\"obner basis techniques, cf. \cref{s:human-1} and \cref{s:human-2}. All the topologies originating from the six geometries in \cref{t:main} are in fact included in the special cases of \cref{s:human-1} and \cref{s:human-2}.

In order to study the realizations of the cases not covered in \cref{s:human-1} and \cref{s:human-2} but listed in \cref{p:second-topological-classification} we employ metric, instead of algebraic, arguments.
As a first step, thanks to the area-minimizing property of the cone over the net, we show that the edges of the net must all satisfy certain lower bounds on their edge lengths. This simple fact can be proved by compactness, but in our case getting (good) {\em effective} lower bounds is of fundamental importance, cf. \cref{s:lower-bound}. The arguments rely on the construction of good competitors. Since, however, the admissibility constraints for competitors in \cref{d:Almgren-minimizing} are difficult to work with, we will dedicate some efforts to show that they can be relaxed in several ways, cf. \cref{s:relaxed}

Thanks to the bounds in \cref{s:lower-bound} and to a number of geometric tools (cf. \cref{s:spher-trig} and \cref{s:path-isometries}), we are then able to devise an efficient computer algorithm which, given an assigned precision $\boldsymbol{\delta}$ and a fixed topology:
\begin{itemize}
    \item Either proves that the topology does not have a geometric realization;
    \item Or, for each edge of the topology, returns an interval of size at most $\boldsymbol{\delta}$. The algorithm proves that, for any geometric realization of the net, the length of any edge belongs to the corresponding interval.
\end{itemize}
The details of the algorithm are explained in \cref{s:secondo-programma} and the theoretical tools employed in it are the lemmata and propositions proved in \cref{s:spher-trig} and \cref{s:path-isometries}. The precision $\boldsymbol{\delta}$ is set (approximately, depending on the net to $10^{-4}$. By running our algorithm for a few hours on a standard workstation over the topologies not covered by \cref{s:human-1} and \cref{s:human-2}, we are then able to show that only a handful of topologies, more precisely $26$, are realizable.
For these we have a suitable approximation, up to the precision $\boldsymbol{\delta}$, of any possible realization.

\subsection{Task (III)}
We are left with $26$ geometric realizations, some found in \cref{s:secondo-programma} and some found in \cref{s:human-2}, which do not belong to \cref{lista}. These must be excluded on ``energetic'' grounds, more precisely we shall show that the associated cones, which are in fact stationary, are not minimizing. The strategy is to fix a finite number of ``approximations'' of such nets $\mathcal{N}$. In the cases found in \cref{s:secondo-programma} we in fact only know the lengths of the edges up to an error $\boldsymbol{\delta}$. But even for the cases found in \cref{s:human-2}, for which we have an {\em exact} knowledge, we introduce suitable approximations because they come as infinite (one-parameter and two-parameters) families. We must therefore prove a sufficiently high area saving on the approximation to guarantee that there is a similar competitor for the hypothetical exact net. \cref{s:competitore} provides the rigorous estimates required to ensure the gain outweighs the approximation error. Akin to what happens with the lower bound of \cref{s:lower-bound}, a compactness argument gives a proof of concept of the latter idea, but in reality we need an explicit and sufficiently accurate estimate.

\medskip

At this stage, we thus have a finite collection $\mathcal{C}$ of cones (which provide a sufficiently good approximation of any cone over any possible existing realization of a candidate topology), and for each element of the collection we shall show the existence of a competitor (in the sense of \cref{d:Almgren-minimizing}) which gains a certain specified amount $\boldsymbol{\eta}$ in area. The competitors are stored in \verb|data/competitors.zip| in \cite{repo} (and can be visualized at \cite{website}) and each consists of a collection of triangles. More rigorously, each stored object is a suitable $3$-dimensional piecewise linear immersion of a finite 2-dimensional simplicial complex, cf. \cref{d:simplicial-complex} and \cref{d:piecewise-linear}. These have been generated with the help of Brakke's surface evolver (with some appropriate modifications). The evolver was intialized with a suitable topological deformation of the corresponding cone $E\in \mathcal{C}$, which has been guessed using real world experiments. We then let it run with some specific choices of parameters (in some cases we needed to try more than one choice) for a suitable amount of time. It is however not necessary to know how Brakke's surface evolver produces the competitor; in particular our proof does not rely on Brakke's code being error-free. Rather, another appropriate code, written by us and described in \cref{s:terzo-programma}, checks the following two relevant conditions for each competitor $E'$ in the database:
\begin{itemize}
\item[(i)] $E'$ gains $\boldsymbol{\eta}$ in area over the respective element $E\in \mathcal{C}$.
\item[(ii)] $E'$ is valid in the sense of \cref{d:Almgren-minimizing}.
\end{itemize}
Point (i) has a straightforward implementation. Point (ii) is instead the most difficult aspect of this part of the proof. Recall that according to \cref{d:Almgren-minimizing}, a competitor for $E$ must be of the form $\Phi (E)$ for a Lipschitz map $\Phi$ satisfying some specific requirements. In \cref{s:triangulations-and-maps}, we derive a list of combinatorial properties for $E'$ which guarantees the existence of a map $\Phi$ as in \cref{d:Almgren-minimizing}. These properties can be described in finite terms, which enables us to check them through a computer code. In fact, strictly speaking, it is not $E'$ which is {\em directly} a valid competitor for $E$, but rather a suitable modification of $E'$, cf. \cref{s:triangulations-and-maps}. Our arguments show, however, that this suitable modification gains at least as much area as $E'$.

\section{Index of Concepts and Notation}

\subsection{Index of Concepts}
\renewcommand{\arraystretch}{1.3}
\begin{longtable}{p{0.60\textwidth} p{0.30\textwidth}}
    Almgren minimizing set & \cref{d:Almgren-minimizing} \\
    Approx. metrized marked spherical graph (AMMSG) & \cref{s:secondo-programma} \\
    Boundary distance / Boundary path & \cref{d:b-distance} \\
    Children (and Nice children) & \cref{s:genera-lista} \\
    Collapse (of a simplicial complex) & \cref{d:collapse} \\
    Combinatorial map & \cref{d:combinatorial} \\
    Extinct (net or graph) & \cref{d:extinction} \\
    Fat fish & \cref{d:pesce-grasso} \\
    Goggles & \cref{d:goggles} \\
    Graph Embedding & \cref{d:graph} \\
    Heppes net & \cref{s:preliminari} \\
    Marked spherical graph & \cref{d:candidate} \\
    Metrized marked spherical graph (MMSG) & \cref{s:secondo-programma} \\
    Once-derived net & \cref{d:once-derived} \\
    Path-isometry & \cref{d:path-isometry} \\
    Pins & \cref{d:stazionarie} \\
    Planar map & \cref{subsec:spherical-planar} \\
    Polygon / Polygonal path (Spherical) & \cref{d:polygonal-path,d:polygon} \\
    Positively curved oriented arc & \cref{d:positive-geod} \\
    Simplicial complex / Piecewise linear map & \cref{d:simplicial-complex,d:piecewise-linear} \\
    Stationary geodesic net (with pins) & \cref{d:stazionarie} \\
    Thin fish & \cref{d:thin-fish} \\
    Twice-derived spherical graph & \cref{d:derived-two} \\
    Wedge & \cref{d:wedge}
\end{longtable}

\subsection{Index of Notation}
\renewcommand{\arraystretch}{1.3}
\begin{longtable}{p{0.60\textwidth} p{0.30\textwidth}}
    $C_\theta(\ell)$ & \cref{s:relaxed} \\
    $\mathrm{Darts}(G)$ & \cref{d:darts} \\
    $E(v,t)$ & \cref{s:path-isometries} \\
    $\exp_t(v)$ & \cref{s:path-isometries} \\
    $F_{\mathcal{N}}(d, d')$ & \cref{def:once-derived-relative-position} \\
    $F_{(G, \mathrm{len})}$ / $F_{(G, \mathcal{I})}$ & \cref{def:mmsg-relative-position,def:ammsg-relative-position} \\
    $\Gamma(X)$ & \cref{s:graphs-empeddings} \\
    $\iota(\sigma, F)$ & \cref{s:secondo-programma} \\
    $\mathrm{next}(d)$ & \cref{d:darts} \\
    $R(v,\alpha)$ & \cref{s:path-isometries} \\
    $\mathrm{rev}(d)$ & \cref{s:graphs-empeddings} \\
    $\mathrm{rot}_\alpha(v)$ & \cref{s:path-isometries} \\
    $T\mathbb{S}^2$ / $T^1\mathbb{S}^2$ & \cref{d:tangent-bundle}
\end{longtable}

\section{Stationary nets with two pins}\label{s:reti-con-pin}

In this first step we use well-known tools in geometric measure theory to infer some basic structural properties for the link of a cone which contains a line $\ell$ and is minimizing with respect to it. The main statement is \cref{p:stationary} below. We start with some relevant definitions.

\begin{definition}\label{d:reti}
A geodesic net on the $2$-dimensional sphere $\Sf\subset \R^3$ consists of
a finite set (possibly empty) of distinct points $V_1, \ldots, V_n\in \Sf$ (called {\em vertices}) and a finite set of distinct
{\em edges} $\sigma_1, \ldots, \sigma_m$ such that
\begin{itemize}
\item[(a)] Each $\sigma_i$ is either an equator or a closed geodesic arc whose endpoints are two distinct vertices $V_i$ and $V_j$;
\item[(b)] Distinct edges $\sigma_i$ can only intersect at their endpoints;
\item[(c)] Each $V_i$ is the endpoint of at least one edge.
\end{itemize}
The support of the net is the subset of $\Sf$ which is the union of the edges. We say that the net is connected if the support is connected.
\end{definition}

This definition allows for the trivial empty net and for the net consisting of an empty set of vertices and a single equator as its only edge. Moreover, two edges can form an angle of $\pi$ at a common endpoint, even if they are the only two edges with that endpoint. This means that the net is not entirely determined by its ``support''. For instance, if the support of the net consists of a single equator, the net could have, as the example just mentioned, an empty set of vertices and a single edge, but also $N$ vertices and $N$ edges with any $N\geq 1$.

\begin{definition}\label{d:stazionarie}
A geodesic net is said to have $k$ pins if there is a subset of the vertices (called {\em pins}) of cardinality $k$ such that:
\begin{itemize}
\item[(a)] Any vertex which is not a pin is the endpoint of exactly three edges;
\item[(b)] A pin is the endpoint of at most three edges.
\end{itemize}
In addition we say that the net is stationary if 
\begin{itemize}
\item[(c)] Any pair of edges with a common endpoint $V$ form an angle which is at least $\frac{2\pi}{3}$.
\end{itemize} 
\end{definition}

Note the following consequence of the above definition: at any non-pinned vertex the three edges of a stationary net must form three angles of $\frac{2\pi}{3}$, and of course this occurs also at a pin which is the endpoint of three edges. Moreover, given condition (c) in \cref{d:reti}, we know that a pin is the endpoint of at least one and at most three edges. Finally, note that since the only vertices which might have degree two are the pins, only at these special points there can be two incoming edges forming a ``flat'' angle $\pi$.
The following proposition explains why we consider this possibility: the pins must be ``marked'' even in the latter case, because they serve a particular purpose.

\begin{proposition}\label{p:stationary}
Assume $E\subset \mathbb R^3$ is a cone containing a line $\ell$ which is minimizing with respect to it. The link of $E$, namely $E\cap \partial B_1$, is a stationary net with two antipodal pins, which are the two points forming the set $\ell \cap \partial B_1$.
\end{proposition}

Note that we are not claiming that the net $E\cap \partial B_1$ is connected. We will however prove this fact in a later stage. The rest of this section is devoted to the proof of \cref{p:stationary} and we will show, in passing, how \cref{c:blow-up} can be derived from \cref{t:main}. The first important tool will be the following general compactness result.

\begin{proposition}\label{p:compactness}
Assume $U\subset \mathbb R^N$ is a bounded open set and $\{\Gamma_j\}$ a sequence of $C^2$ $(k-1)$-dimensional submanifolds of $U$ which converges in the $C^2$ topology to a $(k-1)$-dimensional submanifold $\Gamma$. Assume $\{E_j\}$ is a sequence of relatively closed subsets of $U$ such that
\begin{itemize}
\item $E_j$ is minimizing with respect to $\Gamma_j$ in $U$;
\item $\sup \mathcal{H}^k (E_j) < \infty$;
\item $E_j$ converges locally in the Hausdorff distance to a closed set $E$.
\end{itemize}
Then,
\begin{itemize}
\item $E$ is minimizing with respect to $\Gamma$;
\item $\lim_j \mathcal{H}^k (V\cap E_j) = \mathcal{H}^k (V\cap E)$ for every open or closed $V\subset U$ with $\mathcal{H}^k (E\cap \partial V)=0$.
\end{itemize}
\end{proposition}

The proposition was first proved by Almgren in \cite{Almgren-1976} (even though it does not appear there in this precise form), but for the reader's convenience we will give a proof in \cref{a:allard} using Allard's compactness theorem for stationary integral varifolds (other methods can also be used, see e.g. \cite{DGM}, but Allard's theory seems to be the quickest route).
In fact, two further ingredients will be inferred from Allard's theory: the monotonicity formula (both at the boundary and in the interior, cf. \cref{c:monot} below) and some density lower bounds (cf. \cref{l:lower-bound}). These tools can be applied in our case because minimizing sets (in the sense of Almgren) induce integral stationary varifolds, thanks to the simple \cref{l:varifolds}.

In what follows we denote by $A_\Gamma$ the second fundamental form of $\Gamma$.

\begin{corollary}\label{c:monot}
Let $U\subset \mathbb R^N$ be an open set, $\Gamma \subset U$ a $C^2$ $(k-1)$-dimensional submanifold of $U$, and $E$ a minimizing set in $U$ with respect to $\Gamma$. Fix a point $p\in U$.
If $p\not \in \Gamma$, then
\begin{equation}\label{e:monot}
\frac{d}{dr} \frac{\mathcal{H}^k (E\cap B_r (p))}{r^k} \geq 0 \qquad
\mbox{for } r\in (0, \dist\, (p, \Gamma\cup \partial U))\, .
\end{equation}
If $p\in \Gamma$ and $\|A_\Gamma\|_\infty = \sup_\Gamma |A_\Gamma|$, we have
\begin{equation}\label{e:monotb}
\frac{d}{dr} \left(\exp (C\|A_\Gamma\|_\infty r)\frac{\mathcal{H}^k (E\cap B_r (p))}{r^k}\right) \geq 0 \qquad \mbox{for } r\in (0, \dist\, (p, \partial U))\, ,
\end{equation}
where $C$ is a dimensional constant.

Moreover, if  $\frac{\mathcal{H}^k (E\cap B_{r_0} (p))}{r_0^k}=\frac{\mathcal{H}^k (E\cap B_{r_1} (p))}{r_1^k}$ and
\begin{itemize}
    \item[(i)] either $p\not\in \Gamma$ and $0<r_0<r_1<d(p, \Gamma\cup \partial U)$,
    \item[(ii)] or $\Gamma$ is an affine space including $p$ and $0<r_0<r_1<d(p, \partial U)$,
\end{itemize}
then $E\cap B_{r_1} (p)$ is a cone with vertex $p$ and base $E\cap \partial B_{r_1} (p)$.
\end{corollary}

The first part of the Corollary can be found in \cite{All} and \cite{AllB}. For the second part, the conclusion under assumption (i) can be found in \cite{All}; under assumption (ii) it is folklore but not explicitly stated in \cite{AllB}: in the appendix we show how to derive it quickly from the tools developed in that paper.

 We next register a consequence of the results stated thus far, which is just a generalization to all dimensions and codimensions of the first statement of \cref{c:blow-up}.

\begin{corollary}\label{c:blow-up-2}
Let $U\subset \mathbb R^N$ be open and $\Gamma\subset U$ a $(k-1)$-dimensional $C^2$ submanifold. Assume that $E$ is minimizing with respect to $\Gamma$ in $U$ and fix $p\in \Gamma$. Then, for every sequence $r_j\downarrow 0$, there is a subsequence, not relabeled, such that the following holds:
\begin{itemize}
\item[(a)] the sets $E_{p,r_j} = r_j^{-1} (E-p)$ converge, locally in the Hausdorff distance, to a closed set $E'$ of $\mathbb R^N$ which is a cone with vertex at the origin;
\item[(b)] $E'$ is minimizing with respect to $T_p \Gamma$;
\item[(c)] If $V\cap \Gamma\subset E$ for some neighborhood $V$ of $p$, then $T_p \Gamma \subset E'$.
\end{itemize}
\end{corollary}

\begin{proof}
By compactness (and a diagonal argument), we can extract a subsequence of radii (not relabeled) such that $E_{p, r_j}\to E'$ locally in the Hausdorff distance, where $E'$ is a closed set. We want to apply \cref{p:compactness} to show that $E'$ is minimizing with respect to $T_p\Gamma$. To do so, we just need to prove an upper bound on the area of $E_{p,r_j}$. Observe that, for every fixed $R>0$,
\[
\mathcal{H}^k (E_{p,r_j}\cap B_R) = R^k \frac{\mathcal{H}^k (E\cap B_{R r_j} (p))}{(Rr_j)^k}
\]
is uniformly bounded by \cref{c:monot}.

To prove that $E'$ is a cone, we prove that the mass-ratio in balls centered at the origin is (almost) constant in the radius and then we apply \cref{c:monot}.
Let
\[
\Theta (p) := \lim_{r\downarrow 0} \frac{\mathcal{H}^k (E\cap B_{r} (p))}{\omega_k r^k}\, ,
\]
which exists by \cref{c:monot}. Since $E'$ is rectifiable, except for at most a countable number of exceptional radii $r$, we have $\mathcal{H}^k (E'\cap \partial B_r) = 0$. Thus, \cref{p:compactness} implies that, for every such $r$,
\[
\mathcal{H}^k (E'\cap B_r) = \omega_k r^k \lim_{j\to \infty} \frac{\mathcal{H}^k (E\cap B_{r r_j} (p))}{\omega_k (rr_j)^k} = \omega_k r^k \Theta (p)\, .
\]
Therefore, by the second part of \cref{c:monot}, we infer that $E'$ is a cone with vertex at the origin.

Finally, observe that $\Gamma_{p,r_j}:=r_j^{-1}(\Gamma-p)\subseteq E_{p,r_j}$. Moreover, $\Gamma_{p,r_j}$ converge, locally in the Hausdorff distance, to $T_p\Gamma$. Since $E_{p,r_j}$ converge, locally in the Hausdorff distance to $E'$, we conclude that $T_p\Gamma\subseteq E'$.
\end{proof}

Before coming to the proof of \cref{p:stationary}, we state the following simple lemma.

\begin{lemma}\label{l:invariant}
Let $F\subset \mathbb R^2$ be a one-dimensional rectifiable set with locally finite measure such that $E = F\times \mathbb R \subset \mathbb R^3$ is minimizing with respect to $\ell = \{(0,0)\} \times \mathbb R$. Assume $R>0$ and $\Psi: \mathbb R^2\to \mathbb R^2$ is a Lipschitz map for which there is a closed convex set $K\subset B_R (0)\setminus \{0\}$ contanining $\{x:\Psi (x)\neq x\}$ and its image through $\Psi$.
Then $\mathcal{H}^1 (\Psi (F)) \geq \mathcal{H}^1 (F)$.
\end{lemma}
\begin{proof} We argue by contradiction and assume that $\Psi: \mathbb R^2 \to \mathbb R^2$ is a Lipschitz map satisfying the properties above and violating the inequality.
We next fix a large number $L$ (to be chosen later) and define
\[
\Phi (x_1, x_2, x_3) =
\left\{
\begin{array}{ll}
(\Psi (x_1, x_2), x_3) &\mbox{if $|x_3|\leq L$}\\
((|x_3|-L) (x_1,x_2) + (L+1-|x_3|) \Psi (x_1, x_2), x_3) \quad &\mbox{if $L \leq |x_3| \leq L+1$}\\
(x_1,x_2,x_3) &\mbox{otherwise.}
\end{array}
\right.
\]
$\Phi$ is clearly Lipschitz and is the identity outside of the open set
$U := B_R (0) \times (-L-1, L+1)$ and in a neighborhood of $\ell$. Moreover $W := \overline{\{\Phi (x) : x\neq \Phi (x)\}}$ is by construction contained in $K\times [-L-1, L+1]$ and hence its closure is at a positive distance from $\ell$. By minimimality of $E$, we must have
\begin{equation}\label{e:no-gain}
\Delta := \mathcal{H}^2 (\Phi (E) \cap U) - \mathcal{H}^2 (E\cap U) \geq 0\, .
\end{equation}
Observe that $|\partial_{x_1} \Phi| \leq 1+|\partial_{x_1} \Psi|$ and $|\partial_{x_2} \Phi|\leq 1+|\partial_{x_2} \Psi|$. Moreover, on the open set $U$ we can estimate $|\partial_{x_3} \Psi| \leq 1 + R + |\Psi| \leq 1+2R$. In particular the Lipschitz costant of $\Phi$ can be bounded by a number $M$ which is independent of $L$. Therefore we can compute
\begin{equation}
\Delta \leq 2L (\mathcal{H}^1 (\Psi (F)\cap B_R (0)) - \mathcal{H}^1 (F\cap B_R (0)) + 2 (M^2-1) \mathcal{H}^1 (F\cap B_R (0))\, . \label{e:gain}
\end{equation}
Clearly, given that $M$ is fixed while we can choose $L$ very large, we would then conclude $\Delta < 0$, which is at odds with \eqref{e:no-gain}.
\end{proof}

\begin{proof}[Proof of \cref{p:stationary}]
First of all, observe that, by Taylor's regularity theorem, at any point $p\in E\setminus \ell$, the set $E$ is locally either a smooth submanifold, or it is diffeomorphic to a triple junction, or it is diffeomorphic to the cone over the skeleton of a tetrahedron. However, because $E$ itself is a cone, we also see immediately that any tangent cone to $E$ at $p$ must be invariant under translation along the direction $v=\frac{p}{|p|}$. Since the tetrahedron is not invariant under any translation, we conclude that $E$ at $p$ is either diffeomorphic to a triple junction or it is a smooth submanifold. If it is a smooth submanifold, then $B_r (p) \cap E\cap \partial B_1$ must be a geodesic arc for any $r$ sufficiently small; if it is a triple junction, then $B_r (p) \cap E \cap \partial B_1$ must be the union of three geodesic arcs meeting at $p$ and forming equal angles there (again for $r$ sufficiently small).

Fix now a system of coordinates such that $\ell \cap \partial B_1 = \{e,-e\}$ for some unit vector $e\in\R^3$.
We show that at the point $e$ (and $-e$, but by symmetry we will focus on $e$ only) one of the following three things occur for $r$ sufficiently small:
\begin{itemize}
\item[(a)] Either $\partial B_1 \cap E \cap B_r (e)$ is a single geodesic segment connecting the point $e$ to some point on $\partial B_r (e) \cap \partial B_1$;
\item[(b)] or $\partial B_1 \cap E \cap B_r (e)$ is the union of two geodesic segments forming at $e$ an angle $\alpha \in [\frac{2\pi}{3}, \pi]$;
\item[(c)] or $\partial B_1 \cap E \cap B_r (e)$ is the union of three geodesic segments forming at $e$ three angles of $\frac{2\pi}{3}$.
\end{itemize}
First of all, we use \cref{c:blow-up-2} to infer that, up to subsequences, $E_{e,r_j} := r_j^{-1} (E-e)$ converge locally in the sense of Hausdorff, to some set $E'$ containing $\ell$ which is minimizing with respect to $\ell$ and is a cone with vertex at the origin. Note however that, since $E$ itself is a cone, $E'$ must be invariant along the direction $e$, which is tangent to $\ell$. So, $E' = F\times \ell$, where $F\subset \ell^\perp$ is a $1$-dimensional stationary cone in the punctured plane $(\ell^\perp\setminus \{0\})$. Note that, since $F$ is a cone, for every $x\in \partial B_1 \cap F$ it contains the entire halfline originating in $0$ and passing through $x$. On the other hand, since $F$ has locally finite Hausdorff $1$-dimensional measure there can only be a finite number $N$ of these lines.

We can now apply \cref{l:invariant} to see that no pair of distinct half-lines among those forming $F$ can make an angle which is (strictly) smaller than $\frac{2\pi}{3}$, otherwise we can use the competitor illustrated in \cref{f:angle>=120} to contradict its minimality (this is a very well known elementary argument, which we report for the reader's convenience).

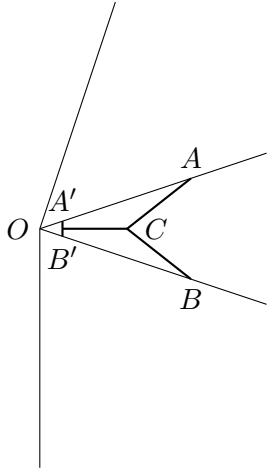
\begin{figure}
\begin{tikzpicture}
\draw (0,0) -- (1,3);
\draw (0,0) -- (3,1);
\draw (0,0) -- (3,-1);
\draw (0,0) -- (0,{-sqrt(10)});
\draw[thick] (2,{2/3}) -- ({2/sqrt(3)}, 0) -- (2,{-2/3});
\draw[thick] (0.3,0) -- ({2/sqrt(3)}, 0);
\draw[thick] (0.3,-0.1) -- (0.3,0.1);
\node[left] at (0,0) {$O$};
\node[above] at (2,{2/3}) {$A$};
\node[below] at (2,{-2/3}) {$B$};
\node[above] at (0.3,0.1) {$A'$};
\node[below] at (0.3,-0.1) {$B'$};
\node[right] at ({2/sqrt(3)+0.1}, 0) {$C$};
\end{tikzpicture}
\caption{Following a well-known argument, we illustrate the competitor disproving the minimality of $F$ if two of its halflines form an angle strictly smaller than $\frac{2\pi}{3}$. W.l.o.g. we can pick the halflines to be a pair of closest ones. We then choose two points $A$ and $B$ on them and connect with the origin $O$ forming a triple junction at some additional point $C$. Because we cannot change the set in a neighborhood of the origin, we ``reconnect'' with the segments $OA$ and $OB$ in a neighborhood of $O$. It is easy to construct a Lipschitz map which deforms the union of the two segments $AA'$ and $BB'$ on the triple junction, while keeping the rest of the set $F$ unchanged. The map can be easily arranged to satisfy the assumptions of \cref{l:invariant} and when $A'$ and $B'$ are sufficiently close to the origin we gain in length.} \label{f:angle>=120}
\end{figure}

The half-lines forming $F$ are, therefore, $0\leq N \leq 3$. Note furthermore that the number $N$ is determined by the density
\[
\Theta (e) = \lim_{r\downarrow 0} \frac{\mathcal{H}^2 (E\cap B_r (e))}{\pi r^2}\, ,
\]
indeed we have $\Theta (e) = \frac{N}{2}$. So, either all tangent cones are single half-lines, or they are all formed by two halflines, or they are all formed by three half-lines. Note that $N$ must be positive because of the condition (3) of \cref{d:Almgren-minimizing}.

We can now apply Allard's interior regularity theorem to infer that, if $r$ is sufficiently small, then
\[
E\cap B_r (e) \setminus B_{r/8} (\ell)
\]
consists of either one, two or three minimal surfaces, with boundaries contained in $\partial (B_r (e) \setminus B_{r/8} (\ell))$, close to either a single half-plane, or to the union of two, or to the union of three (and we can choose all these half-planes so that their boundaries all coincide with the line $\ell$). Since $E$ is a cone, this implies that, for all sufficiently small radii $r$, $E\cap \partial B_1 \cap (B_r (e)\setminus B_{r/8} (e))$ consists of $N=1,2$, or $3$ geodesic segments, each segment joining one endpoint in $\partial B_r (e)$ with an endpoint in $\partial B_{r/4} (e)$. We can now apply the latter conclusion to $r$ and $\frac{r}{2}$: using that the sets $B_r (e) \setminus B_{r/4} (e)$ and $B_{r/2}\setminus B_{r/8} (e)$ overlap, we infer that the same description holds for the set $E\cap (B_r (e)\setminus B_{r/8} (e)) \cap \partial B_1$. Proceeding inductively, we in fact conclude that the same is true in $(B_r (e)\setminus B_{2^{-N} r} (e))\cap \partial B_1$. Letting $N\uparrow \infty$, we conclude that $E\cap \partial B_1 \cap (B_r (e)\setminus \{e\})$ consists of either one, or two, or three geodesic segments $\sigma_i$, with $e$ being a common endpoint to all of them. On the other hand this also implies that $E\cap B_r (e)$ coincides in fact with $E'\cap B_r (e)$ for its tangent cone $E'$, which is then unique. Therefore any pair of distinct geodesic segments $\sigma_i$ and $\sigma_j$ meet at $e$ forming an angle which is at least $\frac{2\pi}{3}$. This completes the proof.
 \end{proof}

\section{Derived nets and topological graphs}\label{s:derivate}

\subsection{Preliminaries on stationary nets with two antipodal pins}\label{s:preliminari}

In this section we prove some useful facts on stationary nets with two antipodal pins. First, we recall that stationary nets with $0$ pins were classified by Heppes in \cite{Heppes}, after a preliminary attempt by Lamarle in \cite{Lamarle}, who missed one of the nets in the list (and did not provide a proof that the list was complete). We will therefore call them ``Heppes nets''. Let now $\mathcal{N}$ be the union of the geodesic arcs forming a Heppes net. If this subset of $\mathbb S^2$ contains two antipodal points $P_1$ and $P_2$, we then obtain easily a stationary net with two antipodal pins. There are in principle three possibilities:
\begin{itemize}
    \item[(i)] Both $P_i$'s are vertices of the original net; then the net is already also a stationary net with two antipodal pins;
    \item[(ii)] One point, say $P_1$, is a vertex, but the other one, i.e. $P_2$, is internal to some segment $\sigma$ of the net $\mathcal{N}$. In this case we add $P_2$ to the collection of vertices and ``split'' the segment $\sigma$ into two segments meeting at $P_2$; naturally, the two new segments ``meet'' at their common endpoint forming an angle of $\pi$;
    \item[(iii)] Both $P_i$'s are internal points of two segments $\sigma_j$ and $\sigma_k$ of $\mathcal{N}$; we then add the pins $P_1$ and $P_2$ to the collection of vertices and subdivide each of the segments $\sigma_j$ and $\sigma_k$ as described in point (ii) above.
\end{itemize}
A simple inspection of the Heppes nets shows that (ii) never occurs. On the other hand, both (i) and (iii) occur in a number of cases. The above considerations motivate the following terminology: a net which is obtained from a Heppes net with the above procedure will be called ``Heppes with two antipodal pins''.

The main goal of this subsection is thus to prove the following proposition. In what follows we will call face of a stationary net $\mathcal{N}$ any connected component of $\mathbb S^2\setminus \mathcal{N}$
\begin{proposition}\label{p:starting}
Let $\mathcal{N}$ be a stationary net with two antipodal pins $P_1$ and $P_2$. Then:
\begin{itemize}
\item[(i)] At least one of the two pins has degree strictly less than $3$ unless $\mathcal{N}$ is Heppes.
\item[(ii)] $\mathcal{N}$ is connected; in particular any connected component of $\mathbb S^2\setminus \mathcal{N}$ is an open disk.
\item[(iii)] If $\mathcal{N}$ has only two vertices, then it consists of one, two, or three half equators connecting the two pins.
\item[(iv)] If $\mathcal{N}$ has at least three vertices, then the two pins cannot belong to the boundary of the same connected component of $\mathbb S^2\setminus \mathcal{N}$ {(in particular an edge $\sigma$ which is incident to $P_1$ cannot have the other endpoint in common with an edge $\tau$ which is incident to $P_2)$}.
\end{itemize}
\end{proposition}

In the proof we will need the following lemmas and convexity will play a crucial role in both.

\begin{remark}\label{r:convexity-tricky}
We warn the reader that a good definition of {\em general} convex subsets of the sphere presents some challenges, as witnessed by the following example. Let $p$ and $-p$ be two antipodal points and consider two distinct equators passing through them. Their complement consists of four open sectors of the sphere, which are convex according to any possible reasonable definition. It seems natural to call their closures convex. However, if we take the closure of two opposite open sectors, we then find two closed convex sets whose intersection consists of a pair of antipodal points. For this reason we restrict the notion of convexity to open sets of the sphere. In the latter case several equivalent definitions can be given, for instance we can require that open convex subsets of $\mathbb S^2$ are intersections of open hemispheres, cf. \cref{s:convessita}.
\end{remark}

\begin{lemma}\label{lem:stationary-in-hemisphere}
    Let $\mathcal N$ be a stationary net with at most two pins.
    Assume that $\mathcal N$ is contained in a (closed) hemisphere. Then, $\mathcal N$ consists of two  vertices connected by either one or two edges (observe that in the latter case the two vertices must be antipodal).
\end{lemma}
\begin{proof}
    Let $\Omega$ be the intersection of all closed hemispheres containing $\mathcal N$. Observe that the cone $\Lambda$ with link $\Omega$ and vertex $0$ is then the Euclidean convex hull of the cone $E$ over $\mathcal{N}$.
Let $p\in \mathbb S^2$ be an extremal point of $\Omega$, namely a point of $\Lambda$ with the property that every line which is going through $p$ but not through the origin intersects $\Lambda$ only at $p$.Then $p\in \mathcal{N}$. Since $p$ cannot be in the interior of an edge, it must be a vertex of the net. Furthermore $p$ cannot be a non-pin vertex and if it is a pin, then it can have at most $2$ edges coming to it. Since $\Omega$ must have at least two extremal points, we conclude that in fact it has exactly two extremal points, which are the two pins $P_1$ and $P_2$ of $\mathcal{N}$. Consequently, $\Omega$ is either the geodesic connecting $P_1$ and $P_2$ or a sector of $\mathbb S^2$, delimited by two half-equators which have $P_1$ and $P_2$ as endpoints. If it is the geodesic, we have proven the statement. If it is a sector, then, since no non-pin vertices can belong to the boundary of $\Omega$, the two half-equators are edges of $\mathcal N$. To prove that there are no other vertices in $\mathcal N$, observe that the other connected components of $\mathcal N$ are Heppes nets and cannot belong to a hemisphere.
\end{proof}

\begin{lemma}\label{lem:face-convexity}
    Let $\mathcal N$ be a connected stationary net with some pins and let $F$ be one of its faces.
    \begin{enumerate}
        \item If the boundary of $F$ does not contain any pins, then $F$ is convex.
        \item If the boundary of $F$ contains one pin, then $F$ is the union of the closures of two open convex sets.
    \end{enumerate}
\end{lemma}
\begin{proof}
    Let us start with (i).
    Since $\mathcal N$ is connected, $F$ is topologically an open disk with connected boundary. Observe that, a priori, it is not obvious that the boundary of $F$ is topologically a loop.
    The \emph{combinatorial} boundary of $F$ is a sequence of directed edges of $\mathcal N$ connecting a its vertices, with internal convex angles (more precisely, all angles are $2\pi/3$).
    Hence, we can apply \cref{l:local-to-global-convexity} and conclude.

    We show how a similar argument works also for (ii).
    Let $P$ be the pin.
    If $P$ has degree one, we continue, past $P$, along the equator containing the only edge ending at $F$, ``traveling in the interior'' of $F$ until we hit $\partial F$ for the first time in $Q$. If $F$ has degree $2$, we travel along the bisector of the internal angle of $P$ until we hit $\partial F$ for the first time in $Q$. In both cases, let $\gamma$ be the line connecting $P$ and $Q$.
    Observe that $\gamma$ is contained in $F$ and it splits $F$ in two connected components which satisfy the assumptions necessary to apply the argument we used in (i). The conclusion follows.
\end{proof}

\begin{proof}[Proof of \cref{p:starting}]
(i) is straightforward: if both pins have degree $3$, then $\mathcal{N}$ is a Heppes net.

\medskip
To show (ii) we provide an argument by contradiction.
If $\mathcal{N}$ had a connected component without pins, the connected component would be a Heppes net and the other connected components would be contained in the faces of the Heppes net, which are convex. By \cref{lem:stationary-in-hemisphere}, we conclude that this cannot be true.

Therefore, $\mathcal{N}$ can consist of at most two connected components, each with at least one of the pins. For the sake of contradiction, assume that $\mathcal{N} = \mathcal{N}_1 \sqcup \mathcal{N}_2$, where both $\mathcal{N}_i$ are connected and both contain one pin. Say that $\mathcal{N}_i$ contains the pin $P_i$.
In particular, $N_2$ belongs to a face of $F$ of $N_1$. If $P_1\not\in \partial F$, then by \cref{lem:face-convexity} we know that $F$ is convex and by \cref{lem:stationary-in-hemisphere} we have the desired contradiction.
Otherwise, $P_1\in \partial F$ and $P_2\in \mathring F$, which leads to a contradiction again thanks to \cref{lem:face-convexity}.

\medskip

In order to prove (iii) and (iv) we proceed as follows.
Assume that the two pins $P_1$ and $P_2$ belong to the boundary of the same connected component of $\mathbb S^2\setminus \mathcal N$.
Thus, there is a path $\gamma \subset \mathcal{N}$ which has nonnegative geodesic curvature and connects $P_1$ and $P_2$. \cref{l:tecnico-2} implies that this path is half of an equator.

Now, rotate this half equator in the counterclockwise direction keeping its two endpoints fixed, until we meet another point of $\mathcal{N}$. Let $\gamma'$ be the rotated meridian. $\gamma'\cap \mathcal{N}$ cannot contain any vertex except for $P_1$ and $P_2$ (because any other vertex will have degree $3$).
If $\gamma=\gamma'$, then we have conclude that the net consists of just one arc connecting $P_1$ and $P_2$.
Otherwise, we repeat the procedure on $\gamma'$ finding $\gamma''$. If $\gamma''=\gamma$, then the net is made of just two arcs connecting the two pins. Otherwise, we have found three arcs connecting the two pins and thus the net is a triple junction.
\end{proof}

\subsection{Once-derived net} Given a stationary net with two pins, any pair of edges which meet at a non-pinned vertex form an angle of $\frac{2\pi}{3}$: only those pairs meeting at a pin might form a different angle, in case the degree of that specific pin is different from three. If we remove the pins with degree smaller than three and the edges which contain them as vertices, we get, instead, a geodesic net with the property that all pairs of edges meeting at some vertex form an angle of $\frac{2\pi}{3}$: this ``reduced net'' is an important stepping stone for a better understanding of the original one. Note however that the resulting net might have $0$, $1$, or $2$ nonconvex faces: in order to have always exactly two nonconvex faces, we ``artificially'' remove also one edge from each pin of degree $3$.

\begin{definition}[Once-derived net]\label{d:once-derived}
    Take a stationary net with two antipodal pins.
    We obtain a \emph{once-derived net} by performing the following operations:
    \begin{itemize}
        \item from the set of vertices, we remove the pins with degree $1$ or $2$;
        \item from the set of edges we remove the edges with at least one endpoint which is a pin removed in the previous step;
        \item from the set of edges we remove one edge connecting two pins of degree $3$, if such an edge exists, otherwise we remove an edge from every pin of degree $3$.
    \end{itemize}
Observe that the once-derived net $\mathcal{N}'$ could very well be the empty set: in the latter case we call it ``extinct''. We will later include more cases in the term {\em extinct}, cf. \cref{r:extinction} and \cref{d:extinction}.
\end{definition}


\begin{proposition}[Properties of once-derived nets]\label{p:derived-1}
    Consider a non-extinct once-derived net $\mathcal{N}'$ derived from a stationary net $\mathcal{N}$ with two antipodal pins. Then $\mathcal{N}$ is connected and has the following properties:
    \begin{itemize}
        \item[(i)] All vertices have degree equal to $2$ or $3$.
        \item[(ii)] All its faces are polygons homeomorphic to open disks.
        \item[(iii)] Exactly two faces are nonconvex.
        \item[(iv)] Every convex angle of any face is $\frac{2\pi}{3}$.
        \item[(v)] Every nonconvex face has either one or two concave angles equal to $\frac{4\pi}{3}$.
        \item[(vi)] The area of a face with $q$ concave angles is\footnote{Note that, when counting the number of sides in a face $F$, every side which is not in common with another face is in fact counted twice for $F$.}
        \begin{equation}\label{e:area-face-net-1}
            \frac\pi3\Big(6-\#\text{sides} + 2q\Big).
        \end{equation}
        So, if $q = 0$ then $\#\text{sides}\le 5$, if $q=1$ then $\#\text{sides}\le 7$, and if $q=2$ then $\#\text{sides}\le 9$.
    \end{itemize}
    The convex faces of $\mathcal N'$ are one of the following types:
    \begin{itemize}
        \item[(vii)] A face consisting of two half equators connecting two antipodal points;
        \item[(viii)] An equilateral triangle;
        \item[(ix)] A parallelogram, determined by the length of one side;
        \item[(x)] A pentagon, determined by the lengths of two sides.
    \end{itemize}
\end{proposition}

\begin{proof}
(i) follows immediately from \cref{p:starting}(iv). Thanks to \cref{p:starting}(ii) we know that $\mathcal N$ is connected and thanks to \cref{p:starting}(iv) we can deduce that also $\mathcal N'$ is connected. Then, (ii) follows.

In order to show (iii), we prove that each pin is associated to exactly one nonconvex face in the once-derived net and that any nonconvex face is associated to exactly one pin.
Observe that if a face of $\mathcal N$ does not have any pin on its boundary, then it is convex and remains the same in $\mathcal N'$.
Furthermore, thanks to \cref{p:starting}(iv) we know that any face of $\mathcal N$ can contain at most one pin in its boundary. This allows us to treat the two pins independently.
So, fix a pin $P$. We shall prove that there is exactly one face of $\mathcal N'$ which contains $P$ inside or on its boundary and is nonconvex. This is sufficient to conclude. Let us consider the various cases depending on the degree of $P$ in $\mathcal N$:
\begin{enumerate}[label=(\alph*)]
    \item If $P$ has degree $1$, then the face of $\mathcal N$ which has $P$ on its boundary is nonconvex and becomes a nonconvex face of $\mathcal N'$.
    \item If $P$ has degree $2$, then the two faces of $\mathcal N$ which have $P$ on their boundary are joined to become a nonconvex face of $\mathcal N'$.
    \item If $P$ has degree $3$, then two of the faces of $\mathcal N$ that have $P$ on their boundary are joined to become a nonconvex face of $\mathcal N'$.
\end{enumerate}
Note that the above argument implies also (v). The removal of a pin of degree $1$ creates a nonconvex face with one nonconvex angle. The removal of a pin of degree $2$ eliminates a convex face and creates a nonconvex face with two nonconvex angles of $\frac{4\pi}{3}$ merging a nonconvex face with the removed convex face. The removal of an edge from a degree $3$ pin creates a nonconvex face with two concave angles of $\frac{4\pi}{3}$.

\medskip

(iv) is clear. (vi) follows from Gauss-Bonnet Formula. The points from (vii) to (x) are already discussed in Heppes' paper~\cite{Heppes}.
\end{proof}

\subsection{Graphs and embeddings}\label{s:graphs-empeddings}
In order to handle the topological properties of once-derived net we will borrow basic concepts from the theory of embeddings of graphs, and some classical facts in topological graph theory will be used to devise our algorithm; as a reference we will use \cite{LZ}. First of all we recall the classical definition of graphs, cf. \cite[Definition 1.3.1]{LZ}.

\begin{definition}[Graph]\label{d:graph-0}
    Let $V$ be a finite set of vertices. An edge is an object that has two (unordered) endpoints in $V$ (which may coincide).
    Two edges with the same endpoints are not assumed to be equal. We say that an edge is \emph{adjacent} to its two endpoints.
A graph is a pair $(V, E)$, where $V$ is a finite set of \emph{vertices} and $E$ is a finite collection of edges.
\end{definition}

We are now ready to introduce the concept of embedding of a graph in a topological space. In order to do so, given a topological space $X$ we let $\Gamma(X)$ be the set of continuous embedded arcs. More precisely:
\begin{itemize}
    \item we first consider the set $Y$ of continuous maps of $[0,1]$ into $X$ which restrict to an embedding of $(0,1)$,
    \item we introduce the equivalence relation for pairs of elements of $Y$ in which one element can be obtained from the other through composition with an increasing continuous map of $[0,1]$ onto itself,
    \item and we define $\Gamma (X)$ as the quotient of $Y$ modulo that relation.
\end{itemize}
Moreover, given an arc in $\Gamma (X)$ and a continuous map $\gamma: [0,1] \to X$ representing its equivalence class, we call ``interior of the arc'' the set $\gamma ((0,1))$. The following definition of graph embedding is equivalent to \cite[Definition 1.3.7]{LZ}.

\begin{definition}[Graph Embedding]\label{d:graph}
    For a graph $G=(V, E)$, an embedding of $G$ in a topological space $X$ is a map $\iota:V \to X$ and $\iota : E \to \Gamma(X)$ such that:
    \begin{itemize}
        \item Distinct edges are mapped into arcs with disjoint interiors;
        \item Distinct vertices are mapped into distinct points;
        \item If $A,B\in V$ are the two vertices adjacent to the edge $e$, then either $\iota (A) = \iota (e) (0)$ and $\iota (B) = \iota (e) (1)$, or viceversa;
        \item If $A$ is any vertex, $e$ any edge, and $t$ any real in $(0,1)$, then $\iota (e) (t) \neq \iota (A)$.
    \end{itemize}
\end{definition}

It is in fact convenient to describe a graph embedding after ``directing the edges''. More precisely,
given a graph $G=(V, E)$ we can associate to it a set $\darts(G)$ which contains two copies of each edge in $E$,
and admits a canonical involution $\rev:\darts(G)\to\darts(G)$ mapping one copy of the edge in the other. We understand the darts as ``orientations'' for the edge: in particular if the two vertices $A$ and $B$ adjacent to an edge are distinct we have a ``dart going from $A$ to $B$'' and a ``dart going from $B$ to $A$''. If $A=B$, i.e. the edge is actually a loop, then we simply have two copies of it. The involution reverses the orientation of the dart. For a description we refer to \cite[Section 1.3.3]{LZ}: note that the map $\rev$ is denoted by $\alpha$ in there.

Given an embedded arc $\gamma$, its ``reverse'' is given by the class of embeddings obtained from composing $\gamma$ with a {\em decreasing} homeomorphism of $[0,1]$ onto itself. We will moreover abuse our notation and denote it $\rev (\gamma)$. The embedding $\iota$ in the above definition can therefore be extended to a map $\iota: \darts (G)\to \Gamma (X)$, with the property that $\rev (\iota (e)) = \iota (\rev (e))$, and from now on we will refer to the latter.

Finally, the union of the images of the curves $\iota (d)$ with $d$ varying over $\darts (G)$, which is a closed subset of $X$, will be called the image of $\iota$ and denoted by $\iota (\darts (G))$. We next define when two embeddings are ``equivalent'', cf. \cite[Definition 1.3.7]{LZ}.

\begin{definition}[Isomorphism between Graph Embeddings]\label{d:iso-graph}
    Let $\iota:G\to X$ and $\iota':G'\to X'$ be two embeddings from two graphs $G,G'$ onto two oriented surfaces $X, X'$.
    An isomorphism between these two embeddings is a pair $(f, F)$, where $f:\darts(G)\to \darts(G')$ is a graph isomorphism commuting with $\rev$, and $F:X\to X'$ is an orientation preserving homeomorphism such that $F(\iota(d))=\iota'(f(d))$ for all $d\in\darts(G)$, cf. \cref{f:comm-diag}.
\end{definition}

\begin{figure}
\begin{tikzpicture}[node distance=3cm, auto]
    \node (TL) {$\darts(G)$};              
    \node (TR) [right of=TL] {$\darts(G')$}; 
    \node (BL) [below of=TL] {$X$};          
    \node (BR) [below of=TR] {$X'$};         

    \draw[->] (TL) -- node[above] {$f$} (TR);

    \draw[->] (TL) -- node[left] {$\iota$} (BL);

    \draw[->] (TR) -- node[right] {$\iota'$} (BR);

    \draw[->] (BL) -- node[above] {$F$} (BR);
\end{tikzpicture}
\caption{The commutative diagram for isomorphic embeddings.}\label{f:comm-diag}
\end{figure}

Clearly, the cases which are relevant to us are when $X$ is either the standard sphere $\mathbb S^2$ or the Euclidean plane $\mathbb R^2$, which are both orientable. In what follows we will often use the terms ``spherical graph'' and ``planar graph'' for, respectively, the embedding of a graph in $\mathbb S^2$ and in $\mathbb R^2$.

\begin{remark}\label{r:mirror-symmetry}
We observe that the isomorphisms given in \cref{d:iso-graph} do not cover all the ones which are pertinent to our problem: for us, two stationary nets with antipodal pins are equivalent also when one is obtained from the other with a mirror symmetry of the two-dimensional sphere. However the latter map is not an orientation-preserving homeomorphism of $\mathbb S^2$. In order to avoid ``duplicating'' our candidate topologies for once-derived nets, this will be taken into account at a certain stage, cf. \cref{r:mirror-symmetry-3}.
\end{remark}

\subsection{Twice-derived spherical graphs}  We are ready to introduce the concept of ``twice-derived spherical graph''. The idea is that it is easier to understand the embedded graph coming from a once-derived net if we get rid of the vertices of degree $2$.

\begin{definition}[Twice-derived spherical graphs]\label{d:derived-two}
    Consider the spherical graph associated to a once-derived net and assume it has at least one vertex of degree $3$. We obtain from it a \emph{twice-derived spherical graph} by performing the following operation (as many times as possible):
    \begin{itemize}
        \item Remove a vertex of degree $2$ and identify its two adjacent edges in a single edge. The embedding of the resulting edge is the concatenation of the embeddings of the two removed edges.
    \end{itemize}
\end{definition}

\begin{remark}\label{r:extinction}
Consider the spherical graph associated to a once-derived net which is not empty and has only vertices of degree $2$. If we attempt the same procedure as in \cref{d:derived-two}, namely we repeat the operation of removing inductively the vertices of degree $2$ and merging the edges, the final result would have no vertices. In fact, since we know that the starting graph is connected, we would end up with a single loop. Strictly speaking this does not respect \cref{d:graph}: a graph with empty set of vertices is for us necessarily the empty graph. The stipulation, in \cref{d:derived-two}, that the spherical graph associated to the once-derived net has at least one vertex of degree $3$ is just meant to ensure that the twice-derived spherical graph is well-defined. In order to avoid differentiating too many cases, we stipulate that the empty graph and the single loop with no vertices are also twice-derived spherical graphs and we call them both {\em extinct}. Likewise, a once-derived net which leads to an extinct twice-dervied spherical graph, will be called extinct.
\end{remark}

Since the cases of an empty graph and of a graph with only vertices of degree $2$ is excluded by \cref{d:derived-two}, they will be treated separately. In fact such cases will be completely classified in \cref{s:human-1}, cf. \cref{t:elementary}.

\begin{proposition}[Properties of twice-derived spherical graphs]\label{p:derived-2}
    A twice-derived spherical graph, if it is not extinct, has the following properties:
    \begin{enumerate}[label=(\roman*)]
        \item It is connected.
        \item All vertices have degree $3$ and there is at least one vertex.
        \item If $v$, $e$, and $f$ denote the number of vertices, edges, and faces of the once-derived net, then they satisfy the relations $3 v =2 e$ and $v-e+f=2$ and\footnote{Note that, when counting the number of sides in a face $F$, every side which is not in common with another face is in fact counted twice for $F$} $\sum_{faces} (6-\#\text{sides})=12$.
        \item Every face has at most $7$ sides.
        \item There are at most two faces with more than $5$ sides.
        \item It has at most $12$ faces and at most $30$ edges.
    \end{enumerate}
\end{proposition}

\begin{proof}
The connectedness follows from the connectedness of the once-derived net, which is proven in \cref{p:starting}(ii), and the fact that the procedure which generates the twice-derived spherical graph from a once-derived keeps it. (ii) is also straightforward because the procedure eliminates all vertices of degree $2$. Concerning (iii), observe that $v-e+f=2$ is Euler's formula for graphs embedded in the spheres. The relation $3v=2e$ is obtained counting the edges in the following way: at each vertex there are 3 edges, but each vertex connects two edges, so $e=\frac{3v}{2}$. The final identity corresponds to $12=6f - 2e$ which can be immediately derived from the other two.

To see (iv) let $\mathcal{N}_1$ be the once-derived spherical graph and $\mathcal{N}_2$ be the twice-derived spherical graph. There is a natural bijection between the faces of the two spherical graphs. We thus fix $F\in \mathcal{F} (\mathcal{N}_2)$ and the corresponding face $G\in \mathcal{F} (\mathcal{N}_1)$. $F$ is obtained from $G$ by removing the vertices of $G$ which have degree $2$ and merging the corresponding edges, in particular the number of edges of $G$ coincides with that of $F$ minus its vertices of degree $2$. It follows from \cref{p:derived-1} that the sides of $G$ are at most $5+2n (G)$, where $n(G)$ is the number of degree $2$ vertices of $G$. Hence the sides of $F$ are at most most $5+n(G)$ and hence at most $7$ because $n(G)\leq 2$ (cf. \cref{p:derived-1}).

Coming to (v), the considerations above show that there are exactly two faces in $\mathcal{F} (\mathcal{N}_2)$ that originate from a nonconvex face of $\mathcal{F} (\mathcal{N}_1)$. So, except for the ones originating from these, all other faces in $\mathcal{F} (\mathcal{N}_2)$ have no more sides than some convex face in $\mathcal{F} (\mathcal{N}_1)$, namely at most $5$ by \cref{p:derived-1}.

Finally, with respect to (vi), note that the number of faces in $\mathcal{F} (\mathcal{N}_2)$ equals the number of faces in $\mathcal{F} (\mathcal{N}_1)$. By \cref{p:derived-1} the area of a face $F\in \mathcal{F} (\mathcal{N}_1)$ is an integer multiple of $\frac{\pi}{3}$ and hence at least $\frac{\pi}{3}$. Since the sum of the area of these faces is the total area of the sphere, namely $4\pi$, their number does not exceed $12$. Using (iii) we get
\[
\sum_{F\in \mathcal{F} (\mathcal{N}_2)} \#\text{sides} (F) = 6 (\#\text{faces}-2) \leq 60\, .
\]
However, given that
\[
\sum_{F\in \mathcal{F} (\mathcal{N}_2)} \#\text{sides} (F) = 2 e\, ,
\]
we conclude the second inequality in (vi).
\end{proof}

\section{Polygonal paths and polygons}\label{s:spher-poly}

In this section study polygons in the sphere. We start with some definitions, being careful about whether self-intersections are allowed.

Consider a Lipschitz curve in the sphere which consists of
a finite number of geodesic arcs $\gamma_i$ on the sphere $\mathbb S^2$, with $i\in \{1, \ldots , N\}$, joining vertices $V_{i-1}$ and $V_i$. The arcs are supposed to have strictly positive lengths and be shorter than $2\pi$. Given these data we can give to the path a canonical orientation: each geodesic arc is oriented so that $V_{i-1}$ is the starting point and we assume further that $\gamma_i$ is an arc-length parametrization of the geodesic arc with $\gamma_i (0) = V_{i-1}$ and $\gamma_i (\lambda_i) = V_i$. Given the latter convention we define the angle of the path at the vertex $V_i$ (with $i\leq N-1$) as the counterclockwise angle (at $V_i$) from the tangent vector $-\dot\gamma_{i} (\lambda_i)$ to the tangent vector $\dot{\gamma}_{i+1} (0)$. The path will be called {\em closed} if $V_N=V_0$. We will identify two such closed paths up to choice of the `` starting point'' $V_0=V_N$. 
Two nonclosed paths are congruent if there is an {\em orientation preserving} isometry mapping one (ordered) collection of geodesic segments into the other, while for closed paths the congruence is understood in the same way up to choosing the starting vertex.

\begin{definition}\label{d:polygonal-path}
A \emph{(spherical) polygonal path} is a curve $\gamma$ as described above.
The number of sides is the number of geodesic arcs in the description.
\end{definition}

Note that we are not requiring that the path does not self-intersect, nor we are requiring that the lengths of the $\gamma_i$'s are smaller than $\pi$ (i.e. that the $\gamma_i$ are the length-minimizing geodesic connecting $V_{i-1}$ with $V_i$).

\begin{definition}\label{d:polygon}
A (spherical) polygonal path will be called a \emph{(spherical) polygon} if the path is closed and it is an embedding of $\mathbb S^1$ into $\mathbb S^2$. In these cases the polygon will be called (sphere) triangle, quadrilateral, pentagon, more in general $N$-gon depending of the number of its sides. The {\em interior} of the polygon is the topological disk bounded by the polygon with the property that the exterior normal to it forms a counterclockwise angle of $\frac{\pi}{2}$ with the tangent to the path. The polygon will be called convex if the interior is a convex set of the sphere.
\end{definition}

\begin{remark}
It follows from \cref{l:local-to-global-convexity} that a spherical polygon is convex if and only if all its internal angles are smaller than or equal to $\pi$.
\end{remark}

\subsection{Triangles}
In this section we focus on spherical triangles. When listing the angles of the triangle we will do so by ordering them in counterclockwise fashion as $A,B$, and $C$. The length of the segments will then be denoted by $a$, $b$, and $c$, with the convention that the angle $A$ is opposite to the side $a$, and so on. This is the standard notation in textbooks on spherical trigonometry, see for instance \cite{Todhunter}.

\begin{theorem}\label{thm:convex_triangles}
    Fix $A,B,C\in (0,\pi)$. There is a convex spherical triangle with angles $A, B, C$ if and only if:
    \begin{itemize}
        \item $\pi<A+B+C<3\pi$,
        \item $\pi + 2\min\{A, B, C\}>A+B+C$.
    \end{itemize}
    Under these assumptions, there exists a unique (up to congruence) convex spherical triangle with angles $A, B, C$.
    The lengths of its sides $a, b, c$ can be found using the (polar of the spherical) cosine rule:
    \begin{equation*}
        \cos(a) = \frac{\cos(A) + \cos(B)\cos(C)}{\sin(B)\sin(C)} \,.
    \end{equation*}
    Viceversa the sides determine the angles via the formula:
    \begin{equation*}
      \cos(A)=\frac{\cos(a)-\cos(b)\cos(c)}{\sin(b)\sin(c)}\, .
    \end{equation*}
    Finally the triangle is also determined by the measure of one side and its two adjacent angles and by the measure of an angle and its adjacent sides.
\end{theorem}

The theorem above is a classical fact in spherical trigonometry known since the nineteenth century (see \cite{Todhunter}), but the reader can also derive it using linear algebra and elementary (planary) trigonometry. It will be useful for us to get a generalization of the theorem above to polygonal paths with three sides. We observe that, according to our definition, the three sides of the polygonal path must then belong to three {\em distinct} equators. We follow the same convention as above for labeling the measures of the angles of the polygonal path and the lengths of its sides

\begin{remark}
Note that when passing from convex triangles to polygonal paths with three sides, or even just to non-convex triangles, the conclusions of \cref{thm:convex_triangles} must be handled with some care. For instance, if we know one angle $A$ and the two measures $b$ and $c$ of the adjacent sides, we have determined, up to congruences, the other two vertices of the triangle. In particular one of the two sides must be different than $\pi$, otherwise the two vertices would coincide. We have thus determined the equator in which the third side must lie. However there are exactly two possibilities for this arc, and hence the knowledge of the triple of numbers $A$, $b$, and $c$ does not determine the polygonal path, but it rather gives exactly two possibilities.
\end{remark}

\begin{theorem}\label{thm:generalized_triangles}
    Fix $A,B,C\in(0,2\pi)\setminus\{\pi\}$.
    If there is an odd number of values in $(A, B, C)$ smaller than $\pi$, let $(A',B',C'):=(A,B,C)$ and otherwise let $(A',B',C'):=(2\pi-A,2\pi-B,2\pi-C)$. Let $A'':=(A'\mod\pi)$ and cyclic.

    There is a polygonal path with three angles $(A, B, C)$ if and only if there is a convex triangle with angles $(A'',B'',C'')$ (see \cref{thm:convex_triangles}).

    Under this assumption, there are exactly two (up to congruence) polygonal paths with $(A,B,C)$ as angles.
    The lengths of their sides can be obtained applying the following set of rules.
    \begin{itemize}
        \item If $(a, b, c)$ is a valid triple of side lengths, then the other triple is $(2\pi-a,2\pi-b,2\pi-c)$.
        \item The valid triples of side lengths for the angles $(A,B,C)$ coincide with the valid triples of side lengths for the angles $(2\pi-A,2\pi-B,2\pi-C)$.
        \item $(a, b, c)$ is a valid triple of side lengths for the angles $(A, B, C)$ (with $B, C<\pi$) if and only if $(2\pi-a, b, c)$ is a valid triple of side lengths for the angles $(A,\pi+B,\pi+C)$.
        \item If $0<A,B,C<\pi$, then a valid triple of side lengths can be obtained as in \cref{thm:convex_triangles}.
    \end{itemize}
\end{theorem}

The above theorem is an exercise which we leave to the reader.

\subsection{Congruence criterion for closed polygonal paths} The following fact is a corollary of the above considerations for spherical triangles and will be particularly useful to analyze the special nets of \cref{s:human-1} and \cref{s:human-2}.

\begin{theorem}\label{thm:congruence}
Assume $V_1V_2\ldots V_N$ and $W_1W_2\ldots W_N$ are the vertices of two closed polygonal paths such that:
\begin{itemize}
\item[(i)] The angles $\hat{V}_i$ are all congruent to the angles $\hat{W}_i$, and all ddiferent than $\pi$;
\item[(ii)] The geodesic arcs $V_iV_{i+1}$ are congruent to the edges $W_iW_{i+1}$ for all $i\in \{1, \ldots, N-3\}$;
\item[(iii)] The last three edges of both paths are all shorter than $\pi$.
\end{itemize}
Then the two polygonal paths are congruent.
\end{theorem}
The proof is given in \cref{s:congruence}.

\subsection{Convex spherical polygons with angles equal $\frac{2\pi}{3}$} The explicit formulas in \cref{p:formulae} below will be of fundamental importance throughout the rest of the paper, as they give a simple way to relate the lenghts of the sides in the four type of poligons which have angles all equal to $\frac{2\pi}{3}$.

\begin{proposition}\label{p:formulae}
A convex spherical triangle with angles all equal to $\frac{2\pi}{3}$ is equilateral and its side measures
\begin{equation}\label{e:triangle_side}
\lambda_T = \arccos (-\tfrac{1}{3})\, .
\end{equation}
In a convex spherical quadrilateral with angles equal to $\frac{2\pi}{3}$ opposite sides are congruent. Moreover, the measures $a$ and $b$ of two adjacent sides satisfy the relation
\begin{equation}\label{e:quadrilateral-formula}
\tan \tfrac{a}{2} \tan \tfrac{b}{2} = \tfrac{1}{2}\, .
\end{equation}
If two adjacent sides are equal (which means that all sides are congruent) then the diagonal measures $\lambda_T$.

A convex spherical pentagon $ABCDE$ with angles equal to $\frac{2\pi}{3}$ is determined by the length of any two of its sides and every side is strictly shorter than $\lambda_T$.
Moreover, if the vertices are labeled consecutively, then
\begin{equation}\label{e:formula-Taylor}
2\cos |DE| - \tfrac13 = f(|AB|, |BC|) \,.
\end{equation}
where
\begin{equation}\label{e:fpentagono}
f(x, y) := \cos x+\cos y +\cos x\cos y -2\sin x\sin y \, .
\end{equation}
\end{proposition}

The sidelength of an equilateral triangle with angles $\frac{2\pi}{3}$ is well known. \eqref{e:formula-Taylor} is instead taken from \cite{Taylor} (where, however, it is not proved) and the other statements appear in \cite{Heppes}, except for \eqref{e:quadrilateral-formula}, which does not seem to appear anywhere. Anyway, we provide proofs for all the statements in \cref{s:app-convex-faces}, based on the trigonometric identities for triangles given in \cref{thm:convex_triangles}.

\section{Sieves}\label{s:sieves}

Eventually, in the following section, we plan to describe an algorithm which will produce a finite list of spherical graphs\footnote{Actually, we will keep track of a bit more information than just a spherical graph, see \cref{d:candidate}.} such that any once-derived net, when interpreted as a spherical graph, belongs to the list.

In order to shorten our list, we show here a number of conditions that any once-derived net must satisfy in addition to those already found in \cref{p:derived-1}.
We will first state all the conditions: each condition will be stated separately in a lemma and will be assigned a name; the proofs will then follow in separate subsections, each bearing the name of the lemma.

We start by introducing the following concepts.
Recalling \cref{p:derived-1}, the area of a face depends only on the number of concave angles and the number of edges of the face.

\begin{definition}[Boundary distance]\label{d:b-distance}
Given two vertices $A$ and $B$ in the boundary of a face $F$ we introduce the following.
\begin{itemize}
    \item A {\em boundary path} from $A$ to $B$ is an ordered list of $N+1$ vertices $A= V_0 , V_1, \ldots , V_N = B$ with the property that, for each $i\in \{0, \ldots , V_{N-1}\}$, $V_i$ and $V_{i+1}$ are connected by an edge $\sigma_i$ of $F$ and the union of the edges $\sigma_i$ form a positively curved oriented arc.
    \item The {\em boundary distance} of $A$ and $B$ in the face $F$ is the minimal $N$ for which there is boundary path as in (a) either from $A$ to $B$ or from $B$ to $A$.
    \item The boundary distance of $A$ and $B$ {\em avoiding a specific set $\mathcal{E}$} of edges is the minimal $N$ for which there is a boundary path with the additional property that none of the edges $\sigma_i$ belongs to $\mathcal{E}$.
\end{itemize}
\end{definition}

We are now ready to state all the sieves. In the next statements, $\mathcal N$ denotes a non-extinct once-derived net.

\begin{lemma}[No convex triangles next to convex pentagons]\label{l:no-T+P}
A convex triangle of $\mathcal{N}$ cannot share an edge with a convex pentagon.
\end{lemma}

\begin{lemma}[No concave biangles]\label{l:no-bi-concavi}
$\mathcal{N}$ cannot contain a concave biangle.
\end{lemma}

\begin{lemma}[Simple boundary]\label{l:S1}
The boundary of any face of $\mathcal{N}$ is topologically an $\mathbb S^1$.
\end{lemma}

The proof of \cref{l:S1} also implies easily the following fact, which is not used as a sieve but will be useful in other contexts.

\begin{corollary}\label{c:convex-faces-intersection}
Two convex faces of the net $\mathcal{N}$ have at most one edge in common.
\end{corollary}

\begin{lemma}[Biangle]\label{l:biangles}
Let $\mathcal{N}$ be such that one face is a biangle $B$. Then
\begin{itemize}
\item[(i)] $B$ is the only biangular face.
\item[(ii)] The two faces adjacent to $B$ are distinct and nonconvex.
\item[(iii)] The convex faces (except for the biangle) are all \emph{strictly} contained in a single hemisphere.
\item[(iv)] Any convex triangle in $\mathcal{N}$ cannot share a side with another convex face.
\end{itemize}
Moreover, the sum of the areas of the convex faces (excluding the biangle), plus the \emph{contribution} of the nonconvex faces is at most $\frac{5\pi}{3}$, where the contribution of a nonconvex face is defined as
\begin{itemize}
\item[(a)] $0$ if the face has a single concave angle;
\item[(b)] $\frac{\pi}{3}\max\{0, 4-d\}$, where $d$ is the boundary distance of the two concave angles avoiding the edges of $B$, if the face has two concave angles.
\end{itemize}
\end{lemma}

\begin{lemma}[Long edge]\label{l:long-edge}
    Assume that there is an edge $\sigma$ of $\mathcal{N}$ which is strictly longer than $\pi$ or with length $\pi$ which is not the edge of a biangle. Then all the conclusions of \cref{l:biangles} hold, specifically:
    \begin{itemize}
        \item[(i)] The two faces adjacent to $\sigma$ are distinct and nonconvex;
        \item[(ii)] The convex faces are all \emph{strictly} contained in a single hemisphere;
        \item[(iii)] The sum of the areas of the convex faces plus the contribution from the nonconvex faces is at most $\frac{5\pi}{3}$, where the contribution of a nonconvex face is defined as
        \begin{itemize}
        \item[(a)] $0$ if the face has a single concave angle;
        \item[(b)] $\frac{\pi}{3}\max\{0, 4-d\}$, where $d$ is the boundary distance of the two concave angles avoiding $\sigma$, if the face has two concave angles.
        \end{itemize}
    \end{itemize}
\end{lemma}

\begin{lemma}[Concave angles are far]\label{l:concave}
    The boundary distance between the two concave angles of a face of $\mathcal{N}$ is at least $s-5$, where $s$ is the number of sides of the face.
\end{lemma}

\begin{proof}[Proof of \cref{l:no-T+P}] This is proved by Heppes in \cite{Heppes}, at any rate it follows immediately from \cref{p:formulae}.
\end{proof}

\begin{proof}[Proof of \cref{l:no-bi-concavi}] If $\mathcal{N}$ contains a concave biangle, then the two vertices of the biangle have degree $2$ and thus the two edges adjacent to them are precisely the edges of the biangle. In particular, the boundary of the biangle is a connected component of $\mathcal{N}$ and, because $\mathcal{N}$ is connected, it is the whole net. But then $\mathcal{N}$ would be an extinct net.
\end{proof}

\begin{proof}[Proof of \cref{c:convex-faces-intersection}] The statement follows immediately from \cref{p:elementary-convexity}(v)
\end{proof}

\begin{proof}[Proof of \cref{l:S1}]
Let $F$ be the face we are considering.
If $F$ is a convex face there is nothing to prove.

Consider the case of a nonconvex face with a single concave angle $A$.
Start from $A$ and follow the bisector of $A$ inside the face $F$ until we hit the first point $B$ of the boundary. The resulting geodesic arc $\eta$ subdivides the face $F$ into parts $F_1$ and $F_2$. Using \cref{l:local-to-global-convexity} we conclude that both are convex. If we consider the equator $\gamma$ containing $\eta$, we conclude that, since in a neighborhood of $A$ there are points of $\gamma$ which do not belong do $F_1$ nor $F_2$, $F_1\cap F_2$ contains $\eta$, but is not $\gamma$. It must therefore be a geodesic arc which contains $\eta$. But note that there are points in $\gamma$ arbitrarily close to $B$ which do not belong to $F_1$ nor to $F_2$. So $F_1\cap F_2=\eta$. In particular, $\partial F_1 \setminus \eta$ and $\partial F_2 \setminus \eta$ are two disjoint embedded open arcs with endpoints $A$ and $B$. Since $\partial F$ is the union of the latter arcs and of the points $A$ and $B$, we conclude that $\partial F$ is topologically $\mathbb S^1$.

Consider now the case of a nonconvex face $F$ with two concave angles $A$ and $B$. In the stationary net that generated $\mathcal{N}$ there is a pin $P$ which is connected to $A$ and $B$ with two edges $s$ and $t$, contained in the face $F$. If $s$ and $t$ form an angle $\pi$ at $P$ then we argue as above decomposing $F$ into two convex faces $F_1$ and $F_2$ whose intersection consists only of $s\cup t$.
Otherwise we take the bisector of the concave angle formed at $P$ by the edges connecting $P$ at $A$ and $B$ and follow it until it hits the first point $C$ in $\partial F$. The three geodesic segments $PA$, $PB$, and $PC$ decompose $F$ into the closures of three convex open sets $F_1$, $F_2$, and $F_3$ and, again, the intersection of any pair of them is the shortest geodesic connecting $P$ with one of the three points $A$, $B$, and $C$. We thus conclude that $\partial F$ can be decomposed as the union of three embedded open arcs with empty intersection plus their endpoints $A$, $B$, and $C$, forming a topological $\mathbb S^1$.
\end{proof}

\begin{proof}[Proof of \cref{l:biangles}] First of all denote by $V$ and $U$ the two vertices of the biangle and observe that it must be a convex biangle by \cref{l:no-bi-concavi}. If one vertex is a pin the other is a pin too and the biangle would be a connected component of $\mathcal{N}$, hence the whole net, which in turn would be elementary. So both $V$ and $U$ have degree $3$. In particular the biangle is also a face of the net $\mathcal{O}$ from which $\mathcal{N}$ is derived. Now rotate the two half meridians forming the biangle keeping $V$ and $U$ fixed, moving them into the complement of the biangle, until they intersect the first points of $\mathcal{O}\setminus \partial B$. If the two meridians meet then the net would be a triple junction, with the two pins lying on the edges, which is impossible. So each rotated meridian must meet a distinct vertex of $\mathcal{O}$ and these two vertices must both be pins. In particular, once they stop, the two rotated meridians form a single equator. The latter subdivides the sphere into two hemispheres. One of them contains the biangle $B$ and the other one contains all the other convex faces of the net $\mathcal{N}$. So the two faces adjacent to $B$ must both be nonconvex.

Next, consider the hemisphere $H$ which contains all the convex faces of $\mathcal{N}$ and fix a pin $P$, which is contained in the boundary of the hemisphere. If $P$ has degree $2$ then both edges adjacent to the pin $P$ either touch $\partial H$ coming from $H$, or they lie in $\partial H$. In particular this implies that $P$ is the vertex of a convex face of $\mathcal{O}$ which was removed when deriving the net $\mathcal{N}$. This face was contained in the hemisphere $H$. We thus conclude that the sum of the areas of the convex faces of $\mathcal{N}$ plus the eventual convex faces of $\mathcal{O}$ which were removed (one for each pin of degree $2$) can be at most $2\pi$.

Consider any of the convex faces $F$ removed from $\mathcal{O}$, with $P$ the corresponding pin of degree $2$. The internal angle of $F$ at $P$ at is at least $\frac{2\pi}{3}$. If $s(F)$ is the number of sides of $F$, then by Gauss-Bonnet the area of $F$ is at least $(6-s(F))\frac{\pi}{3}$. Observe that $F$ is contained in one concave face $G$ of $\mathcal{N}$. The pin $P$ has degree $2$ if and only if this face $G$ has two concave angles. In fact if the pin had degree $1$, then $G$ would have one concave angle. If the pin had degree $3$ then $G$ would need to be the union of two convex polygons. Note however that one would share an edge with a biangle, so would need to be a biangle, implying that the other one would also need to be a biangle. But then the entire net $\mathcal{O}$ would consist of the two pins connected by three meridians and $\mathcal{N}$ would be empty.

Having established that each ``removed'' face $F$ accounted corresponds to a nonconvex faces $G$ of $\mathcal{N}$ with two concave angles containing it, we now show that $(6-s(F))\frac{\pi}{3}$ coincides with our definition of contribution of the face $G$. In fact $F$ has $s(F)-2$ sides in common with $G$. These sides form a boundary path connecting the two concave angles of $G$ and avoiding the side $\sigma$ of $G$ in common with the biangle: $s(F)-2$ is the boundary distance of the two concave angles avoiding $\sigma$. Hence the area of $F$ is $(4-(s(F)-2)) \frac{\pi}{3}$, precisely the ``contribution of the nonconvex face $G$''. Since the contribution of $G$ has been defined as $0$ when there is a concave angle, the argument given so far shows that the sum of the areas of the convex faces of $\mathcal{N}$ plus the contributions of the nonconvex faces is at most the area of the hemisphere $H$, i.e. $2\pi$. Note however that:
\begin{itemize}
\item If one of the two pins has degree $1$, then the union of the faces considered has area strictly smaller than $H$;
\item If both pins have degree $2$, the only way in which the union of the faces considered in the computation could have the same area as $H$ is that the two angles at the pin are exactly $\pi$; but in this case the ``contribution'' of each nonconvex face $G$ is strictly less than the area of the corresponding convex face $F$ which was removed.
\end{itemize}
The sum of the areas of the convex faces plus the contribution of the nonconvex ones must therefore be strictly smaller than $2\pi$. Since however this number is an integer multiple of $\frac{\pi}{3}$, it can be at most $\frac{5\pi}{3}$.

Finally, (iv) follows because a convex triangle cannot share an edge with a convex pentagon and if it shares an edge with a convex quadrilateral, then their union is a set of diameter strictly larger than $\pi$. This verification is left to the reader, after observing that \cref{p:formulae} completely determines the diameter of the union.
\end{proof}

\begin{remark}\label{r:no-pin-3}
We remark here an implicit effect of the proof of \cref{l:biangles}: if the derived-once net is nonextinct and it contains a biangle, then no pin of the net can have degree $3$.
\end{remark}

\begin{proof}[Proof of \cref{l:long-edge}] The proof is achieved using the same idea of the proof of the biangle lemma. If the length of the edge is strictly larger than $\pi$ we choose two interior antipodal points and rotate the corresponding half meridian in both directions keeping its endpoints fixed. If the length of the edge is exactly $\pi$ but it is not an edge of a biangle, then we rotate it in both directions. Like in the proof of \cref{l:biangles} this operation ``sweeps'' an hemisphere in the original net $\mathcal{O}$ from which $\mathcal{N}$ was derived. If $H$ is the complementary hemipshere, then $\partial H$ contains the two pins. The considerations in the proof of \cref{l:biangles} apply now verbatim to this case as well.
\end{proof}

\begin{proof}[Proof of \cref{l:concave}]
If the face $F$ has $s\le 6$ sides, then the statement is trivial.
Thus, we can assume that $F$ is a face with two concave angles and $s\ge 7$ sides.
Its area is $(10-s)\frac{\pi}{3}$. Let $d$ be the boundary distance between the concave angles in $F$. Let $P$ be the pin in the net $\mathcal{O}$ which is contained in the face $F$.

The pin $P$ together with the $d+1$ vertices on the boundary of the face which form the path in $\partial F$ connecting the two concave angles constitute a simple polygon $Q$ strictly contained in $F$. The polygon $Q$ has $d+1$ angles equal to $\frac{2\pi}{3}$, and one (the one at the pin) greater or equal than $\frac{2\pi}{3}$. By Gauss-Bonnet, the area of $Q$ is at least $(6-(d+2))\frac{\pi}{3}$.

Since $Q$ is strictly contained in $F$, the area of $Q$ must be strictly smaller than the area of $F$. Therefore $6-(d+2) < 10-s$ which is equivalent to $s-5\le d$.
\end{proof}

\section{Computer-assisted list of candidate topologies}\label{s:primo-programma}

In this section, we describe the algorithms we use to obtain a finite list of candidates for the twice-derived spherical graphs and for the spherical embedded graphs induced by once-derived nets.
Let us recall in the following definition the special cases which we will neglect in this section.

\begin{definition}\label{d:extinction}
A twice-derived spherical graph is {\em extinct} if it is either empty or a loop with no vertices. A once-derived net is extinct if the corresponding twice-derived spherical graph is extinct.
\end{definition}

The extinct cases include most of the geometries in \cref{t:main}. We will therefore analyze them directly in \cref{s:human-1}.
We can thus focus on the remaining cases, namely the twice-derived spherical graphs in which there is at least one vertex.
The first computer algorithm described in this section, combined with the statements of the previous section, proves:

\begin{proposition}\label{p:first-topological-list}
Any twice-derived spherical graph (as in \cref{d:derived-two}) which is not extinct belongs to the list given in \verb|data/twice_derived_human.txt| in \cite{repo}.
\end{proposition}

Given a once-derived net $\mathcal{N}$ we regard it as a particular embedding in $\mathbb S^2$ of a graph $G$, carrying, however one additional piece of information: for each vertex of degree $2$, we record
which angle is concave.

\begin{definition}[Once-derived marked spherical graph]\label{d:candidate}
 A marked spherical graph is a graph $G$ embedded in $\mathbb S^2$ which has only vertices of degree $3$ and $2$, where for each vertex of degree $2$ we also mark one edge arriving at it.
 Observe that any once-derived net produces a marked spherical graph when we use the concavity of the angle to mark the edge. More precisely, we adopt the convention that the counterclockwise angle from the marked edge to the next one is concave.
 The marked spherical graph produced by a once-derived net is a \emph{once-derived marked spherical graph}.
\end{definition}

\begin{remark}\label{r:mirror-symmetry-3}
Following \cref{d:iso-graph} such two marked spherical graphs are isomorphic if they can be mapped one into the other via an orientation-preserving homeomorphism of $\mathbb S^2$ which also ``respects'' the concave angles, namely if the homeomorphism $\iota$ maps a degree $2$ vertex $A$ of one graph in the vertex $A'$ of the other, then it must map the marked edge arriving at $A$ to the marked edge arriving at $A'$.
Though, for our purposes, we shall also identify mirror symmetric images: in other words we will drop the requirement that $\iota$ is orientation preserving. We warn the readers that the list mentioned in \cref{p:first-topological-list} will instead not distinguish between a twice-derived spherical graph and its mirror image, cf. \cref{r:mirror-symmetry,r:mirror-symmetry-2}.
\end{remark}

The second algorithm described in this section takes as input the list of \cref{p:first-topological-list} and produces a corresponding list of possibilities for once-derived marked spherical graphs, avoiding any ``double counting'' through isomorphic copies (including mirror symmetries, as detailed in the remark above).

\begin{proposition}\label{p:second-topological-classification}
Any once-derived marked spherical graph belongs to the list given in \verb|once_derived_human.txt| in \cite{repo} (see also \cite{website}), up to isomorphism as described in \cref{r:mirror-symmetry-3}.
\end{proposition}

\subsection{High-level illustration of the algorithms}\label{s:hl-top-algorithm}
Here we give a quick illustration of the main ideas, especially suited for readers who already have a working knowledge of concepts in graph theory.
The actual algorithms will be more intricate in order to be computationally feasible (namely ensuring that the algorithms can be run in a reasonable time on common hardware).

Recall that, as shown in \cref{p:derived-2}, any twice-derived spherical graph is a spherical graph where each vertex has degree $3$ and the total number of vertices is at most $30$. In particular, there are finitely many twice-derived spherical graphs.
We want to generate all of them, starting from a cycle and adding one edge at a time.

We generate the list of possible twice-derived spherical graphs by recursively building families of planar graphs, denoted as $\mathcal{C}_k$, where $k$ is the number of edges. While twice-derived spherical graphs are in fact spherically embedded, we can easily regard them as planar (see \cref{subsec:spherical-planar}).
The algorithm builds $\mathcal{C}_{k+1}$ from $\mathcal{C}_k$ as follows. The first element of $\mathcal{C}_{k+1}$ is the cyclic graph with $k+1$ edges. All the other elements are created as ``nice children'' of elements in $\mathcal{C}_k$: the parent graph in $\mathcal{C}_k$ is modified by either attaching a new vertex externally or connecting two external vertices with an edge. Along this process we force the property that graphs in $\mathcal{C}_k$ have at most one leaf.
One can show that this algorithm generates all planar graphs without leaves.

To keep the search feasible (i.e. to prevent a fast growth in the size of $\mathcal{C}_k$), the algorithm applies many pruning criteria based on topological constraint derived from the geometric information (cf., \cref{lem:heuristic-subgraph}).
Moreover, the algorithm discards isomorphic copies so that $C_k$ does not contain two isomorphic planar graphs. In order to perform this step efficiently some care is required (e.g., we employ a polynomial-time isomorphism check and a sophisticated hashing technique for planar graphs).

Once we have produced a list of spherical graphs that contains all twice-derived spherical graphs, we employ a second (simpler) algorithm to obtain from it a list of marked spherical graphs that contains all once-derived marked spherical graphs.
The algorithm iterates through the twice-derived candidates and systematically splits edges to insert between 2 and 4 vertices of degree two. The reconstruction process also keeps track of where the ``concave angles'' can be placed, cf. \cref{d:candidate}). The resulting objects are checked against the properties of once-derived nets which have been proved in the previous sections (cf. \cref{p:starting} and \cref{s:sieves}). This step also accounts for mirror symmetries.

\subsection{Combinatorial maps}
The material contained in this section is well known.
We start with the key observation that equivalence classes of graph embeddings into surfaces can be described through combinatorial data. We recall the notions of darts and darts graph $\darts(G)$, given in \cref{s:graphs-empeddings}.

\begin{definition}[Ordering induced by an embedding]\label{d:darts}
    Let $G=(V, E)$ be a graph and let $\iota:\darts(G)\to \Gamma(X)$ be an embedding of $G$ into an \emph{oriented $2$-dimensional surface} $X$. The embedding $\iota$ canonically induces a function $\next:\darts(G)\to \darts(G)$ as follows.
    Take a dart $e$ starting from $v\in V$. Look at a neighborhood of $\iota(v)$ in $X$ and assume that, if we start from a point on $\iota(e)$ close to $v$ and move counterclockwise around $v$ (which makes sense because $X$ is an oriented surface), the first edge we encounter is $\iota(e')$. Then, we set $\next(e)=e'$.
\end{definition}

The reader can compare with \cite[Section 1.3.3]{LZ}, where the map $\next$ is denoted by $\sigma$. The following lemma summarizes \cite[Lemma 1.3.15 \& Proposition 1.3.16]{LZ} (note that the notation $\varphi$ for the map $\next\circ \rev$ is not completely consistent with the one in \cite{LZ} as in there the same symbol denotes the inverse of our $\varphi$; this change is however immaterial for the conclusion).

\begin{lemma}
    Let $G=(V, E)$ be a graph, and let $\iota:\darts(G)\to \Gamma(X)$ be an embedding of $G$ into an \emph{oriented surface} $X$. Let $\next:\darts(G)\to \darts(G)$ be the associated function. We have:
    \begin{itemize}
        \item For any $e\in \darts(G)$, the starting points of $e$ and $\next(e)$ coincide.
        \item If $e, e'\in\darts(G)$ have the same starting point, there is a $k\in \N$ such that $\next^k(e)=e'$.
        \item Consider the map $\varphi=\next\circ\rev:\darts(G)\to \darts(G)$. The image, through $\iota$, of an orbit of $\varphi$ is (one connected component of) the boundary of a connected component of $X\setminus \iota(\darts(G))$, traversed in the clockwise direction. The orbits of $\varphi$ (as well as their images in $X$ through $\iota$) are called \emph{faces}. A vertex belongs to a face if it is a vertex of a dart forming the face.
    \end{itemize}
\end{lemma}

Since we will only consider embedded graphs and their embeddings $\iota$ into the simply connected spaces $X=\mathbb R^2$ and $X=\mathbb S^2$, to each ``face'' as defined above will correspond a unique connected component $\Omega$ of the complement in $X$ of the image of the graph, with the property that $\partial \Omega$ is its topological boundary. These will also be called faces. When $X=\mathbb S^2$, all of them are topological disks. When $X=\mathbb R^2$, all the bounded ones are topological disks, and there is precisely one unbounded one.

We can now introduce the primary combinatorial object which we will use to give our list of possibilities for twice-derived spherical graphs. The reader can compare with \cite[Definition 1.3.23 \& Remark 1.3.24]{LZ}.

\begin{definition}[Combinatorial Map]\label{d:combinatorial}
    A combinatorial map is a graph $G$ endowed with a function $\next:\darts(G)\to\darts(G)$ such that
    \begin{itemize}
        \item For any $e\in\darts(G)$, the starting points of $e$ and $\next(e)$ coincide.
        \item If $e, e'\in\darts(G)$ have the same starting point, there is $k\in\N$ such that $\next^k(e)=e'$.
    \end{itemize}
\end{definition}

For our purposes, the most relevant fact of the theory is the following one; for the proof we refer the reader to \cite[Construction 1.3.20]{LZ}.

\begin{theorem}\label{t:comb-to-emb}
    Any combinatorial map is induced by an embedding of the underlying graph into an oriented surface such that all faces are simply connected.
\end{theorem}

In our case the relevant combinatorial maps are the ones induced by embeddings in the sphere. On the other hand the topology of the oriented surface given by the above theorem is completely determined by the Euler characteristic of the graph. Another important point in our classification is that we wish to avoid repetitions: two combinatorial maps inducing isomorphic embeddings of the same graph should not appear in our list. There is in fact a simple way of recognizing isomorphic embeddings through their combinatorial maps (cf. again \cite[Definition 1.3.23]{LZ}).

\begin{definition}[Isomorphism of Combinatorial Maps]
    Let $(G,\next)$ and $(G',\next)$ be two combinatorial maps.
    An isomorphism between these two combinatorial maps is a bijection $f:\darts(G)\to\darts(G')$ that commutes with $\rev$ and $\next$.
\end{definition}

For the following theorem we refer to \cite[Section 1.3.3]{LZ} (the theorem is not explicitly stated in the reference, but stated informally in \cite[Bottom of p. 37]{LZ}).

\begin{theorem}
    Let $\iota:\darts(G)\to X$ and $\iota':\darts(G')\to X$ be two embeddings of two graphs into two oriented surfaces such that
    \begin{itemize}
        \item $G$ and $G'$ are connected, and they both have at least one edge.
        \item $X$ and $X'$ are connected.
        \item Every face of $\iota$ and $\iota'$ is simply connected.
    \end{itemize}
    Then, the two embeddings $\iota$ and $\iota'$ are isomorphic if and only if the associated combinatorial maps are isomorphic.
\end{theorem}

From the theorem we infer immediately the following corollary.

\begin{corollary}
    Let $(G, \next)$ be a connected combinatorial map with at least one edge.
    There is a unique (up to isomorphism) embedding of $G$ into a a $2$-dimensional oriented surface, compatible with the combinatorial map, such that its faces are simply connected.
    The $2$-dimensional surface is determined by the Euler characteristic $v-e+f$ of the combinatorial map.
\end{corollary}

\subsection{Spherical and planar maps}\label{subsec:spherical-planar} Even though we are interested in embeddings of graphs in the sphere, the algorithm that proves that all twice-derived spherical graphs belong to a finite list (cf., \cref{p:first-topological-list}) works with embeddings of graphs in the plane. Let us spend some time to explain how to pass from a spherical graph to a planar one and vice versa.

We can pass from a spherical embedding with simply connected faces to a planar embedding by removing from the sphere one point which is not in the image of the embedding. Vice versa, to pass from the planar embedding to the spherical embedding it suffices to take the one-point compactification of the plane.

At the level of combinatorial maps this can be described as follows.

\begin{definition}[Spherical and Planar Maps]
    A \emph{spherical map} is a combinatorial map for a graph with the Euler characteristic of the sphere. A \emph{planar map} is a spherical map with a special face, namely a triple $(G, \next, Y)$, where $(G, \next)$ is a combinatorial map, and $Y$ one of its faces, namely an orbit of $\varphi = \next\circ \rev$.
\end{definition}

Once we fix an embedding corresponding to the underlying spherical map, we get a corresponding planar map by simply choosing a face in the spherical embedding. In geometric terms, this identifies a connected component in the complement in $\mathbb S^2$ of the image of the graph: by removing a point from this connected component, we find a planar embedding. The complement in $\mathbb R^2$ of the image of this embedding consists of one unbounded connected component and finitely many bounded connected components which are topological disks. Moreover the unbounded connected component's first homotopy group is $\mathbb Z$.

Conversely, given a connected graph and a planar embedding of it, if we compactify $\mathbb R^2$ to $\mathbb S^2$ by adding a point at infinity we recover a spherical embedding of the graph. The operations just describe give therefore a correspondence between spherical and planar maps (which clearly is not one-to-one). Note that it is also natural to use the term ``unbounded face'' for the special face of the planar combinatorial map.

\subsection{Generating the list of candidates}\label{s:genera-lista} The algorithm for generating the list of candidates twice-derived spherical graphs is recursive over the topological complexity of the graph. In order to describe the recursive procedure we introduce the following terminology.
Recall that, in graph theory, a ``leaf'' of a graph is a vertex of degree $1$.

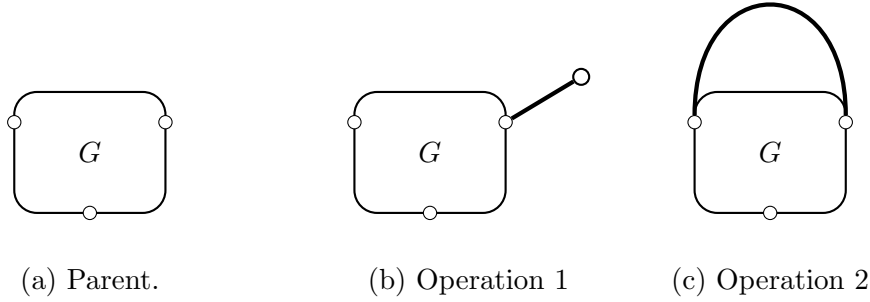
\begin{figure}[htbp]
    \centering
    \begin{tikzpicture}[
        node distance=2cm,
        vertex/.style={circle, draw=black, fill=white, inner sep=0pt, minimum size=5pt},
        newvertex/.style={circle, draw=black, fill=white, inner sep=0pt, minimum size=6pt, thick},
        mapboundary/.style={draw=black, thick, rounded corners=3mm},
        newedge/.style={ultra thick, draw=black},
        labeltext/.style={font=\sffamily\small}
    ]

        \begin{scope}[local bounding box=scope1]
            \draw[mapboundary] (-1,-0.8) rectangle (1,0.8);
            \node at (0,0) {$G$};

            \node[vertex] (v1) at (-1, 0.4) {};
            \node[vertex] (v2) at (1, 0.4) {};
            \node[vertex] (v3) at (0, -0.8) {};

            \node[below=0.5cm of scope1] {(a) Parent.};
        \end{scope}

        \begin{scope}[xshift=4.5cm, local bounding box=scope2]
            \draw[mapboundary] (-1,-0.8) rectangle (1,0.8);
            \node at (0,0) {$G$};
            \node[vertex] (v1) at (-1, 0.4) {};
            \node[vertex] (v2) at (1, 0.4) {};
            \node[vertex] (v3) at (0, -0.8) {};

            \node[newvertex] (new) at (2, 1) {};
            \draw[newedge] (v2) -- (new);

            \node[below=0.5cm of scope2] {(b) Operation 1};
        \end{scope}

        \begin{scope}[xshift=9cm, local bounding box=scope3]
            \draw[mapboundary] (-1,-0.8) rectangle (1,0.8);
            \node at (0,0) {$G$};
            \node[vertex] (u1) at (-1, 0.4) {};
            \node[vertex] (u2) at (1, 0.4) {};
            \node[vertex] (u3) at (0, -0.8) {};

            \draw[newedge] (u1) to[out=90, in=90, looseness=2.5] (u2);

            \node[below=0.5cm of scope3] {(c) Operation 2};
        \end{scope}

    \end{tikzpicture}

    \caption{The operations defining the children of a planar map. (a) A base planar map $(G, \next, Y)$ where the rectangle represents the boundary $Y$. (b) The first operation: attaching a new vertex (thick circle) to the boundary via a new edge (thick line). (c) The second operation: connecting two existing boundary vertices via a new edge (thick curve) drawn through the outer face.}
    \label{fig:operations}
\end{figure}

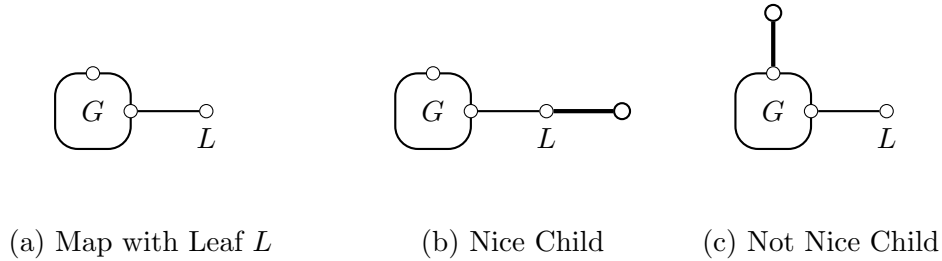
\begin{figure}[htbp]
    \centering
    \begin{tikzpicture}[
        vertex/.style={circle, draw=black, fill=white, inner sep=0pt, minimum size=5pt},
        leafvertex/.style={circle, draw=black, fill=white, inner sep=0pt, minimum size=5pt},
        newvertex/.style={circle, draw=black, fill=white, inner sep=0pt, minimum size=6pt, thick},
        mapboundary/.style={draw=black, thick, rounded corners=3mm},
        newedge/.style={ultra thick, draw=black},
        labeltext/.style={font=\sffamily\small}
    ]

        \begin{scope}[local bounding box=scope1]
            \draw[mapboundary] (-0.5,-0.5) rectangle (0.5,0.5);
            \node at (0,0) {$G$};

            \node[vertex] (root) at (0.5, 0) {};
            \node[vertex] (topvertex) at (0, 0.5) {};

            \node[leafvertex, label=below:$L$] (leaf) at (1.5, 0) {};
            \draw[thick] (root) -- (leaf);

            \node[below=0.8cm of scope1] {(a) Map with Leaf $L$};
        \end{scope}

        \begin{scope}[xshift=4.5cm, local bounding box=scope2]
            \draw[mapboundary] (-0.5,-0.5) rectangle (0.5,0.5);
            \node at (0,0) {$G$};
            \node[vertex] (root) at (0.5, 0) {};
            \node[vertex] (topvertex) at (0, 0.5) {};

            \node[vertex, label=below:$L$] (oldleaf) at (1.5, 0) {};
            \draw[thick] (root) -- (oldleaf);

            \node[newvertex] (new) at (2.5, 0) {};
            \draw[newedge] (oldleaf) -- (new);

            \node[below=0.8cm of scope2] {(b) Nice Child};
        \end{scope}

        \begin{scope}[xshift=9cm, local bounding box=scope3]
            \draw[mapboundary] (-0.5,-0.5) rectangle (0.5,0.5);
            \node at (0,0) {$G$};
            \node[vertex] (root) at (0.5, 0) {};
            \node[vertex] (topvertex) at (0, 0.5) {};

            \node[leafvertex, label=below:$L$] (stilleaf) at (1.5, 0) {};
            \draw[thick] (root) -- (stilleaf);

            \node[newvertex] (new) at (0, 1.3) {};
            \draw[newedge] (topvertex) -- (new);

            \node[below=0.8cm of scope3] {(c) Not Nice Child};
        \end{scope}

    \end{tikzpicture}

    \caption{Illustration of ``Nice'' vs. ``Not Nice'' children. (a) A map $X= (G, \next, Y)$ containing a leaf vertex $L$ and another boundary vertex on top. (b) A {\em Nice Child}: The new edge is attached to $L$. Consequently, $L$ becomes an internal vertex (degree 2) and is no longer a leaf. (c) A {\em Not Nice Child}: The new edge is attached to the other boundary vertex, so $L$ remains a leaf in the child map.}
    \label{fig:nice_children}
\end{figure}

\begin{definition}[Children of a Planar Map]
    Let $(G, \next, Y)$ be a planar map. Its \emph{children} are the planar maps that can be obtained from $(G, \next, Y)$ with one of the following two operations (cf. \cref{fig:operations}):
    \begin{itemize}
        \item Attach a new vertex ``externally'' to a vertex in $Y$ (the special face of the child consists of $Y$ and the two additional darts coming from the added edge).
        \item Connect two vertices (the two vertices can coincide) in $Y$ ``externally'' with a new edge (the special face of the child consists then of a proper subset of $Y$ and one of the darts coming from the added edge).
    \end{itemize}

    The \emph{nice children} are the children such that:
    \begin{itemize}
        \item Any leaf of $G$ is not a leaf of the child (cf. \cref{fig:nice_children}).
    \end{itemize}
A \emph{descendant} of a planar map is either the planar map itself or a child of a descendant of the planar map. A \emph{nice descendant} of a planar map is either the planar map itself or a nice child of a nice descendant of the planar map.
\end{definition}

\vspace{1cm}

The following lemma allows us to consider only nice children during our search.

\begin{lemma}\label{lem:no-leaf-generation}
    Any connected planar map with at least one edge that either has no leaves, or has only one leaf contained in the unbounded face, is a nice descendant of a cycle.
\end{lemma}

\begin{proof} First of all, we observe that the statement is obvious if the map itself is a cycle. Then we prove the general case by induction on the number of edges. If there is only one edge, then the two endpoints of the edge must coincide, and thus the map is a cycle.

Assume now that the statement holds when the number of edges is $k\geq 1$ and fix a planar map $X$ as in the statement with $k+1$ edges. If there is one leaf, we can remove it to get a parent planar map which satisfies the assumption of the lemma, from which we conclude that $X$ is a nice descendant of a cycle from the inductive hypothesis. If there are no leaves, then all vertices have degree at least $2$. Consider the unbounded face. If all the vertices adjacent to it have degree exactly two, then the face is a connected component of the planar map, and thus the entire map is the face itself. In particular it is a cycle. Otherwise there is a vertex $V$ of degree $3$. If we now remove one of the edges in the unbounded face which is adjacent to $V$, the resulting parent planar map $X^p$ has $k$ edges and is still connected. Consider the other vertex $V'$ adjacent to the edge which has been removed. If $V'=V$ then either $V$ had degree $>3$ in $Y$, and hence $X^p$ has no leafs, or it had degree exactly $3$, and then $X^p$ contains only one leaf, which is contained in its unbounded face. If $V'\neq V$, then $V$ has degree at least two in $X^p$, while $V'$, which belongs to the unbounded face of $X^p$, is the only vertex which could be a leaf in $X^p$. In all these cases $X^p$ falls under the assumptions of the lemma, and thus by the induction assumption it is a nice descendant of a cycle, in turn implying the same property for $X$.
\end{proof}

A first algorithm to prove that twice-derived spherical graphs belong to a finite list could go as follows.
We define recursively families $\mathcal F_0, \mathcal F_1, \mathcal F_2,\dots$ of planar maps such that $\mathcal F_k$ contains only planar maps with exactly $k$ edges.
\begin{itemize}
\item We define $\mathcal F_0$ as the empty set.

\item $\mathcal F_{k+1}$ consists of all the nice children of the planar maps in $\mathcal F_k$ and of the cycle of length $k+1$.
\end{itemize}
Because of \cref{lem:no-leaf-generation}, a connected planar map without leaves and with $k$ edges belongs to $\mathcal F_k$. Since we know that a twice-derived spherical graph has at most $30$ edges, we can examine all the generations $\mathcal{F}_k$ with $k\leq 30$ and check which planar maps in these lists satisfy the restrictions of \cref{p:derived-2}. However the cardinality of the families $\mathcal F_0, \mathcal F_1, \mathcal F_2, \dots$ grows quickly (cf. e.g. \cite{oeisA000109}), making the above algorithm unfeasible. Hence, the following criteria help us making $\mathcal{F}_k$ sensibly smaller, since they allow us to throw away elements of $\mathcal{F}_k$ for which we are sure that none of their descendants will satisfy the properties of \cref{p:derived-2}.

\begin{lemma}\label{lem:heuristic-subgraph}
    Let $(G, \next, Y)$ be a planar map. Denote by $\mathcal V(G)$, $\mathcal E(G)$, $\mathcal F(G)$ the vertices, edges, and faces of $G$.
    Denote by $\mathcal{BV}(G)$ the vertices on the boundary and by $\mathcal{IF}(G)$ the \emph{bounded} faces of $G$.
    Assume that a descendant of $(G, \next, Y)$ comes from a twice-derived spherical graph after removing a point from the spherical embedding (cf., \cref{subsec:spherical-planar}).
    The following properties hold.
    \begin{enumerate}
        \item Every $F\in\mathcal{IF}(G)$ has at most $7$ sides.
        \item Let $\gamma:=\#\{F\in\mathcal{IF}(G): |F|>5\}$. Then $0\le \gamma\le 2$.
        \item $\deg_G (w) \leq 3$ for any $w\in \mathcal V(G)$ and $\deg_G(w) = 3$ for $w\in \mathcal{V} (G) \setminus \mathcal{BV}(G)$.
        \item Let $Q_i(j):=6 q_i(j)-j$, where $q_i(j)$ is the smallest size of a collection of positive integers in $[1, 7]$, with sum equal to $j$, with at most $i$ numbers larger than $5$.
        Set $a:= \sum_{v\in \mathcal{BV}(G)} (3-\deg_G(v))$ and
        \[
        \tilde a= \left\{\begin{array}{ll}
        a \quad & \mbox{if $a$ is even},\\
        a+3 \quad &\mbox{if $a$ is odd}.
        \end{array}
        \right.
        \]
        We have
        \[
        \min_{\substack{\tilde a \leq \alpha \leq 60-|\darts (G)| \\ 6\mid \alpha-\tilde a}} Q_{2-\gamma}\big(|Y| + \alpha\big) + \sum_{F\in \mathcal{IF}(G)} (6-\text{sides}(F)) \le 12\, .
        \]
    \end{enumerate}
\end{lemma}
\begin{proof}
Let $(G', \next', Y')$ be the descendant of $(G, \next, Y)$ that is a twice-derived spherical graph (up to arbitrarily choosing a face to obtain a planar graph).

Note that a bounded face of $G$ is necessarily also a bounded face of $G'$. Recalling \cref{p:derived-2} no face of $G'$ can have more than $7$ sides and this proves statement (1). Moreover there can be at most two faces of $G'$ with more than $5$ sides, which proves (2). Next, a vertex in $\mathcal{V} (G)\setminus \mathcal{BV} (G)$ has the same degree in $G$ as in $G'$, which proves the second part of (3). In general the degree of a vertex can only increase in $G'$, which shows the first part of (3).

Coming to (4), let $a'=|\darts(G')\setminus \darts(G)|$.
Since $G'$ is a cubic graph, we have
\begin{equation*}
    a' = a + 3 |\mathcal V(G')\setminus \mathcal V(G)|.
\end{equation*}
Observe also that $a'$ is even because $a'=2|\mathcal E(G')\setminus \mathcal E(G)|$. Thus, we deduce $6\mid a'-\tilde a$.
Moreover, $a'+\darts(G)=\darts(G')\le 60$ by \cref{p:derived-2}, and thus $\tilde a \le a'\le 60-\darts(G)$.

Again by \cref{p:derived-2}, we have
\begin{equation*}
\sum_{\text{$F\in \mathcal F(G')\setminus \mathcal{IF}(G)$}} (6-\text{sides}(F))
+ \sum_{F\in \mathcal{IF} (G)} (6-\text{sides} (F))
= \sum_{F\in \mathcal{F} (G')} (6-sides (F)) = 12\, .
\end{equation*}

Our observations reduce the desired statement to proving
\begin{equation*}
 \sum_{\text{$F\in \mathcal F(G')\setminus \mathcal{IF}(G)$}} (6-\text{sides}(F))
 \ge Q_{2-\gamma}\big(|Y| + a'\big)\, .
\end{equation*}
The sum of the sides of $F$ for $F$ running in $\mathcal F(G')\setminus \mathcal{IF}(G)$ is exactly the number $a'$ of new darts introduced to get $G'$ and the darts in the unbounded face of $G$, in other words
\begin{equation*}
\sum_{\text{$F\in \mathcal F(G')\setminus \mathcal{IF}(G)$}} \text{sides}(F) = |Y| + a'\,.
\end{equation*}
The sides of each face in $\mathcal{F} (G')\setminus \mathcal{IF}(G)$ are at least $1$ and at most $7$ and there are at most $2-\gamma$ such faces with more than $5$ sides. Hence, we deduce that
\begin{equation*}
    |\mathcal F(G')\setminus \mathcal{IF}(G)| \ge q_{2-\gamma}(|Y|+a') ,
\end{equation*}
We finally conclude
\begin{align*}
 \sum_{\text{$F\in \mathcal F(G')\setminus \mathcal{IF}(G)$}} (6-\text{sides}(F))
 &= 6|\mathcal F(G')\setminus \mathcal{IF}(G)| - (|Y|+a') \\
 &\ge 6q_{2-\gamma}(|Y|+a')-(|Y|+a')=Q_{2-\gamma}(|Y|+a') .
\end{align*}
\end{proof}

We can use this Lemma to improve the initial algorithm. We define $\mathcal C_0$ as the empty set and then $\mathcal C_{k+1}$ is defined recursively as the union of the cycle of length $k+1$ (provided $k+1\le 7$) and the set of all nice children of the planar maps in $\mathcal C_k$ satisfying the constraints (1), (2), (3), (4).
Because of \cref{lem:heuristic-subgraph}, we conclude that a twice-derived net with $k$ edges belongs necessarily to $\mathcal C_k$.

However, this new modified algorithm is not yet the one used in our implementation. Without further considerations the families $\mathcal{C}_k$ would contain multiple copies of the same spherical graph, since we are not accounting for isomorphisms of the maps. In theory we could still run the algorithm as described and only deal with the multiple counting when checking which elements of $\mathcal{C}_k$ are twice-derived spherical graphs.  However this would be computationally too costly. In our actual algorithm we therefore avoid these multiple counts in every class $\mathcal{C}_k$. How this is done will be explained in \cref{s:isomorphisms}.

\subsection{Reconstruction of once-derived candidates}\label{s:primo-programma-bis} We describe here how, given the list of \cref{p:first-topological-list}, we produce the list of \cref{p:second-topological-classification}. The goal is just to give the appropriate rules to ``add back'' the vertices of degree $2$ which were removed, recalling that in a once-derived marked spherical graph there are at most four and at least two degree $2$ vertices.
\begin{itemize}
\item[(1)] We fix a graph $G = (V,E)$ and its spherical embedding $\iota$, coming from the list of \cref{p:first-topological-list}.
\item[(2)] We fix the number $2\leq n \leq 4$ of degree two vertices which we wish to add.
\item[(3)] We fix an edge $\sigma\in E$ and let $A$ and $B$ be the two vertices connected by it. We then remove $\sigma$, add a new vertex $C$, add two edges (adjacent to $A$ and $C$ and to $B$ and $C$, respectively), and mark one of them. At the same time we produce an embedding $\iota_{(1)}$ of this new graph $G_{(1)}$ from the embedding $\iota$ of $G$ in the obvious way, namely ``splitting'' the curve connecting $\iota (A)$ and $\iota (B)$ in two pieces.
\item[(4)] We repeat step (3) the necessary number of times to get to $G_{(n)}$ and $\iota_{(n)}$.
\item[(5)] We check whether $G_{(n)}$ and $\iota_{(n)}$ satisfy the conditions stated in \cref{p:derived-1} and in \cref{s:sieves}.
\end{itemize}
The procedure above is repeated over all possible choices of: elements in the list of \cref{p:first-topological-list}, natural numbers $n$ between $2$ and $4$, edges selected in step (3), and markings decided in the same step.
The algorithm described above can be easily expressed in terms of combinatorial maps and we leave the task to the reader.
As is the case for the algorithm of \cref{s:primo-programma}, we would incur into multiple counting of the same marked spherical graph, if we were to ignore the isomorphic copies which might occur in the list. This point is addressed as described in \cref{s:isomorphisms}.

\section{Useful inequalities for convex faces}\label{s:spher-trig}

We prove here some elementary inequalities for the sides of convex triangles, quadrilaterals, and, more importantly, pentagons with angles equal to $\frac{2\pi}{3}$, which are derived from \cref{p:formulae}.  These will be the backbone of the algorithm explained in \cref{s:secondo-programma}, which will allow us to effectively search for approximate metrizations as once derived nets of the marked spherical graphs listed in \cref{p:second-topological-classification}. The inequalities can be thought as quantitative versions of the monotonicity considerations of Heppes in \cite{Heppes}.

\begin{lemma}\label{lem:quadrilateral_consecutive_lowerbound}
    Let $ABCD$ be a convex spherical quadrilateral with angles of $\frac{2\pi}{3}$.
    Then
    \begin{equation*}
        |AB| + |BC| \ge 2\arctan 2\sqrt2\approx 2.46\ldots \,.
    \end{equation*}
\end{lemma}

\begin{lemma}\label{lem:pentagon_consecutive_lowerbound}
    Let $ABCDE$ be a convex spherical pentagon with angles of $\frac{2\pi}{3}$.
    Then
    \begin{equation*}
        |AB| + |BC| \ge 1 \,.
    \end{equation*}
\end{lemma}

\begin{proposition}\label{p:pentagon-monotonicity}
    Let $ABCDE$ be a convex spherical pentagon with angles of $\frac{2\pi}{3}$. Denote by $\lambda_T$ the length of the side of the equilateral triangle with angles $\frac{2\pi}{3}$ and by $f$ the function
\begin{equation}\label{e:fpentagono-2}
f(x, y) := \cos x+\cos y +\cos x\cos y -2\sin x\sin y \,
\end{equation}
(cf. \cref{p:formulae}).
The following holds for any choice of the numbers $x,y,z$ in the range $(0, \lambda_T)$:
    \begin{itemize}
        \item If $|AB|\ge x$ and $|BC|\ge y$, then
        \begin{equation*}
            2\cos(|DE|) - \frac13\le f(x, y) \,.
        \end{equation*}
        \item If $|AB|\le x$ and $|BC|\le y$, then
        \begin{equation*}
            2\cos(|DE|) -\frac13 \ge f(x, y) \,.
        \end{equation*}
        \item If $|AB|\le x$ and $|DE|\ge z$, then
        \begin{equation*}
            2\cos z - \frac13 \ge f(x, |BC|) \,.
        \end{equation*}
        \item If $|AB|\ge x$ and $|DE|\le z$, then
        \begin{equation*}
            2\cos z - \frac13 \le f(x, |BC|) \,.
        \end{equation*}
    \end{itemize}
\end{proposition}

\subsection{Proof of \cref{lem:quadrilateral_consecutive_lowerbound}}
Let $a=|AB|/2$ and $b=|BC|/2$. We have
    \begin{equation*}
        \tan a\cdot \tan b = \frac12 \,,
    \end{equation*}
    with $0<a, b<\tfrac\pi2$.
    Using the tangent addition formula, we get
    \begin{equation*}
        \tan(a+b) = \frac{\tan a + \tan b}{1-\tan a\cdot\tan b} = 2(\tan a + \tan b) \ge 4\sqrt{\tan a\cdot \tan b} = 2\sqrt2\,.
    \end{equation*}
    Therefore,
    \begin{equation*}
        |AB|+|BC| = 2(a+b) \ge 2\arctan 2\sqrt 2 \,.
    \end{equation*}

\subsection{Proof of \cref{lem:pentagon_consecutive_lowerbound}}
    Since $\sqrt2\ge 1$, we can assume that $|AB|, |BC|\le \sqrt2$.

    For $0\le x\le \sqrt2$, we have
    \begin{equation*}
        0\le \sin x \le x \quad\text{ and } \cos x \ge 1-\frac{x^2}2\ge 0 \,.
    \end{equation*}
    Therefore, we have
    \begin{equation*}
        2 \ge 2\cos|DE|-\frac13 = f(|AB|, |BC|)
        \ge 3-|AB|^2-|BC|^2-2|AB|\cdot|BC| = 3-(|AB|+|BC|)^2 \,,
    \end{equation*}
    which is equivalent to the statement.

\subsection{Proof of \cref{p:pentagon-monotonicity}}
The proposition follows immediately from the formula \eqref{e:formula-Taylor} and the following elementary monotonicity property of $f$.

\begin{lemma}\label{lem:pentagon_monotonicity}
The function $f$ of \eqref{e:fpentagono} is decreasing with respect to both its variables in the square $(0, \lambda_T)^2$.
\end{lemma}
\begin{proof}
    Observe that $\tfrac \pi 2 < \lambda_T < \pi$.
    Let us compute the derivative of $f$ with respect to $x$:
    \begin{equation*}
        \partial_x f(x, y) = -\sin x(1+\cos y) - 2\cos x\sin y =
        -\sin x\sin y\Big(\frac{1+\cos y}{\sin y} + \frac{2\cos x}{\sin x}\Big) \,.
    \end{equation*}
    Both $\theta\mapsto \frac{1+\cos \theta}{\sin \theta}$ and $\theta\mapsto \frac{2\cos \theta}{\sin \theta}$ are decreasing functions for $\theta\in(0, \pi)$. Thus, for $x, y\in (0, \lambda_T)$, we deduce
    \begin{equation*}
        \partial_x f(x, y) \le -\sin x\sin y \Big(\frac{1+\cos\lambda_T}{\sin\lambda_T} + \frac{2\cos\lambda_T}{\sin\lambda_T}\Big) = -\sin x\sin y \frac{1+3\cos\lambda_T}{\sin\lambda_T} = 0 \,.
    \end{equation*}
    Hence, $f$ is decreasing with respect to $x$ and, by symmetry, also with respect to $y$.
\end{proof}

\section{Relaxed Almgren competitors}\label{s:relaxed}

Several arguments in the rest of the paper will take advantage of the minimizing assumption of $E$ in \cref{t:main} by exhibiting suitable competitors. However the simpler ``ideal comparison surfaces'' will typically not comply with the constraints required in \cref{d:Almgren-minimizing} and we will need to argue that they can be suitably approximated with competitors which instead comply with them. This section proves three very useful general propositions which will allow us to relax considerably the requirements of \cref{d:Almgren-minimizing}. The first one shows that it suffices to define the map $\Phi$ in \cref{d:Almgren-minimizing} on the cone $E$.

\begin{proposition}\label{p:only-on-E}
Assume that $\ell \subset \mathbb R^3$ is a line through the origin and $\ell\subset E\subset \mathbb R^3$ is the cone over a link $\mathcal{N}\subset \mathbb S^2$ which is a connected geodesic net (cf. \cref{d:reti}). Let $H: E \to \mathbb R^3$ be a Lipschitz map such that:
\begin{itemize}
\item[(i)] $H (x) =x$ outside some ball $B_{r_0}$;
\item[(ii)] If $W:= \{x\in E: H (x)\neq x\}$, then $\overline{W\cup H (W)} \subset \mathbb R^3 \setminus \ell$.
\end{itemize}
Then there is an extension $G$ of $H$ to the whole Euclidean space $\mathbb R^3$ which complies with the requirements of \cref{d:Almgren-minimizing}.
\end{proposition}

The next propositions will give substantial more leeway, because we will allow the set $\overline{W\cup \Phi (W)}$ in \cref{d:Almgren-minimizing} to ``touch'' $\ell$ in a controlled way. For the first proposition we introduce the following terminology.

\begin{definition}
Given $\theta \in (0, \frac{\pi}{2})$ we denote by $C_\theta(\ell)$ the closed ``ice-cream'' cone
\begin{equation*}
    C_\theta(\ell) :=
    \{x\in\R^3:\, \dist\, (x, \ell) \leq |x| \sin\theta\} \,.
\end{equation*}
\end{definition}
We will also use the notation $B_\delta (A)$ for the open set $\{x: \dist\, (x, A) < \delta\}$, namely the open $\delta$-neighborhood of a set $A$.

\begin{proposition}[Ice-cream relaxed condition]\label{p:Almgren-relaxed}
Assume $\ell \subset \mathbb R^3$ is a line through the origin and $\ell\subset E\subset \mathbb R^3$ the cone over a connected geodesic net $\mathcal{N}\subset \mathbb S^2$. Let $F: E \to \mathbb R^3$ be a Lipschitz map such that:
\begin{itemize}
\item[(i)] $F(x)=x$ outside some ball $B_{r_0}$;
\item[(ii)] $F(x) = x$ for $x\in C_\theta(\ell)$ and $F^{-1} (\ell) \subset \ell$.
\end{itemize}
Then, for every $\varepsilon > 0$, there is a Lipschitz map $G: \mathbb R^3 \to \mathbb R^3$ such that:
\begin{itemize}
\item[(a)] $G(x) = x$ outside $B_{r_0}$;
\item[(b)] upon setting $W:= \{x: G (x)\neq x\}$, we have $\overline{W\cup G(W)} \subset \R^3\setminus \ell$;
\item[(c)] $\mathcal{H}^2 (G (E\cap \overline{B}_{r_0})) \leq \mathcal{H}^2 (F(E\cap \overline{B}_{r_0})) + \varepsilon$.
\end{itemize}
\end{proposition}

The next one will require slightly stronger conditions in a structurally different cone.

\begin{definition}\label{d:wedge}
We will call {\em wedge} any open convex cone of $\mathbb R^3$ which is invariant under translations along the direction of $\ell$.
\end{definition}
Note that such a wedge $V$ can be written as $V'+\ell$, where $V'$ is an open convex cone in the plane orthogonal to $\ell$. Furthermore observe that $V'$ is either the whole plane or an open convex subset of the plane bounded by two half-lines emanating from the origin.

\begin{proposition}[Wedge relaxed condition]\label{p:Almgren-relaxed-2}
Assume that $\ell \subset \mathbb R^3$ is a line through the origin, $\ell\subset E\subset \mathbb R^3$ the cone over a connected geodesic net $\mathcal{N}\subset \mathbb S^2$.
Let $F: E\to \mathbb R^3$ be a Lipschitz map such that:
\begin{itemize}
\item[(i)] $F(x)=x$ outside some ball $B_{r_0}$;
\item[(ii)] upon setting $U := \{x:F(x)\neq x\}$, we have $\overline{U \cup F(U)} \subset \R^3\setminus (V\cap B_\delta (\ell))$ for some wedge $V$ and some $\delta>0$;
\end{itemize}
Then, for every $\varepsilon > 0$, there is a Lipschitz map $G: \mathbb R^3 \to \mathbb R^3$ such that:
\begin{itemize}
\item[(a)] $G(x)=x$ outside $B_{r_0}$;
\item[(b)] upon setting $W:= \{x: G (x)\neq x\}$, we have $\overline{W\cup G(W)} \subset \R^3\setminus \ell$;
\item[(c)] $\mathcal{H}^2 (G (E\cap \overline{B}_{r_0})) \leq \mathcal{H}^2 (F(E\cap \overline{B}_{r_0})) + \varepsilon$.
\end{itemize}
\end{proposition}

\subsection{Proof of \texorpdfstring{\cref{p:only-on-E}}{p:only-on-E}}
We will need the following technical lemma, whose proof can be achieved with a regularization procedure.

\begin{lemma}\label{l:approx-1}
Let $\Omega\subset \mathbb R^k$ be a bounded open set and let $H:\overline\Omega \to \mathbb R^N$ be a continuous map whose restriction to $\partial \Omega$ is Lipschitz. Then, for every $\eta>0$, there is a Lipschitz map $G: \overline\Omega \to \mathbb R^N$ such that $G=H$ on $\partial \Omega$ and $\|G-H\|_{C^0} \leq \eta$.
\end{lemma}
\begin{proof}
Let $\tilde H:\overline \Omega\to \R^N$ be a Lipschitz function that coincides with $H$ on $\partial\Omega$.
Let $A:\R^k\to \R^N$ be defined as $A=H-\tilde H$ in $\Omega$ and $A=0$ in $\mathbb R^k \setminus \Omega$. The function $A$ is continuous and compactly supported.
Fix a standard smooth convolution kernel $\varphi$, and, for any $\eps>0$, let $A_\eps:= A*\varphi_\eps$. Observe that $A_\eps\to A$ uniformly as $\eps\to 0$ and furthermore $A_\eps$ is Lipschitz. Define $G_\eps:\overline\Omega\to\R^N$ as
\[
G_\eps(x) = \tilde H(x) + \min\Big\{1, \tfrac{\dist\, (x, \partial\Omega)}\eps
\Big\} A_\eps(x) \,.
\]
For any $\eps>0$, the function $G_\eps$ is Lipschitz and $G_\eps(x)=\tilde H(x)=H(x)$ for $x\in\partial\Omega$. It remains only to bound uniformly $|G_\eps(x)-H(x)|$ for $x\in\overline\Omega$.
If $\dist\, (x,\partial\Omega)\ge \eps$, then
\begin{equation*}
    |H(x)-G_\eps(x)| = |A(x)-A_\eps(x)| \le \norm{A-A_\eps}_{C^0}\,,
\end{equation*}
while if $\dist\, (x,\partial\Omega)<\eps$, then
\begin{equation*}
    |H(x)-G_\eps(x)| \le |A(x)| + |A_\eps(x)| \le \norm{A-A_\eps}_{C^0} +  2\sup_{x:\, \dist\, (x, \partial\Omega)<\eps}|A(x)|\,.
\end{equation*}
Since $A|_{\partial \Omega} \equiv 0$ and $A$ is uniformly continuous, the last two inequalities imply that, if we choose $\eps$ small enough compared to $\eta$, the map $G_\eps$ satisfies the requirements of the statement.
\end{proof}

\begin{proof}[Proof of \cref{p:only-on-E}] First of all, without loss of generality we will assume that $r_0=1$. Next extend $H$ on the complement of $B_1$ simply by setting it equal to the identity. The main task becomes therefore to further extend to $\overline B_1\setminus E$. Let
\[
W' := \{x\in E: x\neq H (x)\}
\]
and choose $\delta$ so that $2\delta < \theta$ and
\[
2\delta < \min \{\dist\, (y, \ell): y\in W' \cup H (W')\}\, .
\]
Impose further that $2\delta$ is strictly smaller than the minimum distance between any pair of vertices of
the net $\mathcal{N}$.

We then extend $H$ as the identity on the closed set
\[
K := (\mathbb R^3 \setminus B_{1-2\delta}) \cup \overline{B}_{2\delta} (0)\cup C_{2\delta} (\ell)\, .
\]
This extension is clearly Lipschitz on $K\cup E$. Now, let $U$ be the open set
\[
U:= B_{1-2\delta} \setminus (K\cup E)\, .
\]
Denote by $P$ and $-P$ the two pins of the net $\mathcal{N}$.
Each connected component $\Omega$ of $U$ corresponds to a connected component $\Lambda$ of $V:= \mathbb S^2\setminus (\mathcal{N}\cap \overline{B}_{2\delta} (P) \cup \overline{B}_{2\delta} (-P))$
 through the relation
\[
\Omega = \{ tx : x\in \Lambda, 2\delta < t <1-2\delta\}\, .
\]
In particular, since each $\Lambda$ is homeomorphic to the open disk of $\mathbb R^2$, each $\Omega$ is homeomorphic to the open $3$-ball. We can then extend the map $H$ {\em continuously} in each connected component $\Omega$ so that the values of the latter extension over each $\Omega$ stay at distance $2\delta$ from $\ell$: this follows from the fact that the second homotopy group of $\mathbb R^3 \setminus \overline{B}_{2\delta} (\ell)$ is trivial. We denote by $H'$ such extension.

We can now use \cref{l:approx-1} to perturb $H'$ in $\Omega$ to a new Lipschitz map $H''$ close to $H'$ in the $C^0$ norm and coinciding with $H$ on $\partial \Omega$: the map $G$ is defined to be $H$ on $\mathbb R^3\setminus \Omega$ and $H''$ on $\overline{\Omega}$. $G$ is then Lipschitz and the closure of the set $W:=\{x: G(x)\neq x\}$ is contained in $B_1\setminus \ell$. If we impose the inequality $\|H-H''\|_{C^0}\leq \delta$ then we also ensure $G (W\cap \Omega) \subset \mathbb R^3 \setminus B_{\delta} (\ell)$. Since however we know that $G (W\setminus \Omega) \ssubset \mathbb R^3 \setminus \ell$, we conclude that $G (W)\ssubset \mathbb R^3\setminus \ell$. Therefore $G$ satisfies all the requirements of the Proposition.
\end{proof}

\subsection{Proof of \texorpdfstring{\cref{p:Almgren-relaxed}}{p:Almgren-relaxed}} Since we can invoke \cref{p:only-on-E}, we really need to define $G$ on $E$. Moreover, without loss of generality we will assume $r_0=1$.

\medskip

{\bf Step 1. Definition of $G$ on $E$.} First of all, since $E$ is a cone and $F$ is the identity on $\partial B_1 \cap E = \mathcal{N}$, we know that $F$ on $E\setminus \overline{B}_1$ as a $1$-homogeneous map, i.e. we have
\[
F (x) = |x| F  \big({\tfrac{x}{|x|}}\big)\qquad \forall x\in E\setminus \overline{B}_1\, .
\]
We then fix $\eta>0$, whose choice will be specified later, and define the following Lipschitz function $\lambda$ on the positive real line:
\[
\lambda (t) =
\left\{\begin{array}{ll}
\textstyle{\frac{1}{\eta}} & \mbox{if $0<t\leq \eta$}\\
\exp (- \eta^{-1} (2\eta - t) \log \eta) \quad& \mbox{if $\eta< t < 2\eta$}\\
1 & \mbox{if $t\geq 2\eta$}
\end{array}
\right.\, .
\]
Note that $1\leq \lambda \leq \frac{1}{\eta}$.

At this point it is notationally convenient to use $d (x)$ for the $\dist\, (x, \ell)$. We then introduce
\[
\Phi (x) := \lambda (d(x)) x\,
\]
and finally define explicitly the map $G$ for $x\in E$, by setting
\[
G (x) = \left\{
\begin{array}{ll}
x \qquad &\mbox{if $d (x) < \eta$}\\ \\
\frac{|x|}{|\Phi (x)|} F (\Phi (x)) \qquad &\mbox{otherwise}
\end{array}
\right.\, .
\]
We show below that this map satisfies all the requirements of the proposition. 

\medskip

{\bf Step 2. Topological conditions}. In this step we check that, upon choosing $\eta$ small enough, $G$ complies with the requirements of \cref{d:Almgren-minimizing}. The fact that (when $d(x)\neq 0$) $x = \frac{|x|}{|\Phi (x)|} \Phi (x)$ will be used repeatedly.
By construction, the map $G$ is the identity both when $d(x)<\eta$ and when $|x|\ge 1$. Moreover, it is continuous and locally Lipschitz. Hence, it is globally Lipschitz.
We also have that the set $W:=\{x: x\neq G (x)\}$ has positive distance from $\ell$ and is contained in $B_1$.

Observe moreover that $
d(\Phi (x)) =
\lambda (d(x)) d(x) \geq d (x)$ for every $x$. In particular $\Phi  (W)$ has positive distance from $\ell$.
Given that $F^{-1} (\ell)\subset \ell$, we immediately conclude that $\overline{F(\Phi (W))} \subset F (\Phi (\overline{W}))$ is a compact subset of $\mathbb R^3\setminus \ell$. Observe next that, since $\lambda \leq \frac{1}{\eta}$,
\[
\inf \{d (G(x)):x\in W\} \geq \eta \inf \{d (F(\Phi (x))): x\in W\} > 0\, .
\]
Next, since $W$ is bounded and $G$ Lipschitz, $G(W)$ is bounded. Summarizing, we have proved that $G$ is Lipschitz, that $\overline{W}\subset \overline{B}_1\setminus \ell $, and that $G(W)\ssubset \mathbb R^3\setminus \ell$.

\medskip

{\bf Step 3. Area estimate}. In order to complete the proof we need to
prove estimate (C) for $\mathcal{H}^2 (G (E\cap \overline{B}_1))$. First of all we know that $G (x) = F (x)$ over $E\cap \overline{B}_1 \setminus B_{2\eta} (\ell)$. In particular $G (E\cap \overline{B}_1 \setminus B_{2\eta} (\ell))$ is contained in the set $F(E\cap \overline{B}_1)$.Thus we can write
\[
\mathcal{H}^2 (G (E\cap \overline{B}_1)) \leq \mathcal{H}^2 (G (E\cap \overline{B}_1\cap B_{2\eta} (\ell))) +
\mathcal{H}^2 (F(E\cap \overline{B}_1))\, .
\]
Observe that $\Phi (C_\theta (\ell))\subset C_\theta (\ell)$. In particular $G$ is the identity on $C_\theta (\ell) \cap E$ and thus
\[
\mathcal{H}^2 ( G(C_\theta (\ell) \cap E \cap \overline{B}_1 \cap B_{2\eta} (\ell))) = \mathcal{H}^2 ( C_\theta (\ell) \cap E\cap \overline{B}_1 \cap B_{2\eta} (\ell)) \leq C_0 \eta
\]
for some constant $C_0$ which depends on $E$ but not on $\eta$. Hence we can write
\begin{equation}\label{e:stimozza-1}
\mathcal{H}^2 (G (E\cap \overline{B}_1)) \leq \mathcal{H}^2 (F(E\cap \overline{B}_1)) + C_0 \eta +
\mathcal{H}^2 (G ((E\setminus C_\theta (\ell)) \cap \overline{B}_1\cap B_{2\eta} (\ell)))\, .
\end{equation}
We now observe further that
\[
(E\setminus C_\theta (\ell)) \cap \overline{B}_1\cap B_{2\eta} (\ell) \subset E \cap B_{3\theta^{-1} \eta}\,
\]
and, in order to estimate $\mathcal{H}^2 (G(E\cap B_{3\theta^{-1} \eta}))$, we introduce the notation $J_\sigma (A)$ for the Jacobian of the restriction of a linear map $A: \mathbb R^3\to \mathbb R^3$ to a $2$-dimensional plane $\sigma$ and define
\[
J = \sup \{J_{T_x E} (dG|_x) : x\in E \cap \overline{B}_{3\theta^{-1} \eta}\}\, .
\]
In particular, using the area formula, we can estimate
\begin{equation}\label{e:stimozza}
\mathcal{H}^2 (G ((E\setminus C_\theta (\ell)) \cap \overline{B}_1\cap B_{2\eta} (\ell)))
\leq \mathcal{H}^2 (G (E \cap \overline{B}_{3\theta^{-1} \eta}))
\leq C_1 J \eta^2\, ,
\end{equation}
where $C_1$ is a constant which depends on the cone $E$ and on $\theta$.

We recall that, given any orthonormal basis $\{e_1, e_2\}$ of the 2-dimensional plane $\sigma$, we have
\begin{equation}\label{e:formulozzo-10}
J_\sigma (A) = |A [e_1] \times A [e_2]| \leq |A [e_1]||A [e_2]|\, .
\end{equation}
Evidently in the region where $G$ is the identity, for every $2$-plane $\sigma$ we have that $J_{\sigma} (dG|_x) = 1$. We are thus interested in the region where $G (x)$ equals
\[
\frac{|x|}{|\Phi (x)|} F (\Phi (x)) = (\lambda (d(x)))^{-1} F ( \lambda (d (x)) x)\, .
\]
There we can compute
\[
dG_x = dF_{\Phi (x)} + \frac{\lambda' (d(x))}{\lambda (d(x))} \left(- \frac{F(\Phi(x))}{\lambda (d(x))} +
dF_{\Phi (x)} [x]\right) \otimes \nabla d (x)\, .
\]
We recall that $d(x)$ is a 1-Lipschitz function and it is differentiable on the complement of $\ell$, with $|\nabla d (x)|=1$. Given an arbitrary $x\in B_{3\theta^{-1} \eta}\setminus \ell \subset \overline{B}_1\setminus \ell$ (so we assume that $\eta\leq \frac{\theta}{3}$) and an arbitrary $2$-dimensional plane, we can choose an orthonormal basis $e_1, e_2$ with the property that $e_1 \perp \nabla d (x)$. In particular \eqref{e:formulozzo-10} gives the estimate
\[
J_\sigma (dG_x) \leq |dF_{\Phi (x)} [e_1]| \left(|dF_{\Phi (x)} [e_2]| + \frac{|\lambda' (d(x))|}{\lambda (d(x))} (|G (x)| + |dF_{\Phi (x)} [x]|)\right)
\]
for the specific $2$-dimensional plane $\sigma$.
Recall that $F$ is Lipschitz, that $|G(x)|\leq \|F\|_{C^0 (B_1)}$ for $x\in \overline{B}_1$, and that
\[
\left|\frac{\lambda'}{\lambda}\right|=| (\log (\lambda))'|\leq \frac{|\log \eta|}{\eta}\, .
\]
We thus achieve the estimate $J (dG_x) \leq C_2 \eta^{-1} |\log \eta|$, where the constant $C_2$ depends on $\Lip\, (F)$, $\|F\|_{C^0 (B_1)}$, and $\theta$ but not on the parameter $\eta$. We can now insert the latter bound in \eqref{e:stimozza} and \eqref{e:stimozza-1} to conclude
\[
\mathcal{H}^2 (G (E\cap \overline{B}_1)) \leq \mathcal{H}^2 (F(E\cap \overline{B}_1)) + C_0\eta (1+ |\log \eta|)\, .
\]
While the constant $C_0$ depends upon $\Lip\, (F)$, $E$, and $\theta$, it does not depend on $\eta$. Clearly if we choose $\eta$ small enough so that $C_0 \eta (1+|\log \eta|) < \varepsilon$ we achieve the estimate (c) in the statement of the proposition, thus completing the proof.

\subsection{Proof of \texorpdfstring{\cref{p:Almgren-relaxed-2}}{p:Almgren-relaxed-2}} Without loss of generality we can assume that $r_0=1$ and (since otherwise $F$ itself already satisfies the assumption of \cref{p:only-on-E}) that $V$ is not the entire space $\mathbb R^3$.

We next fix coordinates so that $\ell = \{(0,0,t) : t\in \mathbb R\}$. In particular $V= V'\times \mathbb R$, where $V'$ is an open convex cone of $\mathbb R^2$ bounded by two half-lines emanating from the origin. We then fix a Lipschitz deformation retraction $r: [0,1] \times \mathbb R^2 \setminus (B_1\cap V') \to \mathbb R^2\setminus B_1$ with the additional roperty that $r (t,x) \neq 0$ for every $t>0$ and define a corresponding Lipschitz retraction $\rho : [0,1]\times \mathbb R^3\setminus (B_1 (\ell) \cap V) \to \mathbb R^3 \setminus B_1 (\ell)$ as
\[
\rho (t, (x', x_3)) = (r (t,x'), x_3)\, .
\]
Finally, for every positive $\eta$ we define further the retraction $\rho_\eta: [0,1]\times \mathbb R^3 \setminus (B_\eta (\ell) \cap V) \to \mathbb R^3 \setminus B_\eta (\ell)$ as
\[
\rho_\eta (t, x) = \eta \rho \left(t, \tfrac{x}{\eta}\right)\, .
\]
Next, for every $x,y\in \mathbb R^3\setminus V$ we let $\gamma_{x,y}\subset \mathbb R^3\setminus V$ be the shortest path connecting $x$ and $y$ in $\mathbb R^3\setminus V$. The path exists because $\mathbb R^3\setminus V$ is closed and arc-connected, but in fact we will show at the end that it is unique and has a simple description. This allows us to define the map $[0,1] \ni t \mapsto \Gamma (t,x,y)$ as the constant speed parametrization of $\gamma_{x,y}$ and in the last step of the proof we will also show that $\Gamma$ is locally Lipschitz in all variables.

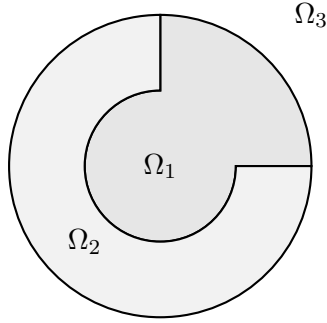
\begin{figure}[hbtp]
\begin{tikzpicture}[scale=1]


    \filldraw[fill=gray!20, draw=black, thick]
        (0,1) -- (0,2) arc (90:0:2cm) -- (1,0) arc (360:90:1cm) -- cycle;

    \filldraw[fill=gray!10, draw=black, thick]
     (0,1) -- (0,2) arc (90:360:2cm) -- (1,0) arc (360:90:1cm) -- cycle;

     \node at (0,0) {$\Omega_1$};
     \node at ({-1},{-1}) {$\Omega_2$};
     \node at (2,2) {$\Omega_3$};
\end{tikzpicture}
\caption{A cross section, orthogonal to $\ell$, of the partition of $\mathbb R^3$ into the regions $\Omega_i$. In this picture the wedge $V$ is given by $\{(x_1,x_2,x_3): x_1, x_2 \geq 0\}$.}\label{f:partition}
\end{figure}

Assuming the latter properties of $\Gamma$, we fix a positive number $\eta < \frac{\delta}{2}$ (whose choice will be specified later) and denote by $d (x)$ the distance of $x$ to $\ell$. We first define the map over $E$ and it is convenient to introduce the following partition of $\mathbb R^3$ into three sets which are clearly invariant under translations along $\ell$ (cf. \cref{f:partition}):
\begin{align*}
\Omega_1 &:= \overline{B}_\eta (\ell) \cup (\overline{B}_{2\eta} (\ell) \cap V)\\
\Omega_2 &:= \overline{B}_{2\eta (\ell)} \setminus (\overline{B}_\eta (\ell) \cup V)\\
\Omega_3 &:= \mathbb R^3\setminus \overline{B}_{2\eta} (\ell)\, .
\end{align*}
We then set
\begin{align*}
G (x) :=
\left\{
\begin{array}{ll}
x \quad &\mbox{if $x\in E\cap \Omega_1$}\\
\rho_\eta \left(\tfrac{d(x)-\eta}{\eta}, \Gamma \left(\tfrac{d(x)-\eta}{\eta}, x, F(x)\right)\right) \quad &
\mbox{if $x\in E\cap \Omega_2$}\\
\rho_\eta (1, F (x)) \quad &\mbox{if $x\in E\cap \Omega_3$}
\end{array}
\right.
\end{align*}
The role of the function $x\mapsto \frac{d(x)-\eta}{\eta}$ is to interpolate between the value $0$ on $\partial B_\eta (\ell)$ and $1$ on $\partial B_{2\eta} (\ell)$: any other interpolation with Lipschitz constant comparable to $\eta^{-1}$ would be suitable.

Observe next that, since $F(x)=x$ on $E\setminus B_1$ and thus $G(x) = \rho_\eta (1, x)$ on $E\cap \Omega_3\setminus B_1$ and $G(x) = \rho_\eta \left(\tfrac{d(x)-\eta}{\eta}, x\right)$ on $E\cap \Omega_2\setminus B_1$. Since however $\rho_\eta (t,x)=x$ for ever $x\in \mathbb R^3\setminus B_\eta (\ell)$ and $\Omega_2\cup \Omega_3 \subset \mathbb R^3\setminus B_\eta (\ell)$, we immediately conclude that $G(x)=x$ on $E\setminus B_1$.

Note next that by construction $W := \{x: G(x)\neq x\}$ is contained in the compact set $K := \overline{B}_1\setminus B_\eta (\ell)$. Next, in order to check that $G(W)\ssubset \mathbb R^3\setminus \ell$ it suffices to show that $G (K) \cap \ell = \emptyset$. By definition of $\rho$, $\rho_\eta (t,x)$ is an element of $\ell$ if and only if $t=0$ and $x\in \ell$. It is therefore obvious that $G (K)\cap \ell = \emptyset$. Thus we can apply \cref{p:only-on-E} and extend $G$ to a map on $\mathbb R^3$ which complies with the requirements (a) and (b). The rest of the proof is devoted to prove that, upon choosing $\eta$ sufficiently small, the map satisfies (c) as well.

\medskip

{\bf Step 1. $G$ is Lipschitz.} We first check the Lipschitz continuity of $G$. Clearly $G$ is Lipschitz in each of the three sets $\Omega_i$
and in fact $G|_{\Omega_i}$ extends to a unique Lipschitz map $G_i$ on $\overline{\Omega}_i$ simply using each of the three explicit formulae.
In order to check that it is globally Lipschitz it suffices to show that the map is continuous on $\overline{\Omega}_i \cap \overline{\Omega}_j$ for every $i\neq j$, which is equivalent to showing that $G_i = G_j$ on the latter set. This is easily checked as follows:
\begin{itemize}
\item $\overline{\Omega}_1\cup \overline{\Omega}_2$ consists of the two pieces $ \partial B_{\eta} (\ell) \setminus V$ and $\partial V \cap \overline{B}_{2\eta} (\ell) \setminus B_\eta (\ell)$. On the first we have
    \[
    G_2 (x) = \rho_\eta (0, \Gamma (0, x, F(x))) = \rho_\eta (0, x) =x = G_1 (x)\, ,
    \]
    while on the second we have $F(x)=x$ and $\rho_\eta (t,x) =x$, and therefore
    \[
    G_2 (x) = \rho_\eta (\tfrac{d(x)-\eta}{\eta}, \Gamma (\tfrac{d(x)-\eta}{\eta}, x, x)) = \rho_\eta (\tfrac{d(x)-\eta}{\eta}, x) = x = G_1 (x)\, ;
    \]
    \item $\overline{\Omega}_1\cap \overline{\Omega}_3 = \partial B_{2\eta} (\ell) \cap \overline{V}$. On this set we have $F(x)=x$ and $\rho_\eta (1, x)=x$ and therefore $G_3 (x) = x = G_1 (x)$.
\item $\overline{\Omega}_2\cap \overline{\Omega}_3= \partial B_{2\eta} (\ell) \setminus V$. On this set we have
\[
G_2 (x) = \rho_\eta (1, \Gamma (1, x, F(x)) = \rho_\eta (1, F(x)) =  G_3 (x)\, .
\]
\end{itemize}

\medskip

{\bf Step 2. Lipschitz constant.} We next claim that the Lipschitz constant of $G$ is bounded by a constant $L$ independent of $\eta$. In order to see this it suffices to show that $|D G|\leq C$ at all points of differentiability in $B_1\cap E$. This is simply achieved using the chain rule and observing that
\begin{itemize}
\item $|D_x \Gamma (t,x,y)| + |D_y \Gamma (t,x,y)| \leq C_0 (|x|+|y|)$;
\item $|\partial_t \Gamma (x,y)|\leq C_0 |x-y|$;
\item $|D_x \rho (t,x)|\leq C_0$;
\item $|\partial_t \rho (t,x)|\leq C_0 \eta$.
\end{itemize}
The last two are in fact obvious conclusions, while the first two estimates will be proved in the last step.
The bound on the derivative becomes then obvious on the regions $\Omega_1$ and $\Omega_3$.

The bound on $\Omega_2$ is more subtle. First define $\Lambda (x) = \Gamma \left(\tfrac{d(x)-\eta}{\eta}, x, F(x) \right)$ and observe that
\begin{align*}
D_x \Lambda (x)
=& \partial_t \Gamma \left(\tfrac{d(x)-\eta}{\eta}, x, F(x) \right) \otimes \tfrac{\nabla d (x)}{\eta}
+ D_x \Gamma \left(\tfrac{d(x)-\eta}{\eta}, x, F(x) \right)\\
&+ D_y \Gamma \left(\tfrac{d(x)-\eta}{\eta}, x, F(x) \right) D_x F (x)\, .
\end{align*}
We can tame the only ``dangerous'' term, namely the first term, using
\[
\left|\partial_t \Gamma \left(\tfrac{d(x)-\eta}{\eta}, x, F(x) \right) \otimes \tfrac{\nabla d (x)}{\eta}\right|
\leq \tfrac{C_0}{\eta} |x-F(x)|\,
\]
and observing that, since the above chain rule is applied when $d(x) \leq 2\eta$, the assumption that $F(x)=x$ on $\ell$ implies $|x-F(x)| \leq d(x) (1+\Lip\, (F)) \leq 2\eta (1 + \Lip\, (F))$. This shows then $|D_x \Lambda|\leq C$. Next we compute
\begin{align*}
D_x G (x) = \partial_t \rho_\eta \left(\tfrac{d(x)-\eta}{\eta}, \Lambda (x)\right) \otimes \tfrac{\nabla d(x)}{\eta} + D_x \rho_\eta \left(\tfrac{d(x)-\eta}{\eta}, \Lambda (x)\right) D_x \Lambda\, .
\end{align*}
This time the dangerous term is tamed using $|\partial_t \rho| \leq C_0 \eta$.

\medskip

{\bf Step 3. Area estimate.} Note that the maps $G$ and $F$ coincide on $\overline{B}_1\setminus (B_{2\eta} (\ell)\cap F^{-1} (B_{\eta} (\ell))$ and that $G(x)= \rho_\eta (1, F(x))$ for $x\in \overline{B}_1\setminus B_{2\eta} (\ell)$ and in particular $G=F$ on the set $\overline{B}_1\setminus (B_{2\eta} (\ell) \cup F^{-1} (B_\eta (\ell)))$. Thus we can estimate
\begin{align*}
\mathcal{H}^2 (G (E\cap \overline{B}_1 \setminus B_{2\eta} (\ell)))
&\leq \mathcal{H}^2 (E\cap \overline{B}_1)
+ \mathcal{H}^2 (\rho_\eta (1, B_{\eta} (\ell) \cap F(E\cap \overline{B}_1)))\\
&\leq \mathcal{H}^2 (E\cap \overline{B}_1)
+ C_0^2 \mathcal{H}^2 (F(E \cap \overline{B}_1)\cap B_{\eta} (\ell))\, .
\end{align*}
Now, $F (E\cap \overline{B}_1)$ is a compact set with finite Hausdorff measure and therefore
\[
\lim_{\eta\downarrow 0} \mathcal{H}^2 (F(E \cap \overline{B}_1)\cap B_{\eta} (\ell)) = 0\, .
\]
Given the above estimates, it is obvious that
\[
\mathcal{H}^2 (G (E\cap \overline{B}_1)) \leq \mathcal{H}^2 (F (E\cap \overline{B}_1)) + \varepsilon
\]
as soon as $\eta$ is sufficiently small.

\medskip
{\bf Step 4. Lipschitz continuity of $\Gamma$.} Consider $x= (x',t)$ and $y = (y', s)$ in $\mathbb R^3\setminus V$. If the segment joining $x$ and $y$ is contained in $\mathbb R^3$, then this is the (unique) path $\gamma_{x,y}$ of least length. This happens if and only if one of the following happens:
\begin{itemize}
    \item $x'$ is the origin;
    \item $x'\neq 0$ and if we let $\alpha^+\subset \mathbb R^2$ be the open half-space containing $V'$ and bounded by the line $\{\lambda x': \lambda \in \mathbb R\}$, then $y'\not\in \alpha^+$.
\end{itemize}
When $x'\neq 0$ and $y'\in \alpha^+$, upon setting $L:= |x'|+|y'|$ and
\[
z:= \left(0, \tfrac{|x'|}{L} t + \tfrac{|y'|}{L} s\right)\, ,
\]
the unique shortest path $\gamma_{x,y}$ is the union of the two segments $[x,z]$ and $[z,y]$.

Given this it is easy to check that $\Gamma$ is locally Lipschitz and that $|D_x \Gamma (t,x,y)| \leq |x|$. It is also easy to check that the length of $\gamma_{x,y}$ is bounded by $C_0 |x-y|$, where $C_0$ depends only upon the opening angle of the cone $V'$, thus giving $|\partial_t \Gamma (t,x,y)|\leq C_0 |x-y|$.

\section{Excluding basic degenerations}\label{s:lower-bound}

The primary goal of this section is to give effective lower bounds on the length of the edges of a once-derived net: while so far we have mostly leveraged the stationarity of the net, the assumption that it comes from the link of an area-minimizing cone as in \cref{t:main} will play here a crucial role. These bounds will be of great help in the computer-assisted search for approximate metrizations of the candidate once-derived topologies (cf. \cref{s:secondo-programma}). We will split the lower bound into three cases, summarized in the following theorem. In order to simplify our statements we will use the notation $\text{len} (\sigma)$ for the length of a geodesic arc $\sigma$ in $\mathbb S^2$.

\begin{proposition}\label{t:effective-short-edge}
Consider the link $\tilde{\mathcal{N}}=E\cap \partial B_1$ of an area-minimizing cone as in \cref{t:main}, let $\mathcal{N}$ be the corresponding once-derived net and $\sigma$ be an edge of $\mathcal{N}$.
\begin{itemize}
    \item[(i)] If $\sigma$ is the side of a convex face of $\mathcal{N}$ or a convex face of $\tilde{\mathcal{N}}$, then $\text{len} (\sigma) \geq \tfrac{1}{100}$.
    \item[(ii)] If $\text{len} (\sigma) < \frac{1}{100}$ then there is an edge $\tau$ of $\mathcal{N}$ adjacent to $\sigma$ such that
\begin{equation*}
        \text{len}(\sigma) \ge \tfrac{1}{10} \text{len}(\tau)\, .
    \end{equation*}
    \item[(iii)] Let $\sigma$ and $\tau$ be the two edges adjacent to a vertex of degree $2$. If both edges are shorter than $\frac{1}{100}$ then there is a third distinct edge $\varsigma$, adjacent to either $\sigma$ or $\tau$, such that
    \begin{equation}\label{e:two-short-edges}
        \text{len}(\sigma) + \text{len}(\tau) \ge \tfrac{1}{10} \text{len}(\varsigma)\, .
    \end{equation}
\end{itemize}
\end{proposition}

We include in this section also another type of lower bound, which will require a more sophisticated analysis to produce a rather tight bound, but which nonetheless is achieved with similar methods. Assume $\tilde{\mathcal{N}}$ is the link of a minimizing cone and $P$ one of the pins. Assume the degree of the pin is strictly less than $3$. Let $F$ be the nonconvex face of $\mathbb S^2\setminus \tilde{\mathcal{N}}$ whose boundary contains $P$. A compactness argument shows that $P$ gets arbitrarily close to an edge of $\partial F$ which is {\em not} adjacent to it, then it in fact gets arbitrarily close to some other vertex. A more precise and effective version of the latter statement is given below in \cref{p:pin-touching} for pins of degree $1$. This proposition will be used in \cref{s:terzo-programma}, when we describe how to discard many stationary links on energetic grounds. 

\begin{proposition}\label{p:pin-touching}
Let $E$ be as in \cref{t:main}, let $\tilde{\mathcal{N}}= E\cap \partial B_1$ be the link of $E$ and $P\in \ell \cap \partial B_1$ a pin of degree $1$ contained in the edge $\sigma$. Let $\mathcal{N}$ be the corresponding once-derived net and $F$ the face of $\mathbb S^2\setminus \mathcal{N}$ which contains the pin. Assume that the distance of $P$ to the closest vertex of $\tilde{\mathcal{N}}$ is at least $0.22$ and, moving from $P$ on the equator $\gamma\supset \sigma$ away from the edge $\sigma$, let $Q$ be the first intersection with $\partial F$.
Assume that the smallest angle formed by $\gamma$ with any edge containing $Q$ is at least $\frac{\pi}{4}$. Then the distance from $P$ to $Q$ is at least $0.01$.
\end{proposition}

\subsection{Minimizing criterion} The crucial ingredient to prove the effective bounds in \cref{t:effective-short-edge} and \cref{p:pin-touching} is the minimizing criterion of \cref{l:homotopy-estimate}, which is shown exhibiting appropriate competitors. We recall that $\# A$ denotes the cardinality of the set $A$. The basic construction is detailed separately in the following lemma because it will be rather useful at a later stage.

\begin{lemma}\label{l:homotopy-construction}
Consider a cone $E$ as in \cref{t:main} and let $\tilde{\mathcal{N}}= E \cap \mathbb S^2$ be its link. Let $\rho:[0,1]\times \tilde{\mathcal{N}}\to \mathbb S^2$ be a Lipschitz homotopy and let $\Omega$ be the smallest closed set such that
\[
\{x: \rho (t,x)\neq x\} \cup \{\rho (t,x): \rho (t,x)\neq x\} \subset \Omega\quad \mbox{for every $t$.}
\]
Assume that
\begin{itemize}
\item[(i)] $\rho(0, x) = x$ for all $x\in \tilde{\mathcal{N}}$;
\item[(ii)] Either there is a wedge $V$ as in \cref{d:wedge} such that $\Omega \subset \mathbb S^2 \setminus V$ or $\Omega \subset \mathbb S^2\setminus \ell$.
\end{itemize}
Then, for every $\delta > 0$ there is $\varepsilon < \delta$ and a Lipschitz map $\Phi: E\cap \overline{B}_1\to \overline{B}_1$ with the following properties:
\begin{itemize}
\item[(a)] $\Phi (x)=x$ for $x\in \tilde{\mathcal{N}}$;
\item[(b)] $\Phi (E \cap \overline{B}_{1-\varepsilon}) = E'\cap \overline{B}_{1-\varepsilon}$, where $E'$ is the cone with base $\rho (1, \tilde{\mathcal{N}})$ and vertex $0$;
\item[(c)] The following inequality holds
\begin{equation}\label{e:competitore-omotopico-1}
\Ha^2(\Phi (E\cap \overline{B}_1\setminus B_{1-\varepsilon})) \leq \int_{\mathbb S^2} \# \rho^{-1}(\{x\}) \,d\Ha^2(x) + \delta \, ;
\end{equation}
\item[(d)] $\Phi$ satisfies the assumptions of \cref{p:Almgren-relaxed} in case (i) and the assumptions of \cref{p:Almgren-relaxed-2} in case (ii).
\end{itemize}
\end{lemma}

\begin{proof}
We let $\varepsilon>0$ be any positive real number and define $\Phi_\varepsilon$ on $E\cap \overline{B}_1$ as follows:
\begin{itemize}
\item $\Phi_\varepsilon (0)=0$;
\item $\Phi_\varepsilon (x) = |x| \rho (1, \frac{x}{|x|})$ for every $x\in E$ with $0<|x|< 1-\varepsilon$;
\item $\Phi_\varepsilon (x) = |x| \rho (\frac{1-|x|}{\varepsilon}, \frac{x}{|x|})$ for every $x\in E$ with $1-\varepsilon \leq |x| \leq 1$.
\end{itemize}
This map clearly is Lipschitz and satisfies (a), (b) and (d). We therefore only need to show that
\begin{equation}\label{eq:local667}
\limsup_{\varepsilon \to 0} \mathcal{H}^2 (\Phi_\varepsilon (E) \cap (\overline{B}_1\setminus B_{1-\varepsilon}))
\leq \int_{\mathbb S^2} \# \rho^{-1} (\{x\})\, d\mathcal{H}^2 (x)\, .
\end{equation}
Let $A_\eps = E \cap (\overline{B}_1\setminus B_{1-\eps})$. Parameterize $x \in A_\eps$ by $(r, \omega) \in [1-\eps, 1] \times \tilde{\mathcal{N}}$ with $x = r\omega$ and $d\mathcal{H}^2 = r\,dr\,d\mathcal{H}^1(\omega)$.
On $A_\varepsilon$, the map is $\Phi_\varepsilon(r, \omega) = r \rho(\frac{1-r}{\varepsilon}, \omega)$.

By the area inequality for Lipschitz maps (and the change of variable $t = \frac{1-r}{\varepsilon}$), we bound the area by integrating the Jacobian $J\Phi_\varepsilon$:
\[
\mathcal{H}^2 (\Phi_\varepsilon (A_\varepsilon)) \leq \int_{\tilde{\mathcal{N}}} \int_{1-\varepsilon}^1 J\Phi_\varepsilon(r, \omega) \,r\,dr\,d\mathcal{H}^1(\omega)
=
o(1)+\int_{\tilde{\mathcal N}} \int_0^1 J\rho(t, \omega)\,dt \,dt\mathcal H^1(\omega)
\,.
\]
Passing to the limit and applying the area formula to $\rho$ yields \cref{eq:local667}.
\end{proof}

\begin{corollary}\label{l:homotopy-estimate}
Consider a cone $E$ as in \cref{t:main} and let $\tilde{\mathcal{N}}= E \cap \mathbb S^2$. Let $\rho$ be as in \cref{l:homotopy-construction}. Then,
\begin{equation}\label{e:competitore-omotopico}
\Ha^1(\tilde{\mathcal{N}})-\Ha^1(\rho(1, \tilde{\mathcal{N}})) \leq 2 \int_{\mathbb S^2} \# \rho^{-1}(\{x\}) \,d\Ha^2(x)\, .
\end{equation}
\end{corollary}
\begin{proof}
Fix any $\delta$ and consider the map $\Phi_\varepsilon$ and the number $\varepsilon = \varepsilon (\delta) < \delta$ given by \cref{l:homotopy-construction}. Since we can apply either \cref{p:Almgren-relaxed} or \cref{p:Almgren-relaxed-2}, we conclude that
\begin{equation}\label{e:competitore-omotopico-2}
\tfrac{1}{2} \mathcal{H}^1 (\tilde{\mathcal{N}}) = \mathcal{H}^2 (E\cap \overline{B}_1)
\leq \mathcal{H}^2 (\Phi (E\cap \overline{B}_1)\, .
\end{equation}
On the other hand, by \eqref{e:competitore-omotopico-1} we get
\begin{align*}
\mathcal{H}^2 (\Phi (E\cap \overline{B}_1)) & = \mathcal{H}^2 (\Phi (E) \cap B_{1-\varepsilon})
+ \mathcal{H}^2 (\Phi (E) \cap \overline{B}_1\setminus \overline{B}_{1-\varepsilon})\\
&\leq \tfrac{(1-\varepsilon)^2}{2} \mathcal{H}^1 (\rho (1, \tilde{\mathcal{N}})) +
\int_{\mathbb S^2} \# \rho^{-1} (\{x\})\, d\mathcal{H}^2 (x) +\delta\, .
\end{align*}
Combining the latter inequality with \eqref{e:competitore-omotopico-2} we then get
\[
 \mathcal{H}^1 (\tilde{\mathcal{N}}) - (1-\varepsilon)^2 \mathcal{H}^1 (\rho (1, \tilde{\mathcal{N}}))
 \leq 2 \int_{\mathbb S^2} \# \rho^{-1} (\{x\})\, d\mathcal{H}^2 (x) + 2\delta\, .
\]
Recalling that $\varepsilon<\delta$ and letting $\delta\downarrow 0$ we achieve \eqref{e:competitore-omotopico}
\end{proof}

\subsection{Suitable homotopies} From the minimizing criterion in \cref{l:homotopy-estimate} and through a suitable choice of appropriate homotopies we will gain the following two ``lower level estimates'', from which in turn we will get \cref{t:effective-short-edge}.

\begin{proposition}\label{prop:effective-short-edge-1}
    Consider the link $\tilde{\mathcal{N}} = E\cap \partial B_1$ of an area-minimizing cone $E$ as in \cref{t:main}, let $\alpha,\beta,\gamma$ be three edges of $\tilde{\mathcal{N}}$ which form a path in this order, and denote by $A$ and $B$ the vertices of $\beta$ in common, respectively, with $\alpha$ and $\gamma$. Assume that:
    \begin{itemize}
        \item[(i)] $\alpha$, $\beta$, and $\gamma$ are in the boundary of a face $F$ of $\mathbb S^2\setminus \tilde{\mathcal{N}}$;
        \item[(ii)] the internal angles of $F$ at $A$ and $B$ measure $\frac{2\pi}{3}$;
       \item[(iii)] let $H$ be the open hemisphere with $\beta\subset \partial H$ containing a neighborhood of $\beta$ in $F$; either $F$ is convex or the eventual pin\footnote{$H$ might contain no pin because both pins might lie in its boundary: in this case the assumption is considered as verified.} contained in $H$ is at distance at least $\frac{\pi}{4}$ from $A$.
    \end{itemize}
    Then
    \begin{equation*}
        \text{len}(\beta) \ge \tfrac{1}{10} \min\{\tfrac{1}{10}, \text{len}(\alpha), \text{len}(\gamma)\}.
    \end{equation*}
\end{proposition}

\begin{proposition}\label{prop:effective-short-edge-2}
 Consider the link $\tilde{\mathcal{N}} = E\cap \partial B_1$ of an area-minimizing cone $E$ as in \cref{t:main}, let $\alpha,\beta,\gamma,\delta$ be four edges of $\tilde{\mathcal{N}}$ which form a path in this order, and denote by $A$, $B$, and $C$ the vertices in common between, respectively, $\alpha$ and $\beta$, $\beta$ and $\gamma$, and $\gamma$ and $\delta$. Assume that:
    \begin{itemize}
        \item[(i)] $\alpha$, $\beta$, $\gamma$, and $\delta$ are in the boundary of a face $F$ of $\mathbb S^2\setminus \tilde{\mathcal{N}}$;
        \item[(ii)] the internal angles of $F$ at $A$, $B$, and $C$ measure $\frac{2\pi}{3}$.
        \item[(iii)] denote by $H_1$ and $H_2$ the open hemispheres with $\beta\subset \partial H_1$, $\gamma \subset \partial H_2$ and containing, respectively a neighborhood of $\beta$ and a neighborhood of $\gamma$ in $F$; either $F$ is convex or the eventual pin contained in $H_1\cap H_2$ is at distance at least $\frac{\pi}{4}$ from $B$.
    \end{itemize}
    Then
    \begin{equation*}
        \text{len}(\beta) + \text{len}(\gamma) \ge \tfrac{1}{10} \min\{\tfrac{1}{5}, \text{len}(\alpha), \text{len}(\delta)\}\, .
    \end{equation*}
\end{proposition}

\begin{figure}[htbp]
    \begin{tikzpicture}[scale=0.5]
        \pgfmathsetmacro{\L}{7}
        \pgfmathsetmacro{\S}{8}
        \coordinate (A) at (-1,0);
        \coordinate (B) at (1, 0);
        \coordinate (A2) at ({-1-\L/2}, {\L*sqrt(3)/2});
        \coordinate (B2) at ({1+\L/2}, {\L*sqrt(3)/2});
        \coordinate (D) at (-1,{\L*sqrt(3)/3});
        \coordinate (C) at (1,{\L*sqrt(3)/3});
        \node[below left] at (A) {$A$};
        \node[below right] at (B) {$B$};
        \node[left] at (A2) {$A_2$};
        \node[right] at (B2) {$B_2$};
        \fill[gray!20] (A) -- (A2) -- (D) to[bend right=20] (C) -- (B2) -- (B) -- (A) -- cycle;
        \draw[very thick] (D) to[bend right=20] (C);
        \draw[very thick] (A) -- (D) -- (A2);
        \draw[very thick] (C) -- (B2);
        \draw[very thick,dashed] (A2) --node[below left]{$\alpha$} (A);
        \draw[very thin] (A) -- (B);
        \draw[very thick,dashed] (B) --node[below right]{$\gamma$} (B2);
        \draw[thin, dashed] (A2) -- (B2);
        \draw[thin, dashed] (B) -- (C);
        \draw[very thin] (A2) -- ({-1-\S/2}, {\S*sqrt(3)/2});
        \draw[very thin] (B2) -- ({1+\S/2}, {\S*sqrt(3)/2});
        \node[below left] at (D) {$D$};
        \node[below right] at (C) {$C$};
        \node[above] at (0,{7/2+0.4}) {$\eta$};
        \node[below] (0,0) {$\beta$};
    \end{tikzpicture}
\caption{The homotopy used to prove \cref{prop:effective-short-edge-1}. All the straight segments drawn in the pictures are length-minimizing geodesics, while the arc $\eta$ is made of points equidistant from the edge $\beta$.}\label{f:competitor-1}
\end{figure}
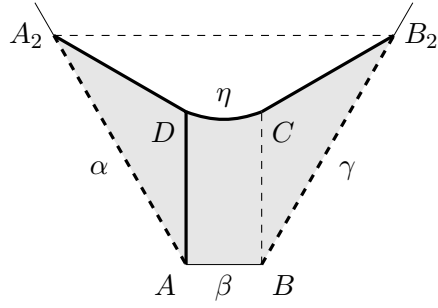

\begin{proof}[Proof of \cref{prop:effective-short-edge-1}]
Assume by contradiction that the edge $\beta$ is strictly shorter than $\frac{1}{100}$, $\frac{\text{len} (\alpha)}{10}$, and $\frac{\text{len} (\gamma)}{10}$. We then mark on $\alpha$ and $\gamma$ the two points $A_2$ and $B_2$ with the property that the lengths of the two segments $A_2A$ and $B_2B$ are $10 \text{len} (\beta)$. Furthermore we consider the isosceles triangles $AA_2D$ and $BB_2C$ with angles $\frac{\pi}{6}$ at $A$, $A_2$, $B$ and $B_2$ with the property that the edges $AD$ and $BC$ are orthogonal to $AB$, cf. \cref{f:competitor-1}. Finally we connect $C$ and $D$ with a special arc $\eta$ which is not geodesic but which rather consists of the points which are equidistant to the edge $AB$. 

We now describe an homotopy $\rho$ to which we wish to apply \cref{l:homotopy-estimate}. The homotopy keeps all points of the net $\tilde{\mathcal{N}}$ fixed except for those in the open arcs $AA_2$ and $BB_2$. The segment $AA_2$ is deformed continuously onto the union of $AD$ and $DA_2$; while $BB_2$ is deformed onto the union of $BA$, $AD$, $\eta$, $CB_2$.
The homotopy is defined in such a way that all the points in the interior of the triangles $AA_2C$, $BB_2C$, and the region $Q$ bounded by $DA\cup AB \cup BC \cup \eta$ have exactly one preimage, while the points which do not belong to their closures or to $\tilde{\mathcal{N}}$ have no preimages.
We also draw the shortest geodesic segment joining $A_2$ and $B_2$ and introduce the quadrilateral $\Omega$, which contains $AA_2D\cup BB_2C \cup Q$.
We claim that
\begin{itemize}
\item[(a)] the following inequality holds
\begin{equation}\label{e:disuguaglianza-area-lunghezza}
|AA_2| + |BB_2| > |A_2D|+|DA|+\text{len} (\eta) + |B_2C| + 2 (\mathcal{H}^2 (AA_2D) + \mathcal{H}^2 (BB_2C)+
\mathcal{H}^2 (Q))\, .
\end{equation}
\end{itemize}
The claims (a) will be proved in \cref{a:homotopy-computations}.

We need however to check that the homotopy satisfies condition (ii) in \cref{l:homotopy-construction} in order to apply the \cref{l:homotopy-estimate}.
If $\Omega$ does not contain any pin, then clearly condition (ii) in \cref{l:homotopy-construction} is satisfied and the proof follows after observing that
\[
\mathcal{H}^1 (\tilde{\mathcal{N}}) - \mathcal{H}^1 (\rho (1, \tilde{\mathcal{N}})) = |AA_2| + |BB_2| - (|A_2D|+|DA|+\text{len} (\eta) + |B_2C|)
\,
\]
while
\[
\int_{\mathbb S^2} \# (\rho^{-1} (x))\, d\mathcal{H}^2 (x) = (\mathcal{H}^2 (AA_2D) + \mathcal{H}^2 (BB_2C)+
\mathcal{H}^2 (Q)) 2 \mathcal{H}^2 (AA_2QM)\, ,
\]
in particular \eqref{e:disuguaglianza-area-lunghezza} would be at odds with \eqref{e:competitore-omotopico}.

We next show that if $\Omega$ contains a pin $P$, then it is either $A$ or $B$ and thus it follows immediately that condition (ii) in \cref{l:homotopy-construction} holds. In order to prove that $P\in \{A,B\}$ note first that, let $F$ be the face in the statement of the proposition. If $F$ is convex, the interior $\omega$ of the segment $A_2B_2$ (which is shorter than $\pi$ and hence length-minimizing) is contained in the interior of $F$. Therefore also the interior of $\Omega$ is contained in the interior of $F$. Since the interior of $F$ does not contain any pin, then $P\in \partial \Omega\setminus \omega \subset \tilde{\mathcal{N}}$. But $P$ is a vertex of $\tilde{\mathcal{N}}$ and the only two vertices of the net contained in $\partial \Omega\setminus \omega$ are $A$ and $B$.
If $F$ is not convex, consider the open hemisphere $H$ of condition (iii) in the statement. Then $\Omega\setminus \beta\subset H$ and, since it can be easily checked that $\Omega$ has diameter strictly smaller than $\frac{\pi}{4}$, we immediately infer $P\in \beta$. As above, since $\beta$ is an edge of $\tilde{\mathcal{N}}$, $P$ must be either $A$ or $B$.
\end{proof}

\begin{proof}[Proof of \cref{prop:effective-short-edge-2}] We argue by contradiction and assume that
\[
\text{len} (\beta) + \text{len} (\gamma) < \min\{ \tfrac{1}{50}, \tfrac{\text{len} (\alpha)}{10}, \tfrac{\text{len} (\delta)}{10}\}\, .
\]
Without loss of generality, we also assume $\text{len} (\beta) \leq \text{len} (\gamma)$.
We then draw the equator $\lambda$ which bisects the internal angle of $F$ at the vertex $B$.
Let $M$ be the point on $\lambda$, in the direction of the face $F$, at distance $5 (\text{len} (\beta) + \text{len} (\gamma))$ from $B$.
Then, let $\bar{\lambda}$ be the equator perpendicular to $\lambda$ at $M$. In \cref{a:homotopy-computations}, we will use elementary geometric arguments to show that:
\begin{itemize}
\item[(a)] $\bar\lambda$ intersects the edge $\alpha$ at an interior point $E$ and the edge $\delta$ at an interior point $D$;
\item[(b)] The pentagon $ABCDE$ is convex;
\item[(c)] $|CD| > |DE| + 2\mathcal{H}^2 (ABCDE)$;
\end{itemize}
cf, \cref{f:competitor-2}.

Now, we describe the homotopy $\rho$ to be used in \cref{l:homotopy-estimate}. First of all $\rho$ moves the points in the the interior of the arc $CD$, while the rest of the net $\tilde{\mathcal{N}}$ is kept fixed. The arc $CD$ is mapped by $\rho (1, \cdot)$ onto the union of the arcs $\gamma$, $\beta$, $AE$, and $DE$. Moreover the homotopy is defined in such a way that each point in the open pentagon $ABCDE$ has exactly one preimage and the points outside the union of its closure and $\tilde{\mathcal{N}}$ have no preimages.  Arguing as in the proof of \cref{prop:effective-short-edge-1} we identify the closed pentagon $ABCDE$ as the region $\Omega$ which satisfies condition (ii) of \cref{l:homotopy-construction}. In particular we claim that, if a pin $P$ belongs to $\Omega$, then $P\in \{A,B,C\}$. In order to show this, we consider first the case in which the face $F$ is convex. Then, since $\Omega$ is convex, its interior and the interior $\omega$ of the segment $DE$ are both contained in the interior of $F$. So $P\in \partial \Omega \setminus \omega \subset \tilde{\mathcal{N}}$ and, because $P$ is a vertex of $\tilde{\mathcal{N}}$ and the only vertices of $\tilde{\mathcal{N}}$ contained in $\partial \Omega \setminus \omega$ are $A$, $B$, and $C$, we must necessarily have $P\in \{A,B,C\}$. If $F$ is not convex, consider however that, again by convexity of the interior of $\Omega$, $\Omega\setminus (\beta \cup \gamma) \subset H_1 \cap H_2$, where $H_1$ and $H_2$ are the hemispheres of assumption (iii) in the proposition. But because we easily see that the diameter of $\Omega$ is strictly smaller than $\frac{\pi}{4}$, this implies that $P\in \beta \cup \gamma$. Again, since $P$ is a vertex of the net $\tilde{\mathcal{N}}$ and the only vertices contained in $\beta\cup \gamma$ are $A,B,$ and $C$, we infer that $P$ is one of them.

    \begin{figure}[htbp]
        \begin{tikzpicture}[scale=3]
            \coordinate (A) at (0,0);
            \coordinate (B) at (-0.4/2, {0.4*sqrt(3)/2});
            \coordinate (C) at (0.1, {sqrt(3)/2});
            \coordinate (M) at (1, {0.4*sqrt(3)/2});
            \coordinate (D) at (1, {sqrt(3)/2});
            \coordinate (E) at (1, 0);
            \node[below left] at (A) {$A$};
            \node[above left] at (C) {$C$};
            \node[left] at (B) {$B$};
            \node[above right] at (M) {$M$};
            \node[above] at (1.8, {0.4*sqrt(3)/2}) {$\lambda$};
            \fill[gray!20] (A) -- (B) -- (C) -- (D) -- (E) -- cycle;
            \draw[very thin] (E) -- (A) --node[below left] {$\beta$} (B) --node[above left] {$\gamma$} (C);
            \draw[very thick, dashed] (C) -- (D);
            \draw[very thick] (D) -- (E);
            \node[below] at (E) {$E$};
            \node[above] at (D) {$D$};
            \draw[very thin, dashed] (B) -- (2, {0.4*sqrt(3)/2});
        \end{tikzpicture}
        \caption{The competitor used to prove \cref{prop:effective-short-edge-2}. The thick dashed segment is removed from the net and replaced with the thick segment. The gray area is the region spanned by the homotopy from the original net to the competitor.}\label{f:competitor-2}
    \end{figure}
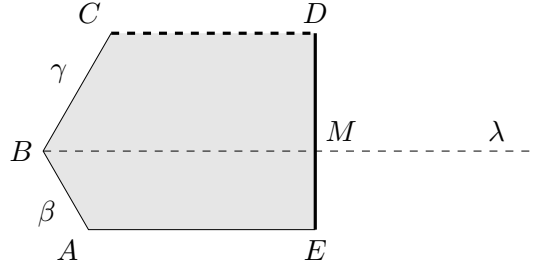

Clearly
\[
\mathcal{H}^1 (\tilde{\mathcal{N}}) - \mathcal{H}^1 (\rho (1, \tilde{\mathcal{N}})) =
|CD| - |ED|\,
\]
while
\[
\int_{\mathbb S^2} \# (\rho^{-1} (x))\, d\mathcal{H}^2 (x) = \mathcal{H}^2 (ABCDE)\, .
\]
But then (c) implies that the inequality \eqref{e:competitore-omotopico} is violated by this explicit $\rho$.
\end{proof}

\subsection{Proof of \cref{t:effective-short-edge}} {\bf Proof of (i).} Fix the convex face $F$ of the assumption. If the face is a biangle or a triangle, the lengths are rigid and the statement holds by \cref{p:formulae}. Now, assume the face $F$ is either a convex quadrilateral or a convex pentagon. First of all observe that, if $F$ is a face of $\mathcal{N}$, then (since it is convex), it is also a face of $\tilde{\mathcal{N}}$. Let $a$ be the length of $\sigma$ and $b,c$ be the lengths of the adjacent sides of $F$. Thanks to \cref{lem:quadrilateral_consecutive_lowerbound,lem:pentagon_consecutive_lowerbound}, we have $a+b, a+c\ge 1$. So, either $a\ge 0.01$ or $b,c\ge 0.99$. In the first case, we have proven the statement. Otherwise, we can apply \cref{prop:effective-short-edge-1} and deduce $a\ge 0.1\min\{0.1, b, c\}\ge 0.01$.

\medskip

{\bf Proof of (ii).} We can assume that
\begin{equation}\label{e:piccolo}
\text{len} (\sigma) < 0.01\, ,
\end{equation}
Let $A$, $B$ be the two vertices of $\sigma$. In $\tilde{\mathcal N}$, both $A$ and $B$ have degree $3$. None of them can be a pin, because otherwise it would be a pin of degree $3$ and hence $\sigma$ would be the edge of a convex face in $\mathbb S^2\setminus \tilde{\mathcal{N}}$: in this case we could apply statement (i) of the Proposition, which is at odds with \eqref{e:piccolo}. In particular:
\begin{itemize}
\item[(p)] if an edge $\tau$ of $\tilde{\mathcal{N}}$ adjacent to $A$ or $B$ does not belong to $\mathcal{N}$, then the other endpoint of $\tau$ is a pin.
\end{itemize}
Let $H$ be an open hemisphere bounded the equator containing $\sigma$ which either does not contain a pin or such that the only pin $P$ in $H$ is at distance at least $\tfrac\pi 4$ from $A$.
Let $\alpha$ be the edge adjacent to $A$ and $\gamma$ the edge adjacent to $B$ which enter into $H$. Observe that if $\alpha$, resp. $\gamma$, does not belong to $\mathcal N$, then, by (p), its length is at least $\tfrac\pi4-0.01$.
Applying \cref{prop:effective-short-edge-1} we obtain
\begin{equation*}
    \text{len}(\sigma) \ge \frac1{10} \min\{\tfrac1{10}, \text{len}(\alpha), \text{len}(\gamma)\} \,.
\end{equation*}
Since $\text{len}(\sigma)<0.01$, it must hold $\min\{\text{len}(\alpha), \text{len}(\gamma)\} < \tfrac1{10}$. Without loss of generality, we assume that $\alpha$ is the shortest between $\alpha$ and $\beta$. So $\text{len}(\sigma)\ge \tfrac1{10}\text{len}(\alpha)$ and $\text{len}(\alpha) < \tfrac1{10}$. Thus, as observed above, the edge $\alpha$ belongs to $\mathcal N$ and this concludes the proof.

\medskip

{\bf Proof of (iii)}. We can, without loss of generality, assume
\begin{equation}\label{e:piccolo2}
\max \{\text{len} (\sigma), \text{len} (\tau)\} < 0.01\, .
\end{equation}
Let us label as $B$ the only vertex in common between $\sigma$ and $\tau$. Observe that:
\begin{itemize}
    \item the other endpoint $A$ of $\sigma$ must be distinct from the other endpoint $C$ of $\tau$, otherwise both $\sigma$ and $\tau$ would be long $\pi$;
    \item Neither $A$, nor $B$, nor $C$ can be a pin, because it would be a pin of degree $3$ and thus either $\sigma$ or $\tau$ would belong to a convex face of $\mathbb S^2\setminus \tilde{\mathcal{N}}$, which, by point (i) of the Proposition, is at odds with \eqref{e:piccolo2}.
\end{itemize}
Recall next that $B$ has degree $2$ in $\mathcal{N}$. Thus there is a third edge $\sigma_B$ adjacent to $B$ in $\tilde{\mathcal{N}}$, which has been removed because it connects $B$ to a pin $P$. Clearly $\text{len} (\sigma_B) < \pi$, otherwise $\sigma_B$ would contain $P$ and its antipodal $P'$, which we know is the other pin of the net.

\begin{figure}[htbp]
        \begin{tikzpicture}[scale=1.2]
        \draw[very thick] ({-sqrt(3)+7},3) --node[left] {$\sigma'_A$} ({-sqrt(3)+7},1);
        \draw[very thick] ({-2*sqrt(3)+7},0) --node[above left] {$\sigma_A$} ({-sqrt(3)+7},1);
        \draw ({-sqrt(3)+7},1) --node[below left] {$\sigma$} (7,0) --node[below right] {$\tau$} ({sqrt(3)+7},1);
        \draw[very thick] ({sqrt(3)+7},1) --node[right] {$\sigma'_C$} ({sqrt(3)+7},3);
        \draw[very thick] ({sqrt(3)+7},1) --node[right] {$\sigma_C$} ({2*sqrt(3)+7},0);
        \node[above] at (7,0) {$B$};
        \node[above right] at ({-sqrt(3)+7},1) {$A$};
        \node[above left] at ({sqrt(3)+7},1) {$C$};
        \draw[very thick, dashed] (7,0) --node[left] {$\sigma_B$} (7,-2);
        \node at ({-sqrt(3)+7},-1) {$G_A$};
        \node at ({sqrt(3)+7},-1) {$G_C$};
        \node at (7,2) {$G$};
        \node[below] at (7,-2) {$P$};
        \end{tikzpicture}
        \caption{The edges and faces used in the proof of \cref{t:effective-short-edge}(iii). The edges $\sigma$ and $\tau$ belong to the once-derived net $\mathcal{N}$. The dashed edge $\sigma_B$ belongs to $\tilde{\mathcal{N}}$ but certainly does not belong to $\mathcal{N}$. All the other (thick but not dashed) edges belong to $\tilde{\mathcal{N}}$ and might or might not be removed by the procedure to derive $\mathcal{N}$.}
        \label{f:short-e-iii}
\end{figure}
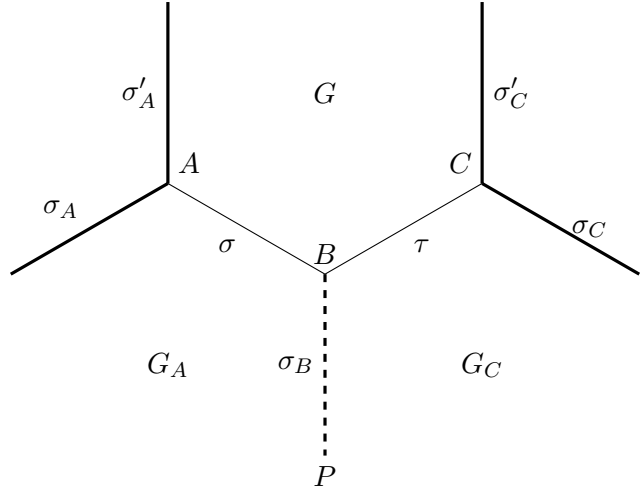

We next label with $\sigma_A$ and $\sigma'_A$ the other two edges of $\tilde{\mathcal{N}}$ adjacent to $A$ and with $\sigma_C$ and $\sigma'_C$ the other two edges adjacent to $C$, cf. \cref{f:short-e-iii}. The labels are placed so that there are three faces $G_A$, $G_C$ and $G$, of $\mathbb S^2\setminus \tilde{\mathcal{N}}$ such that:
\begin{align*}
&\sigma_A \cup \sigma\cup \sigma_B \subset \partial G_A\\
&\sigma_B \cup \tau \cap \sigma_C \subset \partial G_C\\
&\sigma'_A\cup \sigma \cap \tau\cup\sigma'_C \subset \partial G\, ,
\end{align*}
cf. again \cref{f:short-e-iii}.
Note finally that:
\begin{itemize}
\item[(p')] If $\sigma'_A$, resp. $\sigma'_C$, does not belong to $\mathcal{N}$, then it connects $A$, resp. $C$, to the pin $P'$.
\end{itemize}
To complete the proof we analyze two possibilities.

\medskip
\noindent
{\bf Alternative 1: $\text{len} (\sigma_B) \geq \frac{\pi}{4}$}. At least one between $\sigma_A$ and $\sigma_C$ must belong to $\mathcal{N}$, otherwise $\mathcal{N}$ would have a face with three concave angles. Without loss of generality, assume $\sigma_A$ is an edge of $\mathcal N$.
We can apply \cref{prop:effective-short-edge-1} to the three edges $\sigma_B$, $\sigma$, $\sigma_A$ and infer
\[
0.01 > \text{len} (\sigma) \geq \tfrac{1}{10} \min \{\tfrac{1}{10}, \text{len} (\sigma_A), \text{len} (\sigma_B)\}\, .
\]
In particular, $\min \{\tfrac{1}{10}, \text{len} (\sigma_A), \text{len} (\sigma_B)\} < \frac{1}{10}$, and, since $\text{len} (\sigma_B)\geq \frac{\pi}{4}$, we conclude
\[
\min \{\tfrac{1}{10}, \text{len} (\sigma_A), \text{len} (\sigma_B)\} = \text{len} (\sigma_A)\, .
\]
Hence \eqref{e:two-short-edges} holds if we set $\varsigma=\sigma_A$.

\medskip
\noindent
{\bf Alternative 2. $\text{len} (\sigma_B) < \frac{\pi}{4}$}. In this case we have $\dist\, (B, P') > \frac{3\pi}{4}$. Let now $\gamma_A$ be the equator containing $\sigma$ and $H_1$ the hemisphere bounded by $\gamma_A$ in which the segment $\sigma'_A$ enters. Clearly $P$ belongs to the other hemisphere bounded by $\gamma_A$, in particular its antipodal $P'$ is in $H_1$. If we denote by $\gamma_C$ the equator containing $\tau$ and by $H_2$ the hemisphere bounded by $\gamma_C$ in which $\sigma'_C$ enters, we conclude likewise that $P'\in H_2$. We then see that we are in the position to apply \cref{prop:effective-short-edge-2}, where the hemispheres of condition (iii) therein are precisely the ones just identified. We then conclude
\[
0.02 > \text{len} (\sigma) + \text{len} (\tau) \geq \tfrac{1}{10} \{\tfrac{1}{5}, \text{len} (\sigma'_A), \text{len} (\sigma'_C)\}\, .
\]
If we let $\varsigma$ be the shortest segment between $\sigma'_A$ and $\sigma'_C$, we then see that it is necessarily shorter than $\frac{1}{5}$, which means that it must belong to $\mathcal{N}$. In fact (p') implies that, if it were removed, then it would connect either $A$ or $C$ with the pin $P'$, and hence it would be longer than $\frac{3\pi}{4}-0.01$. Hence we conclude the validity of \eqref{e:two-short-edges} in this case as well.

\subsection{Proof of \cref{p:pin-touching}} As in the cases analyzed in \cref{t:effective-short-edge}, we will argue by contradiction showing the existence of an admissible competitor with lower area. In order to do so we start by fixing a system of coordinates (labeled $(x,y,z)$ as usual) in which the point $Q$ is the ``north pole'' $Q= (0,0,1)$ and the equator containing the edge $\sigma$ is the plane $\{y=0\}$. We then intersect $E$ with the plane $\wp = \{z=1\}$: this intersection is a planar net, i.e. a closed set consisting, locally, of a finite number of straight segments connecting a finite number of vertices (the set is not compact, but this plays no role in our argument, which is local around the point $Q$). In fact we can establish a natural map from $E\cap \wp$ to $\tilde{\mathcal{N}}$ by mapping $x$ into $\frac{x}{|x|}$: this maps vertices of $E\cap \wp$ into vertices of $\tilde{\mathcal{N}}$ and, locally, edges into edges. The map has Lipschitz constant $1$ if we endow $\wp$ with the Euclidean distance and $\mathbb S^2$ with the geodesic distance.
The intersection of $\ell$ with $\wp$ consists then of the one point, which we label by $M$, such that $P=\frac{M}{|M|}$. If we let $\bar{\delta}$ be the distance (on the sphere) between the pin $P$ and $Q$ and we assume it is shorter than $0.01$, then the length $\delta$ of the segment $MQ$ is $\tan \bar{\delta}$ at we therefore conclude that it is smaller than $0.011$.

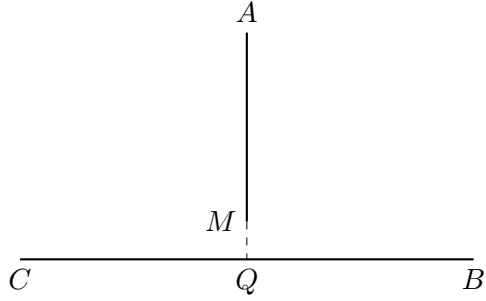
\begin{figure}[hbtp]
\begin{tikzpicture}
\coordinate (B) at (3,0);
\coordinate (C) at (-3,0);
\coordinate (A) at (0,3);
\coordinate (M) at (0,0.5);
\draw[thick] (B) -- (C);
\draw[thick] (A) -- (M);
\draw[very thin, dashed] (M) -- (0,0);
\node[below] at (B) {$B$};
\node[below] at (C) {$C$};
\node[above] at (A) {$A$};
\node[left] at (M) {$M$};
\node[below] at (0,0) {$Q$};
\end{tikzpicture}
\caption{The set $E\cap \wp$ in a neighborhood of $M=\ell \cap \wp$. In this case $\alpha = \frac{\pi}{2}$.}\label{f:9.2-figura-1}
\end{figure}

We next consider the plane containing the equator $\gamma$ and let $\scrs$ be its intersection with $\wp$ and denote by $\scrh$ the intersection of the plane $\{y=0\}$ with $\wp$. Of course $\scrs\cap \scrh = \{Q\}$ and we let $\alpha$ be the (smallest) angle formed by $\scrs$ and $\scrh$. Note that this is the angle formed by the equator $\gamma$ and the edge $\sigma$ at $Q$, therefore $\frac{\pi}{4} \leq \alpha \leq \frac{\pi}{2}$. In particular, after possibly applying a reflection along the plane $\{x=0\}$, we conclude that the point $M$ has coordinates
\[
M = (\delta \cos \alpha, \delta \sin \alpha, 1)\, .
\]
Consider now any vertex $V$ in $E\cap \wp$: having established that the map $V\mapsto \frac{V}{|V|}$ is $1$-Lipschitz and that it maps vertices into vertices, we conclude that the segment $AM$ is at least $0.22$, in particular the segment $AQ$ is at least $0.22-0.011 = 0.209 > 0.2 =: L$. We then mark as $A$, $B$, and $C$ the points
\begin{align}
A &= (L \cos \alpha, L \sin \alpha, 1)\\
B &= (L,0,1)\\
C &= (-L,0,1)\, .
\end{align}
We remark that:
\begin{itemize}
\item The segments $AM$ and $BC$ are contained in $E\cap \wp$;
\item Indeed $\overline{B}_{0.2} (M)\cap E \cap \wp = AM \cup BC$, cf. \cref{f:9.2-figura-1}.
\end{itemize}
Consider now the cone
\[
\mathcal{C} := \{ rx: x\in \overline{B}_{0.2} (M)\cap \wp \quad \mbox{and} \quad 0\leq r\leq 1\}\, .
\]
The set $\Sigma_0:= E\cap \mathcal{C}$ consists of the two triangles $BOC$ and $MOA$, where $O$ denotes the origin: this is the cone with vertex $0$ over the set $\overline{B}_{0.2} (M)\cap E \cap \wp = AM \cup BC$. We construct a competitor $\Sigma$ which consists of $10$ suitable triangles contained in $\mathcal{C}$. We will show that:
\begin{itemize}
\item[(a)] $\mathcal{H}^2 (\Sigma_0) - \mathcal{H}^2 (\Sigma) > 0$;
\item[(b)] while clearly $\mathcal{E} := (E \setminus \mathcal{C}) \cup \Sigma$ will not be an admissible competitor for $E$ in the sense of \cref{d:Almgren-minimizing}, there is a Lipschitz map $F: E\to \mathbb R^3$ which maps $E$ into $\mathcal{E}$ and satisfies the assumptions of \cref{p:Almgren-relaxed-2}.
\end{itemize}
Therefore (a) and (b) contradict the minimizing property of $E$ and hence will complete the proof.

\begin{figure}[hbtp]
\begin{tikzpicture}
\coordinate (B) at (3,0);
\coordinate (C) at (-3,0);
\coordinate (A) at (0,3);
\coordinate (M) at (0,0.2);
\coordinate (D) at ({(3-sqrt(3))/2}, {(3-sqrt(3))/2});
\node[above right] at (D) {$D$};
\draw[thick] (A) -- (D) -- (M);
\draw[thick] (B) -- (D);
\node[below] at (B) {$B$};
\node[below] at (C) {$C$};
\node[above] at (A) {$A$};
\node[above] at (M) {$M$};
\draw[thick] (C) -- (M);
\end{tikzpicture}
\caption{This picture shows the basis of the cone $\Sigma'$ which would be an ideal competitor for $\Sigma_0$. $\Sigma$ contains a scaled copy $r\Sigma'$ of the latter, while the remaining part of $\Sigma$ interpolates between $\Sigma_0\cap \wp$ and the basis of $r\Sigma'$.}\label{f:9.2-figura-2}
\end{figure}
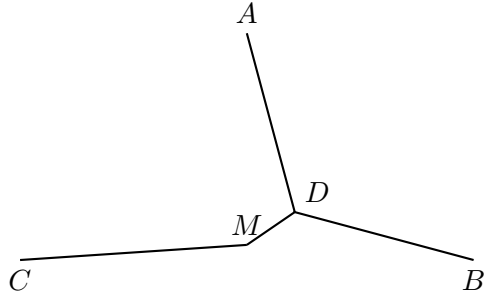

{\bf Description of $\Sigma$.}
We start by describing the triangles forming $\Sigma$. In order to do so we need to introduce:
\begin{itemize}
\item The Torricelli point $D$ of $A$, $B$, and $Q$. This is the unique point $D$ of the plane $\wp$ with the property that $AD$, $BD$, and $QD$ meet at $D$ forming equal angles;
\item The points $A'=r A$, $B' = rB$, $C'= rC$, and $D'=rD$, where $r=0.8$.
\end{itemize}
The surface $\Sigma$ consists then of:
\begin{itemize}
\item The triangles $BCM$ and $CMO$;
\item The cone with basis the net $A'D'\cup B'D'\cup M'D'$, which is formed by the three triangles $A'OD'$, $B'OD'$, and $M'OD'$;
\item The ``transition'' surface formed by the five triangles $AA'D'$, $AMD'$, $MM'D'$, $BMD'$, and $BB'D'$.
\end{itemize}
The geometric idea behind the competitor is that we would like to substitute $\Sigma_0$ in $\mathcal{C}$ with a new cone $\Sigma'$ with vertex $0$ and with base given in the \cref{f:9.2-figura-2}. But transitioning between the two cones directly on the plane $\wp$ would be energetically costly. We therefore insert rather a scaled copy of the cone $r\Sigma'$ ($r=0.8$ is the scaling factor) and interpolate at the intermediate scales. We can afford one horizontal interpolation in the triangle $BCM$ because it is very small, the remaining interpolating surface takes advantage also of the direction orthogonal to $\wp$.

{\bf Proof of (a).} We note that the quantity $\mathcal{H}^2 (\Sigma_0) - \mathcal{H}^2 (\Sigma)$ depends in fact on $x$ and $\alpha$ (indeed it can even be computed through an explicit formula). This gives a function $(\delta,\alpha)\mapsto f (\delta, \alpha)$, defined on the domain $(0,0.011]\cap [\frac{\pi}{4}, \frac{\pi}{2}]$. We will show in \cref{a:homotopy-computations} that $f$ is Lipschitz, but in fact it can be immediately seen, either geometrically or from the explicit formula, that it can be extended to the domain $\mathcal{D} = [0,0.011]\times [\frac{\pi}{4},\frac{\pi}{2}]$. (a) amounts to showing that the function $f$ is positive on the domain. The Lipschitz constant of $f$ will be estimated explicitly in section \cref{a:homotopy-computations}. We will then partition $\mathcal{D}$ in small boxes using a grid fine enough so to guarantee that, if the value of $f$ is above an appropriate threshold $\eta$ at the center of a box of the partition, then the function is in fact positive on the whole $\mathcal{D}$. A computer program will then evaluate the function at all the centers and validate the bound.

\medskip

{\bf Proof of (b).} It remains to show the existence of the map $F$. First of all $F$ is defined to be the identity on $E\setminus \Sigma_0$. The objective is thus to define $F: \Sigma_0 \to \Sigma$ so that $F$ is the identity on the boundaries of the two triangles $COB$ and $AOM$ forming $\Sigma_0$. $F$ is constructed as the composition of two maps $F = F_2\circ F_1$ and from the construction it will be clear that it satisfies the assumptions of \cref{p:Almgren-relaxed-2}. Both maps are in fact explicit maps from a simplicial complex to another simplicial complex.

\smallskip

\begin{figure}
\tdplotsetmaincoords{60}{0}
\begin{tikzpicture}[tdplot_main_coords, scale=2.5, line join=round]
    \coordinate (O) at (0,0,0);
    \coordinate (C) at (-0.5,0,1);
    \coordinate (B) at (0.5,0,1);
    \coordinate (M) at (0,0.3,1);
    \draw[thick] (C) -- (M) -- (O) -- cycle;

    \draw[thick] (B) -- (M) -- (O) -- cycle;

    \draw[fill= gray!20] (C) -- (B) -- (M) -- cycle;

    \draw[thick] (C) -- (B) -- (M) -- cycle; 
    \draw[dashed] (O) -- (M);
    \node[below] at (O) {$O$};
    \node[left] at (C) {$C$};
    \node[right] at (B) {$B$};
    \node[above] at (M) {$M$};

    \node[right] at (-1.5,0,1) {$B$};
    \node[left] at (-2.5,0,1) {$C$};
    \node[below] at (-2,0,0) {$O$};
    \node[above] at (-2,0.3,1) {$M$};
    \fill (-2,0.3,1) circle (0.4pt);
    \draw[fill= gray!20] (-2.5,0,1) -- (-1.5,0,1) -- (-2,0,0.7) -- cycle;
    \draw[thick] (-2.5,0,1) -- (-1.5,0,1) -- (-2,0,0) -- cycle;
    \draw (-2,0,0.7) -- (-2,0,0);
    \node[below left] at (-2,0,0.7) {$\bar M$};
\end{tikzpicture}
    \caption{The map $F_1$ is piecewise affine and maps the triangle $CBO$ in the union of the three triangles $CMB$, $CMO$, and $MBO$.}\label{f:mappaF1}
\end{figure}
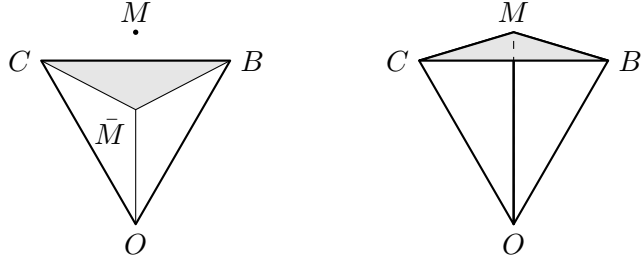

$F_1$ is the identity on $AOM$ and maps the triangle $COB$ onto the union of the triangles $CMO \cup BMO \cup CMB$. This is easily done as follows: we pick a point $\bar M$ interior to the triangle $COB$ and then prescribe $F_1 (C) = C$, $F_1 (B)=B$, $F_1 (O)=O$, and $F_1 (\bar M)=M$. $F_1$ is then defined to be the unique affine extension in each of the three triangles $C\bar M O$, $B \bar M O$, and $CB \bar M$, cf. \cref{f:mappaF1}.

\smallskip

$\Sigma_1 = F_1 (\Sigma_0)$ is therefore the union of four triangles, namely the ``horizontal'' triangle $BCM$ and the three triangles $CMO$, $BMO$ and $AMO$. We next define $F_2$ on $\Sigma_1$. $F_2$ is the identity on the horizontal triangle and on the triangle $CMO$. It thus has to map the two triangles $BMO$ and $AMO$ on the rest of $\Sigma$, which we denote by $\Sigma_2$. Recalling that $D$ is the Torricelli point of $A,B$, and $M$ and that $A'=rA$, $B'=rB$, $M'=rM$, and $D'=rD$, $\Sigma_2$ is the union of
\begin{itemize}
\item The cone with basis the net $A'D'\cup B'D'\cup M'D'$, which is formed by the three triangles $A'OD'$, $B'OD'$, and $M'OD'$;
\item The ``transition'' surface formed by the five triangles $AA'D'$, $AMD'$, $MM'D'$, $BMD'$, and $BB'D'$.
\end{itemize}
Note that $F_2$ is required to be the identity on the boundaries of both triangles $BOM$ and $AOM$, namely on the four segments $MO$, $AO$, $BO$, $BM$, and $AM$.

$F_2$ is in fact defined in very similar fashion to $F_1$. On the triangle $BMO$ we mark the points $B'$ and $M'$ and the midpoint $D_B$ of the segment $B'M'$. We then prescribe that $F_2$ is the identity on all the five points $O$, $M$, $B$, $M'$, and $B'$, while $F_2 (D_B)=D'$. Joining the five fixed points with $D_B$ we subdivide $BMO$ into five triangles: $F_2$ is then he unique affine extension on each of these. The procedure on $AMO$ mirrors the one for $BMO$. More precisely we mark the point $A'$ and let $D_A$ be the midpoint of the segment $A'M'$. We then prescribe $F_2 (A)=A$, $F_2 (A')=A'$, and $F_2 (D_B) = D'$, while of course we already have fixed $F_2 (M')=M'$, $F_2 (M) = M$, and $F_2 (O) = O$. Even in this case we subdivide $AMO$ into the corresponding five triangles and define $F_2$ to be the unique affine extension in each of them.

\section{Path-Isometries}\label{s:path-isometries}

In this section, we describe an algorithm that assigns an orientation-preserving isometry of the sphere to each polygonal path. We then provide a neat formula for determining this isometry and we infer an estimate on the distance between the isometries assigned to two different paths, assuming they share the same internal angles but differ in edge lengths. These tools are central in the computer-assisted search for approximate metrizations of candidate once-derived marked spherical graphs, as detailed in \cref{s:secondo-programma}.

\subsection{Path-isometry Formulas}
As usual, orientation-preserving isometries of the sphere will be called {\em special isometries} and will be denoted (as a group) by $\SO (3)$. $T\mathbb S^2$ will denote the tangent bundle of the sphere $\mathbb S^2$. It will be useful for us to identify $T\mathbb S^2$ with the appropriate subset of pairs $(e, w)$ of vectors in $\mathbb R^3\times \mathbb R^3$, as specified in the following definition.

\begin{definition}[Unit tangent bundle]\label{d:tangent-bundle}
The tangent bundle $T\mathbb S^2$ is the set of pairs $(e,w)\in \mathbb R^3\times \mathbb R^3$ such that
\begin{align}
&|e|=1\label{e:bundle-1}\\
&w\cdot e =0\, .\label{e:bundle-2}
\end{align}
Finally, $\TS$ is the subset of $T\mathbb S^2$ consisting of tangent vectors of length $1$, namely the pairs $(e_1,e_2)$ which satisfy \eqref{e:bundle-1}, \eqref{e:bundle-2}, and $|e_2|=1$.
\end{definition}

We further introduce two families of elements of $\SO (3)$ related to $\TS$

\begin{definition}[Exponentiation and rotation]
If $v= (e_1, e_2) \in T^1 \mathbb S^2$, we let $e_3=e_1\times e_2$. We then denote
\begin{itemize}
\item by $E(v,t)$ the counterclockwise rotation of angle $t$ around the axis directed by $e_3$;
\item by $R (v, \alpha)$ the counterclockwise rotation of angle $\alpha$ around the axis directed by $e_1$.
\end{itemize}
\end{definition}

Note that $R (v, \alpha)$ depends only upon the base point $e_1$ of the element $v$.
We are now ready to define what we will often call \emph{path-isometry}.

\begin{definition}[Path-isometry]\label{d:path-isometry}
    Consider a polygonal path consisting of vertices $V_0, \ldots, V_N$ connected by geodesic arcs $\gamma_i$ (starting at $V_{i-1}$ and ending at $V_i$), parametrized over $[0, \lambda_i]$ by arclength (note that the arcs are not necessarily length-minimizing).
    The associated \emph{path-isometry} is the unique special isometry that maps the initial point $V_0$ into the last point $V_N$ and the initial unit tangent $\dot{\gamma}_1(0)$ into the last unit tangent $\dot{\gamma}_N(\lambda_N)$.
\end{definition}

Because of the definition of $E(v, t)$ and $R(v,\alpha)$, the following formula for the path-isometry is almost tautological.
\begin{lemma}\label{lem:path-isometry-easy}
    Following the same notation for a path as in \cref{d:path-isometry}, assume that the counterclockwise angle at $V_i$ from $-\dot{\gamma}_i(\lambda_i)$ to $\dot{\gamma}_{i+1}(0)$ is $\alpha_i$.
    The \emph{path-isometry} associated to this path is
    \begin{equation}\label{e:path-isometry-easy}
    E(W_{N-1}, \lambda_N) \circ R (W_{N-1}, \alpha_{N-1}-\pi) \circ E (W_{N-2}, \lambda_{N-1}) \circ
    \ldots \circ R (W_1, \alpha_1-\pi) \circ E (W_0, \lambda_1)\, ,
    \end{equation}
    where $W_i := (V_i, \dot{\gamma}_{i+1} (0))\in \TS$.
\end{lemma}
Our goal is to obtain a formula for path-isometries that is more stable than \cref{e:path-isometry-easy} with respect to perturbation in the lengths $\lambda_i$. We state the formula immediately, but defer the proof until the necessary technical tools have been introduced.

\begin{proposition}\label{prop:path-isometry-hard}
    Following the same notation for a path as in \cref{lem:path-isometry-easy}, let $w_0=(V_0, \dot\gamma_1(0))$ and define $w_{i}:=(V_0, R(w_0, \alpha_1 + \cdots+\alpha_i-i\pi)\dot\gamma_1(0))$. Then, the path-isometry associated to this path is
    \begin{equation}\label{e:path-isometry-hard}
        E(w_0, \lambda_1)\circ E(w_1, \lambda_2) \circ \cdots \circ E(w_{N-1}, \lambda_N) \circ R\big(w_0, \alpha_1 + \cdots + \alpha_{N-1} - (N-1)\pi\big) \,.
    \end{equation}
\end{proposition}

First of all, we note that the group $\SO (3)$ acts naturally on $\TS$.

\begin{definition}[Action of $\SO (3)$ on $\TS$]\label{def:action-so-tangent}
The differential of $T\in \SO(3)$ gives a map $T_*:\TS\to\TS$. Using the identification of \cref{d:tangent-bundle}, its action on $(e_1, e_2)\in \TS \subset \mathbb R^3\times \mathbb R^3$ is $(e_1,e_2) = v \mapsto T_* (v)= (T (e_1), T (e_2))$.
\end{definition}

\begin{remark}[Transitivity]\label{r:trans-uniq}
The above action of $\SO(3)$ on $\TS$ is transitive: given any two pairs $(e_1, e_2)$, $(f_1, f_2)$ of orthogonal unit vectors of $\mathbb R^3$ there is an element $T$ of $\SO (3)$ which maps $e_1$ into $f_1$ and $e_2$ into $f_2$. In fact $T$ is unique: if we let $e_3 = e_1\times e_2$ and $f_3=f_1\times f_2$ it follows necessarily that $f_3 = T (e_3)$, thus determining the value of $T$ on a basis of $\mathbb R^3$.
\end{remark}

We next recast the exponential and rotation introduced previously as differentials.

\begin{definition}[Exponential and Rotation]
    Let us define the exponential map
    \[
    \exp:\mathbb S^1\times \TS\to \TS
    \]
    as $\exp_t((e_1,e_2))= (\gamma (t), \gamma'(t))$, where $\gamma:\R\to\mathbb S^2$ is the unit-speed parametrization of the geodesic of $\mathbb S^2$ starting at $\gamma (0)= e_1$ with $\gamma'(0)=e_2$.

    Let us define the rotation map $\rot:\mathbb S^1\times \TS\to\TS$ as $\rot_\alpha(v) = R_*(v)$, where $R\in\SO (3)$ is the (counterclockwise) rotation of $\R^3$ around the (oriented) axis directed along the basepoint $e_1$ (of the element $v= (e_1, e_2)$) with angle $\alpha$.
\end{definition}

The following lemma is a simple exercise left to the reader.

\begin{lemma}[Exponential and Rotation as Linear Maps]\label{lem:local65}
For any $v\in \TS$ and any $t,\alpha\in\mathbb S^1$, we have
$\exp_t(v) = E(v, t)_*v$ and $\rot_\alpha(v) = R(v, \alpha)_*v$.
\end{lemma}

\begin{remark}
    Observe that both $\exp_t:\TS\to\TS$ and $\rot_\alpha:\TS\to\TS$ commute with the action of $\SO(3)$ on $\TS$.
\end{remark}

We are ready to show \cref{prop:path-isometry-hard}.

\begin{proof}[Proof of \cref{prop:path-isometry-hard}]
We prove by induction on $1\le k\le N$, the identity
\begin{equation}\label{eq:local741}
    \big(E(w_0, \lambda_1)\circ\cdots \circ E(w_{k-1}, \lambda_k)\big)_*(w_{k-1})
    =
    (V_k, \dot\gamma_k(\lambda_k))
    \,.
\end{equation}
For $k=1$, the identity boils down to $E(w_0, \lambda_1)_*w_0 = (V_1, \dot\gamma_1(\lambda_1))$, which follows from \cref{lem:local65}.
Now, assuming that \cref{eq:local741} holds for a given $k$, we show it for $k+1$.
Using \cref{lem:local65} and the fact that $\exp$ and $\rot$ commute with the action of $SO(3)$, we obtain
\begin{align*}
    E(w_k, \lambda_{k+1})_*(w_k) =
    \exp_{\lambda_{k+1}}(w_k)
    =
    \exp_{\lambda_{k+1}}(\rot_{\alpha_k-\pi}(w_{k-1}))
    =
    \big(\exp_{\lambda_{k+1}}\circ \rot_{\alpha_k-\pi}\big)(w_{k-1})\,.
\end{align*}
Using the latter identity, we deduce
\begin{align*}
    \big(E(w_0, \lambda_1)\circ\cdots &\circ E(w_{k}, \lambda_{k+1})\big)_*(w_k)
    \\
    &=
    \big(E(w_0, \lambda_1)\circ\cdots \circ E(w_{k-1}, \lambda_{k})\big)_*\big(
        (\exp_{\lambda_{k+1}}\circ \rot_{\alpha_k-\pi})(w_{k-1})
    \big)
    \\
    &=
    (\exp_{\lambda_{k+1}}\circ \rot_{\alpha_k-\pi})
    \Big(
        \big(E(w_0, \lambda_1)\circ\cdots \circ E(w_{k-1}, \lambda_{k})\big)_*
        (w_{k-1})
    \Big)
    \\
    &=
    (\exp_{\lambda_{k+1}}\circ \rot_{\alpha_k-\pi})(V_k, \dot\gamma_k(\lambda_k))
    =
    \exp_{\lambda_{k+1}}(V_k, \dot\gamma_{k+1}(0))
    =(V_{k+1}, \dot\gamma_{k+1}(\lambda_k)) \,,
\end{align*}
which is \cref{eq:local741} when $k$ is replaced by $k+1$.

Applying \cref{eq:local741} with $k=N$, we obtain
\begin{align*}
    \Big(E&(w_0, \lambda_1)\circ \cdots \circ E(w_{N-1}, \lambda_N) \circ R\big(w_0, \alpha_1 + \cdots + \alpha_{N-1} - (N-1)\pi\big)\Big)_*(V_0, \dot\gamma_1(0))
    \\
    &=
    \big(E(w_0, \lambda_1)\circ \cdots \circ E(w_{N-1}, \lambda_N)\big)_* (w_{N-1})
    =(V_N, \dot\gamma_N(\lambda_N))
\end{align*}
which concludes the proof.
\end{proof}

\subsection{Path-isometry stability}
Let us state the stability estimate for path-isometries. As usual the operator norm $\norm{S}_{op}$ of an element $S\in \SO (3)$ is defined as
\[
\norm{S}_{op} := \{\max |S (e)|: e \in \mathbb S^2\}\, .
\]

\begin{proposition}[Path-isometry Stability]\label{p:path-isometry}
Consider two polygonal paths with the following properties:
\begin{itemize}
\item The starting vertex and the unit tangent vector of the first edge at it are the same;
\item The number $N$ of edges is the same;
\item The internal angles at the $i$-th vertex are the same for every $i\in \{1, \ldots , N-1\}$.
\end{itemize}
Let $S$ and $\bar{S}$ be the path-isometries assigned to the respective paths as in \cref{d:path-isometry}.
If we denote by $\lambda_i$ and $\bar \lambda_i$ the lengths of the respective edges, then
\begin{equation}\label{e:path-isometry-est}
\|S-\bar{S}\|_{op} \leq \sum_{i=1}^N |\lambda_i - \bar{\lambda}_i|\, .
\end{equation}
\end{proposition}

The proof of the estimate in \cref{p:path-isometry} is a direct consequence of the formula in \cref{prop:path-isometry-hard} and of the following two simple observations.

\begin{lemma}\label{lem:product-is-lip}
    Let $T_1,T_2,\ldots, T_k\in \SO(3)$ and $\bar T_1, \bar T_2, \ldots, \bar{T}_k\in \SO(3)$ be two $k$-uples of (special) isometries.
    We have
    \begin{equation*}
        \norm{T_1\circ \cdots \circ T_k - \bar{T}_1\circ \cdots \circ \bar{T}_k}_{op}
        \le
        \sum_{i=1}^k \norm{T_i-\bar{T}_i}_{op} \,.
    \end{equation*}
\end{lemma}

\begin{lemma}\label{lem:E-is-lip}
    For any $v\in \TS$ and any $s,t\in\mathbb S^1$, we have
    \begin{equation*}
        \norm{E(v, s) - E(v, t)}_{op} \le \abs{s-t} \,.
    \end{equation*}
\end{lemma}

Before coming to their proofs, we show how \cref{p:path-isometry} follows from them.

\begin{proof}[Proof of \cref{p:path-isometry}]
Let $V_0, \ldots, V_N$ and $\bar{V}_0, \ldots, \bar{V}_N$ be the vertices of the two paths and let $\lambda_i$ and $\bar\lambda_i$ be the length of the geodesic segments connecting $V_{i-1}$ with $V_i$ and $\bar{V}_{i-1}$ with $\bar V_i$.
Observe that the elements $w_i\in \TS$ appearing in \cref{e:path-isometry-hard} coincide for the two paths as the angles coincide. In particular, we have
\begin{equation*}
    \norm{S-\bar S}_{op} = \norm{
        E(w_0, \lambda_1)\circ \cdots \circ E(w_{N-1}, \lambda_N)
        -
        E(w_0, \bar\lambda_1)\circ \cdots \circ E(w_{N-1}, \bar\lambda_N)
    }_{op} \,.
\end{equation*}
The conclusion follows applying first \cref{lem:product-is-lip} and then \cref{lem:E-is-lip}.
\end{proof}

\begin{proof}[Proof of \cref{lem:product-is-lip}]
    We prove it by induction on $k$. For $k=1$, the statement is trivial.

    For $k\ge 2$, we have
    \begin{equation}\begin{aligned}\label{eq:local6657}
        &\norm{T_1\circ \cdots \circ T_k
        -
        \bar{T}_1\circ \cdots \circ \bar{T}_k
        }_{op} \\
        &\le
        \norm{T_1\circ \cdots \circ T_k - \bar{T}_1 \circ \cdots \circ \bar{T}_{k-1} \circ T_k}_{op}
        +
        \norm{\bar{T}_1 \circ \cdots \circ \bar{T}_{k-1} \circ T_k - \bar{T}_1\circ \cdots \circ \bar{T}_k}_{op}
        \\
        & =
        \norm{(T_1 \circ \cdots \circ T_{k-1}
        -
        \bar{T}_1\circ \cdots \circ \bar{T}_{k-1}
        )\circ T_k}_{op}
        +
        \norm{\bar{T}_1\circ \cdots \circ \bar{T}_{k-1}\circ (T_k-\bar{T}_k)}_{op}
        \\
        & =
        \norm{T_1 \circ \cdots \circ T_{k-1}
        -
        \bar{T}_1\circ \cdots \circ \bar{T}_{k-1}}_{op}
        +
        \norm{T_k-\bar{T}_k}_{op} \,.
    \end{aligned}\end{equation}
    In the above computations we have used that the operator norm is a norm, and hence satisfies the triangle inequality, and that it is invariant with respect to the (both left and right) composition with elements of $\SO(3)$.
    Joining \cref{eq:local6657} with the inductive hypothesis, the conclusion follows.
\end{proof}
\begin{proof}[Proof of \cref{lem:E-is-lip}]
    Without loss of generality, $v =(e_1, e_2)\in \TS$ where $(e_1, e_2, e_3:=e_1\times e_2)$ will be used as the basis of $\R^3$.
    For any $\lambda\in\mathbb S^1$, we have
    \begin{equation}\begin{aligned}\label{eq:local123}
        \norm{E(v, s)-E(v, t)}_{op}
        &=
        \le \int_s^t \norm*{\frac{d}{d\lambda}E((e_1, e_2), \lambda)}_{op} \,d\lambda\,.
    \end{aligned}\end{equation}
    Observe that $E((e_1, e_2), \theta)$ is represented by the matrix
    \begin{equation*}
        \begin{pmatrix}
            \cos\theta & \sin\theta & 0 \\
            -\sin\theta & \cos\theta & 0 \\
            0 & 0 & 1
        \end{pmatrix}
    \end{equation*}
    and the derivative with respect to $\theta$ is
    \begin{equation*}
        \begin{pmatrix}
            -\sin\theta & \cos\theta & 0 \\
            -\cos\theta & -\sin\theta & 0 \\
            0 & 0 & 0
        \end{pmatrix} \,.
    \end{equation*}
    The latter is the composition of a projection with an isometry, thus its norm is at most $1$ (and it is, in fact, equal to $1$). This observation, together with \cref{eq:local123}, concludes the proof.
\end{proof}

\section{Computer-assisted restriction for metric realizations}\label{s:secondo-programma}

In this section we will describe an algorithm which accomplishes the following task.
The algorithm takes as input a marked spherical graph $G$ in the list of possibilities obtained in \cref{p:second-topological-classification} for a once-derived net. As an output
\begin{itemize}
    \item[(i)] either it proves that there is no cone $E$ as in \cref{t:main} with a link inducing a once-derived net $\mathcal{N}$ with marked spherical graph $G$;
    \item[(ii)] or it assigns a length to each edge of the graph so that the resulting net $\mathcal N$ embedded in the sphere satisfies \emph{approximately} the constraints we proved to hold for such an $\mathcal{N}$.
\end{itemize}
Furthermore, for all marked combinatorial maps listed in \cref{p:second-topological-classification} but two, the algorithm shows that there is at most one output as in (ii), namely that (up to isometries and up to small perturbations) there is at most one $\mathcal{N}$ ``realizing'' the given combinatorial map.

In order to make the above statement more precise we introduce the following terminology.

\begin{definition}[Metrized Marked Spherical Graph]
    A \emph{metrized} marked spherical graph (MMSG) is a marked spherical graph $G$ (see \cref{d:candidate}) endowed with an assignment $\len:\mathcal E(G)\to (0,2\pi)$ of a positive length to each of its edges.

    Given two metrized marked spherical graphs $(G, \len_1)$ and $(G, \len_2)$ on the same underlying marked spherical graph $G$, their distance is the $L^\infty$ distance between $\len_1$ and $\len_2$.
\end{definition}

\begin{remark}
Observe that a once-derived net $\mathcal{N}$ induces a length assignment $\len_{\mathcal{N}}$ on the associated marked spherical graph $G$. However, a metrized marked spherical graph is not necessarily induced by a once-derived net: we do not even require that it is induced by a geodesic net in the sense of \cref{d:reti}.
\end{remark}

Let us state the Theorem proved by the algorithm.
\begin{theorem}\label{thm:part2}
Consider a set $E$ as in \cref{t:main} and let $\mathcal{N}$ be the net once-derived from its link. If the once-derived marked spherical graph $G$ associated to $\mathcal{N}$ is non-extinct, then it belongs to the following subset of the list give in \cref{p:second-topological-classification}:
    \begin{align*}
    \mathscr{L} := \{&28, 99, 158, 182, 250,
        319, 325, 346, 361, 421,
        468, 484, 520, 614,
        \\ &681,
        712, 736, 779, 827, 909,
        920, 938, 1029, 1096, 1100,
        1122\}\,
    \end{align*}
    Moreover, if $G\in \mathscr{L}\setminus \{827,1029\}$, then there is a metrized marked spherical graph $(G, \len)$ (see \verb|data/metrized_nets/| in \cite{repo} or \cite{website}) with the following property:
    \begin{itemize}
        \item[(U)] let $E$ be as in \cref{t:main} and $\mathcal{N}$ be the net once-derived from its link; if marked spherical graph associated to $\mathcal{N}$ is $G$ and $\len_{\mathcal{N}}$ is the corresponding length assignment for its edges, then the distance from $(G, \len)$ to $(G, \len_{\mathcal{N}})$ is at most $10^{-4}$.
    \end{itemize}
\end{theorem}

For none of the cases mentioned in \cref{thm:part2} we actually prove either existence or uniqueness of a once-derived net $\mathcal{N}$ realizing $G$. Nonetheless, we strongly believe that in all of them a once-derived net actually exists. Algebraic methods might be sufficient to upgrade the statement to an existence and uniqueness one, but we do not pursue this goal since it is not necessary for our scope (and might require a lot of additional work).

Before coming to the description of the algorithm we point out the heuristic reason behind the fact that for most topologies the metrization is absent and, when it is present, in all but two cases (goggles and fat fish families, cf. \cref{s:fat-fish+goggles}) the once-derived net is in fact rigid. Indeed, one can check that the lengths $\sigma_1, \ldots , \sigma_N$ must satisfy $N+1$ equations: this gives an overdetermined system of equations which hence, in most cases, does not have any solution.

\subsection{High-level presentation of the algorithm}\label{ss:algo2presentation}
Let us relax the definition of metrized marked spherical graph by allowing the length to be intervals instead of singletons.
\begin{definition}[Approximately Metrized Marked Spherical Graph]\label{def:ammsg}
    An \emph{approximately} metrized marked spherical graph (AMMSG) is a marked spherical graph $G$ together with an assignment $\inter:\mathcal E(G)\to \{I: I\subseteq[0,2\pi] \text{ is an interval}\}$ of an interval to each of its edges.

    We say that a metrized marked spherical graph $(G, \len)$ is a \emph{realization} of an approximately metrized marked spherical graph $(G, \inter)$ if, for each edge $e\in\mathcal E(G)$, we have $\len(e)\in \inter (e)$.

    We say that an AMMSG $(G, \inter)$ is a refinement of another $(G, \inter')$ if $\inter (e) \subseteq \inter'(e)$ for each edge $e\in\mathcal E(G)$.
\end{definition}

The core routine \texttt{search-valid-realization} of our algorithm takes as input an approximately metrized marked spherical graph $(G,\inter)$ and either
\begin{itemize}
    \item Proves that there is no realization of $(G,\inter)$ which is a net $\mathcal{N}$ once-derived from the link of a minimizing cone $E$ as in \cref{t:main},
    \item Or produces a realization of $(G,\inter)$ which \emph{approximately} satisfies most of the constraints of such an $\mathcal{N}$. Notice that in this case we do not guarantee anything on the metrized marked spherical graph (only that it is a realization of $(G, \inter)$).
\end{itemize}
Furthermore, one can specify some additional constraints (e.g., two edges must have the same length) and restrict the search to metrized marked spherical graphs satisfying these additional constraints.

Let us explain how \texttt{search-valid-metrization} can produce a proof of \cref{thm:part2}.
For each marked combinatorial graph $G$ in the list of \cref{p:second-topological-classification} which does not belong to the list $\mathscr{L}$ of \cref{thm:part2}, we run \texttt{search-valid-metrization} starting from the approximately metrized marked spherical graph $(G, \inter)$ where $\inter (e)$ is the interval $I$ defined as follows.
\begin{itemize}
    \item If $e$ is the edge of a face with $2$ sides, then $I=\{\pi\}$.
    \item If $e$ is the edge of a convex face with $3$ sides, then $I=\{\lambda_T\}$ (recall \cref{p:formulae}).
    \item If $e$ is the edge of a convex face with $4$ sides, then $I=[0.01, \pi-0.001]$ (recall \cref{t:effective-short-edge}).
    \item If $e$ is the edge of a convex face with $5$ sides, then $I=[0.01, \lambda_T-0.001]$ (recall \cref{p:formulae,t:effective-short-edge}).
    \item If $e$ is not the edge of a convex face but we can apply \cref{l:long-edge} to conclude that its length is at most $\pi$, then $I=[0,\pi]$;
    \item If none of the above applies, then $I=[0, 2\pi]$.
\end{itemize}
Then, \texttt{search-valid-metrization} proves that any metrized marked spherical graph on $G$ fails to satisfy at least some of the necessary conditions found so far to be the once-derived net associated to the link of a minimizing cone $E$ as in \cref{t:main}.

It remains to prove rigidity, namely \cref{thm:part2}(U), when $G$ belongs to $\mathscr{L}\setminus \{827,1029\}$. Let $G$ be one such element, and perform the following sequence of operations.
\begin{enumerate}
    \item Run \texttt{search-valid-realization} on $(G, \inter)$ where $\inter$ is the same described above.
    This produces a metrized marked spherical graph $(G, \len)$.
    \item Rerun \texttt{search-valid-realization}, starting again from $(G, \inter)$, but with the additional constraint that we search only metrized marked spherical graphs with distance greater than $10^{-4}$ from $(G, \len)$. This shows that no such metrized marked spherical graph is once-derived from the link of a cone $E$ as in \cref{t:main}.
\end{enumerate}
This procedure works in all cases apart from net $736$, which requires additional considerations. Indeed, if we run this procedure for net $736$, the routine is unable to show the rigidity because there is a $1$-parameter family of ``valid'' metrizations. These are valid only approximately (i.e., they don't come from true once-derived nets) and we are able to show rigidity in the following way.
After the first run of the strategy we can restrict the values of some specific edges to small intervals: in other words, we find rigidity for a specific subnet. We then use this information together with synthetic and algebraic arguments to prove that two specific edges must be congruent, cf. \cref{prop:736}. Finally we run again \texttt{search-valid-realization} subject to this constraint and show rigidity.

The remaining part of this section is devoted to explaining the internal working of the routine \texttt{search-valid-realization}.

\subsection{The core routine: \texttt{search-valid-realization}}
Recall that the input of the routine is an approximately metrized marked spherical graph $(G, \inter)$. The routine consists then of the following three steps.
\begin{itemize}
\item[1.] We run a sequence of checks on $(G, \inter)$. If $(G, \inter)$ admits a realization that is a net once-derived from a link of a minimizing set $E$ as in \cref{t:main}, then all the checks would pass. Therefore, if any of the checks fails, the routine has proven that there is no valid realization. These checks can modify $\inter$ (some intervals may be reduced). If this happens, it is guaranteed that any valid realization before the change is a realization also after.
\item[2.] If $\sum_{e\in \mathcal E(G)} |\inter (e)| < 10^{-5}$, the routine returns the realization of $(G, \inter)$ that is obtained by assigning to $e$ the midpoint of $\inter (e)$ for all $e\in\mathcal E(G)$.
\item[3.] If we reached this step, it means that we are not able to rule out the existence of a valid realization of $(G, \inter)$ and $(G, \inter)$ is not ``precise'' enough for us to be content. Therefore, we create a collection $\{X_1, \dots, X_k\}$ of refinements of $(G, \inter)$ such that any realization of $(G, \inter)$ is a realization of $X_i$ for at least one $1\le i\le k$. Then, we run the routine on each of $X_i$ independently. If any of these subroutines produces a metrized marked spherical graph, we return it. Otherwise, we have proven that no valid realization of $X$ exists.
\end{itemize}
In this sketch, two crucial details are missing: what are the checks we run during Step 1 and how we produce the refinements of Step 3. These decisions are, arguably, the core of the algorithm and slight variations (in the idea or in the implementation) can easily lead to much longer execution times. In fact slight variations can also lead to the checks being less effective on deciding that some realizations are valid. This clearly would not produce a ``wrong'' answer, but just an answer which is less efficient in the road towards a proof of \cref{t:main}.

\subsection{The refinements of Step 3}
The refinements $X_1, \dots, X_k$ of $(G, \inter)$ mentioned in Step 3 are produced as follows.
A subset $\mathcal A$ of the edges of $G$ is selected; the selection process depends only on $G$ and not on $\inter$.
We obtain $k=2^{|\mathcal A|}$ refinements by bisecting each of the intervals assigned to the edges in $\mathcal A$.

Since it is not crucial for the correctness of the algorithm, we give only a brief idea of how $\mathcal A$ is chosen (we refer the interested reader to \verb|embed_or_exclude.cpp| in \cite{repo}).
The set $\mathcal A$ is a small subset of edges such that if we knew the length of all the edges in $\mathcal A$ then we could deduce the lengths of all the remaining edges with a simple sequence of inferences (e.g., thanks to \cref{p:formulae} if we know two sides of a pentagonal convex face, then we immediately know all the sides). Moreover, we select $\mathcal A$ so that the deduction of the lengths of the remaining edges given the lengths of $\mathcal A$ is as numerically stable as possible.

\subsection{The checks of Step 1}\label{subsec:embed-or-exclude}
We consider a metrized marked spherical graph ``valid'' if it is induced by a once-derived net (see \cref{d:once-derived}) which also satisfies the lower bounds on the edge lengths stated in \cref{t:effective-short-edge}.

A check is a process that takes in an AMMSG and tries to prove that there is no valid realization of the given AMMSG. If it fails, it may still refine the given AMMSG by showing that if a valid realization exists then it is a valid realization also of the refinement it produces.

In order to describe the checks, we shall first explain how to obtain a quantitative approximation of the relative position of two vertices (or darts) in an approximately metrized marked spherical graph.

\subsubsection{Approximate Relative Position for AMMSGs}
Let us begin by introducing the concept of relative position for once-derived nets.

\begin{definition}\label{def:once-derived-relative-position}
    Consider a once-derived net $\mathcal N$.
    Identify a dart $d\in \darts(\mathcal N)$ with the element of $\TS$ given by initial vertex and the direction of the dart. Given two darts $d$ and $d'$, let $F_{\mathcal N}(d, d')$ be the (special) isometry that maps $d$ onto $d'$.
\end{definition}
Observe that $F_{\mathcal N}(d, d')$ encodes the relative position of two darts; in particular it is invariant under rotations of the whole net.

We will now describe how to obtain an analogous object for MMSGs and AMMSGs.

\begin{definition}\label{def:mmsg-relative-position}
    Consider a metrized marked spherical graph $(G, \len)$. Let $d_1, d_2, \dots, d_k$ be a path of darts in $G$ (i.e., the second vertex of $d_i$ is the first vertex of $d_{i+1}$).

    We know the assigned lengths $\len(d_1), \dots, \len(d_k)$, and we have the information that the clockwise angle from $d$ to $\next(d)$ is $\frac{2\pi}3$ for all $d\in \darts(G)$, with the exceptions of the marked darts, for which the angle is $\frac{4\pi}3$. Given these pieces of information, after fixing a point $V_0\in \mathbb S^2$ for the starting vertex of $d_1$ and a direction for it, we can build a corresponding polygonal path which following the darts $d_1, d_2, \dots, d_k$, cf. \cref{d:polygonal-path}. Then, as already done in \cref{def:once-derived-relative-position}, we can identify each $d_i$ with an element of $\TS$, and we let
    $F_{(G, \len)}(d_1, [d_2, \dots, d_{k-1}], d_k)$ be the (special) isometry that maps $d_1$ onto $d_k$. Observe that this isometry does not depend on the position chosen initially for $d_1$.
\end{definition}
Observe that, differently from the case of once-derived nets, this time we have to choose a path from $d_1$ to $d_k$ to define their relative position in an MMSG.
Moreover, notice that if $(G,\len)$ is the MMSG associated to a once-derived net $\mathcal N$, then
\begin{equation*}
    F_{\mathcal N}(d_1, d_k) = F_{(G, \len)}(d_1, [d_2, \dots, d_{k-1}], d_k) \, ,
\end{equation*}
in other words it is really in this case that the relative position is {\em independent} of the chosen path.

We can now move to the case of AMMSGs.
\begin{definition}\label{def:ammsg-relative-position}
    Consider an approximately metrized marked spherical graph $(G, \inter)$. Let $(G, \len)$ be the metrized marked spherical graph obtained by setting $\len(e)$ to be the midpoint of $\inter (e)$ for all edges $e\in\mathcal E(G)$.
    Let $d_1, d_2, \dots, d_k$ be a path of darts in $G$ (i.e., the second vertex of $d_i$ is the first vertex of $d_{i+1}$).

    Let $F_{(G, \inter)}(d_1, [d_2, \dots, d_{k-1}], d_k)$ be the subset of $\SO(3)$ given by
    \begin{equation*}
        \Big\{
            T\in \SO(3):\, \norm*{T-F_{(G, \len)}(d_1, [d_2, \dots, d_{k-1}], d_k)}_{op} \le
            \frac{|\inter (d_1)| + \cdots + |\inter (d_{k-1})|}2
        \Big\}
        \,.
    \end{equation*}
\end{definition}

The crucial fact is the following one.

\begin{lemma}\label{lem:approximate-isometry}
    Let $(G, \inter)$ be an approximately metrized marked spherical graph and consider a path $d_1, d_2, \dots, d_k$ of darts in $G$.
If $\mathcal N$ is a once-derived net which is a realization of $(G, \inter)$ (when interpreted as a MMSG), then
    \begin{equation*}
        F_{\mathcal N}(d_1, d_k) \in F_{(G, \inter)}(d_1, [d_2, \dots, d_{k-1}], d_k) \,.
    \end{equation*}
\end{lemma}
\begin{proof}
    After unrolling the definition, this is a restatement of \cref{p:path-isometry}.
\end{proof}

\subsubsection{List of checks}
We are ready to present the checks.
However, we do not provide here a detailed explanation of their implementation. Instead, for each of them, we describe what is the property of the net that they verify and which of the results we proved they are enforcing or applying. We also refer the reader to our \texttt{C++} implementation \verb|metrization_constraints.cpp| in \cite{repo}.
As soon as one of the checks fails, the algorithm terminates as it has proven that no valid realization exists.

Our list of checks may seem arbitrary. We want to emphasize that if any of them is removed, then the list of \emph{plausible} nets produced by the algorithm is larger than the one of \cref{thm:part2}. The reader should interpret the list of checks as a \emph{meta-proof}: we are describing the main ingredients of the minimal proof we could find.
Note also that, even if we do not prove it, we believe that all metrized marked spherical graphs found in \cref{thm:part2} correspond to a true stationary geodesic with two antipodal pins. This implies that our list of checks is exhaustive.

\begin{description}[
    leftmargin=0pt,   
    itemsep=1em,      
    parsep=0pt,       
    font=\bfseries,   
    style=unboxed     
]
\item[Inequalities on convex faces] We apply, on convex faces with $4$ or $5$ sides, the inequalities in \cref{lem:quadrilateral_consecutive_lowerbound}, \cref{lem:pentagon_consecutive_lowerbound}, and \cref{p:pentagon-monotonicity} to restrict the intervals assigned to the edges.
The algorithm keeps cycling the inequalities through all the convex quadrilaterals and pentagons until it cannot produce any further restriction (or one of the intervals assigned becomes empty, proving that no valid realization exists).

\item[Closure of nonconvex faces]
We check whether there is at least one nonconvex face which is ``obstructed from closing''.
Let us show how one can apply \cref{lem:approximate-isometry} to obtain an obstruction.  Let $(G, \inter)$ be the AMMSG we are working with and let $d_1, d_2, \dots, d_k$ be the darts forming the boundary of a nonconvex face. Then, if $(G, \inter)$ admits a valid realization $\mathcal N$, then
\begin{equation*}
    \id = F_{\mathcal N}(d_1, d_1) \in F_{(G, \inter)}(d_1, [d_2, \dots, d_k], d_1) \,.
\end{equation*}
In particular, if the identity does not belong to the family of isometries on the right-hand side (which is computable) then there is no valid realization.

Our algorithm implements a more delicate check which is more numerically stable and allows us to also restrict some of the intervals when we are not immediately able to exclude the existence of a valid realization. We refrain from describing it here and refer the reader to the code.

\item[Pins' antipodality I] In a stationary net with two antipodal pins, the two pins must be antipodal.
In particular, consider a stationary net with two pins $\mathcal N$ whose once-derived version is a realization of the AMMSG we are working on. Consider three darts of $\mathcal N$ whose second vertex is a pin (not necessarily the same pin). Then, the three equators which contain the three darts must all go through the pins.
Applying \cref{lem:approximate-isometry} similarly to how we used it in the previous check, we are able to exploit this observation to obtain an obstruction for the AMMSG.

\item[Angle at pin] In a stationary net with two antipodal pins, if one pin has degree at least two then its two edges form an angle which is at least $\frac{2\pi}{3}$.
Once again, this condition can be transformed in a possible obstruction for the AMMSG using \cref{lem:approximate-isometry}.

\item[Pins are outside convex faces]
Consider a once-derived net $\mathcal N$. Even though we do not know the exact positions of the pins, we do know the \emph{beginning} of the darts that were removed to obtain the once-derived net and thus, as we observed in \emph{Pins' antipodality I}, we can identify a few equators which contain the pins. Assume that, for $\mathcal N$, at least two such equators do not coincide. Thus, the two pins are uniquely identified as the intersection of those two equators. In particular the two intersection points cannot be strictly inside one of the convex faces of $\mathcal N$.
This statement about once-derived nets, thanks to \cref{lem:approximate-isometry}, can be distilled into a possible obstruction for AMMSGs.

\item[Intersection of distinct edges]
Any two edges of a once-derived net which are not adjacent to the same vertex are disjoint. In particular, this must hold for two edges belonging to the boundary of a nonconvex face.
This statement about once-derived nets, thanks to \cref{lem:approximate-isometry}, can be distilled into a possible obstruction for AMMSGs.

\item[Pins' antipodality II]
Consider a once-derived net $\mathcal N$ and let $Q_1$ and $Q_2$ be two degree $2$ vertices belonging to two distinct nonconvex faces $F_1$ and $F_2$.

Starting from $Q_1$ and entering in the face $F_1$ along the bisector of its internal angle at $Q_1$ formed by the two edges at $Q_1$ we must first encounter the pin $P_1$ belonging to $F_1$ (in the stationary net with two antipodal pins that induced $\mathcal N$) and then the boundary of $F_1$.
Let $d_1$ be the amount of distance we have to go along the bisector to meet the boundary of $F_1$. We have $d_1 \ge \text{dist}_{\mathbb S^2}(P_1, Q_1)$.

Define analogously $P_2$ and $d_2$. Then, since the two pins $P_1$ and $P_2$ are antipodal, we have
\begin{equation*}
    \pi = \text{dist}_{\mathbb S^2}(P_1, P_2) \le
    \text{dist}_{\mathbb S^2}(P_1, Q_1)
    +
    \text{dist}_{\mathbb S^2}(Q_1, Q_2)
    +
    \text{dist}_{\mathbb S^2}(P_2, Q_2)
    \le d_1 + d_2 + \text{dist}_{\mathbb S^2}(Q_1, Q_2) \,.
\end{equation*}
Observe that this inequality involves only quantities that are known if we know $\mathcal N$ (we do not need to know the position of the pins).

This statement about once-derived nets, thanks to \cref{lem:approximate-isometry}, can be distilled into a possible obstruction for AMMSGs.

\item[Boundary test]
This last check is probabilistic and is used to rule out some very pathological cases. Because of its probabilistic nature, it is run multiple times to increase the probability of finding an obstruction.

Consider a once-derived net $\mathcal N$.
Take a random oriented geodesic segment $\sigma$ on the sphere whose two endpoints do not belong to $\mathcal N$. For any dart $d$ of $\mathcal N$ (which is an oriented geodesic edge), let $\iota(\sigma, d)$ be the intersection number of $\sigma$ and $d$, that is,
\begin{itemize}
\item $0$ if the two geodesic segment do not intersect;
\item $1$ if they intersect positively, namely the two tangents $\dot\sigma$ and $\dot d$ form a positively oriented pair at the intersection point;
\item $-1$ if they intersect negatively.
\end{itemize}
For a face $F$ of $\mathcal N$, we then define
\[
\iota (\sigma, F) = \sum_{d \in \partial F} \iota (\sigma, d)\, .
\]

We leave the proof of the following lemma to the reader.
\begin{lemma}
    The following holds:
    \begin{itemize}
        \item[(i)] $\iota (\sigma, F)\in \{-1,0,1\}$ for every face $F$ of $\mathcal N$;
        \item[(ii)] $\iota (\sigma, F) \neq 0$ for at most two faces $F$ of $\mathcal N$;
        \item[(iii)] $\sum_F \iota (\sigma, F) =0$ (where $F$ runs over all faces of $\mathcal N$).
    \end{itemize}
\end{lemma}
This statement about a once-derived net (and a geodesic segment $\sigma$) can be distilled into a possible obstruction for AMMSGs thanks to \cref{lem:approximate-isometry}.
\end{description}

\section{Special cases: extinct nets}\label{s:human-1}

In this section we analyze the first batches of special nets. These are the nets which lead to an empty twice-derived spherical graph. Their once-derived nets contain therefore only vertices of degree $2$, cf. \cref{d:once-derived} and \cref{d:extinction}. These are the ones which we called ``extinct'' and they are classified here in \cref{t:elementary}. Any of the nets found in the latter classification is a link of some element of \cref{lista}.

\begin{theorem}\label{t:elementary}
Let $\mathcal{N}$ be the link of an area-minimizing cone $E$ as in \cref{t:main}, let $\mathcal{N}'$ be the corresponding once-derived net, and consider the marked spherical graph $G$ induced by the latter. Assume that $\mathcal{N}'$ is extinct, namely that the vertices of $G$ have degree $2$. Then $G$ is either empty, or it has $2$, $3$, or $4$ vertices. Moreover:
\begin{itemize}
\item[(i)] If $G$ is the empty graph, then $E$ belongs to \cref{lista}(a)\&(b);
\item[(ii)] If $G$ has two vertices, then $E$ belongs to \cref{lista}(c);
\item[(iii)] If $G$ has three vertices, then $E$ belongs to \cref{lista}(f);
\item[(iv)] If $G$ has four vertices, then $E$ belongs to \cref{lista}(d).
\end{itemize}
\end{theorem}

Recalling \cref{p:derived-1}, in any once-derived net there are at most $4$ vertices of degree $2$, while all the others are of degree $3$. Since, under the assumption of \cref{t:elementary}, $\mathcal{N}'$ has only vertices of degree $2$, these vertices are at most $4$. We will analyze the five possible cases (which include the case of the empty net) separately.

\subsection{$G$ is empty.} In this case the only vertices of $\mathcal{N}$ must be the two pins, because the other vertices are never removed. But then each edge of $\mathcal{N}$ is an half-equator joining the two pins. In particular the two pins have the same degree $d\in \{1,2,3\}$. The degree cannot be $3$, because in this case the pins would not have been removed to derive $\mathcal{N}'$. If $d=1$, then $\mathcal{N}$ consists of a single half-equator joining the two pins, and thus $E$ is as in (a) of \cref{lista}. If $d=2$, then $\mathcal{N}$ consists of two half equators joining the two pins and they must form a convex angle $\geq \frac{2\pi}{3}$. This is case (b) of \cref{lista}.

\subsection{$G$ has one vertex.} This is not possible, because then the net would consist of a loop starting from this vertex and getting back to it: the initial net does not have loops and the edges of the once-derived net are edges of the initial one.

\subsection{$G$ has two vertices.} Let us denote them by $A$ and $B$. There must then be two edges of $\mathcal{N}'$ connecting them, because they have degree $2$ and we have already argued that $G$ contains no loops. In particular $A$ and $B$ would need to be antipodal and the two edges would need to be two half equators connecting them and meeting at $\frac{2\pi}{3}$ angles. The complement of $\mathcal{N}'$ consists then of a convex biangle and a concave biangle. The two pins of the original net $\mathcal{N}$ would need to be in the nonconvex face (or on its boundary if they are $A$ or $B$) and would also need to belong to the bisector of one of the two angles. These bisectors are actually the same equator and its intersection with the nonconvex face has length exactly $\pi$. Since the pins are antipodal, we conclude that $A$ and $B$ are the two pins. Since they have not been removed, they both have degree $3$ in $\mathcal{N}$. But then $\mathcal{N}$ consists of the two pins joined by three half-equators meeting at equal angles, namely $E$ is as in (c) of \cref{lista}.

\subsection{$G$ has three vertices.}
\begin{figure}
    \begin{tikzpicture}[scale=0.7]
        \coordinate (A) at (-2,0);
        \coordinate (C) at (2, 0);
    \coordinate (B) at (0, {2*sqrt(3)});
        \coordinate (P) at (0, {-2*sqrt(3)});
        \draw[very thick] (-8,0) -- (A) -- (B) -- (C) -- (8,0);
        \draw[very thin, dashed] (A) -- (C);
        \node[below left] at (A) {$A$};
        \node[below right] at (C) {$C$};
        \node[above left] at (B) {$B$};
        \node[left] at (-1,{-sqrt(3)}) {$\alpha'$};
        \node[right] at (1,{-sqrt(3)}) {$\beta'$};
        \draw[very thick, dashed] (A) -- (P) -- (C);
        \node[below] at (P) {$P'$};
        \node[below] at (0,0) {$\gamma$};
        \draw[very thick, dashed] (B) -- (0, 7);
        \node[above] at (0,7) {$P$};
        \node at (-3,3) {$\mathcal{T}$};
        \node at (-3,-3) {$\mathcal{T}'$};
        \node[left] at (0,5.5) {$\alpha$};
    \end{tikzpicture}
    \caption{The net $\mathcal{N}$ and $\mathcal{N}'$ when $\mathcal{N}'$ is a cycle with three vertices. The thick segments depict the three edges of $\mathcal{N}'$ while the thick dashed segments represent the edges of $\mathcal{N}$ removed in the derivation of $\mathcal{N}'$.}\label{f:estinta-3-vertici}
\end{figure}
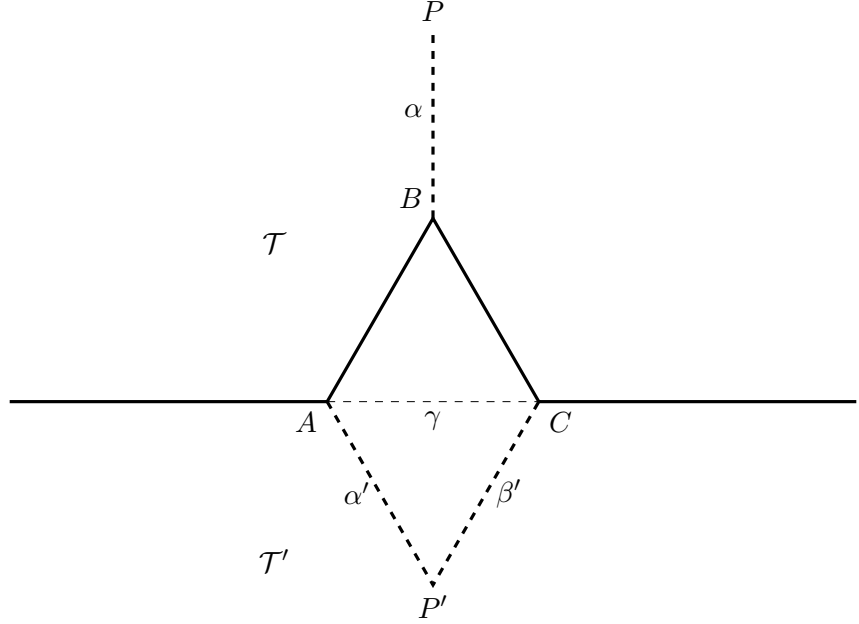

The net $\mathcal{N}'$ forms a cycle which splits the sphere into two triangles, which we denote by $\mathcal{T}$ and $\mathcal{T}'$. We assume, without loss of generality, that the internal angles of $\mathcal{T}$ measure $\frac{2\pi}{3}$, $\frac{2\pi}{3}$, and $\frac{4\pi}{3}$ at three vertices which we denote, respectively, by $A$, $C$, and $B$.
$\mathcal{T}$ is isosceles and using \cref{thm:generalized_triangles} and \cref{thm:convex_triangles} we conclude that the lengths of its edges are determined.
Let $\lambda_T$ and $\lambda_Q$ be, respectively, the side lengths of the equilateral triangle and equilateral quadrilateral with angles equal to $\frac{2\pi}3$. We leave to the reader to check that
\begin{equation*}
    d_{\mathbb S^2}(A, B) = d_{\mathbb S^2}(B, C) = \lambda_Q
    \text{ and }
    d_{\mathbb S^2}(A, C) = 2\pi - \lambda_T \,.
\end{equation*}
Moreover, if $\gamma$ is the equator containing the edge connecting $A$ and $C$, then $\mathcal{T}$ is contained in one of the two hemispheres bounded by $\gamma$, cf. \cref{f:estinta-3-vertici}.

Observe that there must be a pin $P'$ in (or on the boundary of) $\mathcal T'$ connected to the two vertices $A$ and $C$ by two edges $\alpha'$ and $\beta'$ which have been removed from the net $\mathcal{N}$: $P'$, $A$, and $C$ form then a convex triangle congruent to the convex triangle formed by $A$, $B$, and $C$. In particular $P'$ is symmetric to $B$ with respect to $\gamma$. The other pin $P$ is the antipodal of $P'$ and $\mathcal N$ contains an additional edge $\alpha$ joining it to $C$.
We have completely identified the net up to rigid motions and, by \cref{d:thin-fish}, is is the (link of) a thin fish with six edges.

\subsection{$G$ has four vertices.} In this final case the goal is to show that $\mathcal{N}$ is in fact the support of the link of \cref{lista}(d), which we call {\em tetrahedral net}. First of all, the derived-once net $\mathcal{N}'$ splits the sphere into two quadrilaterals, which we denote by $Q$ and $Q'$. Each quadrilateral has two (internal) angles measuring $\frac{2\pi}{3}$ and two (internal) angles measuring $\frac{4\pi}{3}$. In particular we have two possibilities:
\begin{itemize}
\item[(a)] The vertices at which the internal angles of $Q$ measure $\frac{4\pi}{3}$ are consecutive;
\item[(b)] The internal angles of $Q$ at any pair of two consecutive vertices are different.
\end{itemize}
The first case can be excluded by the following lemma.

\begin{lemma}\label{l:no-consecutive-concave}
If $\mathcal{N}$ is as in \cref{t:elementary}(iv), then alternative (a) is ruled out.
\end{lemma}

It is possible to give a direct argument with pen and paper, but we have not been able to find one which is shorter than a couple of pages. On the other hand we can quickly prove it running the code described in \cref{s:secondo-programma}, cf. \verb|processImpossibleQuadrilateralExtinct| in \verb|embed_or_exclude.cpp| in \cite{repo}.

\medskip

From now on, we treat case (b). Therefore, in both quadrilaterals $Q$ and $Q'$ the internal angles alternate in concavity.

{\bf Step 1.} Let us denote the vertices of the face $Q$ as $A$, $B$, $C$, $D$; with $A$ and $C$ having internal angles equal to $\frac{2\pi}{3}$, while $B$ and $D$ having internal angles equal to $\frac{4\pi}{3}$. Both $B$ and $D$ are adjacent to an edge of $\mathcal{N}$ which was removed when deriving $\mathcal{N}'$. If this is the same edge (which would happen when both $B$ and $D$ have degree $3$ in $\mathcal{N}$), then it subdivides $Q$ into two equiangular triangles with angle $\frac{2\pi}{3}$, thereby determining that all of the edges of $\mathcal{N}'$ are long $\lambda_T$. From this we immediately see that $\mathcal{N}$ is the tetrahedral net. Moreover the same consideration can be done on the quadrilateral $Q$'. We can therefore, without loss of generality, assume that:
\begin{itemize}
\item[(c)] $B$ and $D$ are connected to a pin $P$ internal to $Q$ by two edges $\beta$ and $\delta$ of $\mathcal{N}$;
\item[(d)] $A$ and $C$ are connected to a pin $P'$ internal to $Q'$ by two edges $\alpha$ and $\gamma$ of $\mathcal{N}$.
\end{itemize}
Note in addition that, if any of these pairs of edges form a flat angle at one of the two pins, then we again conclude that $\mathcal{N}$ is the tetrahedral net. We can therefore also assume
\begin{itemize}
\item[(e)] the smallest angle formed by $\beta$ and $\delta$ at $P$ is strictly smaller than $\pi$.
\end{itemize}
From now on we want to show that (c), (d), and (e) lead to a contradiction.

\medskip

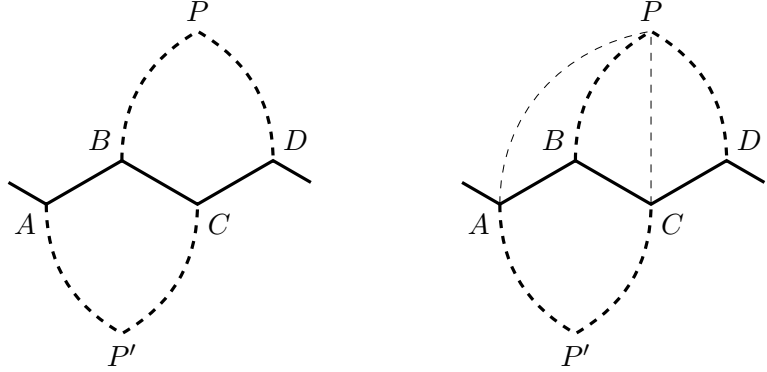
\begin{figure}
        \begin{tikzpicture}[scale=0.5]
         \coordinate (J) at (-2, 0);
         \coordinate (D) at (-3, {sqrt(3)/3});
         \coordinate (C) at (-5, {-sqrt(3)/3});
         \coordinate (B) at (-7, {sqrt(3)/3});
         \coordinate (A) at (-9, {-sqrt(3)/3});
         \coordinate (P) at (-5, 4);
         \coordinate (P') at (-7,-4);
         \coordinate (J') at (-10,0);
         \draw[very thick] (J') -- (A) -- (B) -- (C) -- (D) -- (J);
         \node[below left] at (A) {$A$};
         \node[above left] at (B) {$B$};
         \node[below right] at (C) {$C$};
         \node[above right] at (D) {$D$};
         \node[above] at (P) {$P$};
         \node[below] at (P') {$P'$};
         \draw[very thick, dashed] (B) to [out=90, in=210] (P) to [out=-30, in=90] (D);
          \draw[very thick, dashed] (A) to [out=-90, in=150] (P') to [out=30, in=-90] (C);

         \coordinate (J1) at ({-2+12}, 0);
         \coordinate (D1) at ({-3+12}, {sqrt(3)/3});
         \coordinate (C1) at ({-5+12}, {-sqrt(3)/3});
         \coordinate (B1) at ({-7+12}, {sqrt(3)/3});
         \coordinate (A1) at ({-9+12}, {-sqrt(3)/3});
         \coordinate (P1) at ({-5+12}, 4);
         \coordinate (P'1) at ({-7+12},-4);
         \coordinate (J'1) at ({-10+12},0);
         \draw[very thick] (J'1) -- (A1) -- (B1) -- (C1) -- (D1) -- (J1);
         \node[below left] at (A1) {$A$};
         \node[above left] at (B1) {$B$};
         \node[below right] at (C1) {$C$};
         \node[above right] at (D1) {$D$};
         \node[above] at (P1) {$P$};
         \node[below] at (P'1) {$P'$};
         \draw[very thick, dashed] (B1) to [out=90, in=210] (P1) to [out=-30, in=90] (D1);
          \draw[very thick, dashed] (A1) to [out=-90, in=150] (P'1) to [out=30, in=-90] (C1);
          \draw[very thin, dashed] (C1) -- (P1);
          \draw[very thin, dashed] (A1) to [out =90, in=190] (P1);
        \end{tikzpicture}
        \caption{Final arguments for the rigidity of the tetrahedral case. On the left the quadrilaterals $ABCP'$ and $BPDC$ are used in Step 2 to prove that all edges of $\mathcal{N}'$ are shorter than $\pi$. Thus \cref{thm:congruence} shows that $AB$ is congruent to $CD$ and $BC$ is congruent to $AD$.
On the right, in Step 3 the triangles $CPD$ and $APB$ are shown to be congruent, and likewise the triangles $BPC$ and $ADP$. As a result we conclude that the angle $BPD$ is flat.
}\label{f:estinta-4}
    \end{figure}

{\bf Step 2.}
We next claim that all the edges of $\mathcal{N'}$ are shorter than $\pi$. For the sake of contradiction, assume that, e.g., the edge $AB$ is $\pi$ or longer. Then the quadrilateral $ABDP$ is not convex, cf. \cref{f:estinta-4}, and its area is, by Gauss-Bonnet, at most $\frac{4\pi}{3}$ (recall that the internal angle at $P$ is at most $\frac{4\pi}{3}$). Now, on the edge $AB$ let $\chi$ be the half-equator which has vertices $B$ and $B'$. We can rotate $\chi$ inside $ABPD$ until we hit a point of the boundary which is distinct from $B$ and $B'$. We must be able to rotate until we hit the edge $AP$. Otherwise we would hit some other point on the boundary, but this is not possible because $P$ is the only concave angle of the quadrilateral. So we must have rotated $\chi$ by $\frac{2\pi}{3}$, the angle at $B$, forming a biangle of area $\frac{4\pi}{3}$ contained strictly inside $ABPD$. Since the latter has area at most $\frac{4\pi}{3}$ we have reached a contradiction.

We now orient the closed loop formed by $\mathcal{N}'$ in two ways: one going from $A$ to $B$ to $C$ to $D$, the other going from $B$ to $A$ to $D$ to $C$. In doing so we discover two polygonal paths which have a common edge and equal internal angles: since all the edges are shorter than $\pi$, we can apply \cref{thm:congruence} and conclude that $BC$ is congruent to $DA$. We then repeat apply the argument with the two loops $B$ to $C$ to $D$ to $A$ and $C$ to $B$ to $A$ to $D$, concluding that $AB$ and $CD$ are also congruent.

\medskip

{\bf Step 3.} In this (final) step we will use that the pins $P$ and $P'$ are antipodal. Consider in particular the equator contanining $P'A$: this equator necessarily contains $P$ given that it is antipodal to $P'$. Hence its intersection with $Q$ contains a geodesic segment $PA$ which bisects the internal angle of $Q$ at $A$. Likewise for $C$, cf. \cref{f:estinta-4}. So, the triangles $ABP$ and $CDP$ have the two edges $AB$ and $CD$ congruent and the two angles adjacent to the respective edges congruent (measuring $\frac{\pi}{3}$ and $\frac{2\pi}{3}$). By \cref{thm:generalized_triangles}, the two triangles are congruent. So the angles $A\hat{P}B$ and $C\hat{P}D$ are congruent. An analogous argument implies that $B\hat{P}C$ and $D\hat{P}A$ are also congruent. Looking at the pin $P$ we realize however that $2A\hat{P}B+2B\hat{P}C = 2\pi$, in particular the internal angle at $P$ of the quadrilateral $BADP$ is $\pi$, contradicting (e).

\section{Special cases II: seven difficult families}\label{s:human-2}

Here we will analyze seven marked spherical graphs among the ones listed in \cref{p:second-topological-classification}. The overall conclusions are given in \cref{t:human-classification} and \cref{prop:736}. The link of any element in \cref{lista} generates a marked spherical graph which is either analyzed in \cref{t:elementary} or in \cref{t:human-classification}. However, \cref{t:human-classification} involves also two further families of links which are described in detail in \cref{s:fat-fish+goggles} and depend upon one and two parameters. These are not in \cref{lista}: for these \cref{t:human-classification} must be considered just a step before they will be finally excluded (on energetic grounds) in the next sections. \cref{prop:736} concerns the marked spherical graph numbered as 736 in the list of \cref{p:second-topological-classification}. As it was already mentioned in \cref{s:secondo-programma}, a more subtle interaction of theoretical pen-and-paper considerations with the code described in \cref{s:secondo-programma} allows ultimately to complete the proof of \cref{thm:part2} for the links of minimizing cones $E$ which induce the marked spherical graph 736: this is explained in the proof of \cref{prop:736}.

\subsection{Fat fish and goggles}\label{s:fat-fish+goggles} Our next goal is to analyze metric realizations of the marked spherical graphs numbered as 158, 681, 827, 920, 1029, and 1100, (see \cref{f:158,f:681,f:pesce-grasso-class,f:pesce-magro-class,f:1029,f:falso-cubo} for diagrams depicting them)
in the list of \cref{p:second-topological-classification}. In order to proceed with our analysis we need to introduce two families of stationary geodesic nets with two pins. At a later stage we will show that none of these nets is the link of a minimizing cone, but in this section we will prove that these are the only possible realizations of the marked spherical graphs 827 and 1029.

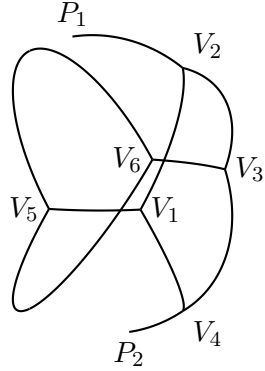
\begin{figure}
\tdplotsetmaincoords{80}{20}
\begin{tikzpicture}[tdplot_main_coords, scale=2]
    \pgfmathsetmacro{\xA}{sqrt(3)/3}
    \pgfmathsetmacro{\yA}{-sqrt(6)/3}
    \pgfmathsetmacro{\zA}{0}
    \pgfmathsetmacro{\xB}{sqrt(3)/3}
    \pgfmathsetmacro{\yB}{0}
    \pgfmathsetmacro{\zB}{sqrt(6)/3}
    \pgfmathsetmacro{\xC}{sqrt(3)/3}
    \pgfmathsetmacro{\yC}{sqrt(6)/3}
    \pgfmathsetmacro{\zC}{0}
    \pgfmathsetmacro{\xD}{sqrt(3)/3}
    \pgfmathsetmacro{\yD}{0}
    \pgfmathsetmacro{\zD}{-sqrt(6)/3}
    \pgfmathsetmacro{\xE}{-1}
    \pgfmathsetmacro{\yE}{0}
    \pgfmathsetmacro{\zE}{0}
    \pgfmathsetmacro{\xF}{0}
    \pgfmathsetmacro{\yF}{-1}
    \pgfmathsetmacro{\zF}{0}
    \pgfmathsetmacro{\xG}{0}
    \pgfmathsetmacro{\yG}{1}
    \pgfmathsetmacro{\zG}{0}
    \pgfmathsetmacro{\xP}{(-0.2)}
    \pgfmathsetmacro{\yP}{0}
    \pgfmathsetmacro{\zP}{sqrt(1-0.2*0.2)}
    \pgfmathsetmacro{\xI}{-0.5}
    \pgfmathsetmacro{\yI}{0}
    \pgfmathsetmacro{\zI}{sqrt(3)/2}
    \pgfmathsetmacro{\xH}{-0.5}
    \pgfmathsetmacro{\yH}{0}
    \pgfmathsetmacro{\zH}{-sqrt(3)/2}

\pgfmathsetmacro{\angle}{acos(\xA*\xB + \yA*\yB + \zA*\zB)}        \pgfmathsetmacro{\angleBC}{acos(\xB*\xC + \yB*\yC + \zB*\zC)}    \pgfmathsetmacro{\angleCD}{acos(\xC*\xD + \yC*\yD + \zC*\zD)} \pgfmathsetmacro{\angleDA}{acos(\xD*\xA + \yD*\yA + \zD*\zA)}    \pgfmathsetmacro{\angleAF}{acos(\xA*\xF + \yA*\yF + \zA*\zF)}    \pgfmathsetmacro{\angleCG}{acos(\xC*\xG + \yC*\yG + \zC*\zG)}
\pgfmathsetmacro{\angleBP}{acos(\xP*\xB + \yP*\yB + \zP*\zB)}
\pgfmathsetmacro{\angleDQ}{acos(-\xP*\xD - \yP*\yD - \zP*\zD)}
\pgfmathsetmacro{\angleFH}{acos(\xF*\xH + \yF*\yH + \zF*\zH}
\pgfmathsetmacro{\angleFI}{acos(\xF*\xI + \yF*\yI + \zF*\zI}
\pgfmathsetmacro{\angleGH}{acos(\xG*\xH + \yG*\yH + \zG*\zH}
\pgfmathsetmacro{\angleGI}{acos(\xG*\xI + \yG*\yI + \zG*\zI}

\draw[thick, samples=100, domain=0:1] plot[variable=\t]
(
{(sin((1-\t)*\angle)/sin(\angle)) * \xA + (sin(\t*\angle)/sin(\angle)) * \xB},
{(sin((1-\t)*\angle)/sin(\angle)) * \yA + (sin(\t*\angle)/sin(\angle)) * \yB},
{(sin((1-\t)*\angle)/sin(\angle)) * \zA + (sin(\t*\angle)/sin(\angle)) * \zB}
);

\draw[thick, samples=100, domain=0:1] plot[variable=\t]
    (
        {(sin((1-\t)*\angleBC)/sin(\angleBC)) * \xB + (sin(\t*\angleBC)/sin(\angleBC)) * \xC},
        {(sin((1-\t)*\angleBC)/sin(\angleBC)) * \yB + (sin(\t*\angleBC)/sin(\angleBC)) * \yC},
        {(sin((1-\t)*\angleBC)/sin(\angleBC)) * \zB + (sin(\t*\angleBC)/sin(\angleBC)) * \zC}
    );

\draw[thick, samples=100, domain=0:1] plot[variable=\t]
    (
        {(sin((1-\t)*\angleCD)/sin(\angleCD)) * \xC + (sin(\t*\angleCD)/sin(\angleCD)) * \xD},
        {(sin((1-\t)*\angleCD)/sin(\angleCD)) * \yC + (sin(\t*\angleCD)/sin(\angleCD)) * \yD},
        {(sin((1-\t)*\angleCD)/sin(\angleCD)) * \zC + (sin(\t*\angleCD)/sin(\angleCD)) * \zD}
    );

\draw[thick, samples=100, domain=0:1] plot[variable=\t]
    (
       {(sin((1-\t)*\angleDA)/sin(\angleDA)) * \xD + (sin(\t*\angleDA)/sin(\angleDA)) * \xA},
        {(sin((1-\t)*\angleDA)/sin(\angleDA)) * \yD + (sin(\t*\angleDA)/sin(\angleDA)) * \yA},
        {(sin((1-\t)*\angleDA)/sin(\angleDA)) * \zD + (sin(\t*\angleDA)/sin(\angleDA)) * \zA}
    );

\draw[thick, samples=100, domain=0:1] plot[variable=\t]
    (
        {(sin((1-\t)*\angleAF)/sin(\angleAF)) * \xA + (sin(\t*\angleAF)/sin(\angleAF)) * \xF},
        {(sin((1-\t)*\angleAF)/sin(\angleAF)) * \yA + (sin(\t*\angleAF)/sin(\angleAF)) * \yF},
        {(sin((1-\t)*\angleAF)/sin(\angleAF)) * \zA + (sin(\t*\angleAF)/sin(\angleAF)) * \zF}
    );

\draw[thick, samples=100, domain=0:1] plot[variable=\t]
    (
        {(sin((1-\t)*\angleCG)/sin(\angleCG)) * \xC + (sin(\t*\angleCG)/sin(\angleCG)) * \xG},
        {(sin((1-\t)*\angleCG)/sin(\angleCG)) * \yC + (sin(\t*\angleCG)/sin(\angleCG)) * \yG},
        {(sin((1-\t)*\angleCG)/sin(\angleCG)) * \zC + (sin(\t*\angleCG)/sin(\angleCG)) * \zG}
    );

\draw[thick, samples=100, domain=0:1] plot[variable=\t]
    (
        {(sin((1-\t)*\angleBP)/sin(\angleBP)) * \xB + (sin(\t*\angleBP)/sin(\angleBP)) * \xP},
        {(sin((1-\t)*\angleBP)/sin(\angleBP)) * \yB + (sin(\t*\angleBP)/sin(\angleBP)) * \yP},
        {(sin((1-\t)*\angleBP)/sin(\angleBP)) * \zB + (sin(\t*\angleBP)/sin(\angleBP)) * \zP}
    );

\draw[thick, samples=100, domain=0:1] plot[variable=\t]
    (
        {(sin((1-\t)*\angleDQ)/sin(\angleDQ)) * \xD - (sin(\t*\angleDQ)/sin(\angleDQ)) * \xP},
        {(sin((1-\t)*\angleDQ)/sin(\angleDQ)) * \yD - (sin(\t*\angleDQ)/sin(\angleDQ)) * \yP},
        {(sin((1-\t)*\angleDQ)/sin(\angleDQ)) * \zD - (sin(\t*\angleDQ)/sin(\angleDQ)) * \zP}
    );

\draw[thick, samples=100, domain=0:1] plot[variable=\t]
    (
        {(sin((1-\t)*\angleFH)/sin(\angleFH)) * \xF + (sin(\t*\angleFH)/sin(\angleFH)) * \xH},
        {(sin((1-\t)*\angleFH)/sin(\angleFH)) * \yF + (sin(\t*\angleFH)/sin(\angleFH)) * \yH},
        {(sin((1-\t)*\angleFH)/sin(\angleFH)) * \zF + (sin(\t*\angleFH)/sin(\angleFH)) * \zH}
    );

\draw[thick, samples=100, domain=0:1] plot[variable=\t]
    (
        {(sin((1-\t)*\angleFI)/sin(\angleFI)) * \xF + (sin(\t*\angleFI)/sin(\angleFI)) * \xI},
        {(sin((1-\t)*\angleFI)/sin(\angleFI)) * \yF + (sin(\t*\angleFI)/sin(\angleFI)) * \yI},
        {(sin((1-\t)*\angleFI)/sin(\angleFI)) * \zF + (sin(\t*\angleFI)/sin(\angleFI)) * \zI}
    );

\draw[thick, samples=100, domain=0:1] plot[variable=\t]
    (
        {(sin((1-\t)*\angleGH)/sin(\angleGH)) * \xG + (sin(\t*\angleGH)/sin(\angleGH)) * \xH},
        {(sin((1-\t)*\angleGH)/sin(\angleGH)) * \yG + (sin(\t*\angleGH)/sin(\angleGH)) * \yH},
        {(sin((1-\t)*\angleGH)/sin(\angleGH)) * \zG + (sin(\t*\angleGH)/sin(\angleGH)) * \zH}
    );

\draw[thick, samples=100, domain=0:1] plot[variable=\t]
    (
        {(sin((1-\t)*\angleGI)/sin(\angleGI)) * \xG + (sin(\t*\angleGI)/sin(\angleGI)) * \xI},
        {(sin((1-\t)*\angleGI)/sin(\angleGI)) * \yG + (sin(\t*\angleGI)/sin(\angleGI)) * \yI},
        {(sin((1-\t)*\angleGI)/sin(\angleGI)) * \zG + (sin(\t*\angleGI)/sin(\angleGI)) * \zI}
    );

\node[above] at (\xP,\yP,\zP) {$P_1$};
\node[below] at (-\xP,-\yP,-\zP) {$P_2$};
\node[right] at (\xA,\yA,\zA) {$V_1$};
\node[right] at (\xC,\yC,\zC) {$V_3$};
\node[above right] at (\xB,\yB,\zB) {$V_2$};
\node[below right] at (\xD,\yD,\zD) {$V_4$};
\node[left] at (\xF,\yF,\zF) {$V_5$};
\node[left] at (\xG,\yG,\zG) {$V_6$};
\end{tikzpicture}
\caption{An element in the fat fish family. The picture should be compared to \cref{f:PlateauB4}: a portion of $\gamma_5$ in \cref{f:PlateauB4} has been removed and replaced by the half-equators which in the picture above join $V_5$ and $V_6$ (another difference is that in \cref{f:PlateauB4} the quadrilateral $V_1V_2V_3V_4$ is made of straight segments, while in the above picture we have drawn the corresponding projections on the sphere).}\label{f:fat-fish}
\end{figure}

\begin{definition}[Fat fish]\label{d:pesce-grasso}
Consider any geodesic net $\mathcal{N}'$ which is an element of the thin fish family, cf. \cref{d:thin-fish}. Let $\gamma_5$ be the (only) geodesic arc in $\mathcal{N}'$ which is longer than $\pi$ and let $V_5$ and $V_6$ be two arbitrary antipodal points in (the interior of) $\gamma_5$. Eliminate from the net the half equator connecting $V_5$ and $V_6$ which is a subset of $\gamma_5$ and add the two half-equators connecting the same vertices which meet the remaining portions of $\gamma_5$ forming angles of $\frac{2\pi}{3}$ at $V_5$ and $V_6$, cf. \cref{f:fat-fish}. If one of the two half-equators joining $V_5$ and $V_6$ intersects one of the two edges $\gamma_6$ and $\gamma_7$, then we discard the net, otherwise it is a {\em fat fish}.
\end{definition}

Once we mod out the symmetries of the sphere, the fat fish family is a $2$-parameter family: the lengths of the arcs $V_5V_1$ and  $V_2 P_1$ both vary in an interval. All the other edges have instead a fixed length.

We next define a second family, which we call ``goggles''. We recall that $\lambda_T$ denotes the length of the sides of an equiangular convex spherical triangle with internal angles equal $\frac{2\pi}{3}$. Note that, given any pair of points $A$ and $C$ on the sphere at distance $\lambda_T$, there are exactly two points $B$ and $D$ with the property that $ABCD$ is a convex equilateral equiangular quadrilateral with angles $\frac{2\pi}{3}$: this can be proved observing that the two diagonals of the quadrilateral $AC$ and $BD$ must measure $\lambda_T$, must be perpendicular, and must intersect at the respective midpoints.

\begin{figure}[htbp]
\begin{tikzpicture}[scale=0.7]
\draw[very thick] (-7,0) -- (-4.5,0) to [bend left] (-2.5,2) -- (-0.5,0) -- (0.5,0) -- (2.5,2) to [bend left] (4.5,0) -- (7,0);
\draw[very thick] (-4.5,0) to [bend right] (-2.5,-2) -- (-0.5,0);
\draw[very thick] (0.5,0) -- (2.5,-2) to [bend right] (4.5,0);
\draw [very thick] (-2.5,2) to [bend left] (0,5) to [bend left] (2.5,2);
\draw [very thick] (-2.5,-2) to [bend right] (0,-5) to [bend right] (2.5,-2);
\draw[very thick, dashed] (-4.5,0) -- (-0.5,0);
\draw[very thick, dashed] (4.5,0) -- (0.5,0);
\node[right] at (7,0) {$\gamma$};
\node[above] at (0,5) {$P_1$};
\node[below] at (0,-5) {$P_2$};
\node[below left] at (-4.5,0) {$V_1$};
\node[below right] at (4.5,0) {$V_4$};
\node[below] at (-0.3,0) {$V_2$};
\node[above] at (0.3,0) {$V_3$};
\node[right] at (-2.3,2) {$V_5$};
\node[left] at (2.3,2) {$V_7$};
\node[right] at (-2.3,-2) {$V_6$};
\node[left] at (2.3,-2) {$V_8$};
\end{tikzpicture}
\caption{A schematic picture of the goggle family. The long straight line represent the equator $\gamma$: the dashed lines are taken away from it, so that the remaining portion consists of two segments connecting $V_2$ and $V_3$ and $V_4$ and $V_1$. The two quadrilaterals $V_1V_5V_2V_7$ and $V_3V_6V_4V_8$ are equilateral and have internal angles $\frac{2\pi}{3}$. The pins $P_1$ and $P_2$ are the two ``poles'', namely the centers of the hemispheres forming the complement of $\gamma$.}\label{f:goggles}
\end{figure}
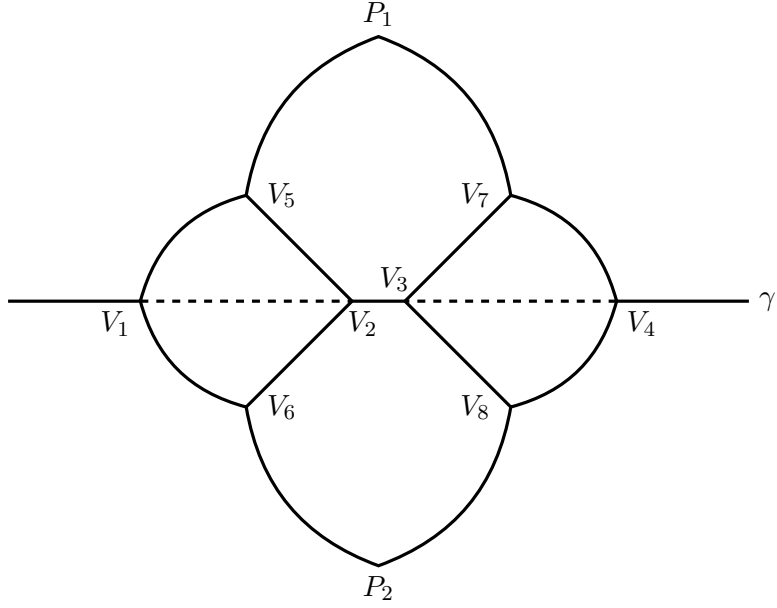

\begin{definition}[Goggles]\label{d:goggles}
Let $\gamma$ be any equator on $\mathbb S^2$ and mark on it four consecutive points $V_1,V_2,V_3,V_4$ with the property that the geodesic distances between $V_1$ and $V_2$ and $V_3$ and $V_4$ are both $\lambda_T$. Remove the respective arcs from $\gamma$, so that what remains consists of two length-minimizing geodesic arcs joining $V_2$ with $V_3$ and $V_4$ with $V_1$. Next add the four vertices $V_5, V_6, V_7, V_8$ such that $V_1V_2V_5V_6$, resp. $V_3,V_4V_7V_8$, are the four vertices of an equilateral quadrilateral with angles $\frac{2\pi}{3}$. The vertices are labeled in such a way that $V_5$ and $V_7$ lie on the same connected component of $\mathbb S^2\setminus \gamma$, cf. \cref{f:goggles}.

Add to the net the four sides of each quadrilateral. Furthermore, let $P_1$ and $P_2$ be the two centers of the two hemispheres forming $\mathbb S^2\setminus \gamma$, labeled so that $P_1$ is in the same hemisphere as $V_5$ and $V_7$, and add the length-minimizing geodesic arcs joining $P_1$ with $V_5$ and $V_7$ and $P_2$ with $V_6$ and $V_8$. We say that the net obtained in this fashion is a {\em goggle} if and only if the convex angle formed by $V_5 P_1$ and $V_7 P_1$ at $P_1$ is at least $\frac{2\pi}{3}$.
\end{definition}

Note that the nets in the goggles family consist of 10 vertices and twelve edges. Each of these nets is symmetric with respect to $\gamma$. There are (up to isometries) two extremal elements in the family. One is given when the convex angle at the pins is exactly $\frac{2\pi}{3}$. If this is the angle internal to the pentagon $V_2 V_3 V_7 P_1 V_5$, then the latter and the symmetric pentagon $V_2 V_3 V_8 P_2 V_6$ are equiangular with angles equal $\frac{2\pi}{3}$. The other extremal case is the one in which the angle at the pin is $\pi$, namely when $P_1$ is the midpoint of the length-minimizing geodesic arc joining $V_5$ and $V_7$ (which we will denote as $V_5V_7$) and $P_2$ is the midpoint of the corresponding geodesic arc $V_6V_8$. In this case the $8$ vertices $V_1V_2V_3V_4V_5V_6V_7V_8$ are in fact the vertices of a cube inscribed in the unit sphere and the line $\ell$ joining $P_1$ and $P_2=-P_1$ goes through the midpoints of two opposite edges of the cube. This extremal net is thus in common with the link of a cone in the cubical family described in \cref{t:main}(e). Note however that this is one case for which we will prove that the corresponding cone is not Almgren minimizing.

\subsection{Classification of six special nets} We are now ready to state the main theorem classifying  the following six special nets.

\begin{theorem}\label{t:human-classification}
Let $\mathcal{N}$ be the link of a cone $E$ as in \cref{t:main} and denote by $\mathcal{N}'$ the corresponding once-derived net. 

Let $id$ be the identifier of the marked spherical graph associated to $\mathcal{N}'$ in the list of \cref{p:second-topological-classification}. Then, for some special values of $id$ we have the following statements:
\begin{itemize}
    \item[$id=158$:] This is impossible.
    \item[$id=681$:] This is impossible.
    \item[$id=827$:] $\mathcal N$ is an element of the fat fish family.
    \item[$id=920$:] $\mathcal N$ is an element of the thin fish family (see \cref{lista}(f)).
    \item[$id=1029$:] $\mathcal N$ is either an element of the cubical family (see \cref{lista}(e)) or the goggles family.
    \item[$id=1100$:] $\mathcal N$ is an element of the cubical family (see \cref{lista}(e)).
\end{itemize}
\end{theorem}

The next subsections will analyze the various cases of \cref{t:human-classification}.

\subsection{1029 is a cube or a goggle} 
\begin{figure}
\begin{tikzpicture}[scale=0.7]
\draw[very thick] (-7,0) -- (-5.5,0) node[above left] {$0$} -- (-3.5,2) node[ above] {$4$} -- (-1.5,0) node[above right] {$1$} -- (1.5,0) node[above left] {$2$} -- (3.5,2) node[above] {$5$} -- (5.5,0) node[above right] {$3$} -- (7,0);
\draw[very thick] (-5.5,0) -- (-3.5,-2) node[below] {$6$} -- (-1.5,0);
\draw[very thick] (1.5,0) -- (3.5,-2) node[below] {$7$} -- (5.5, 0);
\end{tikzpicture}
\caption{A picture of the marked spherical graph 1029 for a once-derived net $\mathcal{N}'$. The picture is on the plane, the graph becomes spherically embedded by adding a point at infinity so that the vertices $0$ and $3$ are connected by an edge. The two quadrilateral $0416$ and $2537$ are convex.}\label{f:1029}
\end{figure}
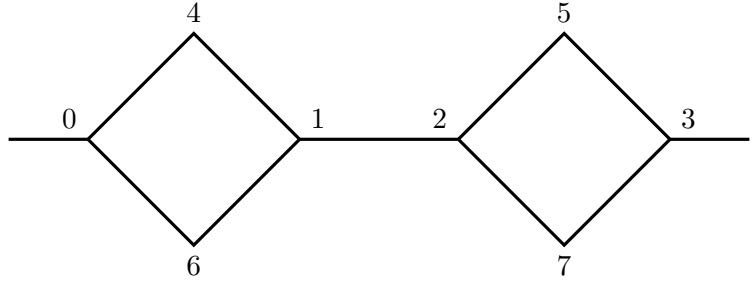

In this step we will prove the following fact. Assume $\mathcal{N}$ and $\mathcal{N}'$ are as in \cref{t:human-classification}, and that the marked spherical graph associated to $\mathcal{N}'$ is numbered 1029. Then $\mathcal{N}$ is necessarily an element of the goggle family or of the cube family. We start by drawing a picture schematizing the properties of the once-derived net $\mathcal{N}'$, cf. \cref{f:1029}.

Since the most important part of the proof will also be used to deal with other cases, we isolate its consequences in the following lemma.

\begin{lemma}\label{l:goggles-rigidity}
Consider a geodesic net such that:
\begin{enumerate}[label=(\roman*)]
    \item Its spherical graph is depicted in \cref{f:1029} (from now on, we identify vertices of the graph with vertices of the net);
    \item Convex angles in \cref{f:1029} correspond to angles of $\tfrac{2\pi}3$ in the net;
    \item The edge $12$ is shorter than $\pi$.
\end{enumerate}
Then, one of the following two statements hold:
\begin{itemize}
    \item[(a)] The two quadrilaterals $0416$ and $2537$ are equilateral and the four vertices $0$, $1$, $2$, and $3$ belong to the same equator.
    \item[(b)] The geodesic segments of length $\pi$ starting from $4$ and $5$ and bisecting the respective concave angles do not intersect.
\end{itemize}
\end{lemma}

In order to apply the lemma to the case at hand we need to check condition (iii): this is, however, a consequence of \cref{l:long-edge}. Moreover, because of the existence of the pins, we can exclude case (b) of the conclusion. Thus, (a) holds and we deduce that $\mathcal N$ is either an element of the goggle family or an element of the cube family.
The rest of this section is devoted to a proof of \cref{l:goggles-rigidity}.

\begin{proof}[Proof of \cref{l:goggles-rigidity}]
Let $Q$ and $Q'$ be the equiangular quadrilaterals $2537$ and $0416$.
Observe that, if any of the two is equilateral, then (a) holds. Indeed, without loss of generality, assume that $Q$ is equilateral. Then, the diagonal $23$ of $Q$ is a bisector of the angles in $2$ and $3$ and thus $0$, $1$, $2$, and $3$ all lie on the same equator $\gamma$. Take now the diagonal $10$ of $Q'$; it must also lie on $\gamma$ and thus the latter diagonal bisects the two angles at $1$ and $0$. So, the convex triangles $106$ and $104$ are both equiangular and equilateral. We conclude that $Q'$ is also equilateral.

\medskip

From now on, we assume that neither $Q$ nor $Q'$ is equilateral.
We draw the two diagonals of the quadrilaterals $Q$ and $Q$', which we denote by $01$ and $23$ and consider:
\begin{itemize}
\item the closed polygonal path $\mathcal{P}$ which start at $0$ and proceeds to $1$, $2$, $3$;
\item the closed polygonal path $\mathcal{P}'$ from $1$ to $0$ to $3$ to $2$.
\end{itemize}
Since $Q$ and $Q'$ are not equilateral, none of the angles of these loops are flat. Furthermore, all edges of the loops, apart from $03$, are shorter than $\pi$ thanks to assumption (iii).
Note moreover that the four angles of the loops are congruent (here we are using that opposite sides in an equiangular convex quadrilateral are congruent).
We can thus apply \cref{thm:congruence} and conclude that the two loops are congruent. This implies that the edges $03$ and $21$ are congruent.

Next, we observe that the two equators containing the two edges $03$ and $21$ are at the same time bisectors of opposite angles of both quadrilaterals. We will then use the following lemma, which can be checked via a direct computation in spherical trigonometry: the explicit computations are given in \cref{a:quadrilateral-diagonal}.

\begin{lemma}\label{l:quadrilateral-diagonal}
Consider two equiangular spherical quadrilaterals $Q_1$ and $Q_2$ with angles of $\frac{2\pi}{3}$. Assume that the two equators bisecting two opposite angles of $Q_1$ form the same angle as the two equators bisecting two opposite angles of $Q_2$. Then $Q_1$ and $Q_2$ are congruent.
\end{lemma}

The lemma implies that the quadrilaterals $Q$ and $Q'$ are congruent. So either $41$ is congruent to $52$ or $41$ is congruent to $27$. If we now look again at the polygonal path $\mathcal{P}$ (where each pair of edges which do not have a vertex in common consists of congruent edges) the two alternatives amount to the angles at $1$ and $2$ in this path being equal or summing to $2\pi$. If the angles were equal, then the diagonals of the two quadrilaterals $Q$ and $Q'$ would lie in the same closed hemisphere $H$ bounded by the equator passing through $1$ and $2$. In particular $3$ and $0$ would belong to the interior of $H$ and so would segment edge $30$, since it is shorter than $\pi$. But then one of the two quadrilaterals $\Omega$ bounded by the polygonal path $\mathcal{P}$ would have area strictly smaller than $2\pi$. But Gauss-Bonnet applied to $\Omega$ shows that its area is exactly $2\pi$, leading to a contradiction.

Summarizing, we have shown that the two convex quadrilaterals are congruent and, moreover:
\begin{itemize}
\item The edges $37$, $16$, $40$, and $52$  are all congruent;
\item The edges $72$, $60$, $41$, and $53$ are all congruent.
\end{itemize}
Let $u_4, u_5\in \TS$ be the unit tangent vectors directed along the bisectors of the concave angles in $4$ and $5$.

Let $O$ be the unique orientation-preserving isometry of the sphere which maps the point $1$ in the point $3$ and the edge $16$ to the the edge $37$ (cf. \cref{f:1029}). The map $O$ maps the geodesic net into itself, and in particular it maps $u_4$ to $u_5$ and viceversa (recall the action of $SO(3)$ on $\TS$ introduced in \cref{def:action-so-tangent}).

Now, let us consider three cases depending on the relative position of $u_4$ and $u_5$.
\begin{enumerate}
    \item Assume that $u_4$ and $u_5$ identify the same equator (via the exponential map) with the same orientation. Then, given that $O_*(u_4)=u_5$ and $O_*(u_5)=u_4$, the points $4$ and $5$ must be antipodal, which is a contradiction because $4$ and $5$ lay on the same side of the equator identified by $12$.
    \item Assume that $u_4$ and $u_5$ identify the same equator (via the exponential map) with opposite orientations. Then, consider the quadrilateral $1254$ where the edge $54$ is the shortest geodesic connecting $5$ to $4$. Observe that the edges of this quadrilateral are non-intersecting because $4$ and $5$ lay on the same side of the equator identified by $12$. The edge $54$ is either the bisector of the concave angles in $5$ and $4$ or it is the bisector of the convex angles in $5$ and $4$.
    \begin{itemize}
        \item If it is the bisector of the concave angles, then the quadrilateral $1254$ is a convex equiangular quadrilateral with angles of $\tfrac{2\pi}3$ and thus $14$ is congruent to $25$; which is a contradiction because it implies that $Q$ and $Q'$ are equilateral.
        \item If it is the bisector of the convex angles, then Gauss-Bonnet implies that its area is $\tfrac{8\pi}3$, which is impossible because the quadrilateral is contained in a hemisphere bounded by the equator identified by $12$.
    \end{itemize}
    \item The equators identified (via the exponential map) by $u_4$ and $u_5$ are distinct. If the conclusion (b) of the statement holds, we do not have to prove anything. Otherwise, there is a point $V$ such that $\exp_s(u_4)=V=\exp_t(u_5)$ with $0 < s, t < \pi$. Since $O$ swaps $u_4$ and $u_5$, we have $O(V)=V$ and in particular $s=t$, which is equivalent to the congruence of $4V$ and $5V$. Applying \cref{thm:congruence} to the pentagon $4125V$ we deduce that $41$ is congruent to $25$ which is a contradiction because it implies that $Q$ and $Q'$ are equilateral.
\end{enumerate}

In particular $O(V)=V$, which in turn implies the congruence of $6V$ and $7V$. Thanks to \cref{thm:congruence}, we deduce that the pentagon $6127V$ is isosceles, from which in turn we infer that $16$ and $27$ are also congruent.
\end{proof}

\subsection{158 and 681 are not possible} 
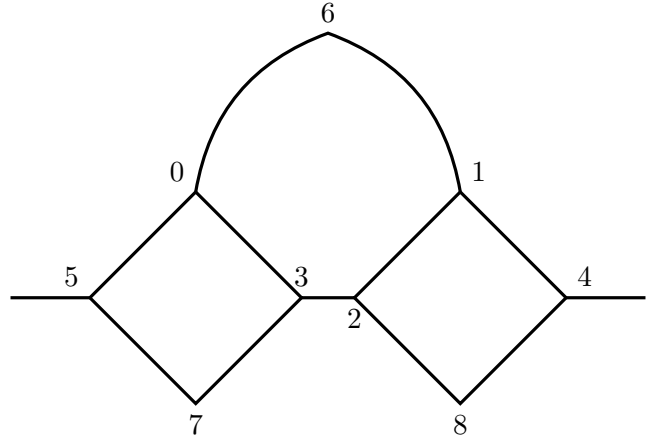
\begin{figure}
\begin{tikzpicture}[scale=0.7]
\draw[very thick] (-6,0) -- (-4.5,0) node[above left] {$5$} -- (-2.5,2) node[above left] {$0$} -- (-0.5,0) node[above] {$3$} -- (0.5,0) node[below] {$2$} -- (2.5,2) node[above right] {$1$} -- (4.5,0) node[above right] {$4$} -- (6,0);
\draw[very thick] (-4.5,0) -- (-2.5,-2) node[below] {$7$} -- (-0.5,0);
\draw[very thick] (0.5,0) -- (2.5,-2) node[below] {$8$} -- (4.5, 0);
\draw[very thick] (-2.5,2) to [bend left] (0, 5) node[above]{$6$} to [bend left] (2.5,2);
\end{tikzpicture}
\caption{A picture of the marked spherical graph 158 for a once-derived net $\mathcal{N}'$. The quadrilaterals $5037$ and $2148$ and the pentagon $30612$ are convex.}\label{f:158}
\end{figure}

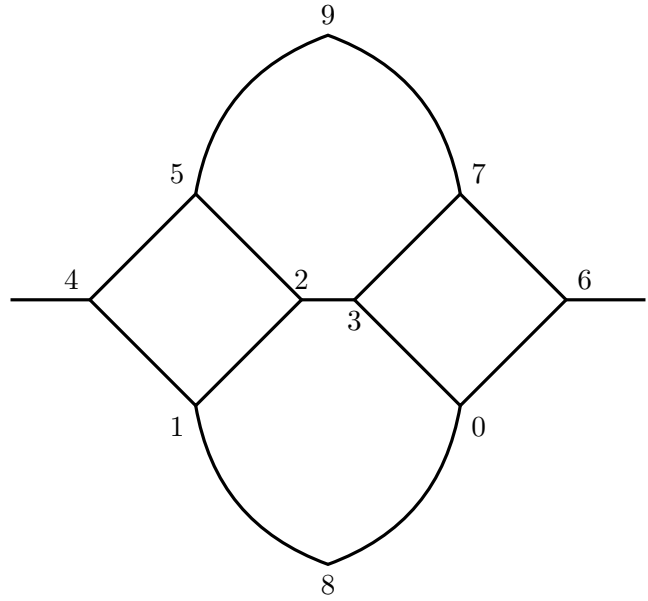
\begin{figure}
\begin{tikzpicture}[scale=0.7]
\draw[very thick] (-6,0) -- (-4.5,0) node[above left] {$4$} -- (-2.5,2) node[above left] {$5$} -- (-0.5,0) node[above] {$2$} -- (0.5,0) node[below] {$3$} -- (2.5,2) node[above right] {$7$} -- (4.5,0) node[above right] {$6$} -- (6,0);
\draw[very thick] (-4.5,0) -- (-2.5,-2) node[below left] {$1$} -- (-0.5,0);
\draw[very thick] (0.5,0) -- (2.5,-2) node[below right] {$0$} -- (4.5, 0);
\draw[very thick] (-2.5,2) to [bend left] (0, 5) node[above]{$9$} to [bend left] (2.5,2);
\draw[very thick] (-2.5,-2) to [bend right] (0, -5) node[below]{$8$} to [bend right] (2.5,-2);
\end{tikzpicture}
\caption{A picture of the marked spherical graph 681 for a once-derived net $\mathcal{N}'$. The quadrilaterals $4521$ and $3760$ and the pentagons $25973$ and $30812$ are convex.}\label{f:681}
\end{figure}

The arguments are very similar and we will therefore carry them in parallel. First of all we draw schematic pictures of both topologies, cf. \cref{f:158} and \cref{f:681}.

Observe that, by \cref{l:long-edge} any once-derived net $\mathcal{N}'$ with either marked spherical graph has all edges which are strictly shorter than $\pi$. We now see that in both cases there is a subnet $\mathcal{N}''$ to which we can apply \cref{l:goggles-rigidity}. In the case of 158 the subnet is obtained by removing the vertex $6$ and the edges $60$ and $61$, while in the case of 681 we obtain it by removing the vertices $9$ and $8$ and the edges $95$, $97$, $80$ and $81$. Since conclusion (b) is false for both nets, \cref{l:goggles-rigidity} guarantees that the two convex quadrilaterals in the respective nets are both equilateral. Note that we have uniquely determined the convex pentagons of both nets, since they have to be isosceles and they have two (nonconsecutive) congruent sides with length which equals the sidelength of the equailateral quadrilaterals. In particular such pentagon is the convex pentagon which appears as two faces of the one element of the goggles family for which the angle at the pins is $\frac{2\pi}{3}$ (this is also the same pentagon of the pentagonal faces of one of Heppes' nets in \cite{Heppes}). We therefore conclude that:
\begin{itemize}
    \item the marked spherical graph 681 has a unique metrization which must be precisely this element of the goggles family; on the other hand this would imply that the vertices $8$ and $9$ are the pins of the net, but the rules for the once-derivation would have removed both, and thus we would have been left with another marked spherical graph (in fact the marked spherical graph 1029);
    \item 158 has a unique metrization which coincides with the latter element of the goggles family with only one pin (and the corresponding two edges) removed; but again this is in contradiction with the derivation rules, which require {\em both} pins to be removed.
\end{itemize}

\subsection{1100 is a cube} Here we analyze the marked spherical graph numbered 1100 and prove the following proposition.

\begin{proposition}\label{p:fake-cube}
Assume $\mathcal{N}$ is the link of a cone $E$ as in \cref{t:main} and $\mathcal{N}'$ a net once-derived from $\mathcal{N}$ whose underlying marked spherical graph is numbered 1100 in the list of \cref{p:second-topological-classification}. Then:
\begin{itemize}
\item[(a)] $E$ is an element in \cref{lista}, namely the cone over the $1$-skeleton of a cube;
\item[(b)] $\ell$ goes through two vertices of the cube, namely the two pins have degree $3$;
\item[(c)] $\mathcal{N}'$ is obtained by removing from $\mathcal{N}$ two edges $\sigma$ and $\sigma'$. The two edges $\sigma$ and $\sigma'$ are adjacent to the pins $P$ and $P'$, and to two other vertices which are {\em not antipodal}.
\end{itemize}
\end{proposition}

We start by drawing the marked spherical graph of the net $\mathcal{N}'$ in \cref{f:falso-cubo} below.  Let us define:
\begin{itemize}
\item the point $F'$ to be the midpoint of the edge $20$;
\item the equator $\gamma$ to be the one passing through $F'$ and orthogonal to the edge $20$;
\item the point $F$ to be the antipodal point to $F'$ on the equator $\gamma$.
\end{itemize}
Note that in the picture the point $F$ is also the midpoint of the edge $57$ and that the whole net is central-symmetric with respect to $F$. These two facts are proven in the following lemma. Before stating it, we recall our convention of denoting by $\lambda_Q$ the length of the edges of an equilateral equiangular spherical quadrilateral with angles $\frac{2\pi}{3}$. A posteriori we must conclude that the length of all the edges of the nets $\mathcal{N}'$ and $\mathcal{N}$ is, in fact, precisely $\lambda_Q$.

\begin{lemma}\label{l:simmetria-uno}
If $\mathcal{N}'$ is as in \cref{p:fake-cube} then the following holds:
\begin{itemize}
\item[(i)] if $\lambda$ is the length of an edge of $\mathcal{N}'$, then $|\lambda-\lambda_Q|< 2\cdot 10^{-4}$;
\item[(ii)] the edges $20$, $63$, and $14$ are congruent;
\item[(iii)] the edges $12$, $40$, $03$ and $26$ are congruent;
\item[(iv)] $\gamma$ intersects $14$ and $63$ orthogonally in the respective midpoints;
\item[(v)] The central symmetry $s$ with respect to the line through $F$ and $F'$ maps $6$ into $4$, $1$ into $3$, $2$ into $0$, $7$ into $5$;
\item[(vi)] $F$ is the midpoint of the edge $57$;
\item[(vii)] If $26$ and $63$ are congruent, then the conclusion of \cref{p:fake-cube} follows.
\item[(viii)] If \cref{p:fake-cube} is false, then the pins of $\mathcal{N}$ are symmetric with respect to the line through $F$ and $F'$. In particular their distance to $F$ is $\frac{\pi}{2}$.
\end{itemize}
\end{lemma}

\begin{figure}
\begin{tikzpicture}
\draw[dashed] (-6,0) -- (6,0);
\node[right] at (6,0) {$\gamma$};
\draw[fill] (0,0) circle (0.07);
\node[below right] at (0,0) {$F$};
\draw[thick] (-4.7,-2) to [out=110, in=-90] (-5,0) to [out=90,in=-110] (-4.7,2) to [out=20, in=160] (-3.3,2) to [out=-70, in=90] (-3,0) to [out=-90, in=70] (-3.3,-2) to [out=-160, in=-20] (-4.7,-2);
\draw[thick] (4.7,-2) to [out=70, in=-90] (5,0) to [out=90,in=-70] (4.7,2) to [out=160, in=30] (3.3,2) to [out=-110, in=90] (3,0) to [out=-90, in=110] (3.3,-2) to [out=-20, in=-160] (4.7,-2);
\draw[thick] (-3.3,-2) to [out=-20, in= -170] (-1.5,-2.2) -- (1.5,2.2) to [out=10, in=160] (3.3,2);
\node[above left] at (-4.7,2) {$2$};
\node[above right] at (4.7,2) {$2$};
\node[below left] at (-4.7,-2) {$0$};
\node[below right] at (4.7,-2) {$0$};
\node[above right] at (-3.3,2) {$6$};
\node[below left] at (3.3,-2) {$4$};
\node[below] at (-3.3,-2) {$3$};
\node[above] at (3.3,2) {$1$};
\node[below right] at (-1.5,-2.2) {$5$};
\node[above left] at (1.5,2.2) {$7$};
\draw[fill] (-5,0) circle (0.07);
\draw[fill] (5,0) circle (0.07);
\node[below left] at (-5,0) {$F'$};
\node[below right] at (5,0) {$F'$};
\end{tikzpicture}
\caption{A diagram describing the marked spherical graph 1100. In the metric realization the angles at $6$ and $7$ facing upwards measure $\frac{4\pi}{3}$, as the angles at $5$ and $1$ facing downwards.} \label{f:falso-cubo}
\end{figure}
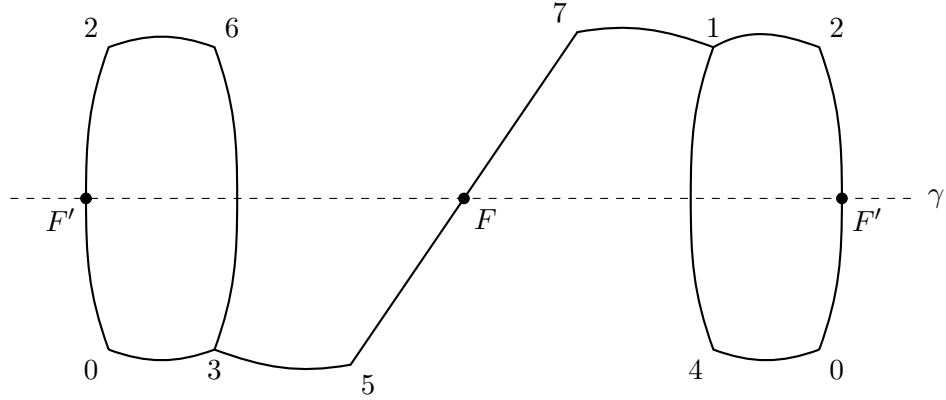

\begin{proof} The statement (i) of the lemma is a consequence of \cref{thm:part2}. The statements (ii), (iii), and (iv) are immediate from the fact that the two quadrilaterals $1204$ and $2603$ are equiangular with angles $\frac{2\pi}{3}$ and share the edge $20$.

Concerning (v), the fact that $s(3)=1$, $s(6)=4$, and $s(2)=0$ follows immediately from (ii) and (iii).
Note that in order to show that $s(5)=7$ we just need to show that the edges $35$ and $71$ are congruent. In order to prove the latter, consider the hexagons with vertices $126357$ and $304175$.
Note that the angles at the vertices (in the given order) form the same sequence $\frac{2\pi}{3}$, $\frac{2\pi}{3}$, $\frac{4\pi}{3}$, $\frac{2\pi}{3}$, $\frac{2\pi}{3}$, $\frac{4\pi}{3}$. Moreover $30$ is congruent to $12$, $04$ to $26$, and $63$ to $41$, thanks to (ii) and (iii). In particular, by \cref{thm:congruence} (which is applicable because by (i) all the segments in the polygonal paths are shorter than $\pi$) $35$ is congruent to $17$. We have therefore shown (v).

As for (vi), since $s(5)=7$ and $s(7)=5$, and the segment $57$ is shorter than $\pi$, $s$ must map it into itself. In particular the midpoint of $57$ is a fixed point of the symmetry $s$. The latter has only two fixed points: $F$ and $F'$. But $F'$ does not belong to $57$.

Under the assumption of (vii), namely that $26$ and $63$ are congruent, we immediately see that $26$, $63$, and $16$ are all congruent and their length is $\lambda_Q$.

Take then the spherical link of the cone over a cube and remove two edges $\sigma$ and $\sigma'$ from this link to obtain a spherical net $\tilde{\mathcal{N}}$.
We impose, however, if $P$ and $A$ are the vertices connected by $\sigma$ and $P'$ and $A'$ the ones connected by $\sigma'$, then $P'$ is antipodal to $P$, but $A'$ is not antipodal to $A$. We see immediately that there is a face $F$ of $\mathbb S^2\setminus \tilde{\mathcal{N}}$ which is an hexagon with the same internal angles as $126357$. Having proven that $36$, $62$, and $21$ have all length $\lambda_Q$ allows us to apply \cref{thm:congruence} to infer that in fact $126357$ is congruent to $F$. This shows that the rigid motion which maps $126357$ onto $F$ maps in fact the net $\mathcal{N}'$ onto $\tilde{\mathcal{N}}$, thus showing \cref{p:fake-cube}.

Finally, coming to (viii), note first that the pins belong to the four equators bisecting the angles in $4$, $5$, $6$, and $7$. At least the ones bisecting $5$ and $6$ must be distinct because of (i).
Therefore the pins $P_1$ and $P_2$ are the only pair of points common to all four bisectors mentioned above. But the symmetry $s$ maps the bisector of $6$ in that of $4$ and the bisector of $7$ in that of $5$. So $s(\{P_1,P_2\})=\{P_1,P_2\}$. But since $s$ is an involution, either $s(P_1) = P_2$, or $s(P_1)=P_1$. The latter option is however not possible, because the only fixed points of $s$ are $F$ and $F'$, which are contained in edges of $\mathcal{N}'$, while the pins cannot.
\end{proof}

We are now ready to close the argument of the proof by contradiction in \cref{p:fake-cube}. We derive six polynomial equations in the sines and cosines of $6$ unknown geometric quantities, which are lengths of edges or measurements of angles.
The list of the geometric quantities is given below in \eqref{e:geom1}--\eqref{e:geom6}, while the six equations are given by \eqref{e:fede-eq-1}, \eqref{e:fede-eq-2}, \eqref{e:fede-eq-3}, \eqref{e:fede-eq-4}, \eqref{e:fede-eq-5}, and \eqref{e:fede-eq-6}.
Using that $\sin^2+\cos^2=1$, we obtain a system of $12$ polynomial equations in $12$ unknowns. Next note that \cref{l:simmetria-uno}(i)  specifies a certain range for the $12$ unknowns. Moreover, if we assume by contradiction that \cref{p:fake-cube} does not hold, this range is further restricted to exclude the cases in which the segments $26$ and $63$ are congruent.
Using a computer algebra software, we will see that there is no solutions of the system in the range.

\subsubsection{Trigonometric identities}
In order to describe the six geometric quantities and the triangles used to relate them we refer to \cref{f:falso-cubo-2} below. The solid segments are all part of the net $\mathcal{N}'$, cf. \cref{f:falso-cubo}. The dashed segments are all geodesics.

While we keep the notation $F$ for the same point in both figures, in the sequel we use $B$, $C$ and $X$ for the three vertices labeled $0$, $3$ and $5$ in \cref{f:falso-cubo}. We extend the edges $20$ and $63$ along the respective equators on the sides of $0$ and $3$ until they meet in the point $A$.

Note that, since the edges $20$ and $63$ intersect the equator $\gamma$ orthogonally, the point $A$ is the pole of the lower hemisphere with boundary $\gamma$, more precisely it is the point equidistant from all points of $\gamma$ which lies in the same connected component of $\mathbb S^2\setminus \gamma$ of the vertices $0$ and $3$.
The dashed segments with extrema $5$ and $6$ are on the bisectors of the corresponding of the angles at the corresponding vertices of $\mathcal{N}'$. Their second common point, labeled $E$, is thus one of the pin, in fact the pin lying in the northern hemisphere (this is justified recalling \cref{l:simmetria-uno}: our net $\mathcal{N}$ is close to the cubical net and so the pin is close to the edge joining $6$ and $7$, which lies on the northern hemisphere). The segment $AE$ intersects $\gamma$ orthogonally at some point $O$.

In order to justify the accuracy of \cref{f:falso-cubo-2} we need to ensure that the segment $AE$ intersects the segment $CX$ in a point $D$ which is interior to both. Recall that we are working under the assumption that the segments $26$ and $63$ have nearly the same length. When they are congruent the point $E$ coincides with $6$, the point $D$ coincides with $C$ and $XF$ is orthogonal to $\gamma$. In particular it follows that in our situation $E$ is close to $6$ and $D$ is close to $C$. Moreover, observe that $E$ cannot be inside the quadrilateral (because it is a pin) and cannot be on its boundary (because then it would coincide with $6$, which would be a pin, and the net $\mathcal{N}$ would be the cubical net). Hence $D$, which is the symmetric of $E$ with respect to the plane containing $\gamma$, is strictly inside $CX$

The six geometric quantities involved in the system are:
\begin{align}
&\mbox{(the size of) the angle $B\hat{A}C$;}\label{e:geom1}\\
&\mbox{(the length of) the segment $AC$;}\label{e:geom2}\\
&\mbox{the angle $A\hat{D}C$;}\label{e:geom3}\\
&\mbox{the segment $AD$;}\label{e:geom4}\\
&\mbox{the segment $EX$;}\label{e:geom5}\\
&\mbox{the angle $D\hat{E}X$.}\label{e:geom6}
\end{align}
By the cosine rule on the triangle $0AC$, we have
\begin{equation}\label{e:fede-eq-1}
\cos \frac{\pi}{3} = - \cos\frac\pi3 \cos B\hat{A}C + \sin \frac{\pi}{3} \sin B\hat{A}C \cos AC\, .
\end{equation}
By the sine rule on the triangle $ACD$, we have
\begin{equation}\label{e:fede-eq-2}
\frac{\sin A\hat{D}C}{\sin AC} = \frac{ \sin \frac{\pi}{3}}{\sin AD }\, .
\end{equation}
By \cref{l:simmetria-uno}(viii), $EF=\tfrac\pi2$ and also $FA=\tfrac\pi2$. In particular, the geodesic connecting $A$ and $E$ is a subset of the equator at distance $\tfrac\pi2$ from $F$ and thus $OF=OF'=\tfrac\pi2$.
Moreover, $0\hat{A}D=OF'$, because $\gamma$ consists is the equator at distance $\tfrac\pi2$ from $A$. Hence, $C\hat{A}D=\frac{\pi}{2}-B\hat{A}C$, in particular $\cos C\hat{A} D = \sin B\hat{A}C$ and $\cos C\hat{A} D = \sin B\hat{A}C$.
By the cosine rule on the triangle $ACD$, we have
\begin{equation}\label{e:fede-eq-3}
\cos \frac{\pi}{3} = -\cos B\hat{A}C \sin A\hat{D}C \cos AD + \sin B\hat{A}C \cos A\hat{D}C\, .
\end{equation}
Observe that $E\hat{D}X=A\hat{D}C$ and $DE=\pi-2\,AD$ so $\sin DE = 2 \sin AD\cos AD$. By the sine rule on the triangle $DEX$, we have
\begin{equation}\label{e:fede-eq-4}
\frac{\sin A\hat{D}C}{\sin EX} = \frac{\sin \frac{\pi}{3} }{2\sin AD \cos AD } \, .
\end{equation}
By the cotangent rule on the triangle $DEX$, we have
\begin{equation}\label{e:fede-eq-5}
(\sin^2 AD - \cos^2 AD) \cos D\hat{E}X
= 2 \cot XE \sin AD \cos AD - \cot A\hat{D}C \sin D\hat{E}X\, .
\end{equation}
Since $X\hat EF=\tfrac\pi2 - D\hat EX$, by the cotangent rule on the triangle $EFX$, we have
\begin{equation}\label{e:fede-eq-6}
\cos EX \sin D\hat{E}X= - \cot \frac{\pi}{3} \cos D\hat{E}X\, .
\end{equation}

\begin{figure}
\begin{tikzpicture}
\node[below left] at (-4.7,-2) {$0=B$};
\node[above left] at (-4.7, 2) {$2$};
\node[below left] at (-5.2, 0) {$F'$};
\draw[dashed] (-6,0) -- (3,0);
\node[right] at (3,0) {$\gamma$};
\draw[thick] (-4.7,-2) [out=120, in=-90] to (-5.2,0) [out=90, in=-120] to (-4.7,2);
\draw[thick] (-4.7,2) [out=20, in=180] to (-3.9,2.2) [out=0, in=160] to (-3.1, 2);
\draw[thick] (-4.7,-2) [out=-20, in=180] to (-3.9,-2.2) [out=0, in=-160] to (-3.1, -2);
\draw[thick, dashdotted] (-4.7,-2) [out =-60, in=160] to (-2, -4);
\draw[thick, dashdotted] (-3.1,-2) [out=-80, in=+130] to (-2,-4);
\node[below] at (-2,-4) {$A$};
\draw[thick] (-3.1,-2) [out=100, in=-90] to (-3.3,0) [out=90, in=-100] to (-3.1,2);
\draw[loosely dashdotted, thick] (-2,-4) -- (-2,3);
\node[above right] at (-2,-3) {$D$};
\draw[thick, dashed] (-3.1,2) -- (-2,3);
\draw[thick] (-3.1,-2) [out=-45, in=170] to (-0.3,-3.5);
\node[below] at (-0.3,-3.5) {$5=X$};
\draw[thick] (-0.3,-3.5) [out=60, in=-100] to (1,0);
\draw[thick, dashed] (-2,3) [out=-80, in=110] to (-0.3,-3.5);
\draw[thick, loosely dashdotted] (-2,3) [out=-20, in=100] to (1,0);
\node[above] at (-2,3) {$E$};
\node[above right] at (1,0) {$F$};
\node[above] at (-3.1,2) {$6$};
\node[above left] at (-3.1,-2) {$3=C$};
\node[below left] at (-2,0) {$O$};
\end{tikzpicture}
\caption{The picture referenced in the concluding argument for \cref{p:fake-cube}.} \label{f:falso-cubo-2}
\end{figure}

\subsubsection{Conclusion via Gr\"obner basis}
The statements that follow are supported by the python script \verb|scripts/almost_cube_computation.py| in \cite{repo} (the output of the script can be found in \verb|scripts/almost_cube_computation_output.txt|), whose execution takes a few seconds.

Observe that the three equations \cref{e:fede-eq-4}, \cref{e:fede-eq-5}, and \cref{e:fede-eq-6} can be reduced to a polynomial system in the sines and cosines of $A\hat D C$, $EX$, $AD$, $D\hat EX$.
Through the computation of a Gr\"obner basis, this system is shown to imply
\begin{equation}\label{eq:local89}
    xy(x - 1)(x + 1)(y - 1)(y + 1)(4x^2y^2 - 4x^2 + 1)(4x^2y^2 - 4x^2 - 4y^2 + 3) = 0,
\end{equation}
where $x=\cos AD$ and $y=\cos A\hat DC$.
Thanks to \cref{l:simmetria-uno}(i), we know that $AD\approx \tfrac{\lambda_T}2$ and $A\hat D C \approx \tfrac{2\pi}3$. Therefore, the only factor of \cref{eq:local89} that can be zero is
\begin{equation*}
4x^2y^2-4x^2+1 = -4\cos^2 AD\sin^2A\hat DC+1=-(2\cos AD\sin A\hat DC-1)(2\cos AD\sin A\hat DC+1)
\end{equation*}
Using again \cref{l:simmetria-uno}(i), we deduce that
\begin{equation}\label{e:fede-eq-7}
    2\cos AD\sin A\hat DC=1 \, .
\end{equation}
Now, observe that the four equations \cref{e:fede-eq-1}, \cref{e:fede-eq-2}, \cref{e:fede-eq-3}, \cref{e:fede-eq-7} can be reduced to a polynomial system in the sines and cosines of $AD, A\hat DC, B\hat AC, AC$. Through the computation of a Gr\"obner basis, this system is shown to imply $\cos B\hat AC=0$. This implies that the quadrilateral $2630$ is equilateral, which is sufficient to conclude thanks to \cref{l:simmetria-uno}(vii).

\subsection{920 is a thin fish and 827 is a fat fish} In this section we conclude the proof of \cref{t:human-classification} analyzing the two remaining marked spherical graphs, those numbered 920 and 827 in the list of \cref{p:second-topological-classification}. The statement is summarized in the following proposition

\begin{proposition}[Classification of thin and fat fish]\label{p:fishes}
The marked spherical graph 920 is only compatible with once-derived nets of the thin fish family, while the marked spherical graph 827 is only compatible with once-derived nets of the fat fish family.
\end{proposition}

\begin{figure}
\begin{tikzpicture}[scale=0.7]
\draw[very thick] (-6,0) -- (-2,0) node[above left] {$1$} -- (0,2) node[ above] {$3$} -- (2,0) node[above right] {$0$} -- (0,-2) node[below] {$2$} -- (-2,0);
\draw[very thick] (2,0) -- (6,0);
\end{tikzpicture}
\caption{A picture of the marked spherical graph 920 for a once-derived net $\mathcal{N}'$. The quadrilateral $1302$ is convex.}\label{f:pesce-magro-class}.
\end{figure}
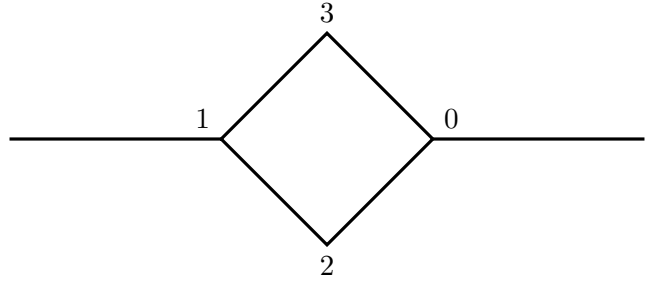

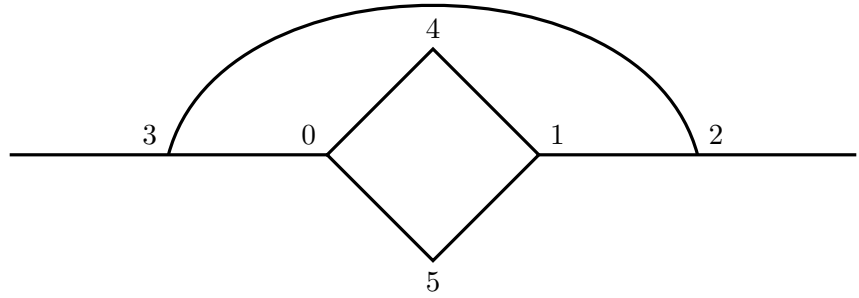
\begin{figure}
\begin{tikzpicture}[scale=0.7]
\draw[very thick] (-8,0) -- (-5,0) node[above left] {$3$} -- (-2,0) node[above left] {$0$} -- (0,2) node[ above] {$4$} -- (2,0) node[above right] {$1$} -- (0,-2) node[below] {$5$} -- (-2,0);
\draw[very thick] (2,0) -- (5,0) node[above right] {$2$} -- (8,0);
\draw[very thick] (-5,0) to [bend left=75] (5,0);
\end{tikzpicture}
\caption{A picture of the marked spherical graph 827 for a once-derived net $\mathcal{N}'$. The quadrilateral $0415$ is convex and the two vertices $2$ and $3$ are joined by two edges and must, therefore, be antipodal points.}\label{f:pesce-grasso-class}.
\end{figure}

\begin{proof}
We examine first the marked spherical graph 920 and we refer to the schematic picture in \cref{f:pesce-magro-class}. Observe that the the triangle $031$ has angles $\frac{2\pi}{3}$, $\frac{2\pi}{3}$ and $\frac{4\pi}{3}$. Note moreover that the two edges $31$ and $30$ must be shorter than $\pi$ because they are in common with a convex quadrilateral. From \cref{thm:generalized_triangles} the lengths of the edges of the triangle are uniquely determined. Note moreover that the edges $31$ and $30$ are congruent. This in particular implies that the convex quadrilateral $Q$ with vertices $3210$ is equilateral. This determines completely the metric of the net and it can be easily checked that it is the metric of the once-derived net which one obtains from any member of the thin fish family for which both pins have degree $1$. 

We analyze next the marked spherical graph 827 and refer to the schematic picture in \cref{f:pesce-grasso-class}. Observe that the vertices $2$ and $3$ are connected by two edges. These two edges must therefore be long exactly $\pi$, while the vertices $2$ and $3$ are antipodal. The edges connecting $2$ and $3$ form therefore a convex biangle $\mathcal B$ of opening $\frac{2\pi}{3}$. Observe in particular that the edges $12$ and $30$ must lie on the bisector $\gamma$ of the angles of $\mathcal B$. If we remove the two edges of the biangle and add the half equator $\gamma\cap \mathcal B$, the resulting net induces the marked spherical graph 920 and therefore we can follow the argument above to prove that the only nets for which the once-derivation procedure would produce the net 827 are the members of the fat fish family for which both pins have degree $2$.
\end{proof}

\subsection{The marked spherical graph 736}\label{s:736} In this section we prove \cref{prop:736}. The proof is computer assisted: first we use the computer program explained in \cref{s:secondo-programma} to show that the lengths of {\em some} specific edges are close to some values; then synthetic arguments, helped by some algebraic computations performed with a computer algebra software, will finish the proof. We recall also that, as it is explained in \cref{s:secondo-programma}, with the additional information gained from \cref{prop:736} an additional run of the code described \cref{s:secondo-programma} will then determine the lengths of {\em all} the edges of any valid metric realization, up to a small error.

\begin{figure}[hbtp]
\begin{tikzpicture}[scale=0.7]
\draw[very thick] (-6,0) -- (-4.5,0) node[above left] {$2$} -- (-2.5,2) node[above left] {$3$} -- (-0.5,0) -- (0.5,0) -- (2.5,2) node[above right] {$9$} -- (4.5,0) node[above right] {$5$} -- (7,0) node[above left] {$6$} -- (9,2) node[above] {$13$} -- (11,0) node[above right] {$7$} -- (12,0);
\draw[very thick] (-4.5,0) -- (-2.5,-2) node[below left] {$1$} -- (-0.5,0);
\draw[very thick] (-2.5,-2) [out=-90,in=150] to (0,-5.2) node[below] {$11$} [out=30,in=-90] to (2.5,-2);
\draw[very thick] (0.5,0) -- (2.5,-2) node[below right] {$4$} -- (4.5, 0);
\node[right] at (0.8,0) {$8$};
\node[left] at (-0.8,0) {$0$};
\draw[very thick] (-2.5, 2) [out=90,in=210] to (0,5.2) node[above] {$12$} [out=-30,in=90] to (2.5,2);
\draw[very thick] (11,0) -- (9,-2) node[below] {$10$} -- (7,0);
\end{tikzpicture}
\caption{A picture of the marked spherical graph 736 for a derived-once net $\mathcal{N}$. Thanks to the algorithm described in \cref{ss:algo2presentation} we have an \emph{almost} rigidity for all edges apart from $56$ and $72$.}\label{f:736}
\end{figure}
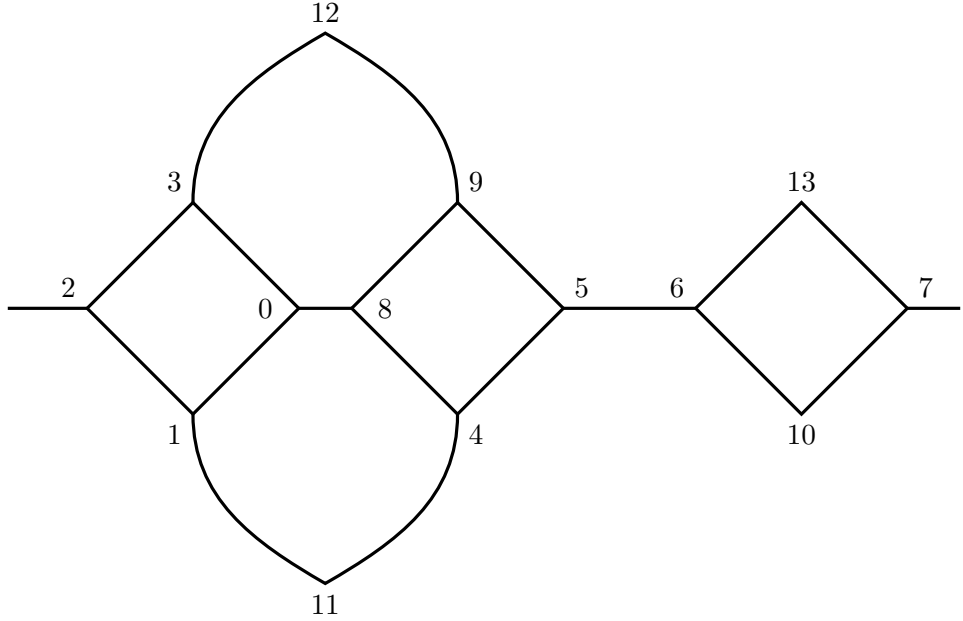

\begin{proposition}\label{prop:736}
Let $\tilde{\mathcal{N}}$ be the link of a minimizing cone as in \cref{t:main} and assume that the once-derived marked spherical graph associated to it is numbered 736, as represented in \cref{f:736}. Then the edges $56$ and $72$ are congruent.
\end{proposition}

\begin{proof}
We denote by $\mathcal N'$ the once-derived net obtained from $\tilde{\mathcal N}$.

\medskip {\noindent \bf Step 1.} We denote by $\lambda_Q$ the length of the edges of a spherical equilateral quadrilateral with angles $\frac{2\pi}{3}$. Consider a pentagon with angles $\frac{2\pi}{3}$ and vertices labeled clockwise as $A,B,C,D$, and $E$. If the length of $DE$ and $BC$ is $\lambda_Q$, \cref{p:formulae} shows that the lengths of the remaining three edges are determined and that the sides $AE$ and $AB$ are congruent. We denote by $\lambda_B$ the length of $CD$ and by $\lambda_U$ the length of $AE$ and $AB$. Upon setting $\boldsymbol{\delta} = 2\cdot 10^{-4}$, we use the algorithm explained in \cref{ss:algo2presentation} (cf. \verb|embed_or_exclude.cpp| in \cite{repo}) to prove:
\begin{itemize}
\item[(a)] The lengths of the twelve edges $32$, $30$, $01$, $12$, $89$, $95$, $54$, $48$, $6(13)$, $(13)7$, $7(10)$, and $(10)6$ belong to the interval $[\lambda_Q-\boldsymbol{\delta}, \lambda_Q+\boldsymbol{\delta}]$.
\item[(b)] The lengths of the four edges $3(12)$, $(12)9$, $4(11)$, and $(11)1$ belong to the interval $[\lambda_H-\boldsymbol{\delta}, \lambda_H+\boldsymbol{\delta}]$;
\item[(c)] The length of the edge $08$ belongs to $[\lambda_B - \boldsymbol{\delta}, \lambda_B + \boldsymbol{\delta}]$.
\end{itemize}
We also record a simple consequence of the above estimates and the smallness of $\boldsymbol{\delta}$:
\begin{itemize}
\item[(d)] The distance between the vertices $5$ and $2$ is strictly less than $\pi$ and the unique shortest geodesic segment joining them does not intersect the 4 convex faces $1032$, $8459$, $089(12)3$, $084(11)1$.
\end{itemize}
This happens because the diagonal $20$ and $85$ would be aligned with the edge $08$ if $\boldsymbol{\delta}$ would be zero, making the union of these three segments a single geodesic arc of length $\lambda_B + 2 \lambda_T > \pi$ (here $\lambda_T$ denotes the length of the diagonal of the equilateral equiangular quadrilateral with angles $\frac{2\pi}{3}$, which equals the length of the sides of the equiangular triangle with angles $\frac{2\pi}{3}$).

\medskip
{\noindent\bf Step 2.} We will use geometric considerations and algebraic computations to show that
\begin{itemize}
\item[(e)] $89$, $45$, $01$, and $23$ are congruent;
\item[(f)] $95$, $84$, $30$, and $21$ are congruent;
\item[(g)] $3(12)$ and $4(11)$ are congruent;
\item[(h)] $9(12)$ and $1(11)$ are congruent.
\end{itemize}
In fact, it suffices to show that $30$ is congruent to $84$, as the remaining congruences are then due to the following facts, recorded in \cref{p:formulae}: opposite edges in an equiangular quadrilateral are congruent; a quadrilateral with internal angles all equal to $\frac{2\pi}{3}$ is determined by the length of one edge; a pentagon with internal angles all equal to $\frac{2\pi}{3}$ is determined by the lengths of two consecutive edges.

In order to show our claim, we denote by $x,y$, and $z$ the values of $\tan \frac{\lambda}{2}$ when $\lambda$ is the length of, respectively, $30$, $08$, and $95$. Recalling \cref{p:formulae}, $\tan \frac{\lambda}{2}$ equals then $\frac{1}{2x}$ and $\frac{1}{2z}$ if $\lambda$ is the length of, respectively, $01$, and $84$. The congruence of $30$ and $84$ is thus equivalent to
\begin{equation}\label{e:relazione-di-congruenza}
2xz =1 \, .
\end{equation}

Using the formulae in \cref{p:formulae} for the edges of the pentagon, some classical trigonometric identities, and computing the Gr\"obner basis of an algebraic system we deduce the following lemma (see \verb|736_computation.py| and \verb|736.sing| in \cite{repo} for the details).
\begin{lemma}
    Let $a$, $b$, $c$ be the length of three consecutive (in this order) sides of a convex spherical pentagon. Then $Q(\tan\tfrac a2, \tan \tfrac b2, \tan \tfrac c2)=0$, where $Q$ is the polynomial
    \begin{equation*}
        Q(x, y, z) =
        x^2y^2z^2 - 2y^2(x^2+z^2) + x^2z^2 - 2(x^2+z^2) - 6xy^2z + 6y(x+z) + 6xz + y^2 - 2 \,.
    \end{equation*}
\end{lemma}
Applying this lemma to the two pentagons $089(12)3$ and $084(11)1$ we infer that $Q (x,y,z) = 0$ and $Q(\frac{1}{2x}, y, \frac{1}{2z})=0$. The second identity can also be recast as a polynomial equation and this enables us to eliminate the variable $y$ and compute the resultant polynomial equation $R(x, z)=0$ in $x$ and $z$. The polynomial $R$ can be factored as
\begin{equation*}
    R(x, z) = 9(2xz-1)^2 S(x, z) \,,
\end{equation*}
where $S(x, z)$ is a polynomial of degree $12$ with integer coefficients (see again \verb|736_computation.py| in \cite{repo} for the details).
Either $2xz-1$, which is what we want, or the factor $S(x, z)$ vanishes. Thanks to (a), we know that $x, z$ are (quantitatively) close to $\tan\tfrac{\lambda_Q}2=\tfrac1{\sqrt2}$. We have $S(\tfrac1{\sqrt2}, \tfrac1{\sqrt2}) = -1944$ and furthermore $S(x', z')\not=0$ whenever $|x'-\tfrac1{\sqrt2}|, |z'-\tfrac1{\sqrt2}| < 0.05$. Thus, thanks to (a), we deduce that $S(x, z)\not=0$ and therefore $2xz=1$ as desired.

\medskip
{\noindent \bf Step 3.}
We prove that either the conclusion of the Proposition holds, or
\begin{itemize}
\item[(i)] The points $5$, $6$, $7$, and $2$ are all on the same equator, and the quadrilateral $(13)7(10)6$ is equilateral.
\end{itemize}
Consider the diagonal of the quadrilateral joining the vertices $6$ and $7$, which we denote by $67$. (i) is equivalent to the edge $67$ forming a flat angle with $56$.
Notice that the angle formed by the polygonal path $5\to 6 \to 7$ at $6$ is congruent to that formed by $2\to 7 \to 6$ at $7$. So either they are both $\pi$, and (i) holds, or they are both distinct from $\pi$. In the second case, consider the closed polygonal paths $2 \to 3 \to 0\to 8 \to 9 \to 5 \to 6 \to 7 \to 2$ and $5\to 4 \to 8 \to 0 \to 1 \to 2 \to 7 \to 6$. Note that they have the same internal angles, all distinct from $\pi$, and the first five edges are congruent. If we know that all the edges are strictly shorter than $\pi$, we then conclude from \cref{thm:congruence} that the last three are also congruent, which however implies the desired conclusion that $56$ and $72$ are congruent.

Note that from (a), (b) and (c) we know that the quadrilateral $6(13)7(10)$ is close to be equilateral, so the polygonal path $5\to 6 \to 7 \to 2$ is close to a geodesic arc connecting $5$ and $2$, which by (d) must be shorter than $\pi$. So the fact that the lengths of $56$ and
$72$ are shorter than $\pi$ is guaranteed.

It remains to show the congruence of the segments $56$ and $72$ under assumption (i). We begin by proving that
\begin{itemize}
\item[(j)] the quadrilaterals $3021$ and $8951$ are also equilateral.
\end{itemize}
Consider a new net in which the vertices $6,10,13$, and $7$ are eliminated and all the edges $56$, $6(13)$, $(13)7$, $7(10)$, $(10)6$, and $72$ are replaced by the geodesic segment joining $5$ and $2$ (and containing $6$ and $7$).
Since the edges $3(12)$ and $9(12)$ are both strictly shorter than $\pi$, we can apply \cref{l:goggles-rigidity} and conclude that (j) holds.

Given all the information (a) -- (j) we are ready to conclude. The three convex quadrilaterals in the net $\mathcal{N}'$ have been proved to be equilateral. In particular, their diagonals $31$, $94$, and $(10)(13)$ identify three concurrent equators all passing through $12$ and $11$. We infer that $12$ and $11$ are the two pins (with degree $3$) and thus the three equators must form angles of $\tfrac{2\pi}3$ with each other. We have completely identified the net (which coincides with the Heppes' net $N_6^{(3)}$ in \cite{Heppes}) and the congruence between $56$ and $72$ follows.
\end{proof}

\section{From triangulated competitors to Almgren minimality}\label{s:triangulations-and-maps}

In \cref{s:terzo-programma} we will complete the proof of \cref{t:main} excluding all the possible nets allowed by \cref{thm:part2} and \cref{t:human-classification} which are not listed in \cref{lista}. In order to achieve this goal we will produce a database of suitable ``triangulated $2$-dimensional sets'' which save area compared to the respective cones. Since the set of possibilities to exclude is infinite, we need a suitable criterion to reduce the problem to the existence of a finite database: this will be accomplished in \cref{s:competitore}, where we show that the existence of a competitor, which saves a certain amount of area, for some specific cone $E$ guarantees the existence of a competitor saving area for every sufficiently close cone $E'$.
The material contained in this section is instead meant to free us from the task of producing a map $F$ satisfying the conditions of \cref{d:Almgren-minimizing} from the cone to the competitor. In particular, \cref{t:comp-to-map} gives us a set of sufficient conditions, to be checked on the (combinatorial data) of the triangulated sets, which guarantees the existence of a map as in \cref{d:Almgren-minimizing}. The reader who is familiar with the language of simplicial complexes can skip directly to the statement of \cref{t:comp-to-map}.

\subsection{Preliminaries on simplicial complexes} We start here by detailing our description of the ``triangulated sets'', which are simply images, under piecewise linear maps, of abstract finite $2$-dimensional simplicial complexes.

\begin{definition}[Simplicial complex]\label{d:simplicial-complex}
An abstract finite simplicial complex $\Sigma$ is a collection of subsets of a finite set $\mathcal{V}(\Sigma)=\{V_1, \ldots, V_N\}$ with the properties that:
\begin{itemize}
\item All singletons are in $\Sigma$;
\item if $S \in \Sigma$, then all the subsets of $S$ belong to $\Sigma$.
\end{itemize}
The singletons will be called vertices, while any element $S\in \Sigma$ of cardinality $k+1$ will be called a $k$-dimensional face of $\Sigma$, and any simplex which is a subset of $S$ will be called a face of $S$. A face of $\Sigma$ is called maximal if it is not a face of a face of $\Sigma$ of higher dimension.
$\Sigma$ is said to be $k$-dimensional if it contains at least a $k$-dimensional face and no $(k+1)$-dimensional faces.
\end{definition}

For one-dimensional and two-dimensional faces it is common to use the terms {\em edges} and {\em triangles}. In particular, since we will only deal with one-dimensional and two-dimensional complexes, we will only encounter vertices, edges, and triangles.

\begin{definition}[Connected simplicial complex]\label{d:connectedness}
The complex is connected if for each pair of vertices $A$ and $B$ there is a collection of vertices $A=V_0, V_1, \ldots, V_L = B$ such that $[V_i,V_{i+1}]$ is an edge of the complex for every $0\leq i \leq L-1$.
\end{definition}

If a face $S$ of $\Sigma$ consists of the vertices $A_0, \ldots , A_k$, we will denote it by $[A_0, A_1, \ldots, A_k]$.

\begin{definition}[Geometric realization]
    For a simplicial complex $\Sigma$, we define its geometric realization $|\Sigma|$ as the subset of $\R^{\mathcal V(\Sigma)}$ defined as
    \begin{equation*}
        |\Sigma| := \bigcup_{S\in \Sigma} \text{conv}(\{e_V:\, V\in S\}) \,.
    \end{equation*}
    We equip $|\Sigma|$ and all its subcomplexes with the restriction of
the Euclidean metric of $\mathbb R^{\mathcal V(\Sigma)}$.

    For any face $S\in\Sigma$, we denote by $|S|$ the subspace of $|\Sigma|$ given by $\text{conv}(\{e_V:\, V\in S\})$.
\end{definition}

We remark that the choice of the metric on $|\Sigma|$ is not particularly important: we will only use maps which are Lipschitz and thus any other biLipschitz equivalent metric (e.g. the path metric as subset of $\mathbb R^{\mathcal V(\Sigma)}$) would serve our purposes.

\begin{definition}[Piecewise linear map]\label{d:piecewise-linear}
A piecewise linear map $\Phi$ from a simplicial complex $\Sigma$ is a map $\Phi : \mathcal{V}(\Sigma) \to \mathbb R^N$ for some $N$. Given a face $S=[V_0, \ldots , V_k]$ we define its image $\Phi (S)$ as the subset of $\mathbb R^N$ consisting of the convex combination of the points $\Phi (V_i)$.
$\Phi (\Sigma)$ will denote the union of all $\Phi (S)$ for $S\in \Sigma$.

We can also interpret $\Phi$ as a map from $|\Sigma|$ into $\R^N$ by extending linearly the map $\Phi(e_V):=\Phi(V)$. With a slight abuse of notation, we will use the same symbol $\Phi$ for the three maps on the domains $\mathcal V(\Sigma)$, $\Sigma$, $|\Sigma|$.

For any face $S$, we denote by $\Phi^\circ(S)$ the image of the \emph{interior} of $S$. The latter consists therefore of the \emph{strict} convex combination of the points $\Phi(V_{i})$, i.e.,
\begin{equation*}
    \Phi^\circ(S) := \{\lambda_0 \Phi(V_0) + \cdots + \lambda_k\Phi(V_k):\, \lambda_0 + \cdots + \lambda_k = 1\text{ and }
    \lambda_0, \lambda_1, \dots, \lambda_k > 0\} \subseteq\R^N \,.
\end{equation*}
\end{definition}

\begin{definition}[Piecewise linear embedding]
A piecewise linear map $\Phi$ from a simplicial complex $\Sigma$ is
called an \emph{embedding} if any of the following equivalent
conditions holds:
\begin{itemize}
\item for every pair of faces $S,S'\in\Sigma$,
\[
\Phi(S)\cap\Phi(S')=\Phi(S\cap S'),
\]
and, for every $k$-dimensional face $S\in\Sigma$, the set $\Phi(S)$
is a $k$-dimensional simplex;
\item
\[
\Phi^\circ(S)\cap\Phi^\circ(S')=\emptyset
\]
whenever $S,S'\in\Sigma$ are distinct faces;
\item $\Phi:|\Sigma|\to\mathbb R^N$ is injective.
\end{itemize}
\end{definition}

If $K\subset \mathbb R^N$ and $\Phi (\Sigma) \subset K$, we will use the shorthand notation $\Phi: \Sigma \to K$. Note that, when $K$ is convex, we only need to check that every vertex of $\Sigma$ is mapped into $K$.

The database of competitors will consist of (computer-generated) images of $2$-dimensional simplicial complexes $\Sigma$ through piecewise linear maps $\Phi$ into $\mathbb R^3$.
Finally, we will need the concept of ``collapse'' of an abstract simplicial complex $\Sigma$ onto a subcomplex $\tilde\Sigma$.

\begin{definition}[Collapse]\label{d:collapse}
Consider a simplicial complex $\Sigma$. A subcomplex $\tilde{\Sigma}$ of $\Sigma$ is obtained from the elementary collapse of a cell if there is a $k$-dimensional simplex $S\in \Sigma$ with the following properties:
\begin{itemize}
\item $S$ is maximal;
\item there is a $(k-1)$-dimensional face $F$ of $S$ which is free, namely it is not a face of any other simplex of $\Sigma$;
\item $\tilde\Sigma$ is obtained from $\Sigma$ after removing $S$ and $F$ from $\Sigma$.
\end{itemize}
We say that $\Sigma$ can be collapsed on a subcomplex $\tilde\Sigma$ if there is a sequence of subcomplexes $\tilde \Sigma = \Sigma_N \subset \Sigma_{N-1} \subset \ldots \subset \Sigma_0 = \Sigma$ with the property that, for every $i\geq 1$, $\Sigma_i$ is obtained from the elementary collapse of a cell of $\Sigma_{i-1}$.
\end{definition}

\subsection{Existence of Almgren maps}
We are finally ready to state the main theorem of this section, which will enable us to produce a ``competitor map'' in the sense of Almgren's minimality (cf. \cref{d:Almgren-minimizing}) from images of $2$-dimensional simplicial complexes through piecewise linear maps.

\begin{theorem}\label{t:comp-to-map}
Let $\mathcal{N}\subset \partial B_1$ be a geodesic net (not necessarily stationary) with two antipodal pins $P$ and $-P$ on the line $\ell$, $E$ be the cone over $\mathcal{N}$, $K\subset B_1$ be a closed convex polyhedron containing the origin in its interior, $\Sigma$ be a connected abstract $2$-dimensional simplicial complex and $\Phi: \Sigma \to \mathbb R^3$ a piecewise linear map. Assume in addition that
\begin{itemize}
\item[(i)] $\Phi (\Sigma)\cap \ell = K \cap \ell$ and there is a one-dimensional subcomplex $\Gamma \subset \Sigma$ such that $\Phi$ is an embedding on $\Gamma$ and $\Phi (\Gamma) = (\ell \cap K) \cup (E\cap \partial K)$;
\item[(ii)] There is a subcomplex $L$ of $\Gamma$ that consists of all faces $F$ of $\Sigma$ such that\footnote{Equivalently $|L|=\Phi^{-1}(\ell)$, considering $\Phi$ as a map on $|\Sigma|$.} $\Phi^\circ(F)\cap \ell\not=\emptyset$;
\item[(iii)] $\Sigma$ can be collapsed on the subcomplex $L$.
\end{itemize}
Then for every $\varepsilon>0$ there is a Lipschitz map $F: \mathbb R^3 \to \mathbb R^3$ such that
\begin{itemize}
\item[(a)] If $W:= \{x: F(x)\neq x\}$ then $\overline{W\cup F(W)}\subset B_1\setminus \ell$;
\item[(b)] the following inequality holds:
\begin{equation}\label{e:quasi-ottimo}
\mathcal{H}^2 (F(E\cap B_1)) \leq \mathcal{H}^2 ((B_1\setminus K)\cap E) + \mathcal{H}^2 (\Phi (\Sigma)) + \varepsilon\, .
\end{equation}
\end{itemize}
\end{theorem}

The rest of this section is devoted to proving \cref{t:comp-to-map}. Hence from now on we fix $\Sigma$, $\Phi$, and $\varepsilon$ as in the statement and our goal is to produce the map $F$.
The proof consists of three parts. First, we modify the pair
$(\Sigma,\Phi)$ near the axial subcomplex $L$.
The resulting complexes retain the collapse property, their images gain only a small amount of area, and provide more control on the combinatorial/geometric structure in a neighborhood of $L$.
This is accomplished in \cref{s:simplicial-complexes}.
Second, in \cref{s:controlled-flows} we show how collapsibility implies the existence of a controlled homotopy from a complex onto a subcomplex it can be collapsed into.
We construct this homotopy for the simplicial complex built in \cref{s:simplicial-complexes} and use it, after some nontrivial modification, to define the map $F$ required in the statement of \cref{t:comp-to-map}.

\subsection{New simplicial complexes and new piecewise linear maps}\label{s:simplicial-complexes}
In this subsection we modify the simplicial complex $\Sigma$ and the map $\Phi$ in two steps.
Let us denote by $\Gamma_{bd}$ the subcomplex of those vertices and edges of $\Gamma$ which are mapped by $\Phi$ into $E\cap \partial K$. Let $Q_+, Q_-$ be the two ``pins'', i.e., the vertices belonging to both $L$ and $\Gamma_{bd}$.

We produce two new simplicial complexes $(\Sigma^1, \Phi^1)$ and $(\Sigma^2, \Phi^2)$.
\begin{enumerate}
    \item $(\Sigma^1, \Phi^1)$ is obtained from $(\Sigma, \Phi)$ as a subdivision that guarantees that all edges with one vertex on the line are short.
    \item $(\Sigma^2, \Phi^2)$ is obtained from $(\Sigma^1, \Phi^1)$ adding a thin collar along the boundary $\Gamma^1_{bd}$ to have a better control of the combinatorial structure around the two pins.
\end{enumerate}

\medskip

{\bf First perturbation: refinement to $\Sigma^1$ and $\Phi^1$.}
The construction depends on a small parameter $\eta>0$ which will be chosen at the end.
Consider the cylinder $B_\eta (\ell)=\{p: \dist (p, \ell) < \eta\}$, namely the tubular neighborhood of $\ell$ of radius $\eta$. Notice that $\partial B_\eta (\ell)$ consists of the points at distance exactly $\eta$ from $\ell$.

Consider the set $\mathcal{E}_\eta$ which contains the edges $\sigma$ of $\Sigma$ such that $\Phi(\sigma)$ intersects $\partial B_\eta (\ell)$.
The set $\mathcal{E}_\eta$ is the same for all sufficiently small values of $\eta$ and contains only edges with exactly one vertex in $L$. In fact consider the set $\mathcal{E}_L$ of edges with exactly one vertex in $L$ and one vertex not in $L$. Any other edge is either an edge of $L$, and hence is mapped in $\ell \subset B_\eta (\ell)$ or an edge $\sigma$ with $\Phi (\sigma)\cap \ell = \emptyset$. Since there are finitely many edges, it suffices to choose $\eta$ sufficiently small to ensure that $\mathcal{E}_\eta = \mathcal{E}_L$. Furthermore, under the latter assumption the intersection $\Phi(\sigma)\cap \partial B_\eta (\ell)$ consists of a single point which belongs to $\Phi^\circ (\sigma)$.

To proceed, let us introduce a standard subdivision for simplicial complexes.
\begin{definition}[Stellar Subdivision]
    Consider a simplicial complex $\Sigma$ and a face $F$ of $\Sigma$. The stellar subdivision of $\Sigma$ at $F$ is the new simplicial complex $\Sigma'$ built as follows:
    \begin{itemize}
        \item The vertices of $\Sigma'$ are the vertices of $\Sigma$ and an additional new vertex $v_F$ (which is, intuitively, \emph{inside} $F$).
        \item The faces of $\Sigma'$ are the faces of $\Sigma$ which are not a superset of $F$ and the following additional family of faces (all containing the new vertex $v_F$):
        \begin{equation*}
            \Big\{\{v_F\}\sqcup B\setminus A:\, \emptyset\not= A\subseteq  F\subseteq B\in\Sigma \Big\}\,.
        \end{equation*}
    \end{itemize}
\end{definition}

Let us build $\Sigma^1$ by performing successive stellar subdivisions on the edges $\sigma$ of $\Sigma$ such that $\Phi(\sigma)$ intersects $\partial B_\eta (\ell)$. The order does not matter for our purposes (even though the actual resulting $\Sigma^1$ depends on the order). Observe that, once the order is specified independently of $\eta$, for every sufficiently $\eta$ the abstract complex $\Sigma^1$ is independent of $\eta$.

Let us build the map $\Phi^1$ so that it coincides with $\Phi$ on the shared vertices, while for each new vertex $v_\sigma$, where $\sigma\in \mathcal{E}_\eta$, we set $\Phi^1(v_\sigma)$ to be the only element of the singleton $\Phi(\sigma)\cap  \partial B_\eta (\ell)$.

Observe that a stellar subdivision of a simplicial complex acts on any subcomplex as a stellar subdivision (or as an empty operation). Therefore, we can define canonically $L^1, \Gamma^1, \Gamma^1_{bd}$ as the subcomplexes of $\Sigma^1$ corresponding to $L, \Gamma, \Gamma_{bd}$ respectively.

Then the following lemma holds.

\begin{lemma}\label{lem:properties_simplicial_one}
    The pair $(\Sigma^1, \Phi^1)$ satisfies all the properties satisfied by $(\Sigma, \Phi)$ in the assumptions of \cref{t:comp-to-map}, the image $\Phi^1(\Sigma^1)$ coincides with the image $\Phi(\Sigma)$, and:
    \begin{itemize}
        \item The number of edges of $\Sigma^1$ adjacent to a vertex in $L^1$ does not depend on $\eta$ (if $\eta$ is sufficiently small).
\item Every vertex of $\Sigma^1$ that does not belong to $L^1$ but is
adjacent to a vertex of $L^1$ is mapped by $\Phi^1$ into
$\partial B_\eta(\ell)$.
    \end{itemize}
\end{lemma}

\begin{proof} 
We use the classical fact that simplicial collapses are preserved
under stellar subdivision: if $K\searrow K_0$ and $\sigma K$ is
obtained from $K$ by a finite sequence of stellar subdivisions, then
\[
\sigma K\searrow \sigma K_0;
\]
see \cite[\S4, Theorem~4, p.~252]{Whitehead1939}. In the present case
all the subdivided edges belong to $\mathcal{E}_L = \mathcal{E}_\eta$ and hence are not faces of
$L$. Therefore the induced subdivision of $L$ is trivial, and we
conclude that
\[
\Sigma^1\searrow L^1=L.
\]
The remaining claims are obvious given the above discussion. Indeed, since each new vertex is placed on the image of the edge being
subdivided and $\Phi^1$ is affine on the resulting simplices, we have
\[
\Phi^1(\Sigma^1)=\Phi(\Sigma).
\]
The induced subdivision of $\Gamma$ preserves the embedding property,
and the natural identification of $|\Sigma^1|$ with $|\Sigma|$
intertwines $\Phi^1$ and $\Phi$. Hence the characterization of $L^1$
as the inverse image of $\ell$ is also preserved.
\end{proof}

\medskip

{\bf Second perturbation: $\Sigma^2$ and $\Phi^2$.}
$\Sigma^2$ is constructed as follows. First of all, it contains all the vertices, segments, and triangles of $\Sigma^1$.

\begin{figure}
\begin{tikzpicture}
\draw (0, 0) -- (-6, 0) -- (0, 2);
\draw (0, 0) -- (6, 0) -- (0, 2) -- (0, 0);
\node[below left] at (-6, 0) {$Q$};
\node[below] at (0,0) {$V$};
\node[below right] at (6,0) {$V'$};
\node[above] at (0,2) {$A(V)$};
\end{tikzpicture}
\caption{The vertices, edges, and triangles added to $\Sigma^1$ to form the simplicial complex $\Sigma^2$. Here $Q$ represents either $Q_+$ or $Q_-$ and belongs to both $L^2$ and $\Gamma^2_{bd}$. The edge $[Q, V]$ belongs to $\Gamma^1_{bd}$, but $V$ does not belong to $\Gamma^2_{bd}$ as it is replaced by $A(V)$. The edge $[V, V']$ belongs to $\Gamma^1_{bd}$ and $V'$ belongs to both $\Gamma^1_{bd}$ and $\Gamma^2_{bd}$.}\label{f:addition}
\end{figure}
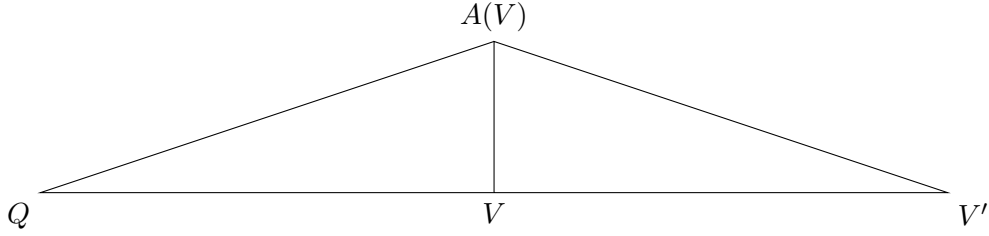

Then, for every $V\in \Gamma^1_{bd}$ adjacent (in $\Gamma^1_{bd}$) to $Q_+$ or $Q_-$, we add a vertex $A(V)$ and the following faces (these additions are depicted in \cref{f:addition}):
\begin{itemize}
    \item The segment $[V, A(V)]$;
    \item For each $V'$ so that $[V, V']$ is an edge in $\Gamma^1_{bd}$, we add the segment $[A(V), V']$ and the triangle $[A(V), V, V']$ (in particular we add the triangle $[A(V), V, Q]$ and the segment $[A(V), Q]$).
\end{itemize}
We let:
\begin{itemize}
    \item $L^2=L^1$;
    \item $\Gamma^2_{bd}$ be the $1$-dimensional simplicial complex of $\Sigma^2$ which coincides with $\Gamma^1_{bd}$ apart from the fact that we replace $V$ with $A(V)$ for all $V$ adjacent (in $\Gamma^1_{bd}$) to $Q_+$ or $Q_-$;
\item $\Gamma^2 = \Gamma^2_{bd} \cup L^2$.
\end{itemize}
For $Q\in\{Q_+, Q_-\}$, let $N(Q)$ be the set of vertices $U\in \Gamma^2_{bd}$ adjacent (in $\Gamma^2_{bd}$) to $Q$. Observe that $N(Q)$ coincides with the set of vertices $A(V)$ for all $V$ adjacent (in $\Gamma^1_{bd}$) to $Q$.

The construction of $\Sigma^2$ guarantees the validity of the following collapsibility result, which is the whole purpose of introducing $\Sigma^2$.

\begin{lemma}\label{lem:properties_simplicial_two_collapsible}
    $\Sigma^2$ can be collapsed onto $L^2$ with a sequence of collapses with the following property.
    For any $U\in N(Q)$ with $Q\in\{Q_+,Q_-\}$, the edge $[Q,U]$ is collapsed onto $Q$.
\end{lemma}
\begin{proof}
    We describe the sequence of collapses explicitly.

    For all $A(V)\in N(Q)$, with $Q\in\{Q_+, Q_-\}$, we perform the following collapses (in this order):
    \begin{enumerate}
        \item For all edges $[V, V']$ in $\Gamma^1_{bd}$ with $V'\not\in\{Q_+, Q_-\}$, collapse the triangle $[A(V), V, V']$ onto $[A(V), V]$ and $[V, V']$.
        \item Collapse the triangle $[A(V), Q, V]$ on $[A(V), Q]$ and $[Q, V]$;
        \item Collapse the segment $[A(V), Q]$ onto $Q$ (as required in the statement).
    \end{enumerate}
    Using the definition of collapse pairs, we thus perform the following collapses:
    \[
\bigl([A(V),V,V'],[A(V),V']\bigr),
\qquad
\bigl([A(V),Q,V],[A(V),V]\bigr),
\qquad
\bigl([A(V),Q],A(V)\bigr).
\]
    At this point, we have collapsed $\Sigma^2$ on $\Sigma^1$, which we already know by \cref{lem:properties_simplicial_one} is collapsible on $L^1=L^2$.
\end{proof}

The map $\Phi^2:\Sigma^2\to\R^3$ is defined as:
\begin{itemize}
    \item $\Phi^2(V) = 0$ for any $V\in L^2$,
    \item $\Phi^2(V) = \Phi^1(V)$ for any $V\in \Sigma^1\setminus L^2$,
    \item $\Phi^2 (A (V)) = \Phi^1 (V)$ for all $A(V)\in N(Q)$ with $Q\in \{Q_+, Q_-\}$,
\end{itemize}

Let us record in the following statement some properties of the image of $\Sigma^2$ via $\Phi^2$.
\begin{lemma}\label{lem:properties_simplicial_two_boundary}
    We have:
    \begin{enumerate}
        \item For all $U\in N(Q)$, where $Q\in\{Q_+, Q_-\}$, the point $\Phi^2(U)$ is at distance exactly $\eta$ from $\ell$.
        \item Let $\tilde \Gamma^2_{bd}$ be the subcomplex obtained by removing $Q_+$ and $Q_-$ (and all the edges adjacent to them) from $\Gamma^2_{bd}$. $\Phi^2$ is an embedding on $\tilde\Gamma^2_{bd}$ with image contained in $E\cap \partial K$.
        \item  We have the decomposition
        \begin{equation*}
            E\cap \partial K
            =
            \Phi^2(\tilde\Gamma^2_{bd})\cup \bigcup_{Q\in\{Q_-, Q_+\},\, V\in N(Q)} [\Phi^1(Q), \Phi^2(V)] \,.
        \end{equation*}
        \item $\Phi^2(L^2) = \{0\}$.
        \item $(\Phi^2)^{-1}(\ell) = |L^2|$; equivalently $L^2$ consists of all faces $F$ of $\Sigma^2$ such that $(\Phi^2)^\circ(F)\cap\ell\not=\emptyset$.
        \item There is a constant $C>0$ (independent from $\eta$) such that
        \begin{equation*}
            \mathcal H^2(\Phi^2(\Sigma^2)) \le \mathcal H^2(\Phi(\Sigma)) + C\eta.
        \end{equation*}
    \end{enumerate}
\end{lemma}
\begin{proof} Parts (1), (2), (3), and (4) are immediate consequences of the construction of $\Phi^2$.

Now, let us show (5).
Consider a face $F$ such that $(\Phi^2)^\circ(F)\cap\ell\not=\emptyset$ and assume, for the sake of contradiction, that $F$ does not belong to $L^2$.

 If $F\subseteq \mathcal V(L^2)$, we know by assumption (ii) of \cref{t:comp-to-map} that $F \in L^2$. Hence, we can assume that at least one vertex of $F$ does not belong to $\mathcal V(L^2)$.

 Let $\tilde F := F\setminus \mathcal V(L^2)$. Since every vertex of $L^2$ is mapped to the origin and $\ell$ is a
linear subspace, we have $(\Phi^2)^\circ(\tilde F)\cap \ell\not=\emptyset$.
    If $\tilde F$ does not contain any vertex in $N(Q)$ (with $Q\in\{Q_+, Q_-\}$), then $\tilde F$ is a face of $\Sigma^1\setminus L^1$, and thanks to \cref{lem:properties_simplicial_one} we have a contradiction.
    Otherwise, let $U\in N(Q)$ be a vertex of $\tilde F$, and consider $V\in\Sigma^1$ such that $U=A(V)$. Let $\hat F:=\tilde F\cup \{V\}\setminus \{A(V)\}$.  
    Since $\Phi^2(U)=\Phi^2(V)$, $(\Phi^2)^\circ(\tilde F) = (\Phi^2)^\circ(\hat F)$. Moreover, by construction of $\Sigma^2$, $\hat F$ is a face of $\Sigma^1$ and thus  we have already proven that $(\Phi^2)^\circ(\hat F)$ is disjoint from $\ell$.

    We prove (6). The faces of $\Sigma^2$ that were not faces of $\Sigma^1$ are mapped into sets with zero area by $\Phi^2$ (because $\Phi^2(A(V))=\Phi^2(V)$). Moreover, for any face $F$ with no vertex from $\mathcal V(L^2)$, we have $\Phi^2(F)=\Phi^1(F)$. Therefore, we have
    \begin{equation*}
        \mathcal H^2(\Phi^2(\Sigma^2)) \le
        \mathcal H^2(\Phi^1(\Sigma^1)) +
        \sum_{\{F\in\Sigma^1:\, F\cap\mathcal V(L^2)\not=\emptyset\}} \mathcal H^2(\Phi^2(F))
    \end{equation*}
    Thanks to \cref{lem:properties_simplicial_one}, $\mathcal H^2(\Phi^1(\Sigma^1))=\mathcal H^2(\Phi(\Sigma))$. Moreover, again thanks to \cref{lem:properties_simplicial_one}, the summation is over a number of faces which is independent of $\eta$ and, for each face $F$ considered in the summation, $\Phi^2(F)\subseteq \overline{B}_\eta (\ell)\cap K$; in particular $\mathcal H^2(\Phi^2(F)) \le C_0\eta$, where $C_0$ depends only on $K$. Thus, the right-hand side can be bounded by $\mathcal H^2(\Phi(\Sigma)) + C\eta$, for some constant $C$ independent of $\eta$.
\end{proof}

\subsection{Controlled homotopies}\label{s:controlled-flows} In this section we wish to produce a controlled homotopy which encodes the collapse of a complex onto a subcomplex. We will then apply this trick on the complex $\Sigma^2$ (which collapses on $L^2$) constructed in the previous subsection.

Let us remark that property (iii) of the following statement will never be fully used (it would have sufficed to know the weaker ``$H(t, x)=x$ for $t$ sufficiently large'') but it is required in a technical step of the inductive proof of \cref{l:controlled-flows}.

\begin{lemma}\label{l:controlled-flows}
Assume $\Sigma$ is a finite connected simplicial complex which can be collapsed onto a subcomplex $\Delta$. Then there is a Lipschitz map $H: [0, \infty) \times |\Sigma| \to |\Sigma|$ such that
\begin{itemize}
\item[(i)] $H (0, |\Sigma|) \subset |\Delta|$;
\item[(ii)] $H (t, \cdot)^{-1} (|\Delta|) = |\Delta|$ for all $t>0$;
\item[(iii)] $H (t,x)=x$ for every $t\geq \dist(x, |\Delta|)$.
\end{itemize}
Furthermore, if $[V, V']\in\Sigma$ is an edge with $V\in \Delta$ that, during the collapsing sequence, is collapsed onto $V$, then we may require that $H(t, |[V, V']|)\subseteq |[V, V']|$ for all $t\ge 0$.
\end{lemma}
\begin{proof}
First, we construct the homotopy $H$ explicitly when $\Delta$ is obtained from the elementary collapse of a single cell, then we show how to concatenate the construction to prove the full statement.

\medskip

{\bf Controlled homotopy for the elementary collapse.} Assume that $\Delta$ is obtained from $\Sigma$ by a single
elementary collapse of a maximal simplex $S=[v_0,\ldots,v_k]$ with free face $F=[v_1,\ldots,v_k]$. So, $|\Sigma|=|\Delta|\cup |S|$.

Let $\lambda_0,\ldots,\lambda_k$ denote the barycentric coordinates on
$S$, i.e., for any point $x\in |S|$ we let $\lambda_i (x)$ be the only real numbers such that $\sum_i \lambda_i (x) =1$ and 
\begin{equation*}
    x = \lambda_0(x)e_{v_0} + \cdots + \lambda_k(x) e_{v_k}\,.
\end{equation*}
Define $\rho:|S|\to [0, 1]$ as
\[
\rho(x):=\min_{1\leq i\leq k}\lambda_i(x) \,
\]
and extend it outside to be identically $0$.
Observe that $|S|\cap|\Delta|=\{x\in |S|:\rho(x)=0\}$.
Set $c:=\max\{1, \Lip\, (\rho)\}$.

We define $H(t,x)=x$ for $x\in|\Delta|$, while for
$x\in |S|$ we define $H(t,x)$ as
\[
H(t,x)
 := x + \max\{0, \rho(x)-ct\}\Big[ke_{v_0} - (e_{v_1} + \cdots + e_{v_k})\Big].
\]
Since $\rho(x)\le \lambda_i(x)$ for all $1\le i\le k$, for $t\ge 0$ the expression for $H(t, x)$ yields an element of $|S|$ (i.e., a convex combination of $e_{v_0}, \dots, e_{v_k}$).
Moreover, since $\rho(x)=0$ on $|S|\cap|\Delta|$, the two definitions agree on $|S|\cap|\Delta|$.

For any $x\in|S|$, we have $\rho(x)=\lambda_i(x)$ for some $1\le i\le k$ and thus $\lambda_i(H(0, x))=0$ which implies $\rho(H(0, x))=0$. It follows that
\[
H(0,|\Sigma|)\subset|\Delta|.
\]
On the other hand, if $x\notin|\Delta|$ and $t>0$, then $\rho(x)>0$ and so
\[
\lambda_i\bigl(H(t,x)\bigr)
=
\lambda_i(x)-\max\{0, \rho(x)-ct\}> \lambda_i(x) - \rho(x)\ge 0
\qquad\text{for every }i=1,\ldots,k.
\]
Consequently $H(t,x)\notin|\Delta|$, and hence
\[
H(t,\cdot)^{-1}(|\Delta|)=|\Delta|
\qquad\text{for every }t>0.
\]

Since $\rho$ vanishes on $|S|\cap|\Delta|$, for any $x\in |S|$, we have
\[
\rho(x)\leq
\Lip\, (\rho)\dist (x,|\Delta|)
\leq
c\,\dist(x, |\Delta|).
\]
Thus, whenever $t\geq\dist(x,\Delta)$, the defining formula gives $H(t,x)=x$.

The Lipschitz continuity of $H$ follows immediately from its definition.
This proves the elementary-collapse case.

\medskip

{\bf General case by induction.} The general case follows then by induction once we can show that, if there is such a map $H_1$ for two complexes $\Delta'\subset \Sigma$ and $H_2$ for $\Delta \subset \Delta'$, then there is a map $H$ for $\Delta\subset \Sigma$. In fact we define $t_0 (x) := \dist\, (H_1 (0,x), \Delta)$ and define first
\[
H' (t,x) :=
\left\{
\begin{array}{ll}
H_2 (t, H_1 (0,x)) \quad & \mbox{if $t< t_0 (x)$}\\
H_1 (t- t_0 (x), x) & \mbox{otherwise}
\end{array}
\right.
\]
and then $H (t,x) = H' (Ct, x)$.
We now check that $H$ satisfies all the conditions up to a choice of a suitably large scaling constant $C$.

\medskip

{\bf Lipschitz continuity in $t$.} We need to estimate $|H' (t,x) - H' (s,x)|$ for $t>s$. Clearly we have
\[
|H' (t,x)-H' (s,x)|\leq \Lip\, (H_1) |t-s| \quad \mbox{if $t > s\geq t_0 (x)$.}
\]
On the other hand, since
\begin{equation}\label{e:identita-banale}
H_2 (t_0 (x), H_1 (0,x)) = H_1 (0,x) = H' (t_0 (x), x)\, ,
\end{equation}
we also have
\[
|H' (t,x) - H'(s,x)|\leq \Lip\, (H_2) |t-s| \quad \mbox{for $s < t\leq t_0 (x)$}\, .
\]
In turn, the two inequalities above give immediately
\[
|H' (t,x)- H'(s,x)|\leq \max \{\Lip\, (H_1), \Lip\, (H_2)\} |t-s| \quad \mbox{when $s< t_0 (x) < t$.}
\]

\medskip

{\bf Lipschitz continuity in $x$.} We need to estimate $|H'(t,x)-H'(t,y)|$ and we can assume, without loss of generality, $t_0 (x) \geq t_0 (y)$. First of all
\begin{equation}\label{e:tutti-sotto}
|H' (t,x) - H' (t,y)|\leq \Lip\, (H_2) {\rm Lip\,} (H_1 (0, \cdot)) |x-y| \qquad
\mbox{if $t \leq t_0 (y) \leq t_0 (x)$.}
\end{equation}
Since $t_0$ is a Lipschitz function itself, we also get
\begin{equation}\label{e:tutti-sopra}
|H' (t,x) - H' (t,y)|\leq \Lip\, (H_1) (1+ \Lip\, (t_0)) |x-y| \qquad
\mbox{if $t \geq t_0 (x) \geq t_0 (y)$}
\end{equation}
It remains to consider the situation $t_0 (y) < t < t_0 (x)$. In this case we estimate
\begin{align*}
 & |H' (t,x) - H' (t,y)|\\
\leq & |H'(t,x)-H' (t_0 (x),x)| + |H' (t_0 (x),x)-H' (t_0 (x), y)| + |H'(t_0 (x), y)- H' (t, y)|\, .
\end{align*}
The first and third terms are estimated with $C_0 |t-t_0 (x)|\leq C_0 |t_0 (y)-t_0 (x)| \leq C_0 \Lip\, (t_0) |x-y|$, using the Lipschitz continuity in $t$ of $H'$ from the previous step, and the Lipschitz continuity of $t_0$. To estimate the second summand we can use \eqref{e:tutti-sopra}.

\medskip

We now check the claims (i), (ii), and (iii) in the Lemma.

\medskip

{\bf (i).} If $0 < t_0 (x)$, then clearly $H (0,x) = H' (0,x) = H_2 (0, H_1 (0,x))\in |\Delta|$. If $t_0 (x) = 0$, then by definition of $t_0$ we have $H_1 (0,x)\in |\Delta|$ and in particular \eqref{e:identita-banale} implies $H (0,x) = H' (0,x) = H_1 (0,x) \in |\Delta|$.

\medskip

{\bf (ii).} We first show $H' (t, \cdot)^{-1} (|\Delta|) \subset |\Delta|$ for positive $t$. Fix $t>0$ and $x\in |\Sigma|\setminus|\Delta|$. We need to show that $H' (t,x)\notin |\Delta|$. If $0<t<t_0(x)$, then $H_1(0,x) \notin |\Delta|$, because
$t_0(x)>0$. Property (ii) for $H_2$ gives
\[
H'(t,x)=H_2(t,H_1(0,x))\notin |\Delta|.
\]
If $t>t_0(x)$, set $s:=t-t_0(x)>0$.
Then
\[
H'(t,x)=H_1(s,x).
\]
If $x\notin |\Delta'|$, property~(ii) for $H_1$ gives
\[
H_1(s,x)\notin |\Delta'|,
\]
and hence $H'(t,x)\notin|\Delta|$. If instead
$x\in|\Delta'|\setminus|\Delta|$, then $\dist(x,|\Delta'|)=0$, so property~(iii) for $H_1$ implies
\[
H_1(s,x)=x\notin|\Delta|.
\]
Thus $H'(t,x)\notin|\Delta|$ in either case.

It remains to analyze the case $t=t_0(x)>0$. Under the latter assumption $H'(t,x)=H_1(0,x)$ and, since $t_0(x)
=
\dist(H_1(0,x),|\Delta|)>0$,
we have $H_1(0,x)\notin|\Delta|$.

Finally, to show $|\Delta|\subset H'(t,\cdot)^{-1}(|\Delta|)$, note that for $x\in|\Delta|$ and any $t\geq 0$ we have $H_1(t,x)=x$, and hence
$H'(t,x)=x$. 

\medskip

{\bf (iii).} Note that, since $H_1 (0,x) =x$ for every $x\in |\Delta'|$, we have
\begin{equation}\label{e:t0-controlled}
t_0(x)
\leq
\Lip\, (H_1)\, \dist(x, |\Delta|).
\end{equation}
We now choose $C \geq \Lip\, (H_1)+1$ and show that the claim holds with the latter choice.

Indeed, let $t\geq\dist(x,\Delta)$.
In particular \eqref{e:t0-controlled} implies
\begin{align*}
Ct-t_0(x)
&\geq
C\dist(x, |\Delta|)
-
\Lip\, (H_1) \dist(x, |\Delta|)\\
&\geq
\dist(x, |\Delta|) \geq
\dist(x, |\Delta'|),
\end{align*}
where the last inequality follows from $\Delta\subset\Delta'$. Note that the last chain of inequalities also implies $Ct \geq t_0 (x)$, so that the second branch in the definition of $H'$ applies, namely
$H(t,x) = H'(Ct,x) =
H_1(Ct-t_0(x),x)$.
Property~(iii) for $H_1$ therefore gives
\[
H(t,x)
=
H_1(Ct-t_0(x),x)
=
x\, .
\]
This proves property~(iii).

The additional property on the image of edges contained in the statement is a direct consequence of the construction of the homotopy. Indeed 
\begin{itemize}
\item before the edge is collapsed, it lies in every remaining subcomplex,
so all preceding elementary homotopies fix it;
\item at the edge collapse, the elementary formula moves it inside itself;
\item every subsequent homotopy fixes its surviving endpoint $V$.
\end{itemize}
This obviously yields $H (t, |[V,V']|) \subset |[V,V']|$. 
.
\end{proof}

\subsection{Construction of the map and proof of \cref{t:comp-to-map}}\label{s:final-map}
We can now construct a map $F$ satisfying the requirements of \cref{t:comp-to-map}.
We will do this in two steps:
\begin{itemize}
    \item Using \cref{l:controlled-flows} on $\Sigma^2$, we will build a function $G:E\cap K\to \R^3$ whose image is \emph{almost} $\Phi^2(\Sigma^2)$ and which satisfies all the assumptions necessary to apply \cref{p:Almgren-relaxed} apart from the Lipschitz continuity in the origin.
    \item We will then perturb $G$ into $\tilde G$ to obtain the Lipschitz continuity in the origin (while keeping the image $\tilde G(E\cap K)$ under control). This allows us to apply \cref{p:Almgren-relaxed} and conclude.
\end{itemize}

Let us recall that $\tilde\Gamma^2_{bd}$ is the $1$-dimensional simplicial complex obtained by removing from $\Gamma^2_{bd}$ the two vertices $Q_+$ and $Q_-$ as well as all the edges adjacent to them.
By \cref{lem:properties_simplicial_two_boundary}, we know that $\Phi^2$ is an embedding from $\tilde \Gamma^2_{bd}$ into a subset $Z$ of $E\cap\partial K$. Let $Z'$ be the union of segments that appear in \cref{lem:properties_simplicial_two_boundary}. Let us recall the decomposition
\begin{equation*}
    E\cap\partial K = Z \cup Z' \,.
\end{equation*}
For notational simplicity, let $\{P_+, P_-\} := \ell\cap \partial K$ (where $P_+$ lies on the same ray of $\ell$ as the pin $P$) and let $M_{\pm}:=\{\Phi^2(V):\, V\in N(Q_{\pm})\}$. Observe that
\begin{equation*}
    Z' = \Bigg(\bigcup_{U\in M_+} [P_+, U]\Bigg) \sqcup \Bigg(\bigcup_{U\in M_-} [P_-, U]\Bigg) \,.
\end{equation*}
We will define $G:E\cap K\to\R^3$ differently in the cone generated by $Z$ and in the cone generated by $Z'$.

\medskip\noindent\textbf{Construction of $G$ on the cone over $Z$.}
Thanks to \cref{lem:properties_simplicial_two_collapsible}, we can apply \cref{l:controlled-flows} on $\Sigma^2$ and obtain a Lipschitz function $H:[0,1]\times |\Sigma^2|\to|\Sigma^2|$ such that:
\begin{itemize}
    \item $H(0, |\Sigma^2|)\subseteq |L^2|$.
    \item $H(t, \emptyparam)^{-1}(|L^2|) = |L^2|$ for all $t>0$.
    \item $H(1, \emptyparam)$ is the identity.
    \item For every $V\in N(Q)$, where $Q\in\{Q_+, Q_-\}$, $H(t, V)\in |[V, Q]|$ for all $t\in [0, 1]$.
\end{itemize}
In fact the \cref{l:controlled-flows} yields an homotopy $\tilde{H}$ which satisfies the first, second, and fourth requirement, while in place of the third we know that there is a $T>0$ such that $\tilde{H} (t,x)=x$ for all $x\in |\Sigma^2|$ and all $t\geq T$. Thus, in order to achieve the third we simply set $H(t,x):= \tilde{H} (T t, x)$. 

Observe that an arbitrary nonzero point in the cone over $Z$ can be expressed (uniquely) as $t\Phi^2(p)$ with $0< t\le 1$ and $p\in |\tilde\Gamma^2_{bd}|$.
We define
\begin{equation*}
    G(t \Phi^2(p)) := \Phi^2(H(t, p)) \,.
\end{equation*}
Since $\Phi^2(H(0,p))=0$, extending it to $G(0) =0$ preserves the above identity no matter what $p$ in $|\tilde \Gamma^2_{bd}|$ we choose.
Since $H(0, |\Sigma^2|)\subseteq |L^2|$ and $\Phi^2(|L^2|)=\{0\}$, the map $G$ is continuous in the origin and maps the origin into the origin. Furthermore, $G$ is locally Lipschitz away from the origin.
By construction, we have that the image of $G$ is contained in $\Phi^2(\Sigma^2)$ and thus, thanks to \cref{lem:properties_simplicial_two_boundary}-(6), 
\[
\mathcal{H}^2 (G(\{tz:0\le t\le 1,\, z\in Z\})) \le \mathcal H^2(\Phi(\Sigma)) + C\eta
\]
for some constant $C>0$ independent from $\eta$.
Moreover, because of \cref{lem:properties_simplicial_two_boundary}-(5),
\begin{equation*}
    G^{-1}(\ell) = G^{-1}(\{0\}) = \{t\Phi^2(p): H(t, p)\in |L^2|\}
    = \{t\Phi^2(p):\, t=0 \text{ or } p\in |L^2|\}
    = \{0\}\,.
\end{equation*}
Finally, for any $U\in M_\pm$ (with $U=\Phi^2(V)$ and $V\in N(Q)$) and for all $0\le t\le 1$, we have
\begin{equation*}
    G(t U) = G(t\Phi^2(V)) = \Phi^2(H(t, |V|)) \in \Phi^2([V, Q]) =[\Phi^2(V), \Phi^2(Q)] = [U, 0],
\end{equation*}
which implies $G([0, U])\subseteq [0, U]$.

\medskip\noindent\textbf{Construction of $G$ on the cone over $Z'$.}
Now, we extend $G$ on the cone over $Z'$ as follows. We define it independently on the cone over $[P_+, U]$ for each $U\in M_+$ (and analogously for $P_-$).
Let $W\in\R^3$ be the midpoint of the segment $[P_+, U]$.
\begin{itemize}
    \item On the cone over $[P_+, W]$, $G$ is the identity.
    \item An arbitrary point in the cone over $[W, U]$ can be expressed as $t(\lambda W+\mu U)$ with $0\le t\le 1$ and $\lambda,\mu\ge 0$ so that $\lambda+\mu=1$. We define
    \begin{equation*}
        G\big(t(\lambda W + \mu U)\big) := \lambda t W + \mu G(tU) \,.
    \end{equation*}
\end{itemize}
The finitely many definitions agree on their
common edges and hence give a locally Lipschitz map away from 0.

The fact that $G([0, U])\subseteq [0, U]$ guarantees that this extension has two crucial properties:
\begin{itemize}
    \item $G(\{tz': 0\le t\le 1,\, z'\in [P_+, U]\}) \subseteq \{tz': 0\le t\le 1,\, z'\in [P_+, U]\}$;
    \item $G^{-1}(\ell)\cap \{tz': 0\le t\le 1,\, z'\in [P_+, U]\} \subseteq\ell$.
\end{itemize}
Thanks to \cref{lem:properties_simplicial_two_boundary}-(1), we deduce that
\begin{equation*}
    G(\{tz':\, 0\le t\le 1,\, z'\in Z'\}) \subseteq E\cap B_\eta (\ell),
\end{equation*}
The latter set consists of the intersection with the tube $B_\eta (\ell)$ of finitely many
planar conical sectors adjacent to the two pin rays and a small region around the origin of radius $\eta$. The intersection of the sectors with the tube contribute $O (\eta)$ to the area, while the small region around the origin contributes only $O (\eta^2)$. So the total area is bounded
bounded by $C\eta$ for some constant $C>0$ (independent from $\eta$).

\medskip\noindent\textbf{Properties of $G$ and construction of $\tilde G$.}
By construction, and thanks to the observations above, we know that the map $G:E\cap K\to \R^3$ is continuous and it is locally Lipschitz away from the origin. Moreover, it satisfies
\begin{enumerate}[label=(\roman*)]
    \item $G(x)=x$ for $x\in E\cap\partial K$.
    \item $G(x)=x$ for $x \in C_\theta(\ell)\cap K$ (with $\theta$ depending on $\eta>0$) and $G^{-1}(\ell) \subset \ell$.
    \item $\mathcal H^2(G(E\cap K)) \le \mathcal H^2(\Phi(\Sigma)) + C\eta$ for some constant $C>0$ independent from $\eta$.
    \item The image of $G$ is a subset of $E\cup \Phi^2(\Sigma^2)$. 
\end{enumerate}
These assumptions are \emph{almost} sufficient to directly apply \cref{p:Almgren-relaxed}. The missing bit is the Lipschitz continuity of $G$ at the origin. We will now use (iv) to perturb $G$ into $\tilde G$ in a way that guarantees the desired Lipschitz continuity.

\medskip

First of all, we note that there is a cone $\mathcal T\subseteq\R^3$ with locally finite area such that, for any small enough $r>0$, $G(E\cap B_r (0))\subset \mathcal T$. Indeed 
let $\mathscr S_0 :=
\{S\in\Sigma^2:0\in\Phi^2(S)\}$.
By \cref{lem:properties_simplicial_two_boundary}-(5), every simplex
in $\mathscr S_0$ has a face contained in $L^2$, and this face is
mapped to $0$. Hence, if we set, 
\[
\mathcal T
:=
E\cup
\bigcup_{S\in\mathscr S_0}
\{t y:t\geq0,\ y\in\Phi^2(S)\}
\]
then $\mathcal{T}$ is a finite union of cones of dimension at most two and therefore is itself a cone with
locally finite $\mathcal H^2$-measure.
Since there are only finitely many simplices not belonging to
$\mathscr S_0$, their images have positive distance from the origin.
Hence $\Phi^2(\Sigma^2)\cap B_r \subset\mathcal T$ for every sufficiently small $r>0$.

\medskip

Choose a small radius $r>0$ such that $G(E\cap B_r (0))\subset \mathcal T\cap B_\eta (0)$ and so that $r < \dist (0, \partial K)$. 
We define $\tilde G:E\cap K\to\R^3$ as
\begin{equation*}
    \tilde G(x) =
    \begin{cases}
        G(x) \quad &\text{if $|x|\ge r$,}\\
        \frac{|x|}r G(r\tfrac x{|x|}) \quad &\text{if $0 < |x|<r$}\\
0 \quad &\text{if $x=0$.}
    \end{cases}
\end{equation*}
The map $\tilde G$ is Lipschitz. Indeed the restriction of $G$ to $E\cap\partial B_r$ is Lipschitz, since
$G$ is locally Lipschitz away from the origin, and since $E$ is a finite
union of planar cones, the positively homogeneous extension of this
boundary map is Lipschitz on $E\cap B_r$. Clearly $\tilde{G}$ retains property (i). To see that property (ii) is retained note that:
\begin{itemize}
\item in the directions where $G$ is the identity, $\tilde{G}$ is also the identity;
\item if $\tilde{G} (x)\in \ell$ and $x\neq 0$, then $G (rx/|x|)\in \ell$ and therefore $x\in \ell$. 
\end{itemize}
Thanks to property (iii) and (iv), we have
\begin{equation*}
    \mathcal H^2(\tilde G(E\cap K)) \le
    \mathcal H^2(G(E\cap K))
    + \mathcal H^2(\mathcal T\cap B_\eta (0))
    \le \mathcal H^2(\Phi(\Sigma)) + C\eta + C'\eta^2,
\end{equation*}
for some constant $C'>0$ independent of $\eta$.

We are now almost ready to apply \cref{p:Almgren-relaxed}. First of all we extend it to a map $\widehat G$ on $E$ by setting it equal to the identity outside of $K$. 
Since $\tilde G$ is the identity on $E\cap\partial K$, such extension is Lipschitz. Moreover it is the identity on $C_\theta (\ell)$ and $\widehat{G}^{-1} (\ell)\subset \ell$. 
Furthermore,
\begin{align*}
\mathcal H^2\bigl(\widehat G(E\cap B_1)\bigr)
&\leq
\mathcal H^2\bigl((B_1\setminus K)\cap E\bigr)
+
\mathcal H^2\bigl(\widetilde G(E\cap K)\bigr)\\
&\leq
\mathcal H^2\bigl((B_1\setminus K)\cap E\bigr)
+
\mathcal H^2\bigl(\Phi(\Sigma)\bigr)
+
C\eta+C'\eta^2.
\end{align*}
Choose $\eta>0$ sufficiently small that $C\eta+C'\eta^2<\frac{\varepsilon}{2}$.
Applying \cref{p:Almgren-relaxed} to $\tilde{G}$ with error
$\varepsilon/2$ gives a Lipschitz map
$F:\mathbb R^3\to\mathbb R^3$ satisfying (a) and (b) in the statement of the theorem. 

\section{Quantitative stability of competitors}\label{s:competitore}

In this section we prove two propositions which, in combination with the theory developed in \cref{s:triangulations-and-maps}, will allow us to use the code described in \cref{s:terzo-programma} to produce effective competitors for the remaining topologies. In \cref{s:triangulations-and-maps} we addressed the fact that the code produces competitors which are triangulated surfaces, while \cref{d:Almgren-minimizing} requires additionally the existence of suitable Lipschitz maps. In this section we address a different issue, namely, that the competitors produced by the code do not exactly match the desired boundary conditions, in other words they ``span'' a net which is only an approximation of the one which we want to exclude. The first proposition will deal with the fat fish family, while the second one will be used for the other nets. 

\subsection{Fat fish} First, recall \cref{d:pesce-grasso} and \cref{f:fat-fish}. In order to parametrize the fat fish family, we denote by $\xi$ the length of the edge $P_1V_2$ and by $\zeta$ the length of the edge $V_5V_1$ (cf. \cref{d:pesce-grasso}). Moreover, without loss of generality, we will assume that $P_1V_2$ is shorter than (or congruent to) $P_2 V_4$ and that $V_5V_1$ is shorter than (or congruent to) $V_3V_6$.
Recall (cf. \cref{d:pesce-grasso}) that an element of the fat fish family is required to be an embedding of the corresponding spherical graph. In particular $\zeta>0$ and the edges $P_1V_2$ and $P_2V_4$ cannot intersect the two half-equators connecting $V_5$ and $V_6$. A direct computation shows that this implies $\xi>0$. These and other useful pieces of information are summarized in the following lemma.

\begin{lemma}\label{l:fat-fish-classification}
The pins of any element $\mathcal{N}$ of the fat fish family have degree one. Moreover there is a rigid motion $R$ of the sphere such that:
\begin{itemize}
    \item the pins of $R(\mathcal{N})$ are the north and south poles $(0,0,1)$ and $(0,0,-1)$;
    \item the edges adjacent to the pins lie in the halfplane $\{x_2=0, x_1\geq 0\}$;
    \item the edge adjacent to $(0,0,1)$ is shorter than (or congruent to) the one adjacent to $(0,0,-1)$ (and its length is $\xi$);
    \item between the two edges adjacent to the two edges of length $\pi$, the one in the half-space $\{x_2<0\}$ is no longer than the one lying in $\{x_2>0\}$ (and its length is $\zeta$).
\end{itemize}
We call this the {\em canonical position} of $\mathcal{N}$. Furthermore:
\begin{itemize}
\item[(a)] Two elements in canonical position always coincide in a neighborhood of the pins;
\item[(b)] Given two elements in canonical position, the four vertices forming the convex quadrilateral face coincide if and only if their $\xi$-coordinates coincide;
\item[(c)] Two elements in canonical position are the same net if and only if they have the same coordinates $(\xi, \zeta)$.
\end{itemize}
\end{lemma}

We are now ready to state the stability proposition.

\begin{proposition}\label{p:stability-fat}
Consider two elements $\mathcal{N}$ and $\mathcal{N}'$ of the fat fish family with corresponding cones $E$ and $E'$, and coordinates $(\xi,\zeta)$ and $(\xi',\zeta')$. Assume that $|\xi-\xi'|+|\zeta-\zeta'|$ is strictly smaller than the distance of any pin of $\mathcal{N}$ to any point of $\mathcal{N}$ not belonging to an edge adjacent to the pin in question. Let $F: \mathbb R^3\to \mathbb R^3$ be a Lipschitz map such that
\begin{itemize}
\item[(i)] $F$ is the identity outside $B_1$ and on a neighborhood of $\ell$;
\item[(ii)] $\overline{\{F(x):x\neq F(x)\}}\cap \ell = \emptyset$.
\end{itemize}
Then, for every positive $\varepsilon$, there is a Lipschitz map $G$ which satisfies (i) and (ii) and such that
\begin{align}
\mathcal{H}^2 (G (E\cap \overline{B}_1))
&\leq  \mathcal{H}^2 (F(E'\cap\overline{B}_1)) + 4 (|\xi-\xi'|+|\zeta-\zeta'|) + 4|\xi-\xi'| \Big( \tfrac{2\sqrt{3}}{3} +
\sqrt{\tfrac{3-\sqrt{6}}{6}}\Big) + \varepsilon \nonumber\\
&\leq  \mathcal{H}^2 (F(E'\cap\overline{B}_1)) + 10 |\xi-\xi'| + 4|\zeta-\zeta'| + \varepsilon\, .\nonumber
\end{align}
\end{proposition}

\subsection{Adjusting triangulated competitors} In this section we deal with the other cases. However, rather than assuming that we have a ``mapped competitor'' for the approximating net, we assume that we have a triangulated competitor. 

We will need a suitable ``controlled change of topology'', since the approximating net might have a slightly different topology than the real one.

\begin{definition}\label{d:topology-change}
Let $\mathcal{N}$ and $\tilde{\mathcal{N}}$ be two geodesic nets with the same antipodal pins, not necessarily stationary, and let $G$ and $\tilde{G}$ be the corresponding embedded spherical graphs. 
We say that there is a controlled change of topology from
$\tilde{\mathcal N}$ to $\mathcal N$ if the graph $G$ can be obtained from $\tilde{G}$ by contracting at most two edges, each incident to a different pin.

Note that the contraction eliminates at most one vertex adjacent to each pin of $\tilde{G}$. 
Such a vertex will be called the ``special vertex associated with the pin''.
We introduce the following maps $\tau: \mathcal{V} (\tilde{G}) \to \mathcal{V} (G)$ and $\iota: \mathcal{V} (G) \to \mathcal{V} (\tilde{G})$:
\begin{itemize}
\item $\tau$ maps each pin $\tilde{P}$ of $\tilde{G}$ into the corresponding pin $P$ of $G$ and the special vertex associated with $\tilde{P}$ (if it exists) in $P$;
\item After contracting the edges (and identifying the corresponding special vertices with the corresponding pins), $\tau$ is an isomorphism of embedded graphs.
\item $\iota$ maps any vertex $V$ of $G$ which is not a pin into its unique preimage $\tau^{-1} (V)$. For a pin $P$ of $G$, if $\tau^{-1} (P)$ is a singleton $\{\tilde{P}\}$, then $\iota (P)=\tilde{P}$; otherwise we let $\iota (P) =\tilde{S}$, the special vertex associated with $\tilde{P}$. 
\end{itemize}
Therefore $\tau$ is a surjective map, $\iota$ an injective map, and $\tau\circ \iota$ is the identity. 
\end{definition}

We are now ready to state the second proposition. 

\begin{proposition}\label{p:stability}
Let $\mathcal{N}'$ be a stationary geodesic net with two antipodal pins $\pm P$, let $\ell$ be the line through the pins, and let $E'$ be its cone. Let $\mathcal{N}$ be a geodesic net (not necessarily stationary) with the same pins such that
\begin{itemize}
\item[(a)] There is a controlled topological change from $\mathcal{N}'$ to $\mathcal{N}$ and we let $\tau$ and $\iota$ be the corresponding maps in \cref{d:topology-change}. 
\item[(b)] All the edges of $\mathcal{N}$ and all the edges of $\mathcal{N}'$ are shorter than $\pi$;
\item[(c)] For every pair $(U,V)$ of adjacent vertices of $\mathcal{N}$ which are not pins, there is an open half-space\footnote{An open half-space of $\mathbb R^3$ is any set $H_+$ of the form
$H_+ :=\{x\in \mathbb R^3: \langle x, v\rangle > 0\}$
for some nonzero vector $v\in \mathbb R^3$.} $H_+$ with $\partial H_+ \supset \ell$ and $U,V, \iota (U), \iota (V) \in H_+$;
\item[(d)] No vertex $V$ of $\mathcal{N}$ is adjacent to both pins.
\end{itemize}
Define, for every edge $\sigma = [U,V]$ of $\mathcal{N}$, the function
\begin{equation}\label{e:maximization-D}
D (\sigma) := \max \bigl\{
\mathcal{H}^2 ([U,V,\iota (V)]) +
\mathcal{H}^2 ([U, \iota (U), \iota (V)]), 
\mathcal{H}^2 ([V,U,\iota (U)]) +
\mathcal{H}^2 ([V, \iota (V), \iota (U)])
\bigr\}\, ,
\end{equation}
where $[A,B,C]$ denotes the triangle of $\mathbb R^3$ with vertices $A,B$, and $C$.
Let $K$ and $K'$ be the convex hulls of the vertices of $\mathcal N$ and $\mathcal N'$ respectively, while $K_c$ is the convex hull of the union of their vertices.

We have the following two conclusions:
\begin{itemize}
    \item Assume that $\Sigma$ is a $2$-dimensional simplicial complex and $\Phi:\Sigma\to\mathbb\R^3$ is a piecewise linear map satisfying assumptions (i), (ii), (iii) of \cref{t:comp-to-map} with respect to $(\mathcal N, E, K)$ (i.e., a valid competitor for $E$). 
    Then, there is a simplicial complex $\Sigma'$ and a piecewise linear map $\Phi':\Sigma'\to\mathbb\R^3$ satisfying the assumptions (i), (ii), (iii) of \cref{t:comp-to-map} with respect to $(\mathcal N', E', K_c)$ (i.e., a valid competitor for $E'$) such that
    \begin{equation}
        \mathcal{H}^2 (\Phi' (\Sigma')) 
        \leq \mathcal{H}^2 (\Phi (\Sigma)) + \sum_{\sigma\in\mathcal E(\mathcal N)} D(\sigma) \label{e:stab-estimate}\, .
    \end{equation}
    \item Assume that $\Sigma'$ is a $2$-dimensional simplicial complex and $\Phi':\Sigma'\to\mathbb\R^3$ is a piecewise linear map satisfying assumptions (i), (ii), (iii) of \cref{t:comp-to-map} with respect to $(\mathcal N', E', K')$ (i.e., a valid competitor for $E'$). 
    Then, there is a simplicial complex $\Sigma$ and a piecewise linear map $\Phi:\Sigma\to\mathbb\R^3$ satisfying the assumptions (i), (ii), (iii) of \cref{t:comp-to-map} with respect to $(\mathcal N, E, K_c)$ (i.e., a valid competitor for $E$) such that
    \begin{equation}
        \mathcal{H}^2 (\Phi (\Sigma)) 
        \leq \mathcal{H}^2 (\Phi' (\Sigma')) + \sum_{\sigma\in\mathcal E(\mathcal N)} D(\sigma) \, .
    \end{equation}
\end{itemize}
\end{proposition}

\subsection{Proof of \texorpdfstring{\cref{p:stability-fat}}{p:stability-fat}}
\cref{l:fat-fish-classification} has a straightforward proof, which we leave to the reader.

We next come to the proof of \cref{p:stability-fat}.
Fix two elements of the fat fish family $\mathcal{N}$ and $\mathcal{N}'$. Without loss of generality we can assume that $\mathcal{N}$ and $\mathcal{N}'$ are both in canonical position (this only requires conjugating the map $F$ with isometries, which does not change any of the assumptions on $F$). We will produce a Lipschitz homotopy $\rho: [0,1] \times \mathcal{N} \to \mathbb S^2$ with $\rho (1, \mathcal{N}) = \mathcal{N}'$ which satisfies requirement (ii) in \cref{l:homotopy-construction} and the inequality
\begin{equation}\label{e:homotopy-estimate-200}
\int_{\mathbb S^2} \# (\rho)^{-1} (\{x\})\, d\mathcal{H}^2 (x) \leq 4 (|\xi-\xi'|+|\zeta-\zeta'|) + 4|\xi-\xi'| \Big( \tfrac{2\sqrt{3}}{3} +
\sqrt{\tfrac{3-\sqrt{6}}{6}}\Big) \,.
\end{equation}
In order to define the homotopy $\rho$ it is convenient to use the labels $V_i$ for the non-pin vertices of $\mathcal{N}$ and corresponding labels $V_i'$ for the non-pin vertices of $\mathcal{N}'$, where we refer to \cref{f:fat-fish}.

The homotopy $\rho$ is obtained by concatenating two homotopies.
The first one, $\rho^1$, adjusts the coordinate $\xi$ to $\xi'$ and it is defined as follows:
\begin{itemize}
\item Denote by $R_\theta$ the rotation of angle $\theta$ around the axis $x_2$. If $x$ belongs to any edge of $\mathcal{N}$ which is not adjacent to the pins, then $\rho^1 (t,x) = R_{t(\xi'-\xi)} (x)$. Note that, by \cref{l:fat-fish-classification}, $\rho^1 (1, V_i) = V'_i$ for $i\in \{1,2,3,4\}$.
\item Let $W_2$ and $W_4$ be the two points adjacent to the pins at distance $a=\tfrac{1}{2}\min \{\xi, \xi'\} $. $\rho^1$ is then the identity on $P_1W_2$ and $P_2W_4$. On $W_2V_2$ we extend $\rho^1$ so as to take values on the equator $\mathbb S^2 \cap \{x_2=0\}$. Likewise on the other edge adjacent to the pin $P_2$. Moreover we extend it so that it is arclength-affine along the relevant equatorial arc at every $t$.
\end{itemize}
Note that $\mathcal{N}'' := \rho^1 (1, \mathcal{N})$ is also an element of the fat fish family, with coordinates $(\xi', \zeta)$. We now wish to find a second homotopy $\rho^2$ which maps $\mathcal{N}''$ onto $\mathcal{N}'$. Denote now by $V_i''$ the vertices $\rho^1 (1, V_i)$ and note that $V_i''=V_i'$ for $i\in \{1, \ldots , 4\}$. So we define $\rho^2$ to be the identity on the edges $P_1V_2''$, $P_2V_4''$, $V_1''V_2''$, $V_2''V_3''$, $V_3''V_4''$, and $V_4''V_1''$. For points in the two half-equators connecting $V_5''$ and $V_6''$ we let $S_\theta$ be the rotation of angle $\theta$ around the axis $x_3$ and define $\rho^2 (t,x) = S_{t(\zeta'-\zeta)} (x)$. On the edges $V_1''V_5''$ and $V_3''V_6''$ we then extend the map so that it takes values on the equator $\mathbb S^2 \cap \{x_3=0\}$, for instance we can take it arclength-affine along the relevant equatorial arc at every time $t$. 

Note that the concatenation $\rho$ of the homotopies $\rho^1$ and $\rho^2$ satisfies (ii) in \cref{l:homotopy-construction}:
\begin{itemize}
    \item by construction it is the identity in a neighborhood of the pins;
    \item by construction the points on the two edges adjacent to the pins which are displaced by the homotopy are anyway all mapped to a positive distance from the pins;
    \item all the other points, belonging to the other edges, are moved at most at a distance $|\xi-\xi'|+|\zeta-\zeta'|$ from themselves, and hence they all stay at a positive distance from the pins.
\end{itemize}

Therefore it remains to show \eqref{e:homotopy-estimate-200}, for which we will invoke the following lemma.

\begin{lemma}\label{l:rotation-lemma}
Denote by $O_\theta$ the rotation around an axis $r$ of $\mathbb R^3$ of angle $\theta$ and let $\sigma$ be a geodesic segment on the sphere $\mathbb S^2$ of length at most $\pi$. Consider the map $\tau: [0,1]\times \sigma \to \mathbb S^2$ given by $\tau (t,x) = O_{t\theta} (x)$. Then
\begin{equation}\label{e:omotopia-rotante}
\int_{\mathbb S^2} \# \tau^{-1} (\{x\})\, d\mathcal{H}^2 (x) \leq 2 |\theta| \sin \left(\tfrac{\text{len} (\sigma)}{2}\right)\, .
\end{equation}
\end{lemma}

We will prove \cref{l:rotation-lemma} below, while for the moment we assume its validity. Note first that the interpolated portions of the pin-adjacent edges remain in
$\mathbb S^2\cap\{x_2=0\}$, while the interpolations on
$V_1''V_5''$ and $V_3''V_6''$ remain in
$\mathbb S^2\cap\{x_3=0\}$. Their images are therefore
one-dimensional, so their two-dimensional Jacobians vanish
almost everywhere. This observation puts us in the position of estimating the left-hand side of \eqref{e:homotopy-estimate-200} for all the other edges.

Note that the two edges connecting $V_5$ and $V_6$ have length $\pi$ and thus contribute each $2 (|\xi-\xi'|+|\zeta-\zeta'|)$, while the four edges $V_1V_2$, $V_2V_3$, $V_3V_4$, and $V_4V_1$ have all length $\lambda_Q$, which is the sidelength of a spherical equilateral quadrilateral with angles $\frac{2\pi}{3}$. Since $\sin \frac{\lambda_Q}{2} = \frac{1}{\sqrt{3}}$ we can then estimate
\[
\int_{\mathbb S^2} \# \rho^{-1} (\{x\})\, d\mathcal{H}^2 (x)
\leq 4 (|\xi-\xi'|+|\zeta-\zeta'|) + \tfrac{8}{\sqrt{3}} |\xi-\xi'|
+ 2 |\xi-\xi'| \left(\sin \tfrac{\text{len} (V_1V_5)}{2}+ \sin \tfrac{\text{len} (V_3V_6)}{2}\right)\, .
\]
Now, observe that $\text{len} (V_1V_5) + \text{len} (V_3V_6) = \pi - \lambda_T$, where $\lambda_T = \arccos (-\frac{1}{3})$ is the sidelength of a spherical equilateral triangle with angles $\frac{2\pi}{3}$. We can thus estimate
\[
\sin \tfrac{\text{len} (V_1V_5)}{2}+ \sin \tfrac{\text{len} (V_3V_6)}{2}
\leq \max_{0\leq x \leq \pi-\lambda_T} \left(\sin \tfrac{x}{2} + \sin \tfrac{\pi-\lambda_T-x}{2}\right)\, .
\]
The latter maximum can be easily computed to be at $x= \frac{\pi-\lambda_T}{2}$, thus giving the bound
\[
\sin \tfrac{\text{len} (V_1V_5)}{2}+ \sin \tfrac{\text{len} (V_3V_6)}{2}
\leq 2 \sin \tfrac{\pi-\lambda_T}{4} = 2 \sqrt{\tfrac{3-\sqrt{6}}{6}}\, .
\]
Putting all these inequalities together gives us the estimate \eqref{e:homotopy-estimate-200}.

We now apply \cref{l:homotopy-construction} to $F$ and the homotopy
$\rho$, and using \eqref{e:homotopy-estimate-200}, for every
$\varepsilon>0$ we obtain a Lipschitz map $G$ satisfying conditions
(i) and (ii) of \cref{p:stability-fat} such that
\begin{align*}
\mathcal{H}^2\bigl(G(E\cap\overline{B}_1)\bigr)
&\leq
\mathcal{H}^2\bigl(F(E'\cap\overline{B}_1)\bigr)
+\int_{\mathbb S^2}
\#\rho^{-1}(\{x\})\,d\mathcal{H}^2(x)
+\varepsilon\\
&\leq
\mathcal{H}^2\bigl(F(E'\cap\overline{B}_1)\bigr)
+4\bigl(|\xi-\xi'|+|\zeta-\zeta'|\bigr)\\
&\qquad
+4|\xi-\xi'|
\left(
\tfrac{2\sqrt{3}}{3}
+\sqrt{\tfrac{3-\sqrt{6}}{6}}
\right)
+\varepsilon.
\end{align*}
This proves the first claimed inequality. Moreover,
\[
4
+4\left(
\tfrac{2\sqrt{3}}{3}
+\sqrt{\tfrac{3-\sqrt{6}}{6}}
\right)
<10.
\]
Consequently,
\[
\mathcal{H}^2\bigl(G(E\cap\overline{B}_1)\bigr)
\leq
\mathcal{H}^2\bigl(F(E'\cap\overline{B}_1)\bigr)
+10|\xi-\xi'|
+4|\zeta-\zeta'|
+\varepsilon.
\]
This completes the proof of \cref{p:stability-fat} and we are thus left with the proof of \cref{l:rotation-lemma}.

\begin{proof}[Proof of \cref{l:rotation-lemma}]
Let us parametrize $\sigma$ over $[0,\alpha]$ at unit speed, where clearly $\alpha$ is the length of $\sigma$. We can then write $\tau (t,s) =
O_{t\theta} (\sigma (s))$: if we consider it as a map from the planar domain $[0,1] \times [0,\alpha]$ into $\mathbb R^3$ we can use the area formula to compute
\[
\int_{\mathbb S^2} \# \tau^{-1} (\{x\})\, d\mathcal{H}^2 (x)
= \int_0^1 \int_0^\alpha \left|\tfrac{\partial \tau}{\partial t} \times \tfrac{\partial \tau}{\partial s}\right|\, ds\, dt\, .
\]
Consider the unit vector $v$ orienting the axis of the rotation $O$ so that $O_\theta$ rotates counterclockwise around $v$ (i.e. $\frac{\partial}{\partial t} (O_{t\theta} (x)) = \theta v \times O_{t\theta} (x)$). We then have
\[
\tfrac{\partial \tau}{\partial t} (t,s) = \theta v \times \tau (t,s) = \theta v \times O_{t\theta} (\sigma (s))\, .
\]
On the other hand
\[
\tfrac{\partial \tau}{\partial s} (t,s) = O_{t\theta} (\dot{\sigma} (s))\, .
\]
Since $O_{t\theta} (v) =v$ and rotations of $\mathbb R^3$ commute with the cross product, we can write
\[
\tfrac{\partial \tau}{\partial t} \times \tfrac{\partial \tau}{\partial s} = \theta
O_{t\theta} ((v \times \sigma (s)) \times \dot{\sigma} (s))\,
\]
and thus
\[
\int_0^1 \int_0^\alpha \left|\tfrac{\partial \tau}{\partial t} \times \tfrac{\partial \tau}{\partial s}\right|\, ds\, dt = |\theta|
\int_0^\alpha \left|(v \times \sigma (s)) \times \dot{\sigma} (s)\right|\, ds\, .
\]
Observe next that
\[
(v\times \sigma) \times \dot{\sigma} = \sigma (v\cdot \dot{\sigma}) -
v (\sigma\cdot \dot{\sigma}) = \sigma (v\cdot \dot{\sigma})\, ,
\]
given that $\dot\sigma$ is orthogonal to $\sigma$. Therefore
\[
\int_{\mathbb S^2} \# \tau^{-1} (\{x\})\, d\mathcal{H}^2 (x)
= |\theta| \int_0^\alpha |v\cdot \dot\sigma (s)|\, ds\, .
\]
Consider now the plane containing the geodesic $\sigma$ and observe that we can choose an orthonormal basis $\{w_1,w_2\}$ for it with the property that $w_1$ is orthogonal to $v$. $\sigma (s)$ is now equal to 
$w_1\cos(s-\beta)+w_2\sin(s-\beta)$ for some angle $\beta$.
Then $\dot\sigma(s) = -w_1\sin(s-\beta)+w_2\cos(s-\beta)$, 
and hence $|v\cdot\dot\sigma(s)| =|v\cdot w_2|\,|\cos(s-\beta)|
\leq|\cos(s-\beta)|$. Changing variables in the integral we can therefore estimate
\[
\int_{\mathbb S^2} \# \tau^{-1} (\{x\})\, d\mathcal{H}^2 (x)
\leq |\theta| \int_{-\beta}^{\alpha-\beta} |\cos (s')|\, ds' \, .
\]
Maximizing the right-hand side in $\beta$ (and using $0\leq \alpha \leq \pi$) we easily see
\[
\int_{\mathbb S^2} \# \tau^{-1} (\{x\})\, d\mathcal{H}^2 (x)
\leq |\theta| \int_{-\frac{\alpha}{2}}^{\frac{\alpha}{2}} |\cos (s')|\, ds'
= 2 |\theta| \sin \tfrac{\alpha}{2}\, .
\]
\end{proof}

\subsection{Proof of \texorpdfstring{\cref{p:stability}}{p:stability}}\label{s:stability} 
We show only the first conclusion; the proof of the second is entirely analogous.

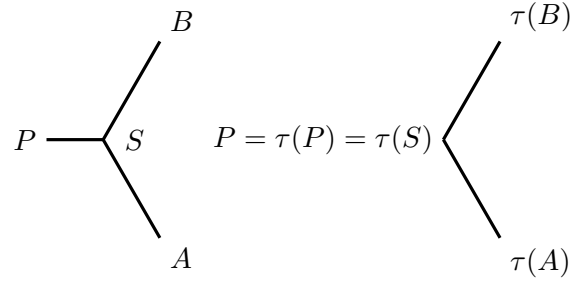
\begin{figure}
\begin{tikzpicture}[scale=1.5]
\draw[very thick] (-0.5, 0) -- (0,0) -- (0.5, {sqrt(3)/2});
\draw[very thick] (0,0) -- (0.5,{-sqrt(3)/2});
\draw[very thick] (3,0) -- (3.5,{sqrt(3)/2});
\draw[very thick] (3,0) -- (3.5,{-sqrt(3)/2});
\node[left] at (-0.5,0) {$P$};
\node[right] at (0.1,0) {$S$};
\node[below right] at (0.5,{-sqrt(3)/2}) {$A$};
\node[above right] at (0.5,{sqrt(3)/2}) {$B$};
\node[left] at (3,0) {$P= \tau (P) = \tau (S)$};
\node[below right] at (3.5,{-sqrt(3)/2}) {$\tau (A)$};
\node[above right] at (3.5,{sqrt(3)/2}) {$\tau (B)$};
\end{tikzpicture}
\caption{The change in topology at the pin of case (a1).}\label{f:topology-change-1} \end{figure}

\begin{figure}
\begin{tikzpicture}[scale=1.5]
\draw[very thick] (-1,0) -- (-0.3, -0.1) -- (0,0) -- (0.5, {sqrt(3)/2});
\draw[very thick] (0,0) -- (0.5,{-sqrt(3)/2});
\draw[very thick] (2,0) -- (3,0) -- (3.5,{sqrt(3)/2});
\draw[very thick] (3,0) -- (3.5,{-sqrt(3)/2});
\node[below] at (-0.3,-0.1) {$P$};
\node[right] at (0.1,0) {$S$};
\node[left] at (-1,0) {$C$};
\node[below right] at (0.5,{-sqrt(3)/2}) {$A$};
\node[above right] at (0.5,{sqrt(3)/2}) {$B$};
\node[below right] at (3.5,{-sqrt(3)/2}) {$\tau (A)$};
\node[above right] at (3.5,{sqrt(3)/2}) {$\tau (B)$};
\node[right] at (3.1,0) {$P=\tau (P) = \tau (S)$};
\node[left] at (2,0) {$\tau (C)$};
\end{tikzpicture}
\caption{The change in topology at the pin of case (a2).}\label{f:topology-change-2}   
\end{figure}
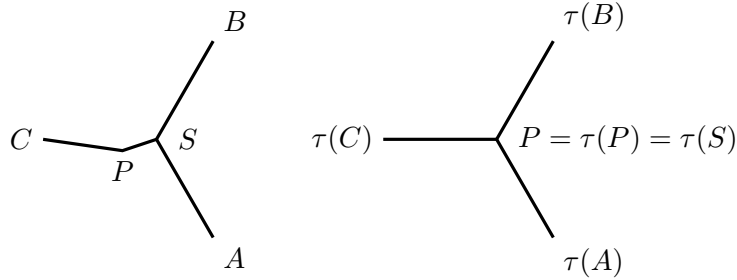

We argue that there are exactly three possible topological behaviors for $\mathcal{N}'$ at a pin $P$ compared to $\mathcal{N}$. Since the nets have the same pins, we will use the notation $\deg_{\mathcal{N}'} (P)$ and $\deg_{\mathcal{N}} (P)$ for the degrees in the corresponding nets. The first possibility is that
\begin{itemize}
\item[(a0)] No edge incident to $P$ is contracted. In this case $\deg_{\mathcal{N}'} (P) = \deg_{\mathcal{N}} (\tau (P))$ and $\tau$ maps the neighbors of $P$ bijectively onto the neighbors of $P=\tau (P)$.
\end{itemize}
Otherwise there is an edge $\sigma=[P,S]$ which is contracted. Since $S$ is a non-pin 
vertex of the net
$\mathcal N'$, it has degree three. Moreover, every pin of
$\mathcal N$ has degree at most three. Hence
\[
\deg_{\mathcal{N}} (P)=\deg_{\mathcal{N}'} (P)+1\leq3,
\]
and therefore $\deg(P)\in\{1,2\}$. We thus have exactly two possibilities:
\begin{itemize}
\item[(a1)] $\deg_{\mathcal{N}} (P) = 1$ and $\deg_{\mathcal{N}} (P') =2$. Then $S$ is the only vertex adjacent to $P$, $\tau (S) = P$ and the two neighbors of $S$ distinct from $P$ are mapped by $\tau$ bijectively onto the two neighbors of $P$, cf. \cref{f:topology-change-1}.
\item[(a2)] $\deg_{\mathcal{N}'} (P) =2$ and $\deg_{\mathcal{N}} (P)=3$. Then $S$ is one of the two vertices adjacent to $P$. The other two neighbors of $S$ are mapped by $\tau$ bijectively onto two of the three neighbors of $P$. The other neighbor of $P$ is mapped to the remaining neighbor of $P$, cf. \cref{f:topology-change-2}.
\end{itemize}

\medskip

{\bf Definition of $\Sigma'$ and $\Phi'$.} Recall the subcomplexes $\Gamma$ and $L$ of $\Sigma$ as in \cref{t:comp-to-map}, and denote by $\Gamma_{bd}$ the subcomplex of those vertices and edges of $\Gamma$ which are mapped by $\Phi$ into $E \cap \partial K$. Observe that $L$ and $\Gamma_{bd}$ have no edges in common and the vertices they have in common are mapped onto the two pins, which we denote by $Q_+$ and $Q_-$ (with the convention that $\Phi (Q_\pm) = \pm P$). First of all, we impose that $\Sigma'$ contains $\Sigma$ as a subcomplex. Next, we add the following simplices to $\Sigma'$. As it will become clear, it is actually convenient to introduce immediately also the map $\Phi'$. Being piecewise linear, it will be sufficient to define it on all the vertices. On the ones common to $\Sigma$ it coincides with $\Phi$. 

\smallskip

{\bf New vertices.} For every vertex $V\in \Gamma_{bd}$ which is not a pin we add $V'$ and we prescribe that $\Phi' (V') = \iota (\Phi (V))$. If at the pin $P = \Phi (Q_+) \in \mathcal{N}'$ there is a change of topology, then $\mathcal{N}'$ has an additional special vertex $S_+$ in $\mathcal{N}'$ which is mapped by $\tau$ to $P = \Phi (Q_+)$. In this case we add one vertex to $\Sigma'$ as well, we denote it by $Z_+$ and we set $\Phi' (Z_+) = S_+$. The same is applied at the pin $-P= \Phi (Q_-)$ and, if there is a change of topology, the corresponding vertices are denoted by $S_-$ and $Z_-$.

\smallskip

{\bf New edges.} For every $V\in \Gamma_{bd}$ which is not a pin we add the edge $[V,V']$. 

Next, for every edge $[V_1, V_2]$ in $\Gamma_{bd}$ such that neither $V_1$ nor $V_2$ is a pin, we add the edge $[V'_1, V'_2]$ and one of the two diagonals $[V_1, V_2']$ and $[V_1', V_2]$: which diagonal we choose is irrelevant.

If at the pin $P = \Phi (Q_+)$ there is no change in topology, then, for every $V\in \Gamma_{bd}$ adjacent to it we add the edge $[V', Q_+]$. If there is a change in topology we add $[Z_+,Q_+]$. Furthermore we include the following additional edges, depending on the two model alternatives (a1) and (a2). 
\begin{itemize}
\item In alternative (a1), we let $A$ and $B$ be the two vertices of $\Gamma_{bd}$ adjacent to $Q_+$ and we add the edges $[Z_+, A']$, $[Q_+,A']$, $[Z_+,B']$, and $[Q_+, B']$. 
\item In the alternative (a2), we let $A$, $B$, and $C$ be the three vertices of $\Gamma_{bd}$ adjacent to $Q_+$. Then we denote by $A'$, $B'$, and $C'$ the corresponding vertices added to form $\Sigma'$, and we ensure that the labels are assigned so that $\Phi' (A')$ and $\Phi' (B')$ are adjacent to $S_+$ in $\mathcal{N}'$, while $\Phi' (C')$ is adjacent to $P=\Phi' (Q_+)$. Then we add the edges $[A', Z_+], [B', Z_+], [A', Q_+]$, $[B',Q_+]$, and $[C',Q_+]$. 
\end{itemize}
We proceed in the same way at the pin $-P = \Phi (Q_-)$, with $Q_-$ and $Z_-$ replacing $Q_+$ and $Z_+$ in the above constructions.

\smallskip

{\bf New triangles.} 
Consider any pair $A, B\in \Gamma_{bd}$ of adjacent vertices which are not pins and not special vertices. We label them so that $[A,B']$ is the ``diagonal edge'' added to $\Sigma'$ in the previous step. Then we add the two triangles $[A,A',B']$ and $[A,B,B']$. 

At the pin $Q_+$ for every vertex $V\in \Gamma_{bd}$ adjacent to $Q_+$ we add the triangle $[Q_+,V,V']$. If there is a change of topology we add the triangles $[Q_+,Z_+,A']$ and $[Q_+,Z_+,B']$. 
We proceed in the same way at the pin $Q_+$, with $Q_-$ and $Z_-$ replacing $Q_+$ and $Z_+$ in the above constructions.

\smallskip

{\bf Definition of $\Gamma'_{bd}$ and $\Gamma'$.} $\Gamma'_{bd}$ is a one-dimensional simplicial complex. Its vertices are $Q_+$, $Q_-$ and all the new vertices added to $\Sigma'$. Moreover, given an edge $[A,B]$ of $\Sigma'$ connecting two vertices of $\Gamma'_{bd}$ we include it in the complex $\Gamma'_{bd}$ if the points $\Phi' (A)$ and $\Phi' (B)$ are adjacent in the net $\mathcal{N}'$. 

\smallskip

The subcomplex $\Gamma'$ is finally defined as the union of the subcomplexes $L$ and $\Gamma'_{bd}$.

\medskip

{\bf Verification that \cref{t:comp-to-map} is applicable.} First of all, $\Phi'$ maps the vertices of $\Sigma$ in the convex hull $K$ of the vertices of $\mathcal{N}$, and the new vertices in the vertices of $\mathcal{N}'$. Hence, by piecewise linearity, $\Phi'$ maps $\Sigma'$ into $K_c$, which is the convex hull of $\mathcal{V} (\mathcal{N}) \cup \mathcal{V} (\mathcal{N}')$. The fact that $\Phi'$ is an embedding of $\Gamma'_{bd}$ onto $\partial K' \cap E'$ is obvious from the construction. The fact that it is an embedding of $L$ onto $\ell \cap K_c$ is also obvious, since it coincides with $\Phi$ and $\ell \cap K_c=\ell \cap K$. Finally, the only faces in common between $\Gamma'_{bd}$ and $L$ are the points $Q_+$ and $Q_-$ and these are mapped onto $P$ and $-P$: since these two points form the intersection of $\partial K'\cap E'$ with $\ell \cap K_c$, we conclude immediately that $\Phi'$ is an embedding of $\Gamma'$ onto $(\partial K'\cap E')\cup (\ell\cap K')$ and that $\Phi' (\Sigma')\cap \ell = K' \cap \ell$.

Having verified \cref{t:comp-to-map}(i) we proceed to verify \cref{t:comp-to-map}(ii). Consider a face $F$ with the property that $(\Phi')^\circ (F) \cap \ell \neq \emptyset$. If $F$ is a face of $\Sigma$, then $F$ is by assumption a face of $L$. If $F$ is not a face of $\Sigma$ and it is a vertex, then we have already proved above that it cannot be mapped into $\ell$. If $F$ is a triangle $T$ we then distinguish two cases.

First of all consider a $T$ which does not have one of the pins $Q_\pm$ among its vertices.
Then this triangle is of type $[A,B,B']$ or $[A,A',B']$, where both $A$ and $B$ are neither a pin nor a special vertex of $\mathcal{N}$ and $A', B'$ are the corresponding vertices added to $\Sigma'$ following the description given above. 
In particular we can apply condition (c) in \cref{p:stability} to the four points $U = \Phi' (A) = \Phi (A)$, $V= \Phi' (B) = \Phi (B)$, $U' = \Phi' (A')$, and $V'= \Phi' (B')$; by definition $U' = \iota (U)$ and $V' = \iota (V)$. Condition (c) in \cref{p:stability} provides an open half-space $H_+$ with $\ell \subset \partial H_+$ and $U,V,U', V'\in H_+$. Since the three vertices of $\Phi' (T)$ are a subset of $\{U,V,U', V'\}$ we can use therefore convexity to establish that the closed convex set $\Phi' (T)$ cannot intersect $\ell$. 

Consider next a triangle $T$ which has a pin among its vertices, and assume without loss of generality that it is $Q_+$. Then $T$ can be of two types:
\begin{itemize}
\item $T=[Q_+, A, A']$, where $A$ is a vertex adjacent to $Q_+$ in $\Gamma_{bd}$ which is not a special vertex. Then the segment $[\Phi' (A), \Phi' (A')]$ does not intersect $\ell$, since it belongs to a triangle which has already been examined: thanks to (d) $A$ must be adjacent to a vertex which is neither a pin nor a special vertex. In particular the triangle $\Phi' (T)$ intersects $\ell$ only at the point $\Phi' (Q_+) = P$.
\item $T=[Q_+,Z_+,A']$. But then $S_+:= \Phi' (Z_+)$ and $W:= \Phi' (A')$ are two adjacent vertices of the net $\mathcal{N}'$ which are not pins. In particular the edge connecting them in the net is a length-minimizing geodesic of the sphere strictly shorter than $\pi$, which does not contain any pin. We conclude that the corresponding segment $[S_+,W]$ does not intersect $\ell$. In particular, for the triangle $\Phi' (T)$, which is the convex hull of $[S_+,W]$ and $P$, we conclude that its only intersection with $\ell$ is the point $\Phi' (Q_+)= P$.
\end{itemize}
We have thus verified that, for every new triangle $T$, either $\Phi'$ maps it into $\mathbb R^3\setminus \ell$, or the only point of $T$ which is mapped in $\ell$ is a vertex. In particular, since every edge added to form $\Sigma'$ is the side of some new triangle $T$, this verifies condition (ii) in \cref{t:comp-to-map} also for the new edges. 

We finally check condition (iii) which is about the collapsibility of $\Sigma'$ onto $L$. We start noticing that $([A,A',B'], [A',B'])$ is a collapse pair for any edge $\sigma = [A, B]$ of $\Gamma_{bd}$ such that $A$ and $B$ are not pins. We thus collapse all of these triangles and arrive at a second simplicial complex $\Sigma''$. Furthermore, for every such edge, $([A,B,B'], [A,B'])$ is now a collapse pair and, after collapsing all of these triangles we arrive at a third complex $\Sigma'''$. At this point, we can collapse the pair $([A,A'],A')$ for every $A$ which is not adjacent to a pin in $\Gamma_{bd}$. After performing all these collapses we arrive at a fourth complex $\Sigma^{(4)}$. All the triangles of $\Sigma^{(4)}$ which are not triangles of $\Sigma$ have now a pin among their vertices, while every edge and every vertex of $\Sigma^{(4)}$ which are not in $\Sigma$ are contained in one such triangle. Consider first the triangles which have $Q_+$ as a vertex. Among them, for those of the form $[Q_+,A,A']$, the pair $([Q_+,A,A'],[A,A'])$ is a collapse pair, because, thanks to condition (d) in \cref{p:stability}, the vertex $A$ cannot be adjacent to $Q_-$ and in particular any other triangle of $\Sigma'$ which had $[A,A']$ as a side has already been collapsed. We therefore collapse all such triangles and all analogous triangles at $Q_-$. We are now left with the triangles of the form $[Q_+,Z_+,A']$ (which exist only if there was a change of topology at $\Phi (Q_+) = P$). Since $([Q_+,Z_+,A'],[Z_+,A'])$ is a collapse pair, we can collapse all these triangles and likewise all the analogous triangles with one vertex in $Q_-$. 
Compared to $\Sigma$, the resulting new complex $\Sigma^{(5)}$ has no additional triangles. If there are additional edges of the form $[Q_+,V]$ and $[Q_-,V]$ where $V$ is a vertex of $\Gamma'_{bd}$, the pairs $([Q_\pm,V], V)$ are all collapse pairs, and after collapsing them we arrive at the complex $\Sigma$.

\medskip

{\bf Area estimate.} We now need to verify \cref{e:stab-estimate}. We count the triangles $\Phi' (T)$ added to $\Phi (\Sigma)$. Each triangle $\Phi' (T)$ is associated to an edge $\sigma\in \mathcal{E} (\mathcal{N})$ in the following way:
\begin{itemize}
\item If $T=[A,B,B']$ or $[A, A', B']$ and neither $A$ nor $B$ is a pin in $\Gamma_{bd}$, then $\sigma=[\Phi (A) , \Phi (B)]$;
\item otherwise $T=[Q_\pm, A, A']$ or $[Q_\pm,Z_\pm, A']$ for some $A\in \Gamma_{bd}$ adjacent to the pin $Q_\pm$; in this case we set $\sigma = [\Phi (Q_\pm), \Phi (A)]$.
\end{itemize}
Under this count, for any $\sigma\in \mathcal{E} (\mathcal{N})$ which is not incident to a pin, there are exactly two triangles, while for an edge incident to a pin, there may be two triangles or just one. 

If $\sigma\in \mathcal{E} (\mathcal{N})$ is not adjacent to a pin then the sum of the areas of the two added triangles is obviously bounded by $D (\sigma)$. Assume next that $\sigma = [U,P]$ and that there is no change of topology at $P$. Then we have $\iota (P)=P$ and it follows immediately that 
\[
D (\sigma ) = \mathcal{H}^2 ([U, \iota (U), P])
\]
which is the area of the only triangle added. In the other case, $\iota (P)$ is a special vertex $S$ and we have added at most the two triangles $[P, U, \iota (U)]$ and $[P, S, \iota (U)]= [P, \iota (P), \iota (U)]$: the sum of their areas is still bounded by $D(\sigma)$.

\section{Area-Decreasing Competitors and proof of \texorpdfstring{\cref{t:main}}{main-theorem}}\label{s:terzo-programma}
In this section we employ all the results proven up to now to conclude the proof of \cref{t:main}. Also this final part is computer-assisted.
Let us begin by summarizing what needs to be proven and what we already know.

Consider an Almgren minimizing cone $E$ with respect to a line $\ell$ as in \cref{t:main}, let $\mathcal N$ be its link and denote by $\mathcal N'$ the corresponding \emph{once-derived} net. We want to show that $E$ belongs to the \cref{lista}.

If $\mathcal N'$ is extinct, we have shown in \cref{t:elementary} that the cone $E$ belongs to the \cref{lista}. So, if $\mathcal N'$ is extinct there is nothing left to prove.

On the other hand, if $\mathcal N'$ is not extinct, we have proven in \cref{thm:part2} that its topology is encoded by one of the marked spherical graphs of \cref{p:second-topological-classification} with identifiers
\begin{align*}
    \mathscr{L} := \{&28, 99, 158, 182, 250,
    319, 325, 346, 361, 421,
    468, 484, 520, 614,
    \\ &681,
    712, 736, 779, 827, 909,
    920, 938, 1029, 1096, 1100,
    1122\} \,.
\end{align*}
Moreover, we have also proven \emph{approximate metric rigidity} for all these topologies apart from $827$ and $1029$.
Thanks to \cref{t:human-classification}, if the once-derived topology is $827$ then $\mathcal N$ is an element of the fat fish family (see \cref{d:pesce-grasso}), while if the once-derived topology is $1029$ then $\mathcal N$ is an element of either the cubical family or the goggles family. If $\mathcal N$ is an element of the cubical family, which belongs to the \cref{lista}, we do not have anything to prove.

In \cref{t:human-classification}, we have also proven that the once-derived topology cannot be $158$ or $681$ and that if the once-derived topology is one of $920$ or $1100$ then $E$ belongs to the \cref{lista}.

Summing up what we said, we are left to show the following proposition.

\begin{proposition}\label{p:final}
None of the cones whose links belong to the following list of stationary nets with two antipodal pins is Almgren minimizing with respect to the line $\ell$ passing through the pins:
\begin{itemize}
    \item The nets whose once-derived topology belongs to
    \begin{align*}
        \mathscr{L}_{rigid} := \{
        182, 250, 319, 325, 361, 
        421, 468, 484, 520, 614,
        712, 736,  909, 938,  1122\} \,.
    \end{align*}
    \item The nets whose once-derived topology belongs to
    \begin{align*}
        \mathscr{L}_{flexible} := \{28, 99, 346, 779, 1096\} \,.
    \end{align*}
    \item Any element of the goggles family (corresponding to the once-derived topology $1029$).
    \item Any element of the fat fish family (corresponding to the once-derived topology $827$).
\end{itemize}
\end{proposition}

In the above schematic plan, we have split the once-derived topologies to rule out into two seemingly arbitrary lists $\mathscr{L}_{rigid}$ and $\mathscr{L}_{flexible}$. 
Let us explain why we came up with this partition.
Any of the \emph{approximately rigid} once-derived topologies we are considering is associated to an Approximately Metrized Marked Spherical Graph (cf. \cref{def:ammsg}). Our presentation will be informal and will assume that the reader understands the content of \cref{subsec:embed-or-exclude}. For the details about the omitted arguments, we refer the reader directly to \verb|find_competitor.cpp| in \cite{repo}.

Observe that if $(G, \mathcal I)$ is an AMMSG containing the once-derived net $\mathcal N'$, one could approximately compute (up to isometries of the whole net) the positions of the two pins that were removed from $\mathcal N$ to obtain $\mathcal N'$: it is sufficient to intersect all the equators containing edges adjacent to a pin (recall \cref{lem:approximate-isometry}). There are two cases:
\begin{itemize}
    \item If the equators do not all (approximately) coincide, then we can intersect them and (approximately) identify the position of the pins and so (approximately) determine the location of all the vertices of $\mathcal N$. This is the case for the topologies in $\mathscr{L}_{rigid}$. One exception is the topology $736$, whose approximate metrization allows us to (approximately) identify the two pins even if all the edges adjacent to a pin belong to the same equator. 
    \item If the equators are all (approximately) the same, then the position of the pins is not (approximately) identified. This is the case for the topologies in $\mathscr{L}_{flexible}$. In this case, to (approximately) determine the location of the pins of $\mathcal N$ we must introduce an additional parameter which determines the length of one (arbitrarily chosen) edge adjacent to a pin.
\end{itemize}

Finally, let us remark that we employ \cref{p:pin-touching} to conclude that, if the respective cone is area-minimizing, the parameters that define a net in $1096$ and in the fat fish family must be (quantitatively) away from the extreme points of their domains (this is the only place where we use \cref{p:pin-touching}). Being able to state a-priori that every pin $P$ cannot be too close to an edge that is not adjacent to an edge adjacent to $P$ allows us to avoid some pathological cases that would make the rest of our strategy fail.

\subsection{Ambiguity in the topology of $\mathcal N$}
Let us understand whether the topology of $\mathcal N'$ determines the topology of $\mathcal N$. Let us assume that we are given an AMMSG $(G, \mathcal I)$ and, if necessary, the additional parameter that is necessary to determine the (approximate) location of the pins. We refer the reader to \cref{d:once-derived} to refresh the exact method to obtain $\mathcal N'$ from $\mathcal N$.

If the position of a pin is distinct from the position of any vertex of $\mathcal N'$ (observe that we know the positions only approximately, so we can only determine whether they are surely distinct, not whether they are surely the same), then there is no ambiguity and the topology of $\mathcal N$ is uniquely determined from the topology of $\mathcal N'$.
If we cannot rule out the possibility that a pin of $\mathcal N$ coincides with a vertex of $\mathcal N'$ using the data of $(G,\mathcal I)$, then the topology of $\mathcal N$ is not uniquely determined. The two possible ambiguous cases are represented in \cref{f:topology-change-1,f:topology-change-2}.

\subsection{Reducing all cases to the same subproblem}
It remains to prove a list of statements all of the form ``no net in this family generates an Almgren-minimizing cone''. We will now explain how to reduce this list to a (much longer) list of more structured subproblems.

We have observed that the nets whose once-derived topology belongs to $\mathcal L_{flexible}$ are (approximately) parametrized by a single real parameter (and their topology is defined up to the ambiguities described in the previous section). Similarly, the family of goggles is parametrized by a single parameter. The fat fish family instead is parametrized by two real parameters.

Observe that a geodesic net with antipodal pins can be encoded by a triple $(G, \{p_1, p_2\}, e)$, where $G$ is a spherical graph, $\{p_1, p_2\}$ are two vertices of $G$ (the pins), and $e:\mathcal V(G)\to\mathbb S^2$ is an embedding of the vertices into the unit sphere such that any two adjacent vertices are not antipodal\footnote{If two vertices are antipodal, then the triple $(G, \{p_1, p_2\}, e)$ does not contain enough data to reconstruct the edge connecting the two vertices. In this case, one can add an artificial vertex in the middle of the edge and the resulting geodesic net, which coincides with the original one as a set, can be fully encoded by this framework.} and $e(p_1)=-e(p_2)$.

\begin{definition}[Approximate Geodesic Net with two pins]
    An approximate geodesic net with two pins is a triple $(G, \{p_1, p_2\}, \tilde e)$, where $G$ is a spherical graph, $\{p_1, p_2\}$ are two vertices of $G$ (the pins), and $\tilde e:\mathcal V(G)\to\mathcal P(\mathbb S^2)$ is a map from the vertices of $G$ into the subsets of the unit sphere.

    A realization $(G, \{p_1, p_2\}, e)$ of an approximate geodesic net $(G, \{p_1, p_2\}, \tilde e)$ is a geodesic net with antipodal pins on the same graph (with the same pins) such that $e(v)\in \tilde e(v)$ for all $v\in \mathcal V(G)$.
\end{definition}

Now, assume that the once-derived net associated to $\mathcal N$ has a topology in $\mathcal L_{rigid}$. Then, applying \cref{lem:approximate-isometry}, we can construct an Approximate Geodesic Net $(G, \{p_1, p_2\}, \tilde e)$ with two pins  such that, up to isometry of $\mathbb S^2$, the net $\mathcal N$ is a realization of $(G, \{p_1, p_2\}, \tilde e)$ and $\tilde e(v)$ is a \emph{small} set for all $v\in\mathcal V(G)$.
What we have just said is not correct, because the topology of $\mathcal N$ cannot always be inferred from the topology of its once-derived net $\mathcal N'$. We solve this annoying issue as follows. We choose $G$ to be the more \emph{complicated} topology that is consistent with the once-derived topology. Then, $\mathcal N$ is either a realization of $(G, \{p_1, p_2\}, \tilde e)$ or it is approximated by a realization of $(G, \{p_1, p_2\}, \tilde e)$. More precisely, in the latter case, for any $\eps>0$ there is a realization $\mathcal{\tilde N}$ of $(G, \{p_1, p_2\}, \tilde e)$ such that we can apply \cref{p:stability} to $\mathcal N$ and $\mathcal{\tilde N}$ with the sum $\sum_{\sigma\in\mathcal E(\mathcal{\tilde  N})}D(\sigma)$ appearing in \cref{e:stab-estimate} being smaller than $\eps$.
In particular, if we can find a competitor (in the sense of \cref{t:comp-to-map}) for all realizations of $(G, \{p_1, p_2\}, \tilde e)$ that decreases the area by a fixed amount, then \cref{p:stability,t:comp-to-map} guarantee that $\mathcal N$ is not Almgren minimizing.
Given that we will show the existence of such a competitor, we are allowed to assume that $G$ represents the \emph{most complicated} topology, and from now on we will not discuss this matter anymore.

Analogously, let us instead assume that the once-derived net associated to $\mathcal N$ has a topology in $\mathcal L_{flexible}$ and let $I$ be a \emph{small} interval such that the length of an (arbitrarily fixed) edge adjacent to a pin belongs to $I$. Then, using only the information on the topology of $\mathcal N'$ together with the interval $I$, we can construct an Approximate Geodesic Net with two pins satisfying the same constraints of the previous paragraph.

The same can be done also if $\mathcal N$ is a goggle and we know that the parameter identifying it belongs to a small interval $I$. 
If it is a fat fish, and we know that the two parameters defining it belong to two small intervals $I, J$, the argument is equivalent but we apply \cref{p:stability-fat} instead of \cref{p:stability}.

In all cases, if we allow ourselves to discretize the zero, one, or two-dimensional space of parameters associated to a topology in sufficiently small regions\footnote{We refer the reader to the code for details on how we do this.}, we can reduce the proof of \cref{t:main} to many instances of the following problem:
\begin{center}
    Given an Approximate Geodesic Net with two pins $(G, \{p_1, p_2\}, \tilde e)$, prove that the cone over any of its realizations is not Almgren minimizing.
\end{center} 

\subsection{Solving the subproblem}
We describe how to treat one instance of the above-mentioned problem.
Consider a geodesic net with antipodal pins $\mathcal N := (G, \{p_1, p_2\}, e)$ chosen so that, for all $v\in\mathcal V(G)$, the quantity
\begin{equation*}
    \sup_{p\in \tilde e(v)} |p-e(v)|
\end{equation*}
is \emph{small}. In our implementation, $\tilde e(v)$ is always a parallelepiped (intersected with $\mathbb S^2$) and we choose $e(v)$ to be its center (with an exception at the pins, to guarantee antipodality).
Here again, there is an exception related to the fact that the topology is ambiguous when the location of the pin falls very close to the location of a vertex of the once-derived net. Let us focus on the pin $p_1$, the analogous argument shall be repeated verbatim for the other pin $p_2$. If there is a vertex $v\not\in\{p_1, p_2\}$ such that the minimum distance between $\tilde e(p_1)$ and $\tilde e(v)$ is below a given threshold, then $\mathcal N$ will be induced by $(\tilde G, \{p_1, p_2\}, e)$ where the vertices $p_1$ and $v$ are collapsed onto a single vertex in $\tilde G$. So, in this case, the topology of $\mathcal N$ and the topology of the realizations of $(G, \{p_1, p_2\}, \tilde e)$ do not match but there is a controlled topological change (cf. \cref{d:topology-change}) from any $\mathcal N$ to any realization. Treating explicitly this case would lead to some very cumbersome notation without yielding any further understanding, so we will assume from now on that $\mathcal N$ has the topology of $G$ and we will just relieve the worried reader by making it clear that the only place where this would actually play a role is in the application of \cref{p:stability} which conveniently supports nets which are obtained one from the other via a controlled topological change.\footnote{More delicate, and more important, is the matter of \emph{why} we decide to opt for this complicated choice. The main reason is that, eventually, we will want to construct a competitor (in the sense of \cref{t:comp-to-map}) for $\mathcal N$ and to do so we will employ Brakke's surface evolver, which tends to behave erratically when one vertex is very close to a pin. Furthermore, if we did not employ this strategy, then proving \cref{p:stability}(c) would become much more cumbersome, or impossible, when a vertex falls very close to a pin.}

Take any realization $(G, \{p_1, p_2\}, e')$ of $(G, \{p_1, p_2\}, \tilde e)$. Observe that there exists an isometry $T$ (depending on the realization) such that
\begin{itemize}
    \item $\norm{T-\id}_{op} \le \sup_{p\in \tilde e(p_1)} |p-e(p_1)| =: \lambda$,
    \item $T(e'(p_1)) = e(p_1)$ and $T(e'(p_2)) = e(p_2)$.
\end{itemize}
In particular, the geodesic net with antipodal pins encoded by $(G, \{p_1, p_2\}, T\circ e')$ is a realization of the approximate geodesic net $(G, \{p_1, p_2\}, f)$ where $f$ is defined as
\begin{itemize}
    \item $f(p_1)=\{e(p_1)\}$ and $f(p_2)=\{e(p_2)\}$,
    \item $f(v) = B(\tilde e(v), \lambda)$ for all $v\in\mathcal V(G)\setminus\{p_1, p_2\}$.
\end{itemize}
Since being Almgren minimizing is invariant under isometries, if we are able to prove that none of the cones over a realization of $(G, \{p_1, p_2\}, f)$ is Almgren minimizing, then we have solved the problem.

What we have achieved in the last paragraph is to restrict our study on a family of nets with the two pins fixed in the same position as the pins of the special net $(G, \{p_1, p_2\}, e)$ that we had previously chosen. The price that we had to pay was the enlargement of $\tilde e(v)$ to $f(v) = B(\tilde e(v), \lambda)$ for a positive $\lambda>0$ (which is as small as the precision with which we initially knew the position of the pins).
This trick is necessary because later we will want to apply \cref{p:stability}.

Let us assume that for any pair $(U, V)$ of adjacent vertices of $G$ that are not pins, there is a half-space $H_+$ with $e(p_1), e(p_2) \in \partial H_+$ and $e(U), e(V)\in H_+$ and $f(U), f(V)\subset H_+$.
This assumption is mild given that all geodesic edges have length strictly less than $\pi$ and we assume to have an instance such that $f(V)$ is a rather \emph{small} set containing $e(V)$ for any $V\in\mathcal V(G)$.

Let $K$ be the convex hull of $\{e(v):\, v\in\mathcal V(G)\}$. Let $E$ be the cone over the geodesic net $(G, \{p_1, p_2\}, e)$.

Now, imagine that we are able to create an abstract $2$-dimensional simplicial complex $\Sigma$ and a piecewise linear map $\Phi:\Sigma\to \R^3$ satisfying the assumptions of \cref{t:comp-to-map}. Thanks to \cref{t:comp-to-map}, if
\begin{equation}\label{eq:local349}
    \mathcal H^2(\Phi(\Sigma)) < \mathcal H^2(E\cap K),
\end{equation}
then $E$ is not Almgren minimizing. We will soon notice that we need a slightly stronger inequality than \cref{eq:local349} for our purposes.

Fix a realization $(G, \{p_1, p_2\}, e')$ of $(G, \{p_1, p_2\}, f)$ and let $K'$ be the convex hull of $\{e'(v):\, v\in\mathcal V(G)\}$. Let $E'$ be the cone over the geodesic net $(G, \{p_1, p_2\}, e')$.

By applying\footnote{For simplicity, we describe only the details of the strategy for nets that are not in the fat fish family. For nets in the fat fish family, the argument is in fact simpler because one can directly apply \cref{p:stability-fat} and skip some of the computations described.} \cref{p:stability,t:comp-to-map}, for any $\eps>0$, there is a map $F:\R^3\to\R^3$ which is the identity outside of $B_1$, satisfies the requirements of \cref{d:Almgren-minimizing}, and
\begin{equation*}
    \mathcal H^2(F(E'\cap B_1)) \le \mathcal H^2(E'\cap B_1\setminus K') + \mathcal H^2(\Phi(\Sigma)) + \eps + \sum_{\sigma\in \mathcal E(G)} D_{e, e'}(\sigma) ,
\end{equation*}
where $D_{e, e'}((U, V))$ is defined as
\begin{align*}
    D_{e, e'}((U, V)) := 
    \max\{&\mathcal H^2([e(U), e(V), e'(V)])
    +
    \mathcal H^2([e(U), e'(U), e'(V)]), \\
    &\mathcal H^2([e(V), e(U), e'(U)])
    +
    \mathcal H^2([e(V), e'(V), e'(U)])\} \,.
\end{align*}

Therefore, if we could show that 
\begin{equation}\label{eq:local742}
    \mathcal H^2(\Phi(\Sigma)) + \sum_{\sigma\in \mathcal E(G)} D_{e, e'}(\sigma)
    <
    \mathcal H^2(E'\cap K')
\end{equation}
we would deduce immediately that $E'$ is not Almgren minimizing as desired.
First of all, observe that this condition is \emph{almost} equivalent to \cref{eq:local349}. More importantly, this condition can be transformed in a statement which does not involve $E'$ or $K'$ at all arguing as follows.

Define
\begin{equation*}
    D_{e}(\sigma) := \sup_{e':\,e'(v)\in f(v)\forall v\in\mathcal V(G)} D_{e, e'}(\sigma),
\end{equation*}
and observe that
\begin{equation*}
    \mathcal H^2(E'\cap K')
    =
    \sum_{(U, V)\in\mathcal E(G)} \mathcal H^2([0_{\R^3}, e'(U), e'(V)])
    \ge 
    \sum_{(U, V)\in\mathcal E(G)} \inf_{p\in f(U), q\in f(V)} \mathcal H^2([0_{\R^3}, p, q]) \,.
\end{equation*}
So, \cref{eq:local742} is implied by
\begin{equation*}
    \mathcal H^2(\Phi(\Sigma))
    <
    \sum_{(U, V)\in\mathcal E(G)} \inf_{p\in f(U), q\in f(V)} \mathcal H^2([0_{\R^3}, p, q])
    -\sum_{\sigma\in \mathcal E(G)}  D_{e}(\sigma) \,. 
\end{equation*}
In this latter inequality, all terms are computable (or at least approximately computable) and thus we have reduced the subproblem to finding a good enough competitor $(\Sigma, \Phi)$ for one particular antipodal net $(G, \{p_1, p_2\}, e)$.

This is, in fact, the core of the problem. As we will describe next, we are able to find a \emph{good enough} competitor in all the necessary cases and this concludes the proof.
The competitors for all the instances of the subproblem can be found in \verb|data/competitors.zip| of \cite{repo} or in the website \cite{website}.

\subsection{Constructing a Competitor}
Even though the proof itself does not depend on \emph{how} we have constructed these competitors, as far as we can verify that they respect all the required constraints, we think that it is important to document how they were constructed.

We are given a geodesic net with antipodal pins $\mathcal N$ and we want to construct a competitor, in the sense of \cref{t:comp-to-map}, for the cone over $\mathcal N$.
We may assume that the geodesic net $\mathcal N$ is not pathological in any way, in particular we may assume that no edge of the net is very short (we have spent quite some effort in the previous subsections to obtain these assumptions).

The construction proceeds in two steps: first we build a competitor that has a larger area but has a topology different from the topology of the cone, then we run Brakke's (constrained) mean curvature flow. At the end, we are just left to verify that the result is a valid competitor with small enough area, which are straight-forward checks described in \cref{t:comp-to-map}.

\subsubsection*{Cone-popping}
To alter the topology of the cone, we perform an operation that we call \emph{popping the cone} (this terminology is taken from Brakke, who introduced a pop command in his code; however our cone-popping differs from his).
The standard cone over a net $\mathcal{N}$ consists of triangles connecting every edge of $\mathcal{N}$ to the origin. Popping the cone replaces the portion of the surface inside a small ball $B(0, r)$ centered at the origin with a different simplicial structure.

\begin{definition}[Popped cone]\label{d:cone-popping}
    Let $\mathcal{N}$ be a geodesic net on $\partial B_1$ with two antipodal pins $P$ and $-P$.

    Assume that we are provided with a pair of combinatorial data $(\mathcal Z, \gamma)$:
    \begin{itemize}
        \item $\mathcal{Z}$ is a set of faces of $\mathcal{N}$, which we refer to as the \emph{excluded faces}.
        \item $\gamma$ is a $1$-dimensional subgraph of $\mathcal N$ that contains both pins.
        \item The $2$-dimensional simplicial complex associated to $\mathbb S^2\setminus \mathcal Z$ can be collapsed on its $1$-dimensional subcomplex $\gamma$.
    \end{itemize}

    The \emph{popped cone} associated to $(\mathcal Z, \gamma)$ is an abstract $2$-dimensional simplicial complex $\Sigma$ mapped into $\R^3$ via a piecewise linear map $\Phi: \Sigma \to \mathbb{R}^3$.
    The geometric realization $\Phi(\Sigma)$ is an immersed simplicial complex constructed as the union of three components, depending on a parameter $0<r<1$:
    \begin{itemize}
        \item \textbf{The outer cone.} For every edge $[u, v]$ of $\mathcal{N}$, $\Phi(\Sigma)$ contains a quadrilateral (triangulated into two faces) that connects the boundary vertices $\Phi(u)$ and $\Phi(v)$ to their scaled counterparts $r\Phi(u)$ and $r\Phi(v)$.
        \item \textbf{The capped faces.} For every face $F$ of $\mathcal{N}$ that is \emph{not} in the excluded set $\mathcal{Z}$, we add a (not necessarily flat) $2$-dimensional (triangulated) polygonal face with vertices $\{r\Phi(v):\, v\in F\}$. If the face contains a pin, the triangulation is pivoted at that pin (otherwise, any triangulation generates a valid popped cone).
        \item \textbf{The inner cone.} For every edge $[u, v]$ in $\gamma$, we add the triangle $[0_{\R^3}, r\Phi(u), r\Phi(v)]$ (note that this is a small triangle if we choose $r$ small).
    \end{itemize}
\end{definition}

The following lemma motivates our choice, even though it ultimately does not play a role in the actual argument for \cref{p:final}. 

\begin{lemma}[Popped cones are valid competitors]\label{l:popped}
    Let $\mathcal N$ be a geodesic net with two antipodal pins $P, -P$ on the line $\ell$. Let $K$ be the convex hull of the vertices of $\mathcal N$.

    Given $(\mathcal Z, \gamma)$, as in \cref{d:cone-popping}, the associated popped cone $(\Sigma, \Phi)$ is a valid competitor for the cone over $\mathcal N$ in $K$ in the sense of \cref{t:comp-to-map}.
\end{lemma}
\begin{proof}
    Let $\iota:\mathcal V(\mathcal N)\to \mathcal V(\Sigma)$ be the map such that $\Phi(\iota(V)) = V$. Moreover, for any vertex $V\in\mathcal V(\mathcal N)$, let (with a slight abuse of notation) $v=\iota(V)$ and let $v'\in\mathcal V(\Sigma)$ be the vertex such that $\Phi(v')=rV$.
    Let $o\in\mathcal V(\Sigma)$ be the vertex mapped into the origin by $\Phi$.

    We must verify the conditions (i), (ii), (iii) of \cref{t:comp-to-map}.

    Let us identify $L$ as the subcomplex connecting, in this order, the vertices $p, p', o, -p', -p$ (with a slight abuse of notation we denote by $p$ and $-p$ the vertices of $\Sigma$ mapped onto $P$ and $-P$).

    Given the definition of $L$ and the fact that $(\Sigma, \Phi)$ and $\text{cone}(\mathcal N)\cap K$ coincide in $K\setminus rK$, it is a tedious but straightforward task to verify that $(\Sigma, \Phi)$ satisfies \cref{t:comp-to-map}(i). We skip this check.

    To check the validity of (ii) we shall verify that, for any face $F$ of $\Sigma$ not in $L$, $\Phi^\circ(F)$ is disjoint from the line $\ell$.

    Let us consider a few cases:
    \begin{itemize}
        \item If $F$ is part of the \textbf{outer cone}, then it is tedious but straightforward to check that $\Phi^\circ(F)$ is disjoint from $\ell$.
        \item Similarly, if $F$ is part of the \textbf{inner cone}, then it is tedious but straightforward to check that $\Phi^\circ(F)$ is disjoint from $\ell$.
        \item The hard case is when $F$ is part of the \textbf{capped faces}. Here we must distinguish between the case in which $F$ is part of a face that contained the pin or not. If the original face of $\mathcal N$ that is being capped did not contain the pin, then it was convex. Therefore $\Phi(F)$ must belong (independently of the triangulation that was chosen) to the convex cone generated by the face that is being capped and since this face does did not contain either $P$ or $-P$, its cone is disjoint (up to the origin) from $\ell$.
        On the other hand, if the original face of $\mathcal N$ that is being capped contained a pin then we know, thanks to the special triangulation we have chosen, that either $F$ contains $p'$ or $-p'$ or it is an edge connecting $[u', v']$ where $[U, V]$ is an edge of $\mathcal N$. In both cases, it is easy to see that $\Phi^\circ(F)$ is disjoint from $\ell$.
    \end{itemize}

    Now, we conclude by showing a collapsing sequence that proves (iii).
    \begin{itemize}
        \item Consider a face $[u, v, v', u']$ (which is made of two triangles) of $\Sigma$, where $[U, V]$ is an edge of $\mathcal N$. Then, collapse it (in two steps) onto $[u, u'], [u', v'], [v, v']$. Do this for all edges $[U, V]$ of $\mathcal N$.
        \item For each vertex $V\in\mathcal N$ different from $P$ and $-P$, collapse $[v, v']$ onto $v'$.
        \item Let $\mathcal A$ be the subcomplex of $\Sigma$ whose vertices are $\{v':\, V\in\mathcal N\}$ and its faces are the \textbf{capped faces}. Let $\mathcal B$ be the $1$-dimensional subcomplex of $\mathcal A$ whose vertices are $\{v': V\in \mathcal V(\gamma)\}$ and whose edges are the ones induced naturally by $\gamma$.
        We have that $\Phi$ is an embedding restricted on $\mathcal A$ and maps it onto a subset of $\R^3$ equivalent (in the right sense) to $\mathbb S^2\setminus \mathcal Z$, where $\mathcal Z$ is the set of excluded faces. Therefore, by assumption, we deduce that we can collapse $\mathcal A$ onto $\mathcal B$ and that is exactly what we decide to do.
        \item For each edge $[U, V]$ in $\gamma$, we collapse $[o, u', v']$ onto $[o, u']$ and $[o, v']$.
        \item Finally, for each vertex $V$ in $\gamma$ that is not $P$ or $-P$, we collapse $[v', o]$ onto $o$. 
    \end{itemize} 
\end{proof}

\subsection*{Brakke's Surface Evolver}
Once we have obtained a popped cone, we evolve it using Brakke's surface evolver \cite{Brakke}. 
More precisely, we set the line $\ell$ as a parametrized boundary, and we \emph{fix} the vertices and the edges on $\partial K\cap \Sigma$ (where $\Sigma$ is the cone and $K$ is the convex hull of the vertices of the net).

Brakke's surface evolver runs a discrete mean curvature flow which is designed to let the surface converge towards a constrained minimum of the area. Let us remark that, during the evolution, Brakke's surface evolver is allowed to change the topological structure of the surface but should, in theory, guarantee that the resulting surface is a valid competitor for the initial surface, given that we started from a valid one by \cref{l:popped}.
We do not rely on the correctness of Brakke's surface evolver as we do simply check that the final surface satisfies the assumptions of \cref{t:comp-to-map} (and so, we do not rely on \cref{l:popped} either). 

A large number of ad-hoc tricks are necessary to make the algorithm converge for all our surfaces. Let us mention just a few difficulties that shall be taken care of.
\begin{itemize}
    \item Brakke's surface evolver has some bugs that lead to bad behaviors in otherwise non-pathological situations. We fixed some of these bugs; our version of Brakke's surface evolver can be found here \verb|https://gist.github.com/dario2994/25b6b268ffd76eeef8e8d36904a12c09|.
    \item If initially one of the edges of the net $\mathcal N$ is very short, the surface evolver has the (reasonable) tendency to collapse the two vertices of the boundary into one. This operation is not valid as it breaks the property that the final surface is a valid competitor for the first one. We avoid this issue by proving a-priori that no edges are short in \cref{t:effective-short-edge}.
    \item The fat fish family presents two idiosyncratic difficulties. It is a $2$-parameter family and to obtain a proof that none of its nets generates a minimizing cone we need to discretize the family into approximately $50\,000$ nets that shall be treated independently. So, we need to run a very fast flow otherwise the task becomes computationally unfeasible. Unluckily, some of these nets are very close to being minimizers (the competitor has area that is just $0.1\%$ better than the initial cone; which is also why we need to discretize the family so much to obtain a proof) and to find it we need a refined flow. Therefore, for each of these nets, we run first a very rough routine which is fast and, if it fails to find a competitor, we then run a much more refined routine (in the sense that the time-evolution steps are smaller and the grid of triangles given to the computer to describe the competitor is finer).
    Moreover, sometimes the optimal competitor produced by the surface evolver exhibits some self-intersections on the line $\ell$ (we do not fully understand why this issue does not appear for the other nets we treat). Any such self-intersection makes the resulting surface an invalid competitor for the initial cone (cf. \cref{t:comp-to-map}-(ii)). To tackle this problem, we consider three different initial \emph{popped cones} and we verify that, for each of the $\approx 50\,000$ nets, at least one of them leads to a valid competitor.
\end{itemize}

\section{Minimality: known cases and open problems}\label{s:davvero-minimi}

In this section we will discuss the minimizing properties of the cones in \cref{lista}.

\subsection{The four easy cases} The cones in (a), (b), (c), and (d) are all minimizing.

\begin{proposition}\label{p:minimizing}
All the elements listed in the cases (a), (b), (c), and (d) are Almgren minimizing with respect to $\ell$.
\end{proposition}
\begin{proof}
First of all, for (c) and (d) the statement follows immediately from \cite{Taylor}. Indeed, consider a cone $E$ as in (c) or (d). According to \cite{Taylor} the inequality \eqref{e:Almgren-minimizing} is satisfied for every bounded open set $V$ and every Lipschitz map $\Phi: \mathbb R^3 \to \mathbb R^3$ such that $\{x:\Phi (x)\neq x\}$ is compact. On the other hand the Almgren-minimizing property with respect to $\ell$ requires that the same inequality holds for a {\em strict} subset of such maps $\Phi$.

\medskip

Consider next a cone $E$ as in case (a) and (b). In case (a) $E$ consists of a half-plane $\tau$ with boundary $\ell$, while in case (b) it consists of two half-planes $\tau$ and $\tau'$ with boundary $\ell$ and forming an angle $\theta \in [\frac{2\pi}{3}, \pi]$. Without loss of generality we can then assume that 
\begin{align*}
\ell & = \{(0,0,s): s \in \mathbb R\}\\
\tau & = \{(r,0,s): s \in \mathbb R, r\geq 0\}\\
\tau' &= \{(r\cos \theta, r\sin \theta, s): s \in \mathbb R, r\geq 0\}\, .
\end{align*}
Assume that $E$ is not Almgren minimizing with respect to $\ell$ and fix a map $\Phi$ which contradicts \eqref{e:Almgren-minimizing} but satisfies the requirements of \cref{d:Almgren-minimizing}. By assumptions there are $R>\delta >0$ such that 
\[
{\rm spt}\, (\Phi - {\rm Id}) \subset (B_R (\ell) \setminus \overline{B}_\delta (\ell)) \cap \{|x_3| < R\} \, .
\]
In particular the ``slice'' 
\[
\Gamma (s) := \{x_3=s\} \cap \Phi \bigl(\overline{B}_R (\ell) \cap \{|x_3|\leq R\} \cap E\bigr)
\]
is a subset of the horizontal plane $\alpha (s) := \{x_3=s\}$ which, for every $s\in (-R, R)$, contains
\begin{itemize}
\item The two points $p_0 = (0,0,s)$ and $q= (R,0,s)$ in case (a);
\item The three points $p_0$, $q$, and $q'= (R\cos \theta, R \sin \theta, s)$ in case (b).
\end{itemize}
From \cref{l:connesso} below we infer that these points belong to the same connected component of $\Gamma (s)$, which we denote by $\Gamma^c (s)$: in case (a) this is obvious, while in case (b) we apply the lemma separately to each of the two rectangles forming $E\cap \overline{B}_R (\ell) \cap \{|x_3|\leq R\}$, which are the intersections of the latter sets with the two halfplanes $\tau'$ and $\tau$. We next apply the coarea formula to get 
\begin{align*}
\int_{-R}^R \mathcal{H}^1 (\Gamma (s))\, ds & \leq 
\mathcal{H}^2 \bigl(\Phi (E\cap \overline{B}_R (\ell) \cap \{|x_3|\leq R\})\bigr)\\
& < \mathcal{H}^2 \bigl(E \cap \overline{B}_R (\ell) \cap \{|x_3|\leq R\}\bigr)
= 2 c_0 R^2\, ,
\end{align*}
where $c_0=1$ in case (a) and $c_0=2$ in case (b). In particular we infer the existence of an $s$ such that 
\[
\mathcal{H}^1 (\Gamma^c (s)) \leq \mathcal{H}^1 (\Gamma (s)) < c_0 R\, .
\]
In both cases $\mathcal{H}^1 (\Gamma^c (s))$ (which is, respectively, strictly less than $R$ in case (a) and strictly less than $2R$ in case (b)) gives an upper bound for the length of the solution of the corresponding Steiner problem: the minimum of the length over connected sets containing, respectively, $p_0$ and $q$ in case (a) and $p_0$, $q$, and $q'$ in case (b). It is however elementary and well known to that in case (a) the minimizer is a straight segment with endpoints $p_0$ and $q$ while in case (b) it is the union of two straight segments with one common endpoint $p_0$ and the other endpoints being $q$ and $q'$. With this we reach a contradiction.
\end{proof}

\begin{lemma}\label{l:connesso}
Let \(R>0\) and let
\[
Q:=\bigl\{(x_1,x_2,x_3)\in\mathbb{R}^3:
x_2=0,\ 0\leq x_1\leq R,\ |x_3|\leq R\bigr\}.
\]
Suppose that \(\Phi\colon Q\to\mathbb{R}^3\) is continuous and agrees
with the identity in a neighborhood of \(\partial Q\) in \(Q\). For
each \(s\in(-R,R)\), set $
\Gamma(s):=\Phi(Q)\cap\{x_3=s\}$.
Then the points $a_s = (0,0,s)$ and $b_s = (R,0,s)$
belong to the same connected component of \(\Gamma(s)\).
\end{lemma}

\begin{proof} We first suppose that \(\Phi\) is smooth. Let $
\pi_3\colon\mathbb{R}^3\to\mathbb{R}$
denote the projection onto the third coordinate, and set $
f:=\pi_3\circ\Phi$.
Since \(\Phi\) agrees with the identity near \(\partial Q\), we have $
f(x_1,0,x_3)=x_3$ on a neighborhood of $\partial Q$.
By Sard's theorem, there is a sequence of regular values
\(t_k\in(-R,R)\) of \(f\) such that \(t_k\to s\). For every \(k\), the
level set $L_k:=f^{-1}(t_k)$
is a compact one-dimensional manifold. Moreover, since \(f=x_3\) near
\(\partial Q\), $\partial L_k=L_k\cap\partial Q=\{a_{t_k},b_{t_k}\}$.
Every compact one-dimensional manifold is a disjoint union of circles
and compact intervals. It follows that the connected component \(A_k\)
of \(L_k\) containing \(a_{t_k}\) is (topologically) a compact interval whose other
endpoint must be \(b_{t_k}\). 

Since \(Q\) is compact, after passing to a subsequence the sets \(A_k\)
converge in the Hausdorff topology to a nonempty compact set \(A\subset
Q\). As the Hausdorff limit of compact connected sets, $A$ is connected, and moreover
$a_s,b_s\in A$.
Furthermore, if \(x\in A\), we can choose \(x_k\in A_k\) such that
\(x_k\to x\). Consequently,
\[
f(x)=\lim_{k\to\infty}f(x_k)
    =\lim_{k\to\infty}t_k
    =s.
\]
Thus \(A\subset f^{-1}(s)\). It follows that \(\Phi(A)\) is a connected
subset of \(\Gamma(s)\) containing \(a_s\) and \(b_s\). This proves the
claim when \(\Phi\) is smooth.

We now consider a merely continuous map \(\Phi\) and approximate it uniformly with
a sequence of smooth maps $\Phi_j$ which are the identity in a fixed neighborhood of the boundary of the rectangle $Q$.

By the smooth case, for every \(j\) there is a compact connected set
\(C_j\subset Q\) such that
\[
a_s,b_s\in C_j
\qquad\text{and}\qquad
\pi_3\bigl(\Phi_j(x)\bigr)=s
\quad\text{for every }x\in C_j.
\]
After passing to a subsequence, \(C_j\) converges in the Hausdorff
topology to a compact connected set \(C\subset Q\). Since every \(C_j\)
contains \(a_s\) and \(b_s\), so does \(C\). On the other hand the
uniform convergence of \(\Phi_j\) to \(\Phi\) immediately gives
$\Phi(C)\subset\Phi(Q)\cap\{x_3=s\}=\Gamma(s)$.
The set \(\Phi(C)\) is connected and, because \(\Phi\) fixes
\(\partial Q\), contains \(a_s\) and \(b_s\). Hence these two points
belong to the same connected component of \(\Gamma(s)\).
\end{proof}

\subsection{The cube family} In case (e) we know the following.

\begin{proposition}\label{p:cubi}
There is a positive $\varepsilon>0$ with the following property.
Let $E$ be the cone over the $1$-skeleton of a cube with the eight vertices placed at the points $(i \frac{\sqrt{3}}{3},  j \frac{\sqrt{3}}{3}, k \frac{\sqrt{3}}{3})$ with $(i,j,k)\in \{-1,1\}^3$. Let $\ell$ be the line going through zero and through the point 
\[
q (\lambda) = \lambda (\tfrac{\sqrt{3}}{3}, \tfrac{\sqrt{3}}{3}, \tfrac{\sqrt{3}}{3}) + (1-\lambda) (- \tfrac{\sqrt{3}}{3}, \tfrac{\sqrt{3}}{3}, \tfrac{\sqrt{3}}{3})
\]
with $\lambda \in [\frac{1}{2},1]$. If $|\lambda - \frac{1}{2}|< \varepsilon$, then $E$ is not Almgren minimizing with respect to $\ell$.
\end{proposition}
\begin{proof}
Observe that, when $\lambda = \frac{1}{2}$, then the link of the cone $E$, together with the two points $\{q(\lambda), -q (\lambda)\}$, forms a stationary geodesic net with antipodal pins which is already listed as an element of the goggles family, cf. \cref{d:goggles}. We have already shown in the proof of \cref{t:main} that in this case $E$ is not Almgren minimizing with respect to $\ell$. But in fact our argument provides (with the assistance of a computer) a simplicial complex $\Sigma$ and a piecewise linear map $\Phi$ which satisfy the assumptions of \cref{p:stability}. But then we can apply the latter proposition to conclude the proof of the claim that $E$ is not Almgren minimizing with respect to $\ell$ when $\lambda$ is sufficiently close to $\frac{1}{2}$.
\end{proof}

\begin{figure}

\begin{tikzpicture}[
    scale=0.85,
    x={(1.085cm,-0.030cm)},
    y={(0.600cm,0.280cm)},
    z={(0cm,1.15cm)},
    line join=round,
    line cap=round
]

\def\a{2}
\def\d{0.58}


\coordinate (ppp) at ( \a, \a, \a);
\coordinate (mpp) at (-\a, \a, \a);
\coordinate (pmp) at ( \a,-\a, \a);
\coordinate (mmp) at (-\a,-\a, \a);

\coordinate (ppm) at ( \a, \a,-\a);
\coordinate (mpm) at (-\a, \a,-\a);
\coordinate (pmm) at ( \a,-\a,-\a);
\coordinate (mmm) at (-\a,-\a,-\a);


\coordinate (Cpp) at (0, \d, \d);
\coordinate (Cmp) at (0,-\d, \d);
\coordinate (Cmm) at (0,-\d,-\d);
\coordinate (Cpm) at (0, \d,-\d);


\draw[black] (Cmp) -- (Cpp) -- (ppp) -- (pmp) -- cycle;
\draw[black] (Cmp) -- (Cpp) -- (mpp) -- (mmp) -- cycle;

\draw[black] (Cmm) -- (Cpm) -- (ppm) -- (pmm) -- cycle;
\draw[black] (Cmm) -- (Cpm) -- (mpm) -- (mmm) -- cycle;

\draw[black] (Cpm) -- (Cpp) -- (ppp) -- (ppm) -- cycle;
\draw[black] (Cpm) -- (Cpp) -- (mpp) -- (mpm) -- cycle;

\draw[black] (Cmm) -- (Cmp) -- (pmp) -- (pmm) -- cycle;
\draw[black] (Cmm) -- (Cmp) -- (mmp) -- (mmm) -- cycle;


\draw[black] (mpp) -- (ppp) -- (Cpp) -- cycle;
\draw[black] (mmp) -- (pmp) -- (Cmp) -- cycle;
\draw[black] (mpm) -- (ppm) -- (Cpm) -- cycle;
\draw[black] (mmm) -- (pmm) -- (Cmm) -- cycle;


\draw[black!55,thick] (mpp)--(ppp);
\draw[black!55,thick] (mmp)--(pmp);
\draw[black!55,thick] (mpm)--(ppm);
\draw[black!55,thick] (mmm)--(pmm);

\draw[black!55,thick] (ppp)--(pmp);
\draw[black!55,thick] (mpp)--(mmp);
\draw[black!55,thick] (ppm)--(pmm);
\draw[black!55,thick] (mpm)--(mmm);

\draw[black!55,thick] (ppp)--(ppm);
\draw[black!55,thick] (pmp)--(pmm);
\draw[black!55,thick] (mpp)--(mpm);
\draw[black!55,thick] (mmp)--(mmm);


\filldraw[
    fill=gray!45,
    draw=black,
    very thick
]
(Cpp) -- (Cmp) -- (Cmm) -- (Cpm) -- cycle;


\draw[red,very thick]
    (0,-1.25*\a,-1.25*\a)
    --
    (0, 1.25*\a, 1.25*\a);

\end{tikzpicture}
\caption{Taylor's competitor for the cubical cone when $\ell$ bisects two edges of the $1$-skeleton.}\label{f:cube-competitor}
\end{figure}
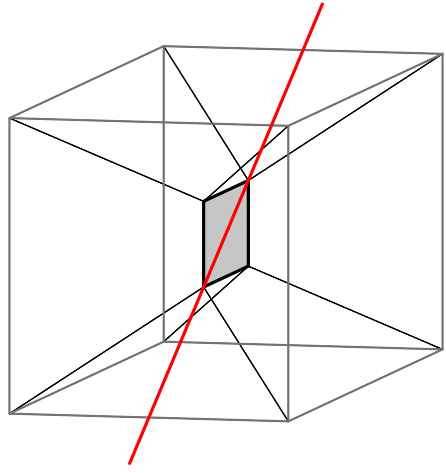

\begin{remark}
The interested reader can easily avoid the assistance of a computer, because one can use the explicit competitor already furnished by Taylor in \cite{Taylor} to exclude the Almgren minimality of $E$ when $\ell$ is replaced by the empty set. For the reader's convenience we describe it quickly, a visual illustration of the competitor is moreover given in \cref{f:cube-competitor}. We let $V (i,j,k)$ the eight vertices of the cube and introduce further the four vertices $(0, \delta, \delta)$, $(0, -\delta, \delta)$, $(0, -\delta, -\delta)$, and $(0, \delta, -\delta)$, where $\delta$ is a small positive number. 
 The competitor consists of the ``central square'' (the one which has the latter four points as vertices), four triangles, and eight trapezoids. For each $\sigma$ which is an edge of the central square there are two trapezoids which have $\sigma$ as shortest basis and one of the two parallel edges of the cube closest to $\sigma$ as longest basis. To construct the four triangles take the four remaining edges of the cube and for each add the closest vertex of the central square. An obvious triangulation of this competitor gives an embedded simplicial complex satisfying the assumptions of \cref{p:stability} (as convex set $K$ we take the convex hull of the vertices $V(i,j,k)$). Computing the area of the competitor as a function of $\delta$, one sees that the competitor has strictly less area than the cone over the skeleton of the cube if $\delta$ is positive and sufficiently small: we omit the computations and rather refer to \cite{Taylor}.
\end{remark}

Brakke's surface evolver suggests that $E$ is not area minimizing if the angle formed by $\ell$ and the plane $\{x_1=0\}$ is smaller than $0.11$ radians. The tools developed in this paper certainly allow to give an effective bound for the number $\varepsilon$ in \cref{p:cubi} and with some effort one might be able to even prove the latter statement rigorously. However we do not pursue this issue here. Note that for $\lambda=1$ the angle between $\ell$ and the same plane is $\frac{1}{2} \arccos (\frac{1}{3}) \sim 0.61$ radians. 

\appendix

\section{Allard's theory and minimizing sets}\label{a:allard}

In this section we give the proofs of \cref{p:compactness}, \cref{c:monot}, and \cref{l:lower-bound}. The cornerstone of all the arguments is the following simple lemma.

\begin{lemma}\label{l:varifolds}
Let $U\subset \mathbb R^N$ be an open set, $\Gamma \subset U$ a $C^2$ $(k-1)$-dimensional submanifold of $U$, and $E$ a minimizing set in $U$ with respect to $\Gamma$. Then the varifold induced by $E$ (which is rectifiable by assumption) is stationary in the sense of Allard in $U\setminus \Gamma$.
\end{lemma}
\begin{proof}
The stationarity in the sense of Allard is defined as follows (see e.g. \cite[Definition 1.3]{DL-All}). Having fixed a vector field $X\in C^\infty_c (U\setminus \Gamma)$ we consider the one-parameter family of diffeomorphisms generated by it, namely the map $\Phi : U\times \mathbb R \to U$ such that $\Phi (x,0)=x$ for every $x$ and $\frac{d\Phi}{dt} (x,t) = X (\Phi (x,t))$. The stationarity in $U\setminus \Gamma$ means that, for every such $X$, after choosing $U'\ssubset U$ containing its support, the first variation
\[
\delta E (X) := \left. \frac{d}{dt}\right|_{t=0} \mathcal{H}^k (U' \cap \Phi (E,t))
\]
vanishes.
Note however that, because of the minimizing property, the point $t=0$ must be a minimum of the function $f(t) = \mathcal{H}^k (U' \cap \Phi (E,t))$: indeed $\Phi (\cdot, t)$ is Lipschitz and it is the identity on $\Gamma$ and on the complement of $U'$, therefore the minimizing property must imply that $\mathcal{H}^k (U' \cap \Phi (E,t)) \geq \mathcal{H}^k (U'\cap E)$. Finally, we recall that the differentiability of the function $f$ is a consequence of the area formula.
\end{proof}

Thanks to the lemma we can apply Allard's regularity theory in our setting.
An important byproduct of Allard's theory is that, if $p\in E$ (irrespectively of whether $p\in \Gamma$ or not), then we have a lower bound on $\mathcal{H}^k (E\cap B_r (p))$ which depends only upon the curvature of $\Gamma$. This is a consequence of the monotonicity formulae and the following lower bound for the density.

\begin{lemma}\label{l:lower-bound}
Assume $U$, $\Gamma$, and $E$ are as in \cref{l:varifolds}. Let $\omega_k$ be the Lebesgue measure of the unit ball in $\mathbb R^k$. Then:
\begin{equation}
\lim_{r\downarrow 0} \frac{\mathcal{H}^k (E\cap B_r (p))}{\omega_k r^k} \geq
\left\{
\begin{array}{ll}
1 \qquad &\mbox{if $p\in E\setminus \Gamma$;}\\
{\textstyle{\frac{1}{2}}} & \mbox{if $p\in E\cap \Gamma$.}
\end{array}
\right.
\end{equation}
\end{lemma}

\begin{proof}
Given that the varifold induced by $E$ is stationary, the interior case ($p\not \in \Gamma$) can be found e.g. in \cite[Proposition 2.2(ii)]{DL-All}. The boundary case ($p\in \Gamma$) is in \cite[Theorem 3.5, conclusion (2)]{AllB}.
\end{proof}

But in fact Allard's theory is much more powerful: his papers \cite{All} and \cite{AllB} imply the existence of a geometric constant $\delta>0$ such that:
\begin{itemize}
    \item If $p\in E$, $B_{2r} (p)\subset U\setminus \Gamma$ and $\mathcal{H}^k (E\cap B_{2r} (p)) < (1+\delta) \omega_k (2r)^k$, then $E\cap B_r (p)$ is a $C^{1,\alpha}$ submanifold of $B_r (p)$ without boundary \cite[Section 8, Regularity Theorem]{All};
    \item If $p\in E\cap \Gamma$, $B_{2r} (p) \subset U$ and $\mathcal{H}^k (E\cap B_{2r} (p)) < (\frac{1}{2}+\delta) \omega_k (2r)^k$, then $E\cap B_r (p)$ is a $C^{1,\alpha}$ submanifold of $B_r (p)$ with boundary $\Gamma \cap B_r (p)$ see \cite[Section 4, Regularity Theorem]{AllB}.
\end{itemize}
Clearly, the lower bounds of \cref{l:lower-bound} can also be inferred as immediate consequences of the latter regularity results.

We can now tackle the proof of \cref{c:monot}.

\begin{proof}[Proof of \cref{c:monot}] We will use here that the varifold induced by $E$ is stationary. Therefore,
\eqref{e:monot} is a direct consequence of Allard's interior monotonicity formula, see \cite[Theorem 5.1]{All}, \cite[Equation 4.15]{Simon3}, or \cite[Theorem 2.1]{DL-All}. \eqref{e:monotb} is a direct consequence of his boundary monotonicity formula, see \cite[Theorem (2) of Section 3.4]{AllB} or \cite[Equation (31)]{Theodora}. Concerning the second part of the corollary, the case $p\not\in \Gamma$ is contained in \cite[Theorem (2) of Section 5.2]{All} and can also be found in the proof of \cite[Theorem 4.5.3]{Simon3}. The second case, namely $p\in \Gamma$ and $A_\Gamma \equiv 0$ can be inferred from the first case using the Reflection Principle of \cite[Section 3.2]{AllB}. In fact, since $A_\Gamma \equiv 0$, $\Gamma$ must be a $(k-1)$-dimensional plane. The reflection principle of \cite[Section 3.2]{AllB} implies that the varifold $W$ induced by $E$ and its reflection along $\Gamma$ is stationary in some ball $B_{r_0+\varepsilon} (p)$: for $W$ $p$ is then an interior point and we can apply the first case of the statement to conclude that the support of $W$ is a cone with vertex $p$. This easily implies that $E$ itself is a cone.
\end{proof}

We finally come to \cref{p:compactness}. Part of the proof will be based on the following lemma.

\begin{lemma}\label{l:squash}
Let $U\subseteq \mathbb R^N$ be an open set and $\Sigma\subset U$ a $k$-dimensional $C^2$ submanifold without boundary (in $U$). For every open $U'\ssubset U$ there is a number $\bar\delta (U', \Sigma)>0$ with the following property. For any $0<\delta<\bar\delta$, there is a Lipschitz map
$\Phi : U \to U$ such that:
\begin{itemize}
\item The Lipschitz constant of $\Phi$ is at most a geometric constant $C_0$;
\item $\Phi$ is the identity outside some open $V\ssubset U$;
\item $\Phi$ maps $B_\delta (\Sigma \cap U')$ in $\Sigma\cap U'$.
\end{itemize}
\end{lemma}

\begin{proof}
Using the $C^2$ regularity of $\Sigma$ and the existence of a normal neighborhood in some open set $U''\ssubset U$ which contains $U'$, the lemma is a simple exercise. We fix two parameters $1 >> \eta >> \delta > 0$. For a point $p\in B_\eta (\Sigma)\cap U'$. If the distance $d(p)$ of $p$ to $\Sigma$ is smaller than $\delta$ we define the value of a map $\Psi$ to be $\Psi (p) = q (p)$, the orthogonal projection of $p$ onto $\Sigma$. If $d(p)= \eta$ we define it as $\Psi (p)=p$. For intermediate values of the distance $d$ we interpolate and map $p$ in
\[
\Psi (p) := \frac{\eta -d (p)}{\eta-\delta} q (p) + \frac{d(p)-\delta}{\eta-\delta} p\, .
\]
We next fix a smooth cut-off function which is identically equal to $1$ on $U'$, compactly supported in $U''$ and ranging in the interval $[0,1]$. We define $\Phi$ in $U''$ as
\[
\Phi (p) = \varphi (p) \Psi (p) + (1-\varphi (p)) p\, .
\]
We will show below that $|D\Psi|$ can be bounded by a geometric constant $C_0$ choosing the ratio $\frac{\delta}{\eta}$ and the number $\eta$ sufficiently small. On the other hand we can estimate
\[
|D\Phi (p)| \leq |\nabla \varphi (p)| |\Psi (p) - p| + \varphi (p) |D\Psi (p)| + (1-\varphi (p))\, .
\]
Since $|\Psi (p) - p|\leq \eta$, it suffices to choose $\eta$ and $\frac{\delta}{\eta}$ sufficiently small to estimate the Lipschitz constant of $\Phi$.

\medskip

It remains to estimate $D\Psi$. We compute
\[
|D\Psi (p)| \leq |\nabla (p)| \frac{|p-q(p)|}{\eta-\delta} + \max \{|D q (p)|, |{\rm Id}|\}
\]
in the range where $D\Psi$ is not the identity, which is contained in an $\eta$ tubular neighborhood of some compact subset $K$ of $\Sigma$. In particular, we can assume that, if $A$ is the second fundamental form of $\Sigma$, $|A|$ is bounded over $K$ by some fixed number $M$. We then use the well-known estimates
\[
|\nabla d (p)| \leq 1 + C_1 M |p-q(p)| \leq 1 + C_1 M \eta
\]
and
\[
|\nabla q(p)| \leq |{\rm Id}| + C_1 M |p-q(p)|\leq C_1 (1+ M \eta)\, ,
\]
where $C_1$ is a geometric constant. The claimed estimate for $D\Psi$ follows easily.
\end{proof}

\begin{proof}[Proof of \cref{p:compactness}] We consider the Radon measures $\mu_j := \mathcal{H}^k \res E_j$ and, up to extraction of a subsequence, not relabeled, we assume that $\mu_j\rightharpoonup^* \mu$ for some Radon measure $\mu$ on $U$.
Using \cref{l:varifolds} and Allard's compactness theorem for stationary integral varifolds, cf. \cite[Theorem 6.4]{All}, we know that, in every open set $U'\ssubset U\setminus \Gamma$, we have a rectifiable set $R$ and an integer-valued Borel function $\Theta$ such that $\mu \res U' = \Theta \mathcal{H}^k \res R$. Clearly the same conclusion holds on the whole set $U\setminus \Gamma$ for a single integer rectifiable $R$ with (locally) finite Hausdorff measure and an integer-valued Borel function. Since $V$ is stationary in $U\setminus \Gamma$, we can assume that $R$ is the intersection of the support of $\mu$ with $U\setminus \Gamma$ (hence a relatively closed set) and that $\Theta$ is the upper semicontinuous function
\[
\Theta (p) := \lim_{r\downarrow 0} \frac{\mu (B_r (p))}{\omega_k r^k}\, .
\]
Moreover, thanks to the lower density bounds, we know that $E_j \cap U'$ converge locally in the Hausdorff topology to $R\cap U'$. Indeed:
\begin{itemize}
\item If $p\in R \cap U'$ we must have that $\mu (B_r (p)) > 0$ for evert $r>0$, in particular there must be at least one sequence of points $\{p_j\}$ with $p_j \in E_j$ and $p_j \in B_r (p)$. Otherwise there is a subsequence of $E_j$, not relabeled, with $E_j \cap B_r (p) = \emptyset$, leading to $\mu (B_r (p)) \leq \liminf_j \mathcal{H}^k (E_j \cap B_r (p)) = 0$.
\item We now take $p\in U'$ which is the limit of some (subsequence of a) sequence $\{p_j\}$ with $p_j \in E_j$ and we wish to show that $p\in {\rm spt}\, (\mu)$, namely that $\mu (B_\rho (p))>0$ for every $\rho< \dist\, (p, U')$. We distinguish two cases.

In the first case we assume $p\in U'\setminus \Gamma$. Then, for any $\rho>0$ sufficiently small so that $\dist\, (p, \Gamma) \geq 2\rho$, we must have
\[
\mu (\overline{B}_\rho (p)) \geq \limsup_{j\to \infty} \mathcal{H}^k (E\cap \overline{B}_\rho (p)) \geq \limsup_{j\to\infty} \mathcal{H}^k (E_j \cap \overline{B}_{\rho/2} (p_j))
\geq C_0^{-1} \omega_k \rho^k\, ,
\]
by \cref{l:lower-bound}. This implies $\mu (B_\rho (p)) >0$ for every $\rho$, namely that $p$ belongs to the support of $\mu$.

In the second case we assume that $p\in \Gamma \cap U'$ and fix $\rho$ such that $2\rho < \dist\, (p, \partial U)$ .Then we let $q_j\in \Gamma_j$ be the closest point to $p_j$ in $\Gamma_j$. Set $\varepsilon_j := |p_j - q_j|$, which clearly is converging to $0$. Moreover $q_j \to p$, so $B_\rho (q_j) \subset U'$ for $j$ sufficiently large. If $\varepsilon_j = 0$ then \cref{l:lower-bound} and the monotonicity formula implies immediately that
\[
\mathcal{H}^k (E_j \cap B_\rho (q_j)) \geq c_0 \rho^k\, .
\]
If $\varepsilon_j > 0$, note however that
\[
\mathcal{H}^k (E_j \cap B_{2\varepsilon_j} (q_j))
\geq \mathcal{H}^k (E_j \cap B_{\varepsilon_j} (p_j))
\geq C_0 \varepsilon_j^k
\]
because we can use the interior case of \cref{l:lower-bound} and the monotonicity formula for radii up to $\varepsilon_j$, given that $B_{\varepsilon_j} (p_j)\subset U'\setminus \Gamma$. But now we can use the boundary case of the monotonicity formula and $2\varepsilon_j \leq \rho/2$ (which holds for $j$ large enough) to bound
\begin{equation}\label{e:lower-bound-at-boundary}
\mathcal{H}^k (E_j \cap B_{\rho/2} (q_j)) \geq C_1 \rho^k\, .
\end{equation}
We now have established the validity of \eqref{e:lower-bound-at-boundary} irrespectively of whether $q_j = p_j$ or not. Since $q_j \to q$ we can now argue as in the interior case and conclude that $\mu (\overline{B}_\rho (p)) > 0$ passing into the limit in $j$ and using \eqref{e:lower-bound-at-boundary} together with $q_j \to q$ to bound
\[
\mathcal{H}^k (E_j \cap \overline{B}_\rho (p)) \geq C_1 \rho^k\, .
\]
\end{itemize}

\medskip

In the next step we show that the density $\Theta$ equals $1$ almost everywhere. It suffices to show that the density is equal to one at every point where the rectifiable set $R$ has a weak tangent plane $\pi$. Fix such a point and assume, without loss of generality, that it is the origin. Moreover define the rescaling maps $\lambda_r (x) := \frac{x}{r}$ We consider the rescaled integral varifolds with rectifiable sets given by $\lambda_r (R)$ and densities $\Theta \circ \lambda_r^{-1}$. These varifolds converge to the varifold induced by the plane $\pi$ with constant density equal to $\Theta (0)$. Again by the compactness of integral varifold we conclude that in fact $\Theta (0)$ is an integer, while the density lower bound implies that it is positive. Our aim is thus to exclude that $\Theta (0) \geq 2$. Using a simple diagonal argument, the measure given by $\nu = \Theta (0) \mathcal{H}^k \res \pi$ is in fact the weak$^*$ limit of $\nu_j := \mathcal{H}^k \res (\lambda_{r_j} (E_j))$ for an appropriately chosen sequence $r_j\downarrow 0$. Moreover, arguing as in the previous step, we can see that the sets $\lambda_{r_j} (E_j)$ must in fact be converging locally in the sense of Hausdorff to the plane $\pi$. We can now apply \cref{l:squash} and, for every fixed $\varepsilon >0$, find a Lipschitz map $\Phi$ with the following properties.
\begin{itemize}
\item $\Lip\, (\Phi)\leq C_0$;
\item $\Phi$ is the identity outside $B_1$;
\item $\Phi$ maps $B_\eta (\pi) \cap B_{1-\varepsilon}$ onto $\pi \cap B_{1-\varepsilon}$.
\end{itemize}
We now wish to estimate $\mathcal{H}^k (B_1 \cap \Phi (\lambda_{r_j} (E_j)))$. Given the properties of the map $\Phi$ we can easily come up with the upper bound
\[
\mathcal{H}^k (B_1 \cap \Phi (\lambda_{r_j} (E_j)))
\leq \omega_k (1-\varepsilon)^k + C_0^k \mathcal{H}^k ((B_1\setminus B_{1-\varepsilon}) \cap \lambda_{r_j} (E_j))
\]
and using the properties of the weak$^*$ convergence of measures we infer
\[
\limsup_{j\to \infty} \mathcal{H}^k (B_1 \cap \Phi (\lambda_{r_j} (E_j)))
\leq \omega_k (1-\varepsilon)^k + C_0^k \Theta (0) (1-(1-\varepsilon)^k)\, .
\]
On the other hand we also know that
\[
\lim_{j\to \infty} \mathcal{H}^k (B_1 \cap \lambda_{r_j} (E_j))
= \Theta (0) \omega_k\, .
\]
Given the minimizing property of $E_j$, which is invariant under dilations, we must have
\[
\Theta (0) \leq (1-\varepsilon)^k + C_0^k \Theta (0) (1-(1-\varepsilon)^k)\, .
\]
We now let $\varepsilon \downarrow 0$ to conclude that $\Theta (0) \leq 1$.

\medskip

We can at this point set $E:= {\rm spt}\, (\mu)$. In order to complete the proof of the proposition it suffices to show the following two things:
\begin{itemize}
\item[(a)] $\mu (\Gamma) =0$;
\item[(b)] $E$ is minimizing with respect to $\Gamma$.
\end{itemize}
In order to show both properties we will use the following fact, a consequence of the $C^2$ convergence of $\Gamma_j$ to $\Gamma$. There is a sequence of $C^1$ maps $\Psi_j$ and a sequence of sets $\Omega_j \uparrow U$ with the following property:
\begin{itemize}
\item $\Psi_j$ is a diffeomorphism of $\Omega_j$ with itself;
\item $\Psi_j -{\rm Id}$ and $\Psi_j^{-1} - {\rm Id}$ converge to $0$ in the $C^1$ norm;
\item $\Psi_j$ is a diffeomorphism of $\Omega_j\cap \Gamma_j$ with $\Omega_j \cap \Gamma$.
\end{itemize}
It is then easy to see that:
\begin{itemize}
\item[(i)] $\mathcal{H}^k \res \Psi_j (E_j \cap \Omega_j)$ converges weakly$^*$ to $\mu$;
\item[(ii)] for every fixed Lipschitz map $\Phi$ which is the identity on $\Gamma\cup (U\setminus U')$ for some $U'\ssubset U$ with $\mu (\partial U')=0$, we must have
\begin{equation}\label{e:almost-minimality}
\liminf_{j\to \infty} (\mathcal{H}^k (U' \cap \Phi (\Psi_j (E_j \cap \Omega_k))
\geq \mu (U')\, .
\end{equation}
\end{itemize}
Assume therefore that there is a $U'\ssubset U$ such that $\mu (U'\cap \Gamma) >0$. By varying $U'$ slightly we can assume that $\mu (\partial U')=0$. For every fixed $U''\ssubset U'$ and for every fixed $\varepsilon >0$ we can find a $\delta>0$ with the following property. For any $\eta<\delta$ there is a Lipschitz map $\Phi$ with the following properties:
\begin{itemize}
\item $\Lip\, \Phi \leq C_0$;
\item $\Phi$ is the identity on the complement of $U'$, on the complement of $B_\varepsilon (\Gamma \cap U')$ and on
$B_{\kappa^2 \delta} (\Gamma)$;
\item $\Phi$ shrinks the normal neighborhood of $\Gamma \cap U'$ by a factor $C_0 \kappa$ in $B_\delta (\Gamma\cap U') \setminus B_{\kappa^2 \delta} (\Gamma)$.
\end{itemize}
Since $\Gamma$ has dimension $k-1$, the tangent space to any $k$-dimensional set at a point $x\in B_\delta (\Gamma\cap U') \setminus B_{\kappa^2 \delta} (\Gamma)$ contains at least one vector which is shrinked by the factor $\kappa$. Moreover as long as $\kappa>0$ we can arrange that $\Phi$ satisfies the condition in \cref{d:Almgren-minimizing} to be an admissible competitor  we can thus estimate
\begin{align*}
0 \leq & \mathcal{H}^k (U'\cap \Phi (\Psi_j (E_j\cap \Omega_j))) -
\mathcal{H}^k (U' \cap \Psi_j (E_j \cap \Omega_j))\\
\leq & C_0^k \mathcal{H}^k (U'' \cap (\overline{B}_\eta (\Gamma)\setminus B_\delta (\Gamma))
\cap \Psi_j (E_j \cap \Omega_j))+
C_0^k \mathcal{H}^k ((U'\setminus U'')\cap \Psi_j (E_j\cap \Omega_j))\\
&- (1-C_0 \kappa) \mathcal{H}^k ((B_\delta (\Gamma)\setminus B_{\kappa^2 \delta} (\Gamma))\cap U''\cap \Psi_j (E_j\cap \Omega_j))\, .
\end{align*}
Next, by slightly varying the sets we can assume:
\[
\mu (\partial U'') = \mu (\partial (B_\eta (\Gamma)))= \mu (\partial (B_\delta (\Gamma)))=0\, .
\]
We now first let $\kappa\downarrow 0$ and then let $j\uparrow \infty$ to get
\[
0 \leq C_0^k \mu (U'' \cap (\overline{B}_\eta (\Gamma)\setminus \Gamma))
+ C_0^k \mu (U'\setminus U'') - \mu (\Gamma \cap U'')\, .
\]
Now, letting $\eta\downarrow 0$ and $U''\uparrow U'$, we get in the limit
\[
0 \leq - \mu (\Gamma \cap U'')\, ,
\]
which clearly is a contradiction.

\medskip

In order to prove (b) we will take advantage of Allard's regularity theory. Given that the density is $1$ a.e. on the rectifiable set $E$, we actually know that there is a set $U'\subset U$, containing $\Gamma$ with the property that $\mathcal{H}^k (E\cap (U\setminus U')) = 0$ and that $E\cap U'$ is a classical embedded minimal surface in $U'$, i.e. an analytic submanifold.

Assume now, by contradiction, that $E$ is not minimizing with respect to $\Gamma$ and therefore that there is a map $\Psi$ satisfying the condition of \cref{d:Almgren-minimizing} and with the property that
\[
\mathcal{H}^k (\Psi (E)\cap W) < \mathcal{H}^k (E\cap W)\, .
\]
Without loss of generality we can assume that $\mathcal{H}^k (E\cap \partial W)=0$. We fix an open set $U''\ssubset W\cap U'$ and use \cref{l:squash} to find a Lipschiz map $\Psi$ which maps some $\delta$-neighborhood of $E\cap U''$ onto $E$, is the identity outside $U'$, and has Lipschitz constant at most $C_0$. We can then estimate that
\begin{align*}
\mathcal{H}^k (\Psi \circ \Phi (\Phi_j (E_j \cap \Omega_j)) \cap W)
\leq C_0^k \Lip\, (\Phi) \mathcal{H}^k (E_j \cap (W\setminus U''))
+ \mathcal{H}^k (\Psi (E)\cap W)\, .
\end{align*}
Letting $j\uparrow \infty$ and using the minimizing condition we infer
\[
0 \leq C_0^k \Lip\, (\Phi) \mathcal{H}^k (E\cap (W\setminus U''))
+ \mathcal{H}^k (\Psi (E) \cap W) - \mathcal{H}^k (E\cap W)\, .
\]
On the other hand we can let $U''\uparrow U'$ and conclude that
\[
0 \leq C_0^k \Lip\, (\Phi) \mathcal{H}^k (E\cap (W\setminus U'))
+ \mathcal{H}^k (\Psi (E) \cap W) - \mathcal{H}^k (E\cap W)\, .
\]
Since $W\setminus U' \subset U\setminus U'$ and $\mathcal{H}^k (E\cap (U\setminus U'))=0$, we thus conclude
\[
0 \leq  \mathcal{H}^k (\Psi (E) \cap W) - \mathcal{H}^k (E\cap W)\, .
\]
This gives a contradiction and hence completes the proof.
\end{proof}

\section{Convex spherical open sets and positively curved arcs}\label{s:convessita}

First of all, given an open set $\Omega\subset \mathbb S^2$ we say that it is convex if and only if the open cone with link $\Omega$ and vertex the origin is a convex subset of $\mathbb R^3$. As it can be easily seen using separating planes, this is equivalent to being the intersection of a collection of open hemispheres of $\mathbb S^2$. For open sets which is not the whole sphere this is also equivalent to not containing any pair of antipodal points and always containing the length-minimizing geodesic arc connecting any two points.

The following criterion for recognizing open convex sets plays a crucial role in several places and will be proved in this section (the statement is folklore but we have not been able to find it in the literature; in the Euclidean space it is a classical result proved independently by Nakajima and Tietze, see e.g. \cite{Tietze}).

\begin{lemma}\label{l:local-to-global-convexity}
Let $\Omega\subseteq\mathbb S^2$ be a connected open domain, which is {\em weakly locally convex}, i.e.:
\begin{itemize}
\item[(WLC)] For every $p\in \partial \Omega$ there is a neighborhood $U$ such that $U\cap \Omega$ has finitely many connected components which are all convex.
\end{itemize}
Then $\Omega$ is convex. In particular $\overline{\Omega}$ is convex and, if $\Omega\neq \mathbb S^2$, then $\partial \Omega$ is a Jordan curve.
\end{lemma}

\begin{figure}
\begin{tikzpicture}
  \draw (0,0) circle (2.5);
  \draw[fill] (0,-2.5) circle (0.07);
  \node[below] at (0,-2.5) {$\gamma (0) = \eta (0)$};
  \draw[->, >={Latex[length=3mm,width=2mm]}] (2.5,0) to (2.5,0.5);
  \draw[->, >={Latex[length=3mm,width=2mm]}] (0,0) to (-0.5,-0.75);
  \draw[thick] (0,-2.5) to [out=0, in=-120] (1,-2) to [out=60, in=-90] (1.4,-0.5) to [out=90, in=-30]  (0.3,1);
 \draw[thick] (-1,-1.5) to [out=-90, in=180] (0,-2.5);
\draw[thick] (-1,-1.5) -- (2,3);
\draw[->, >={Latex[length=3mm,width=2mm]}] (-2.5,0) to (-1.5,0);
\node[above] at (-1.5,0.2) {$\dot\eta\times \ddot\eta$};
  \node[below] at (0, 2) {$H$};
\end{tikzpicture}
\caption{An example of a positively curved oriented arc $\gamma$, in a stereographic projection of the sphere. Note that the arc is not necessarily embedded (even though it is locally embedded) and not necessarily differentiable. $\eta$ parametrizes an equator with the property that, locally at $\gamma (0)$, the curve $\gamma$ is locally a graph over the equator and all ``on the left side'' of it (assuming $\gamma$ and $\eta$ are compatibly oriented). The property is only local, e.g. portions of the curve away from $\gamma (0)$ are allowed to escape the hemisphere $H$.}\label{f:def-pos-curved}
\end{figure}
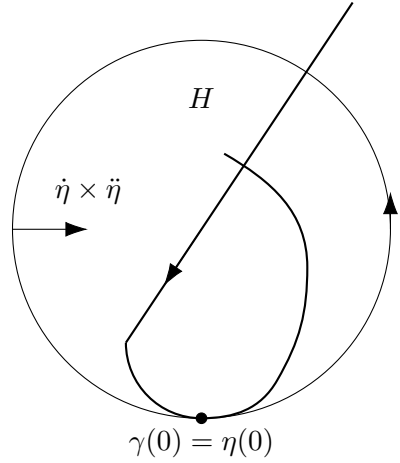

Closely connected to the notion of convexity is the following one of positively curved oriented arcs, cf. \cref{f:def-pos-curved}.

\begin{definition}\label{d:positive-geod}
Let $\gamma : [a,b] \to \mathbb S^2$ be a continuous locally injective map. We will say that $\gamma$ is a positively curved oriented arc if for every $T_0\in [a,b]$ there are an equator $\alpha\subset \mathbb S^2$ and a positive number $\varepsilon$ with the following properties:
\begin{itemize}
\item[(a)] $\alpha$ goes through $\gamma (t_0)$;
\item[(b)] Upon setting $I:= [t_0-\varepsilon, t_0+\varepsilon] \cap [a,b]$ the arc $\gamma (I)$ is contained in a closed hemisphere $H$ bounded by $\alpha$;
\item[(c)] $\gamma (I)$ is a normal graph over $\alpha$, namely the nearest point projection $\mathbf{p}:\gamma (I)\to \alpha$ is well defined and injective;
\item[(d)] If $\eta: (-\pi, \pi)$ is the arclength parametrization of $\alpha$ with $\dot \eta \times \ddot\eta$ pointing towards $H$ and $\eta (0) = \gamma (t_0)$, then $\gamma^{-1} \circ \mathbf{p}^{-1} \circ \eta$ is a strictly increasing function from the interval $\eta^{-1} (\mathbf{p} (\gamma (I))$ to $I$.
\end{itemize}
\end{definition}

Whether an oriented arc is positively curved is invariant under parametrizations of the arc which are orientation preserving. The property can be easily extended to oriented loops. The following lemma plays an important role in our arguments.

\begin{lemma}\label{l:tecnico-2}
Let $\gamma$ be a positively curved oriented arc which is embedded and whose endpoints are antipodal. Then $\gamma$ is necessarily half of an equator.
\end{lemma}

The rest of this section is devoted to the proofs of \cref{l:local-to-global-convexity} and \cref{l:tecnico-2}.

\subsection{Elementary remarks} First of all we record the following remarks, which are simple consequences of standard facts in convex geometry

\begin{proposition}\label{p:elementary-convexity}
$\phantom{6}$
\begin{itemize}
\item[(i)] A connected open set $\Omega\subsetneq \mathbb S^2$ is convex if and only it is contained in at least one open hemisphere and its boundary (oriented so that the unit normal points inwards) is a positively curved oriented loop.
\item[(ii)] If $\Omega \subsetneq \mathbb S^2$ is open and convex than its boundary is a Jordan curve.
\item[(iii)] If $\gamma$ is a positively curved oriented arc, then it has finite length and the normal graph in \cref{d:positive-geod}(c) is Lipschitz.
\item[(iv)] Consider a positively curved oriented arc $\gamma:[0, L]\to\mathbb S^2$ that is parametrized by arc-length. For any $t_0 \in [0, L]$, the following statements hold:
\begin{itemize}
\item $\gamma$ is left and right differentiable at $t_0$.
\item The properties (c) and (d) in \cref{d:positive-geod} are valid for a sufficiently small $\varepsilon$ and any equator passing through $\gamma (t_0)$ and tangent to $\dot\gamma (t_0^-)$ or $\dot\gamma (t_0^+)$.
\item When $t_0\in (0, L)$ the oriented angle between $v^-= \dot \gamma (t_0^-)$ and $v^+ = \dot\gamma (t_0^+)$ is strictly smaller than $\pi$ and the properties (c) and (d) in \cref{d:positive-geod} are valid for a sufficiently small $\varepsilon$ and any equator passing through $\gamma (t_0)$ and tangent to any convex combination of $v^-$ and $v^+$.
\end{itemize}
\item[(v)] If $A$ and $B$ are two disjoint open convex sets, then $\overline{A}\cap \overline{B}$ is either an equator, or a geodesic arc long at most $\pi$, or two antipodal points, or a singleton, or the empty set.
\end{itemize}
\end{proposition}

\begin{proof}
{\bf Proof of (i).} The statement follows immediately from a general theorem of Karcher, cf. \cite[Proposition 1]{Karcher} on open connected sets of Riemannian manifolds: the latter states that an open connected set is convex if and only if it does not intersect its own cut locus and has a support element at every boundary point. The notion of support element is equivalent, in our context, to the existence of the equator $\alpha$ as in \cref{d:positive-geod}.

\medskip

{\bf Proof of (ii).} The statement follows immediately from the fact that the open cone with link $\Omega$ and vertex the origin is a convex set: this implies that $\Omega$ is topologically a disk and that its boundary is locally Lipschitz, in particular $\partial \Omega$ is an embedded loop.

\medskip

{\bf Proof of (iii).} By compactness, without loss of generality we can assume that the arc is embedded. Consider the cone $E$ with vertex the origin and link $\gamma$. Fix any $t_0$. Then we see immediately that the latter is, in a neighborhood of $\gamma (t_0)$ the graph of a convex function $\varphi$ over the plane $\pi$ containing the equator $\alpha$. It is therefore a Lipschitz graph in a slightly smaller neighborhood and the claim follows at once.

\medskip

{\bf Proof of (iv).} Consider $\varphi$ as in (iii). Let $D$ be the subdifferential of $\varphi$ at $\gamma (t_0)$, or equivalently consider the supporting planes of the graph at the point. It follows immediately from the invariance under dilation that any supporting plane must contain the vector $\gamma (t_0)$. This constrains immediately the subdifferential $D$ to be either a point (in which case $\varphi$ is differential at $\gamma (t_0)$) or a closed segment and the claim follows easily.

\medskip

{\bf Proof of (v).} We denote by $A'$ and $B'$ the open convex cones with links $A$ and $B$ and the origin as vertex. We note that $\overline{A}\cap \overline{B} = \overline{A}'\cap \overline{B}'\cap \mathbb S^2$. Since $\overline{A}'$ and $\overline{B}'$ are convex, so is $D = \overline{A}'\cap \overline{B}'$. On other hand $D$ has also empty interior, so it is either empty, or of dimension $1$ (and hence contained in an affine line $l$) or of dimension $2$ and contained in some $2$-dimensional affine plane $\alpha$. Since it is also a cone, namely $\lambda D = D$ for every $\lambda > 0$, the plane $\alpha$ and the line $l$ must is in fact be linear. If $D$ is not empty we have therefore the following four possibilities:
\begin{itemize}
\item $D$ coincides with the plane $\alpha$ and thus $D\cap \mathbb S^2$ is an equator.
\item $D$ is $2$-dimensional, but a strict subset of the plane; then $D\cap \mathbb S^2$ is a geodesic arc of length at most $\pi$.
\item $D$ is the $1$-dimensional line $l$, and hence $D\cap \mathbb S^2$ consists of two antipodal poins.
\item $D$ is $1$-dimensional, but a strict subset of the line $l$; then $D$ is necessarily a half-line starting at the origin, hence $D\cap \mathbb S^2$ is a singleton.
\end{itemize}
\end{proof}

\subsection{Proof of \cref{l:local-to-global-convexity}} {\bf Step 1. Locally convex case.} We first prove the lemma under the following stronger assumption:
\begin{itemize}
\item[(LC)] for every $q\in \partial \Omega$ there is a neighborhood $U$ of $q$ and an open hemisphere $H$ such that $q\in \partial H$ and $U\cap \Omega \subset H$.
\end{itemize}
Fix indeed a point $q\in \partial \Omega$ and an hemisphere $H$ as in (LC). Let $\Omega'= \Omega\cap H$. Observe that the property (LC) still holds for $\Omega'$ in place of $\Omega$ and moreover $\Omega'$ does not contain any pair of antipodal points. We can therefore apply \cite[Proposition 1]{Karcher} and conclude that $\Omega'$ is convex. Let now $H'':= \mathbb S^2\setminus \overline{H}$. If $\Omega'' := H''\cap \Omega$ is empty, then we are done. Assume otherwise that $\Omega''$ is nonempty. Then $H''$ is also convex. Consider that $\alpha := \overline{\Omega}' \cap \overline{\Omega}''$ is a subset of the equator $\gamma:= \overline{H}\cap \overline{H}''$. Since $\Omega$ is connected, it must contain a point of $\gamma$ and thus $\alpha$ is necessarily one-dimensional. On the other hand $\overline{\Omega}''$ cannot contain the point $q$ because $U\cap \Omega \subset H$ for some neighborhood $U$ of $q$. From \cref{p:elementary-convexity}(v) we conclude that $\alpha$ is an arc of $\gamma$, with length at most $\pi$.  By convexity $\beta := \overline{\Omega}'\cap \gamma$ is also an arc and it must contain $\alpha$ strictly. So there is an endpoint $\bar q$ of $\alpha$ which is in the interior of $\beta$ By (LC) there must be one hemisphere $H'''$ which contains $B_r (\bar q)\cap \Omega$. On the other hand, given that $\beta$ contains $\bar q$ in the interior and $\Omega'$ is convex, we easily see that there is $\rho>0$ such that $B_\rho (\bar q) \cap \Omega' = H \cap \Omega'$. So $H'''$ must necessarily coincide with $H$. But by definition this would imply $q\not\in \overline{\Omega}''$, which is a contradiction.

\medskip

{\bf Step 2. Regularization.} In the general case we will first approximate $\Omega$ with a new connected open set which satisfies the local convexity property (LC). First of all consider the function $f:\mathbb S^2\to \mathbb R$ given by $f (p) = \dist\, (p, \partial \Omega)$, where the distance is the geodesic distance on $\mathbb S^2$. Observe that, for every $p\in \partial \Omega$ there is a neighborhood $U$ of $p$ with the property that:
\begin{itemize}
\item $f$ is concave on any connected component of $U\cap \Omega$.
\end{itemize}
In particular, by a simple compactness argument, there is a positive number $\varepsilon_0>0$ such that, for every positive $\varepsilon<\varepsilon_0$, the set $\Omega_\varepsilon:= \{f> \varepsilon\}$ satisfies (LC). We can now argue as in Step 1 and conclude that any connected component of $\Omega_\varepsilon$ is an open convex set. It is elementary to construct a sequence $\Lambda_k\uparrow \Omega$ of connected open sets with the property that $\Lambda_k \subset \Omega_{\varepsilon_k}$ for some sequence $\varepsilon_k\downarrow 0$. In particular, if $\Omega'_{\varepsilon_k}$ is the connected component containing $\Lambda_k$, we then conclude $\Omega'_{\varepsilon_k}\uparrow \Omega$. We now show that $\Omega$ is convex. For ever point $q\not \in \Omega$ we clearly have $q\not \in \Omega'_{\varepsilon_k}$. Thus, we can find an open hemisphere $H_k$ which contains $\Omega_{\varepsilon_k}$ but does not contain $q$. Up to subsequences we can assume that the hemispheres $H_k$ converge to some hemisphere $H^q$. We see immediately that $\Omega \subset H^q$ and $q\notin H^q$. The intersection of all $H^q$'s must then be $\Omega$, this showing that the latter must be convex.

\subsection{Exit lemma} The main tool to prove \cref{l:tecnico-2} is the following.

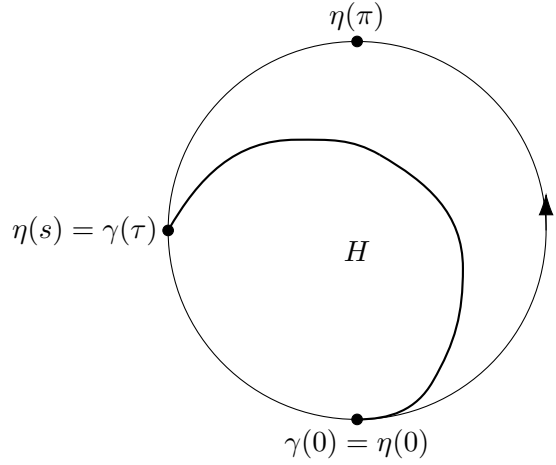
\begin{figure}
\begin{tikzpicture}
  \draw (0,0) circle (2.5);
  \draw[fill] (0,-2.5) circle (0.07);
  \node[below] at (0,-2.5) {$\gamma (0) = \eta (0)$};
  \draw[->, >={Latex[length=3mm,width=2mm]}] (2.5,0) to (2.5,0.5);
  \draw[thick] (0,-2.5) to [out=0, in=-120] (1,-2) to [out=60, in=-90] (1.4,-0.5) to [out=90, in=-30] (0.3,1);
  \draw[thick] (0.3,1) to [out=150, in=0] (-0.8,1.2) to [out=180, in=60] (-2.5,0);
  \draw[fill] (-2.5,0) circle (0.07);
  \draw[fill] (0,2.5) circle (0.07);
  \node[below] at (0,0) {$H$};
  \node[left] at (-2.5,0) {$\eta(s) = \gamma (\tau)$};
  \node[above] at (0,2.5) {$\eta (\pi)$};
\end{tikzpicture}
\caption{A picture illustrating the Exit \cref{l:exit}, using a stereographic projection of the half-sphere $H$. The boundary of $\partial H$ is oriented in such a way that the normal points outward. The ``exit point'' $\gamma (\tau)$ of the curve $\gamma$ is ``past'' the antipodal point $\eta (\pi)$ of the ``entry point'' $\eta (0)=\gamma (0)$.}\label{f:exit-lemma}
\end{figure}

\begin{lemma}[Exit lemma]\label{l:exit}
Let $H\subset \mathbb R^2$ be an open hemisphere and let $\gamma: [0, L]\to \mathbb S^2$ be an arc-length parametrization of positively curved oriented arc of $\mathbb S^2$. Assume that $\gamma (0)\in \partial H$, $\gamma ((0, \delta))\subset H$ for some positive $\delta$, but $\gamma ((0, L])$ is not contained in $H$.  Let $\eta: [0, 2\pi) \to \partial H$ be a parametrization by arclength of $\partial H$ such that $\eta (0) = \gamma (0)$ and $\dot\eta (0)\times \ddot\eta (0)$ points towards $H$. Let $\tau>0$ be the smallest positive real such that $\gamma (\tau)\in \partial H$. Then $\gamma (\tau) = \eta (s)$ for some $s\geq \pi$. Moreover, if $s=\pi$ then $\gamma ([0, \tau])$ is half of an equator passing through the antipodal points $\gamma (0)$ and $\gamma (\tau)$. See \cref{f:exit-lemma}.
\end{lemma}

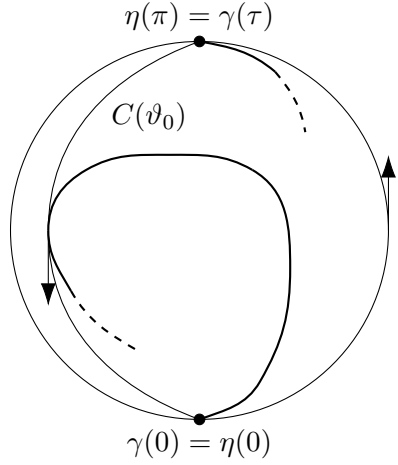
\begin{figure}
\begin{tikzpicture}
  \draw (0,0) circle (2.5);
  \draw[fill] (0,-2.5) circle (0.07);
  \node[below] at (0,-2.5) {$\gamma (0) = \eta (0)$};
  \draw[->, >={Latex[length=3mm,width=2mm]}] (2.5,0) to (2.5,1);
  \draw[->, >={Latex[length=3mm,width=2mm]}] (-2,0) to (-2,-1);
  \draw[fill] (0,2.5) circle (0.07);
  \draw (0,-2.5) to [out=160,in=-90] (-2,0) to [out=90,in=-160] (0,2.5);
  \node[right] at (-1.3,1.5) {$C (\vartheta_0)$};
  \node[above] at (0,2.5) {$\eta (\pi) = \gamma (\tau)$};
  \draw[thick] (0,-2.5) to [out=20, in=-120] (0.8,-2.0) to [out=60, in=-90] (1.2,-0.5) to [out=90, in=-30] (0.7,0.8);
  \draw[thick] (0.7,0.8) to [out=150, in=0] (-0.4,1) to [out=-180, in=30] (-1.5,0.8) to [out=-150, in=90] (-2,0);
  \draw[thick] (-2,0) to [out=-90, in=120] (-1.7,-0.8);
  \draw[thick, dashed] (-1.7,-0.8) to [out=-60, in=150] (-0.8,-1.6);
  \draw[thick] (0,2.5) to [out=-10, in=140] (1, 2.1);
  \draw[thick, dashed] ((1,2.1) to [out=-50, in=100] (1.4,1.3);
\end{tikzpicture}
\caption{A picture of $C (\vartheta_0)$ in the hemisphere $H$, under the assumption that $\gamma (\tau) = \eta (s)$ for $s\leq \pi$. The existence of the point $\gamma (t_0) = \tilde{\eta} (s_0) \in C (\vartheta_0)$ forces a self-intersection.} \label{f:exit-lemma-3}
\end{figure}

\begin{proof} Fix $\gamma$, $\eta$, $\tau$, and $s$ as in the statement. Assume moreover that $s\leq \pi$.
First of all observe that $\dot\gamma (0^+)\neq - \dot\eta (0)$ and $\dot\gamma (\tau^-)\neq -\dot\eta (s)$ (recalling \cref{p:elementary-convexity} $\gamma$ is left and right differentiable at every point). We choose coordinates in $\mathbb R^3$ such that $H = \{v\in \mathbb S^2: \langle v, e_3\rangle >0\}$ and $\gamma (0) = e_1$. In particular $\eta (0) = e_1$ while $\eta (\frac{\pi}{2}) = e_2$ For $\vartheta>0$ we let $e (\vartheta) = \cos \vartheta e_3 + \sin \vartheta e_2$ and $C (\vartheta) = \{v\in \mathbb S^2 : \langle v, e (\vartheta)\rangle =0 \;\mbox{and}\; \langle v, e_3\rangle \geq 0\}$. For small enough $\vartheta$ the intersection of $\gamma ([0, \tau])$ with $C (\vartheta)$ is contained in $\{\gamma (0), \gamma (\tau)\}$, cf. \cref{f:exit-lemma-2}.

\begin{figure}
\begin{tikzpicture}
  \draw (0,0) circle (2.5);
  \draw[fill] (0,-2.5) circle (0.07);
  \node[below] at (0,-2.5) {$\gamma (0) = \eta (0)$};
  \draw[->, >={Latex[length=3mm,width=2mm]}] (2.5,0) to (2.5,0.5);
  \draw[fill] (0,2.5) circle (0.07);
  \draw (0,-2.5) to [out=160,in=-90] (-2,0) to [out=90,in=-160] (0,2.5);
  \node[right] at (-2,0) {$C (\vartheta)$};
  \node[above] at (0,2.5) {$\eta (\pi) = \gamma (\tau)$};
  \draw[thick] (0,-2.5) to [out=20, in=-120] (0.5,-2.2);
  \draw[thick,dashed] (0.5,-2.2) to [out=60, in=-100] (0.9,-1.3);
  \draw[thick] (0,2.5) to [out=-10, in=140] (1, 2.1);
  \draw[thick, dashed] ((1,2.1) to [out=-50, in=100] (1.4,1.3);
\end{tikzpicture}
\caption{A stereographic projection of the hemisphere $H$. We are assuming that $\gamma (\tau) = \eta (s)$ for $s\leq \pi$. Given the angle formed by $\dot{\gamma} (0)$ and $\dot\eta (0)$ and, in case $s=\pi$, the one formed by $\dot\gamma (s)$ and $\dot\eta (\pi)$, for a sufficiently small $\vartheta$ the half-equator $C (\vartheta)$ does not intersect $\gamma ((0, s))$.}\label{f:exit-lemma-2}
\end{figure}
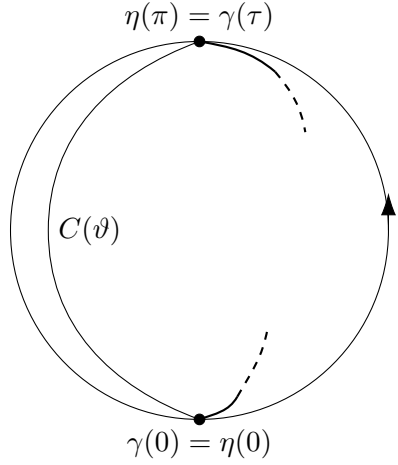

Let $\vartheta_0$ be the largest $\vartheta$ with the property that $\gamma ([0, \tau])$ is contained in the closure of the open convex set $\Omega$ bounded by $C (\vartheta) \cup \eta ([0, \pi])$. Let $\tilde{\eta}$ be a parametrization by arc-length of $\partial \Omega$ which makes it a positively curved oriented loop as in \cref{d:positive-geod} such that $\eta (0) = \tilde{\eta} (0)$ and $\eta (t) = \tilde{\eta} (t)$ for all $t\in [0, \pi]$.

Note that $\vartheta$ must be strictly smaller than $\pi$. Hence we cannot have $\dot\gamma (0^+) = \dot {\tilde{\eta}} (2\pi^-)$, given that the latter vector points towards $\mathbb S^2\setminus H$. On the other hand,  if $\dot\gamma (0^+) = - \dot{\tilde{\eta}} (2\pi^-)$, then we would in fact conclude that $\gamma ([0, \tau])$ lies in $C (\vartheta)$. Likewise, if it is the case that $s=\pi$, $\dot\gamma (\tau^-)$ must be different from $\dot{\tilde{\eta}} (\pi^+)$, while it can be $-\dot{\tilde{\eta}} (\pi^+)$ only if $\gamma ([0, \tau])$ coincides with $C (\vartheta)$. In all the remaining hypothetical cases ($s<\pi$ and $s=\pi$ but $\gamma ([0, \tau])$ is not half of an equator) we conclude that there must be at least one third point $t\in (0, \tau)$ such that $\gamma (t)\in C (\vartheta)$, namely $\gamma (t) = \tilde{\eta} (t')$ for some $t'\in (\pi, 2\pi)$. This would however force $\gamma$ to have a self-intersection, which is a contradiction, cf. \cref{f:exit-lemma-3}.
\end{proof}

\subsection{Proof of \cref{l:tecnico-2}} Fix an arclength parametrization $\gamma$ over $[0, L]$ of an arc as in the statement. Recalling \cref{p:elementary-convexity} $\gamma$ is left and right differentiable at every point and in particular it is differentiable at the two endpoints. Let $n$ be the tangent unit vector to $\mathbb S^2$ which is the counterclockwise rotation of $\dot\gamma (0)$ by $\frac{\pi}{2}$ and let $H$ be the open hemisphere $\{p\in \mathbb S^2: \langle p, n\rangle < 0\}$. If $\gamma$ is entirely contained in $\overline{H}$ the conclusion follows immediately from the Exit \cref{l:exit}. Otherwise let $[0, t_1]$ be a maximal interval with the property that $\gamma ([0, t_1])\subset \overline{H}$. Parametrize $\partial H$ with $\eta$ as in the Exit Lemma and observe that we must necessarily have $\gamma (t_1) = \eta (s_1)$ for some $s_1\in (\pi, 2\pi)$. Consider now a maximal interval $[t_1, t_2]$ with the property that $\gamma ([t_1, t_2])\subset \mathbb S^2 \setminus H$ and again it must be that $\gamma (t_2) = \eta (s_2)$ for some $s_2\in (0, 2\pi]$. We can again apply the exit lemma to the open hemisphere $\mathbb S^2\setminus \overline{H}$. Consider however that a parametrization as in \cref{d:positive-geod} for the latter is now $-\eta$. In particular we conclude that $s_2 \in (s_1, 2\pi)\cup (0, s_1-\pi]$. We now distinguish two cases.

\medskip

\begin{figure}
\begin{tikzpicture}
  \draw (0,0) circle (2.5);
  \draw[fill] (0,-2.5) circle (0.07);
  \node[below] at (0,-2.5) {$\eta (0)$};
  \draw[fill] (0,2.5) circle (0.07);
  \node[above] at (0,2.5) {$\eta (\pi)$};
  \draw[thick] (0,-2.5) to [out=0, in=-90] (1.5,-0.5) to [out=90, in= -30] (0,2) to [out=150,in=-10] (-1,{sqrt(5.25)});
  \draw[fill] (-1,{sqrt(5.25)}) circle (0.07);
  \node[above left] at (-1,{sqrt(5.25)}) {$\eta (s_1)$};
  \draw[thick, dashed] (-1,{sqrt(5.25)}) to [out=-130,in=90] (-2,0) to [out=-90,in=140] (-0.5,{-sqrt(6)});
  \node[below left] at (-0.5,{-sqrt(6)}) {$\eta (s_2)$};
  \draw[fill] (-0.5,{-sqrt(6)}) circle (0.07);
  \draw[thick] (-0.5,{-sqrt(6)}) to [out=50, in=-90] (0.5,0) to [out=90, in=-10] (-1.5,2);
  \draw[fill] (-1.5,2) circle (0.07);
  \node[left] at (-1.5,2.1) {$\eta (s_3)$};
  \draw[thick, dashed] (-1.5,2) to [out=-130, in=70] (-2.2,1);
  \draw[thick, dotted] (-2.2,1) to [out=-110,in=90] (-2.4,0);
  \draw (7,0) circle (2.5);
  \draw[fill] (7,-2.5) circle (0.07);
  \node[below] at (7,-2.5) {$\eta (0)$};
  \draw[fill] (7,2.5) circle (0.07);
  \node[above right] at (7,2.5) {$\eta (\pi)$};
  \draw[thick] (7,-2.5) to [out=0, in=-90] (8.5,-1) to [out=90, in= -30] (7,1.5) to [out=150,in=-10] (5.5,2);
  \draw[fill] (5.5,2) circle (0.07);
  \node[left] at (5.3,1.9) {$\eta (s_1)$};
  \draw[fill] (8.5,-2) circle (0.07);
  \node[right] at (8.5,-2) {$\eta (s_1-\pi)$};
  \draw[thick, dashed] (5.5,2) to [out=-120, in = 120] (5.5,-0.5) to [out=-60, in=160] (8, {-sqrt(5.25)});
  \draw[fill] (8,{-sqrt(5.25}) circle (0.07);
  \node[below right] at (8, {-sqrt(5.25)}) {$\eta (s_2)$};
  \draw[fill] (6,{sqrt(5.25)}) circle (0.07);
  \node[above] at (6,{sqrt(5.6)}) {$\eta (s_2+\pi)$};
  \draw[thick] (8,{-sqrt(5.25)}) to [out=20, in=-90] (9.2,0) to [out=90, in=-25] (8,1.8) to [out=155,in=0] (5.7,{sqrt(6.25-1.3*1.3)});
  \draw[fill] (5.7,{sqrt(6.25-1.3*1.3)}) circle (0.07);
  \node[left] at (5.7,{sqrt(5.25)}) {$\eta (s_3)$};
  \draw[thick, dashed] (5.7,{sqrt(6.25-1.3*1.3)}) to [out=-90, in=120] (6.2,0.4);
  \draw[thick, dotted] (6.2,0.4) to [out=-60, in=130] (7.2,-1);
\end{tikzpicture}
\caption{This picture illustrates the argument by contradiction that $\gamma$ is contained in the closure of the hemisphere $H$ in the proof of \cref{l:tecnico-2}. On the left we see an illustration of Case 1, while on the right we see an illustration of Case 2. The circle serves as a projection of both $\overline{H}$ and $\mathbb S^2\setminus H$. The continuous arcs are those contained in $\overline{H}$, while the dashed arcs are those contained in $\mathbb S^2\setminus H$. } \label{f:cases-1-2}
\end{figure}
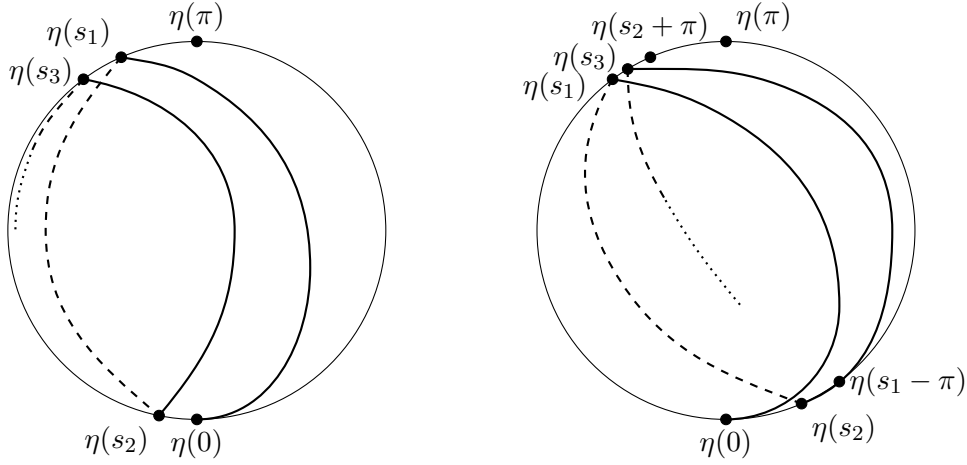

{\bf Case 1} $s_2\in (s_1, 2\pi)$. Clearly then $\gamma ([s_1, L])$ cannot be contained in $\overline{H}$ because, since $\gamma (L) = \eta (\pi)$, this would violate the injectivity of $\gamma$. Moreover, for the same reason, if we let $[t_2, t_3]$ be the maximal interval such that $\gamma ([t_2, t_3])\subset \overline{H}$, we must have $\gamma (t_3)=\eta (s_3)$ for some $s_3\in (s_1, s_2)$. But then it cannot be that $\gamma ([t_3, L])$ is entirely contained in $\overline{H}$ and arguing in the same way we would conclude the existence of a maximal interval $[t_3, t_4]$ with $\gamma ([t_3, t_4])\subset \mathbb S^2\setminus H$ and $\gamma (t_4) = \eta (s_4)$ for some $s_4\in (s_3, s_2)$. Cf. the left picture in \cref{f:cases-1-2}. This process can be iterated indefinitely leading to the conclusion that $\gamma$ cannot reach the antipodal of $\gamma (0)$.

\medskip

{\bf Case 2} $s_2\in (0, s_1-\pi]$. In this case we conclude from the Exit Lemma that, if $[t_2, t_3]$ is the the maximal interval such that $\gamma ([t_2, t_3])\subset \overline{H}$, then $\gamma (t_3)=\eta (s_3)$ for some $s_3\in [s_2+\pi, 2\pi)$. On the other hand, since $\gamma$ is injective, we must also have $s_3<s_1$ and thus $s_3\in [s_2+\pi, s_1)$. We can now keep iterating as above and prove inductively that:
\begin{itemize}
\item When we are looking at the a maximal intervals $[t_{2k}, t_{2k+1}]$ with $\gamma ([t_{2k}, t_{2k+1}])\subset \overline{H}$ and $[t_{2k+1}, t_{2k+2}]$ with $\gamma ([t_{2k+1}, t_{2k+2}]\subset \mathbb S^2 \setminus H$, we have $s_{2k}\geq s_{2k-2}>0$ and $s_{2k+1} > s_{2k-1}$ by injectivity of $\gamma$.
\item The Exit Lemma then implies that $s_{2k+1}\geq s_{2k}+\pi > s_0 + \pi = \pi$, while  $s_{2k+2} \leq s_{2k+1} - \pi \leq s_1 - \pi < 2\pi-\pi = \pi$.
\end{itemize}
Cf. the right picture in \cref{f:cases-1-2}. In particular we conclude again that $\gamma$ cannot reach the point $\eta (\pi)$, namely the antipodal of $\gamma (0)$.

\section{Proof of \texorpdfstring{\cref{p:formulae}}{p:formulae}}\label{s:app-convex-faces} The two formulae \eqref{e:quadrilateral-formula} and \eqref{e:formula-Taylor} will be proven at the end and the argument for \eqref{e:quadrilateral-formula} will show first that opposite sides of the quadrilateral are congruent.

\medskip

The statement about equiangular triangles is an immediate consequences of the formulae in \cref{thm:convex_triangles}. Next, if two adjacent sides of a convex quadrilateral with angles equal $\frac{2\pi}{3}$ are equal, it follows from the congruence between opposite dies that all four sides are congruent. Draw therefore a diagonal of the quadrilateral. Such diagonal is then also a bisector of the angles at the vertices that it connects. If we label the vertices consecutively as $ABCD$ and assume that the diagonal connects $A$ and $C$, we then have that $ABC$ is a convex triangle with angles at $A$ and $C$ congruent and equal to $\frac{\pi}{3}$. It then follows immediately from \cref{thm:convex_triangles} that the diagonal $AC$ measures $\lambda_T$.

The fact that the lengths of two adjacent sides of a convex pentagon with angles $\frac{2\pi}{3}$ determine the pentagon is a very particular case of the more general \cref{thm:congruence}.
If we know the length of two sides which are not adjacent, we can use formula \eqref{e:formula-Taylor} to determine the length of a third side (cf. the function \verb|pentagon2Constraint| in \verb|metrization_constraints.cpp| in \cite{repo}). This side must then be adjacent to one of the other two and hence we have determined the pentagon thanks to \cref{thm:congruence}.


\medskip

We now come to \eqref{e:quadrilateral-formula} and we will prove at the same time that opposite sides in the quadrilateral are congruent. Let $ABCD$ be the quadrilateral with vertices labeled clockwise. We let $E$ and $E'$ be the two antipodal points which are the intersections of the two equators containing $AD$ and $BC$ (note that the convexity excludes that this equator could be the same). We label these two intersections in such a way that the interiors of the triangles $EAB$ and $E'CD$ both are disjoint from the quadrilateral. Note that, in these triangles, the (interior) angles at $E$ and $E'$ are congruent and that the angles at $A$, $B$, $C$, and $D$ all measure $\frac{\pi}{3}$. By \cref{thm:convex_triangles} the two triangles are congruent (which shows the congruence of the sides $AB$ and $CD$) and isosceles. Now let $\gamma$ be the equator bisecting the angles in $E$ and $E'$ and let $H$ be its intersection of the equator with the side $AB$. Since $ABE$ is isosceles, the triangle $EHA$ is rectangle in $H$. Moreover, if we denote by $a$ the length of $AB$ and by $b$ the length of $AD$ we have:
\begin{align*}
|AH| &=\frac{a}{2}\\
|EA| &= \frac{\pi}{2} - \frac{b}{2}\, .
\end{align*}
Let $\alpha$ be the angle at $E$ of the triangle $EHA$. We first apply the law of sines to get
\[
\frac{\sin \tfrac{a}{2}}{\sin \alpha} = \frac{\sin (\tfrac{\pi}{2} - \tfrac{b}{2})}{\sin \tfrac{\pi}{2}}\, ,
\]
which we rewrite as
\begin{equation}\label{e:furba1}
\frac{\sin \tfrac{a}{2}}{\cos \tfrac{b}{2}} = \sin \alpha\, .
\end{equation}
We next use \cref{thm:convex_triangles} a first time to get
\[
\cos |EH| = \cos \tfrac{a}{2} = \frac{ \cos \alpha + \cos \tfrac{\pi}{3} \cos \tfrac{\pi}{2}}{\sin \tfrac{\pi}{3}\sin \tfrac{\pi}{2}}\, ,
\]
which we rewrite as
\begin{equation}\label{e:furba2}
\cos \tfrac{a}{2} = \frac{\cos \alpha}{\sin \tfrac{\pi}{3}}\, ,
\end{equation}
and a second time to get
\[
\cos |EA| = \cos (\tfrac{\pi}{2} - \tfrac{b}{2}) =
\frac{\cos \frac{\pi}{2} + \cos \tfrac{\pi}{3} \cos \alpha}{\sin \tfrac{\pi}{3} \sin \alpha}\, ,
\]
which we rewrite as
\begin{equation}\label{e:furba3}
\sin \tfrac{b}{2} = \frac{\cos \tfrac{\pi}{3} \cos \alpha}{\sin \tfrac{\pi}{3} \sin \alpha}\, .
\end{equation}
We now first divide \eqref{e:furba3} by \eqref{e:furba2} in order to get
\begin{equation}\label{e:furba4}
\frac{\sin \tfrac{b}{2}}{\cos \tfrac{a}{2}} = \frac{\cos \tfrac{\pi}{3}}{\sin \alpha}\,
\end{equation}
and then multiply the latter by \eqref{e:furba1} to get
\[
\tan \tfrac{a}{2} \tan \tfrac{b}{2} = \cos \tfrac{\pi}{3} = \tfrac{1}{2}\, .
\]

\medskip

We then come to \eqref{e:formula-Taylor}. We label the vertices of the pentagon $P=ABCDE$ counterclockwise and we set $b=|AB|$, $a=|BC|$, and $c=|DE|$. We then first prolong the edges $AB$ and $DC$ ``after'' $B$ and $C$ until they meet in a point $F$, forming a (convex) quadrilateral $Q=AFDE$. Hence we prolong $AE$ and $DF$ on both sides until they meet in the (antipodal) points $G$ (on the side of $E$ and $D$) and $G'$. Clearly, in this way we find a convex region $\Omega$ delimited by two half-equators joining the antipodal points $G$ and $G'$ and containing $Q$ (which in turn contains the initial pentagon $P$). We denote by
\begin{itemize}
\item $\gamma$ the internal angles at $G$ and $G'$ of $\Omega$;
\item $\alpha$ the internal angle at $F$ of $Q$.
\end{itemize}
Note moreover that the two triangles $GDE$ and $G'FA$ have, respectively:
\begin{itemize}
\item Internal angles at $E$, $D$, and $A$ measuring $\frac{\pi}{3}$;
\item Internal angle at $G$ and $G'$ measuring $\gamma$;
\item Internal angle at $F$ measuring $\pi-\alpha$.
\end{itemize}
Applying \cref{thm:convex_triangles} to the triangle $GDE$ we can compute
\begin{equation}\label{e:furba11}
\cos \gamma = \cos |ED| \sin^2 \frac{\pi}{3} - \cos^2 \frac{\pi}{3}
= \tfrac{3}{4} \cos c - \tfrac{1}{4}\, .
\end{equation}
If we denote by $\bar b$ the length of the side $AF$, we can do the same to the triangle $G'FA$ to derive
\begin{align}
\cos \gamma &= \cos |AF| \cos \tfrac{\pi}{3} \sin (\pi-\alpha) - \cos \tfrac{\pi}{3} \cos (\pi-\alpha)\nonumber\\
&=\tfrac{\sqrt{3}}{2} \cos \bar b \sin \alpha + \tfrac{1}{2} \cos \alpha\, .\label{e:furba12}
\end{align}
In particular we can use \eqref{e:furba11} and \eqref{e:furba12} to eliminate $\gamma$ and derive
\begin{equation}\label{e:furba21}
\tfrac{3}{4} \cos c - \tfrac{1}{4} = \tfrac{\sqrt{3}}{2} \cos \bar b \sin \alpha + \tfrac{1}{2} \cos \alpha\, .
\end{equation}
We now use \cref{thm:convex_triangles} on the triangle $BCF$, which is isosceles and has angles at $B$ and $C$ measuring $\frac{\pi}{3}$. We then can get
\begin{equation}\label{e:furba22}
\cos \alpha = \cos a \sin^2 \tfrac{\pi}{3} - \cos^2 \tfrac{\pi}{3} = \tfrac{4}{3} \cos a - \tfrac{1}{4}\, .
\end{equation}
Substituting the latter in \eqref{e:furba21} we then get
\begin{equation}\label{e:furba23}
\tfrac{3}{4} \cos c - \tfrac{3}{8} \cos a - \tfrac{1}{8}
= \tfrac{\sqrt{3}}{2} \cos \bar b \sin \alpha\, .
\end{equation}
we next wish to compute $\cos \bar b$. We note that $\bar b = |AB| + |BF| = b + |BF|$. Setting $b'=|BF|$ we can use the formula for the sum of cosines and substitute it in \eqref{e:furba23} to find
\begin{equation}\label{e:furba24}
\tfrac{3}{4} \cos c - \tfrac{3}{8} \cos a - \tfrac{1}{8} =
\tfrac{\sqrt{3}}{2} \cos b \cos b' \sin \alpha - \tfrac{\sqrt{3}}{2} \sin b \sin b' \sin \alpha\, .
\end{equation}
Applying the law of sines to the triangle $CBF$ we can first find
\[
\sin b' \sin \alpha = \sin a \sin \tfrac{\pi}{3} = \tfrac{\sqrt{3}}{2} \sin a\, ,
\]
which we can insert in \eqref{e:furba24} to get rid of $\sin b \sin \alpha$:
\begin{equation}\label{e:furba25}
\tfrac{3}{4} \cos c - \tfrac{3}{8} \cos a - \tfrac{1}{8} + \tfrac{3}{4} \sin b \sin a =
\tfrac{\sqrt{3}}{2} \cos b \cos b' \sin \alpha \, .
\end{equation}
Now we can again apply \cref{thm:convex_triangles} to deduce
\[
\cos b' \sin \alpha \sin \tfrac{\pi}{3} = \cos \tfrac{\pi}{3} + \cos \alpha \cos \tfrac{\pi}{3}\, .
\]
This can be rewritten, using \eqref{e:furba22}, as
\[
\cos b' \sin \alpha = \tfrac{1}{\sqrt{3}} (1+\cos \alpha)
= \tfrac{1}{\sqrt{3}} (\tfrac{3}{4} \cos a + \tfrac{3}{4})\, .
\]
We then complete the proof using the latter identity to get rid of the $\cos b' \sin \alpha$ term in \eqref{e:furba25} getting to
\[
\tfrac{3}{4} \cos c - \tfrac{3}{8} \cos a - \tfrac{1}{8} + \tfrac{3}{4} \sin b \sin a
= \tfrac{3}{8} \cos a \cos b + \tfrac{3}{8}\, .
\]
This is easily seen to be equivalent to \eqref{e:formula-Taylor}. 

Finally, it remains to show that all sides of a pentagon are strictly shorter than $\lambda_T$. Thanks to \eqref{e:formula-Taylor} we know
\begin{equation*}
    2\cos(DE) -\frac13 = f(AB, BC)
\end{equation*}
and thus $DE < \lambda_T=\arccos(-\tfrac13)$ is equivalent to 
\begin{equation}\label{eq:goal654}
    f(AB, BC) + 1 > 0.
\end{equation}
Thanks to the trigonometric identity
\begin{equation*}
    f(x, y) + 1 = 4\cos^2\frac x2\cos^2\frac y2\Big(1 - 2\tan\frac x2\tan\frac y2\Big) ,
\end{equation*}
we deduce that \cref{eq:goal654} is equivalent to
\begin{equation*}
    \tan \frac{AB}2\tan \frac{BC} 2 < \frac12 \,.
\end{equation*}
To prove the latter, we use what we learned for quadrilaterals. Assume, by contradiction, that $\tan \frac{AB}2\tan \frac{BC} 2 \ge \frac12$; then find $A', C'$ on $AB$, $BC$ respectively such that $\tan \frac{A'B}2\tan \frac{BC'} 2 = \frac12$. By convexity, the triangle $A'BC'$ is strictly contained in the pentagon. Moreover, recalling \cref{e:quadrilateral-formula}, we know that $A'BC'$ is exactly half of a convex spherical quadrilateral with angles equal to $\frac{2\pi}3$. We find the contradiction by looking at the areas:
\begin{equation*}
    \frac12 \frac{2\pi}3=\text{Area}(A'BC') < \text{Area}(ABCDE) = \frac{\pi}3.
\end{equation*}

\section{Proof of \texorpdfstring{\cref{thm:congruence}}{thm:congruence}}\label{s:congruence}

Observe that we can assume $N\geq 4$, because the case $N=3$ is already covered by \cref{thm:convex_triangles}. Clearly, it suffices to show that the lengths of the three remaining arcs are determined. Consider now
\begin{itemize}
\item the equator $\alpha$ containing $V_{N-2}$ and $V_{N-1}$ (which is uniquely determined because the length of $[V_{N-2}, V_{N-1}]$ is strictly shorter than $\pi$);
\item the equator $\beta$ containing $V_1$ and $V_N$.
\end{itemize}
The two equators $\alpha$ and $\beta$ must be distinct, otherwise the four points would be all on the same equator and hence the three arcs $[V_{N-2},V_{N-1}]$, $[V_{N-1},V_N]$, and $[V_N,V_1]$ would also be contained in that equator (given that they are all shorter than $\pi$). This would mean that the angles at $V_{N-1}$ and $V_N$ are equal to $0$ or to $\pi$, which is excluded by assumption (the angles at the vertices of a closed polygonal path are supposed to be distinct from $0$ and $\pi$). For the same reason it cannot be that $V_{N-1}$ belongs to $\beta$ or that $V_N$ belongs to $\alpha$. In particular:
\begin{itemize}
\item $\mathbb S^2 \setminus (\alpha\cup \beta)$ consists of four convex open connected sets;
\item one of them, which we call $\Omega$, is determined by the property that both $V_{N-1}$ and $V_{N-2}$ are in its boundary;
\item the arc $[V_{N-1},V_N]$ is contained in $\overline{\Omega}$ (because the latter is convex and the arc is length-minimizing).
\end{itemize}
Observe next that another connected component, which we call $\Omega'$, is determined by the following property:
\begin{itemize}
\item $\partial \Omega'$ contains $V_1$, $V_{N-2}$, a portion of $[V_1,V_N]$ and a portion of $[V_{N-2}, V_{N-1}]$.
\end{itemize}
Let $A$ and $B$ be the two antipodal points at the intersection of $\alpha$ and $\beta$. In particular we mark $A$ so that if we travel from $V_1$ towards $V_N$ along the arc $[V_1, V_N]$, then we get closer to $A$.
We now distinguish eight cases depending on:
\begin{itemize}
\item[(i)] whether traveling along the arc $[V_{N-2}, V_{N-1}]$ from $V_{N-2}$ to $V_{N-1}$ we move towards $A$ (we will then say that the arcs point in the same direction) or $B$;
\item[(ii)] the relative position of $\Omega$ and $\Omega'$.
\end{itemize}
First of all we see that (i) is in fact determined by the rest of the polygonal path. More precisely, consider the original situation in which we had another polygonal path $W_1, \ldots, W_N$ with the same angles and the first $N-3$ segments which are congruent and move the second path with an orientation preserving isometry $S$ so that we make the first $N-2$ vertices coincide with $O (V_i)=W_i$ for all $i\in \{1, \ldots , N-2\}$. Then we easily see that, starting from $V_1$ and traveling towards $O(W_N)$, we in fact move in the same direction as traveling towards $V_N$, while starting from $V_{N-2}$ and traveling towards $O(W_{N-1})$ we move towards the same direction as traveling towards $V_{N-1}$. This holds because the angles at $V_1$ and $V_{N-2}$ are the same as those at $W_1$ and $W_{N-2}$ and because we are assuming that $N\geq 4$.

\begin{figure}
\begin{tikzpicture}[scale=0.9]
\draw (0,0) + (-45:{2*sqrt(2)}) arc (-45:45:{2*sqrt(2)}) node[above] {$A$};
\draw (4,0) + (135:{2*sqrt(2)}) arc (135:225:{2*sqrt(2)});
\node[left] at ({4-2*sqrt(2)},0) {$\alpha$};
\node[right] at ({2*sqrt(2)},0) {$\beta$};
\draw[very thick] (0,0) + (-30:{2*sqrt(2)}) node[below right] {$V_1$} arc (-30:15:{2*sqrt(2)}) node[above right] {$V_N$};
\draw[very thick] (4,0) + (150:{2*sqrt(2)}) node[above left] {$V_{N-1}$} arc (150:200:{2*sqrt(2)}) node[below left] {$V_{N-2}$};
\draw[very thick] (0,0) + (15:{2*sqrt(2)}) -- ($(4,0) + (150:{2*sqrt(2)})$);

\draw (4,0) + (-45:{2*sqrt(2)}) arc (-45:45:{2*sqrt(2)});
\draw (8,0) + (135:{2*sqrt(2)}) arc (135:225:{2*sqrt(2)});
\node[below] at (6,1.8) {$A$};
\node[left] at ({8-2*sqrt(2)},0) {$\alpha$};
\node[right] at ({4+2*sqrt(2)},0) {$\beta$};
\draw (4,4) + (-45:{2*sqrt(2)}) arc (-45:0:{2*sqrt(2)});
\draw (8,4) + (180:{2*sqrt(2)}) arc (180:225:{2*sqrt(2)});
\draw[very thick] (4,0) + (-30:{2*sqrt(2)}) node[below right] {$V_1$} arc (-30:45:{2*sqrt(2)});
\draw[very thick] (8,4) + (200:{2*sqrt(2)}) node[left] {$V_N$} arc (200:225:{2*sqrt(2)});
\draw[very thick] (8,0) + (135:{2*sqrt(2)}) arc (135:200:{2*sqrt(2)}) node[below left] {$V_{N-2}$};
\draw[very thick] (4,4) + (-45:{2*sqrt(2)}) arc (-45:-10:{2*sqrt(2)}) node[right] {$V_{N-1}$};
\draw[very thick] (4,4) + (-10:{2*sqrt(2)}) -- ($(8,4) + (200:{2*sqrt(2)})$);

\draw (8,0) + (-45:{2*sqrt(2)}) arc (-45:45:{2*sqrt(2)});
\draw (12,0) + (135:{2*sqrt(2)}) arc (135:225:{2*sqrt(2)});
\node[below] at (10,1.8) {$A$};
\node[left] at ({12-2*sqrt(2)},0) {$\alpha$};
\node[right] at ({8+2*sqrt(2)},0) {$\beta$};
\draw[very thick] (8,0) + (-30:{2*sqrt(2)}) node[below right] {$V_1$} arc (-30:15:{2*sqrt(2)}) node[below left] {$V_N$};
\draw (8,4) + (-45:{2*sqrt(2)}) arc (-45:0:{2*sqrt(2)});
\draw[very thick] (12,0) + (135:{2*sqrt(2)}) arc (135:200:{2*sqrt(2)}) node[below left] {$V_{N-2}$};
\draw[very thick] (8,4) + (-45:{2*sqrt(2)}) arc (-45:-10:{2*sqrt(2)}) node[above right] {$V_{N-1}$};
\draw[very thick] (8,4) + (-10:{2*sqrt(2)}) to [bend left] ($(8,0) + (15:{2*sqrt(2)})$);

\draw (12,0) + (-45:{2*sqrt(2)}) arc (-45:45:{2*sqrt(2)});
\draw (16,0) + (135:{2*sqrt(2)}) arc (135:225:{2*sqrt(2)});
\node[above right] at (14,2) {$A$};
\node[left] at ({16-2*sqrt(2)},0) {$\alpha$};
\node[right] at ({12+2*sqrt(2)},0) {$\beta$};
\draw (16,4) + (180:{2*sqrt(2)}) arc (180:225:{2*sqrt(2)});
\draw[very thick] (16,0) + (160:{2*sqrt(2)}) node[below right] {$V_{N-1}$} arc (160:200:{2*sqrt(2)}) node[below left] {$V_{N-2}$};
\draw[very thick] (12,0) + (-30:{2*sqrt(2)}) node[below right] {$V_1$} arc (-30:45:{2*sqrt(2)});
\draw[very thick] (16,4) + (200:{2*sqrt(2)}) node[above right] {$V_N$} arc (200:225:{2*sqrt(2)});
\draw[very thick] (16,4) + (200:{2*sqrt(2)}) to [bend right] ($(16,0) + (160:{2*sqrt(2)})$);
\end{tikzpicture}
\caption{Given the relative orientation of the arcs $[V_{N-2}, V_{N-1}]$ and $[V_1, V_N]$ the relative position of $\Omega$ and $\Omega$' is determined by the concavity of the angles in $V_{N-2}$ and $V_{N-1}$. The picture illustrates the situation when the arcs $[V_{N-2}, V_{N-1}]$ and $[V_1, V_N]$ point in the same direction. From the far left to the right we see: $\Omega=\Omega'$, when both angles are convex; $\Omega\cap \Omega'=
\{A,B\}$ when both angles are concave; $\Omega \cap \Omega' \supset [V_1, V_N]$ when $\hat{V}_{N-1}$ is convex and $\hat{V}_N$ is concave; $\Omega \cap \Omega' \supset [V_{N-2},V_{N-1}]$ when $\hat{V}_{N-1}$ is concave and $\hat{V}_N$ is convex.}\label{f:congruence-2}
\end{figure}
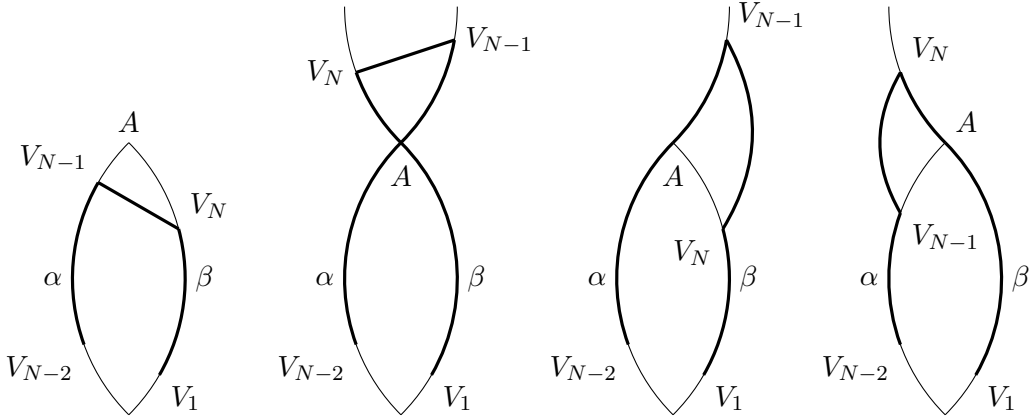

The two cases of (i) are entirely analogous and here we just focus on one of them, namely when the arcs $[V_{N-2}, V_{N-1}]$ and $[V_1,V_N]$ point in the same direction.
We then observe that the relative position of $\Omega$ and $\Omega'$ is completely determined by the concavity/convexity of the two angles at the vertices $V_{N-1}$ and $V_{N-2}$. This is illustrated in \cref{f:congruence-2}.

At this point, having established in which of the four situations we are, we consider the convex triangle with vertices $A$, $V_{N-1}$, and $V_N$. We immediately see that the angles of this triangle are completely determined by the angles of the polygonal path at $V_{N-1}$ and $V_N$. In particular its sides are also completely determined. This in turn determines the length of the arc $[V_{N-1}, V_N]$ and the relative position with respect to $A$ of the points $V_{N-1}$ and $V_N$ on the equators $\alpha$ and $\beta$. In particular we have completely determined the path and hence proved the claim of the theorem.

\section{Recognizing isomorphic planar maps}\label{s:isomorphisms}

Let us present an algorithm which detects whether two combinatorial maps are isomorphic (cf. \cite{LZ,NC}). The algorithm will be used in order to avoid ``multiple counts'' in both the list of twice-derived and once-derived candidates.

\medskip

We first introduce the darts' graph of a graph.

\begin{definition}[Darts' graph]
    Let $(G, \next)$ be a combinatorial map. Its associated Darts' graph is the oriented graph with $\darts(G)$ as vertices and such that there is a directed edge from $d_1$ to $d_2$ if either $\rev(d_1)=d_2$ or $\next(d_1)=d_2$.
    The edges of the Darts' graph keep track of whether they are $\rev$-edges or $\next$-edges.
\end{definition}

\begin{remark}[Isomorphism test]\label{r:iso-test}
    We are given two connected combinatorial maps $(G, \next)$ and $(G', \next')$ with at least one edge.
    Let us describe an algorithm, that runs in polynomial time, which determines whether they are isomorphic.
    Let us fix $d_0\in\darts(G)$ and let us search for an isomorphism $f:\darts(G)\to\darts(G')$ such that $f(d_0)=d'_0$, where $d'_0\in\darts(G')$ is an arbitrarily chosen dart. We will repeat this process for each $d'_0\in\darts(G')$.

    Now, run a depth-first search on the darts' graph associated to $G$, starting from $d_0$, which gives priority to $\rev$-edges over $\next$-edges. Let $d_0, d_1, \dots, d_{k-1}$ be the elements of $\darts(G)$ in the order they were visited by this search.
    Do the same for $\darts(G')$ starting from $d'_0$, obtaining $d'_0, d'_1, \dots, d'_{k-1}$.
    Then, the isomorphism we are looking for shall satisfy $f(d_i)=d'_i$ for all $i=0, \dots, k-1$. In particular, we have identified the candidate isomorphism and we need only to check whether it works.
\end{remark}

Note that the isomorphism test works for embeddings in $\Sf$, as the isomorphism of the embedding is equivalent to the isomorphism of the associated combinatorial maps, but it does not work for planar maps because of the special unbounded face.
To make up for this issue, it is sufficient to ensure that the unbounded face of $(G, \next)$ is mapped into the unbounded face $(G', \next)$.

\medskip

For our purposes (see \cref{s:hl-top-algorithm}), we need to perform the following operation: given a family of planar maps, we want to select a subset containing exactly one representative of each isomorphism class. If we do this na\"ively, the time complexity would be quadratic in the size of the family.

To speed up the process, we come up with a \emph{signature} for combinatorial maps and for planar maps. This is a $\mathbb Z$ invariant, namely a function which to each planar map $X$ associate an integer, with the property that it takes the same value on any other planar map $X'$ isomorphic to $X$. Having a signature allows us to solve the above-mentioned problem as follows: first we split the family in buckets indexed by the value of the signature, then we perform the isomorphism test on each pair in the same bucket.

We next describe the signature functions we employ for combinatorial maps.
We choose two arbitrary integers $p_{\rev}$ and $p_{\next}$.
For a combinatorial map $X=(G, \next)$, let $M_X$ be the matrix\footnote{This is the adjacency matrix of the darts' graph, up to the choice of weights.} such that
\begin{equation}\label{e:adjac-matrix}
M_X (d, d') = \left\{
\begin{array}{ll}
p_{\rev}+p_{\next} \quad & \mbox{if $d'= \next(d)=\rev(d)$}\\
p_{\rev} & \mbox{if $d'=\rev(d)$}\\
p_{\next} & \mbox{if $d'=\next(d)$}\\
0 &\mbox{otherwise}
\end{array}
\right.
\end{equation}
Assume that $X$ is isomorphic to $\bar X$ via the isomorphism $\sigma:\darts(X)\to\darts(\bar X)$. Then $M_X (d, d') = M_{\bar X} (\sigma (d), \sigma (d'))$. It is then immediate to check that the following function $\textrm{Sign}$, taking values in $\mathbb Z$, is a signature:
\begin{equation}\label{e:signature}
{\rm Sign}\, (X) := \sum_{d\in \darts(X)} M_X^{256}(d, \rev(d)) \oplus M_X^{256}(d, \next(d)) \,
\end{equation}
where $M^k_X$ denotes the $k$-th power of the matrix $M_X$, and $\oplus$ denotes the bitwise exclusive OR\footnote{Any commutative operation that is sufficiently ``distinct'' from sum and product of integers would work. The choice of bitwise exclusive OR is standard.}.

Let us now present the signature for planar maps. Observe that the one given above would be a valid signature, but it would not distinguish two planar maps which are equal as combinatorial maps but have different unbounded faces.
The signature for planar maps is the sum of the signature in \eqref{e:signature} and
\begin{equation*}
\sum_{d\in \text{unbounded}(X)} \bigoplus_{e\in\darts(X)} M_X^{256}(d, e) \, .
\end{equation*}
We leave to the reader to check that two isomorphic planar maps have the same signature.

\medskip

The algorithms described above are used verbatim to avoid multiple count while producing the list of \cref{p:first-topological-list}. It must be slightly modified in order to account for the special features of once-derived marked spherical graphs of \cref{p:second-topological-classification}: we just need to add the requirement that, when the algorithm checks whether a candidate is in fact an isomorphism, it also checks the following additional condition: if a degree two vertex $A$ is mapped into $A'$, then the special edge of \cref{d:candidate} adjacent to $A$ must be mapped into the special edge adjacent to $A'$.

\begin{remark}\label{r:mirror-symmetry-2}
Finally, even if we avoid isomorphic planar combinatorial maps we are, for our purposes, still double-counting the spherical graphs in question, because we do not take into account mirror symmetry, cf. \cref{r:mirror-symmetry}. We in fact ignore this issue in the list of \cref{p:first-topological-list} because it is not as computationally costly as ignoring {\em all} isomorphisms: keeping a copy of the mirror symmetric in each list just doubles the cardinality of the $\mathcal{C}_k$ of \cref{s:genera-lista}. Eventually this redundancy is instead accounted for when listing the {\em once-derived} candidate topologies in \cref{p:second-topological-classification}: the list produced by the latter proposition does not contain mirror symmetric copies.
\end{remark}

\section{Computations for the proofs of Propositions \texorpdfstring{\ref{p:pin-touching}}{pin-touching}, \texorpdfstring{\ref{prop:effective-short-edge-1}}{short-edge-1}, and
\texorpdfstring{\ref{prop:effective-short-edge-2}}{short-edge-2}}

\subsection{Computations for \texorpdfstring{\cref{p:pin-touching}}{pin-touching}}\label{a:homotopy-computations} The aim of this section is to complete the proof of \cref{p:pin-touching} giving the argument for the claim (a).

First of all we recall the function
\[
f (\delta, \alpha) := \mathcal{H}^2 (\Sigma_0) - \mathcal{H}^2 (\Sigma)\, ,
\]
which we have defined when the parameters $(\delta,\alpha)$ belong to $(0,0.011]\times [\frac{\pi}{4}, \frac{\pi}{2}]$. However, $\mathcal{H}^2 (\Sigma_0)$ is the sum of the areas of two triangles and $\mathcal{H}^2 (\Sigma)$ is the sum of the areas of 10 triangles. While $D'$ is defined as the Torricelli point of $A$, $Q$, and $B$, for the other vertices we know the coordinates as explicit smooth functions of $\delta$ and $\alpha$ which extend smoothly to $\delta =0$. We will moreover show below that the Torricelli point is a Lipschitz function of the three other points. We thus have a Lipschitz extension of $f$ to $\mathcal{D} = [0,0.011]\times [\frac{\pi}{4}, \frac{\pi}{2}]$.

We next give an upper bound for
\begin{align*}
L_\alpha & := \max_{\mathcal{D}} \left|\tfrac{\partial f}{\partial \alpha}\right|\\
L_\delta & := \max_{\mathcal{D}} \left|\tfrac{\partial f}{\partial \delta}\right|\, .
\end{align*}
We first observe that, given two triangles $XYZ$ and $XYZ'$ we have the obvious estimate
\begin{equation}\label{e:triangle-Lipschitz}
|\mathcal{H}^2 (XYZ)-\mathcal{H}^2 (XYZ')|\leq \frac{1}{2} |(X-Y) \times (Z-Z')|
\leq \frac{1}{2} |XY||Z-Z'|\, .
\end{equation}
Hence the area of a triangle is Lipschitz in each vertex, with a Lipschitz constant bounded by the diameter of the triangle divided by $2$.

In the case of $\alpha$, only the areas of the 10 triangles forming $\Sigma$ depend on it. Moreover only the points $A$, $A'$, $M$, $M'$, and $D'$ depend on $\alpha$, while the others are fixed.
The Lipschitz constants of $A$, $A'$, and $M'$ as functions of $\alpha$ can be bounded by $0.2$, $0.2\cdot 0.8$, $\delta \leq 0.011$, and $0.8 \delta \leq 0.008$ respectively. The point $D'$ is $rD$, where $D$ is the Torricelli point of $A$, $Q$, and $B$, where only $A$ depends on $\alpha$. We will show below that the Torricelli point has a Lipschitz dependence on each of the three vertices, with Lipschitz constant bounded by $\frac{2}{\sqrt{3}}$, cf. \cref{l:Torricelli} below. Therefore the Lipschitz constant of $D'$ with respect to $\alpha$ is bounded by $0.8 \cdot 0.2 \cdot 2/\sqrt{3}$. The diamaters of all 10 triangles are bounded by $\sqrt{1+0.2^2}$. Therefore it is easy to see that the Lipschitz bound for the dependence in $\alpha$ is
\begin{equation}\label{e:Lip_alpha}
L_\alpha \leq \tfrac{\sqrt{1+0.2^2}}{2} \left(0.008 \cdot 2 + 0.011\cdot 6 + 0.16 \cdot 2 + 0.2\cdot 3 + \tfrac{2}{\sqrt{3}}\cdot 0.16 \cdot 8\right) \leq 1.3\, .
\end{equation}
Now, consider that only the points $M$ and $M'$ depend on $\delta$. So, of the $12$ triangles which enter in the computations of $f$, seven have vertices which depend on $\delta$ and are $AMO$, $CMO$, $M'OD'$, $BCM$, $AMD'$, $BMD'$, and $MM'D'$. In the first six only one vertex, $M$, or $M'$, depends $\delta$, with Lipschitz constant respectively $1$ and $r=0.8$. Moreover their diameters can all be bounded by $\sqrt{1+0.2^2}$. For the final triangle two vertices depend on $\delta$, while again the diameter can be bounded by the same number. Therefore we have the bound
\begin{equation}\label{e:Lip_delta}
L_\delta \leq \frac{6 +0.8\cdot 2}{2} \cdot \sqrt{1+0.2^2} \leq 4\, .
\end{equation}
In order to prove that the function $f$ is positive we will then run a machine calculation of the values of $f$ on the centers of the rectangles of a grid with size $s_\alpha = 0.001$ in $\alpha$ and size $s_\delta = 0.0001$ in $\delta$. Given the Lipschitz bound it suffices to show that the minimum value is sufficiently (accounting for precision in floating point computations) over
the threshold
\[
\frac{L_\delta\cdot s_\delta}{2} + \frac{L_\alpha \cdot s_\alpha}{2} = 0.00085\, .
\]
The code, which can be found in the file \verb|pin-lower-bound.cpp| in \cite{repo}, computes the coordinates of the various points for $\delta$ and $\alpha$ varying in the centers of the grid and then uses the usual vector algebra formula for evaluating the areas of the various triangle in the expression of $f$. The code returns the minimum value of $f$ found in the search.

All the coordinates have been explicitly given except for the ones of the Torricelli point $D$. For the latter, one possibility would be to use an efficient general routine like the Weiszfeld algorithm. In this case we would need to certify error bounds. But given the structure of our problem we actually have explicit trigonometric expressions.
In fact elementary geometric considerations show that $D$ lies on the bisector of the angle $AQB$. Therefore its coordinates take the form $(|DQ| \cos \frac{\alpha}{2}, |DQ| \sin \tfrac{\alpha}{2}, 1)$. Moreover, again by elementary geometric considerations, $|DQ| = L (\cos \tfrac{\alpha}{2} - \sin \frac{\alpha}{2} \tan \frac{\pi}{6})$. After elementary manipulations we then get
\begin{equation}\label{e:D(alpha)}
    D(\alpha) = \left( L \left[ \tfrac{1}{2} + \tfrac{1}{\sqrt{3}} \sin\left(\tfrac{\pi}{3} - \alpha\right) \right], \, L \left[ -\tfrac{1}{2\sqrt{3}} + \tfrac{1}{\sqrt{3}} \cos\left(\tfrac{\pi}{3} - \alpha\right) \right], \, 1 \right)\, .
\end{equation}
The maximum value found in the grid search is $0.00095565$. The difference with the threshold is $0.00010565> 10^{-4}$, well above the accuracy of the floating point computations. In fact the maximum is found at (an approximation of) $\delta=0.011$ and $\alpha =\frac{\pi}{2}$, which should be the maximum point of the function on the domain $\mathcal{D}$, a claim which can be likely proved with careful analytical considerations.

\medskip

We conclude the argument with a geometric proof that the Torricelli point has the claimed Lipschitz dependence. Note that in our specific case the Lipschitz bound in $\alpha$ could also be proved directly by simply differentiating the explicit formula \eqref{e:D(alpha)}.

\begin{lemma}\label{l:Torricelli}
Fix two points $P$ and $R$ on the Euclidean plane, and vary a third point $S$. Consider the Torricelli point $T$ of the triangle $PRS$. Then $\textrm{Lip}\, (T) \leq \frac{2}{\sqrt{3}}$.
\end{lemma}
\begin{proof} When one of the angles of the triangle is $\frac{2\pi}{3}$, the Torricelli point is the obtuse angle of the triangle and so $T(S)$ is constant. Otherwise, the Torricelli point forms a triangle of degree $\frac{2\pi}{3}$ at $T$ with the vertices $P$ and $R$, so it lies on a circle $\beta$ which contains the points $P$ and $R$ and is centered at a point $S'$ which forms an equilateral triangle with $P$ and $R$, which is opposite to $S$. And it is also well known from the Fermat-Torricelli construction that $T$ lies on the segment connecting $S$ with $S'$. Hence $T(S)$ is the unique intersection of the circle $\beta$ with the segment $SS'$ and the Lipschitz constant is bounded by an elemetary computation using this fact.
\end{proof}
\subsection{Computations for \texorpdfstring{\cref{prop:effective-short-edge-1}}{short-edge-1}} In this section we complete the proof of \cref{prop:effective-short-edge-1} giving the arguments for the claims (a) and (b).

For claim (a) consider first the equators containing the edges $\alpha$ and $\gamma$. They subdivide the sphere in four convex open sets and there is exactly one of them, say $\Lambda$, which contains the convex triangles $A2DA$ and $BCB_2$. The geodesic segment $A_2B_2$ is contained in $\Lambda$ and for this reason the quadrilateral $AA_2B_2B$ has internal angles at $A_2$ and $B_2$ which are strictly smaller than $\pi$. Given that its internal angles at $A$ and $B$ are $\frac{2\pi}{3}$, we conclude that $AA_2B_2B$ is convex. Since $D$ and $C$ are internal to this quadrilateral, the convex triangles $A_2DA$ and $B_2CB$ are obviously conained in $AA_2B_2B$. Observe also that the length-minimizing geodesic segment $CD$ connecting $C$ and $D$ is contained in $AA_2B_2B$. The four geodesic arcs $AB$, $BC$, $CD$, and $DA$ bound a quadrilateral $ABCD$ contained in $AA_2B_2B$. The distance between an interior point of $CD$ and the edge $AB$ is strictly larger than that of $C$ and $D$ to $AB$. Since the points of the arc $\eta$ are equidistant from $AB$, the arc is contained in the quadrilateral $ABCD$. Therefore the region $Q$ is contained in the quadrilateral $AA_2B_2B$. This completes the proof of claim (a).

It remains to show (b). We set $x= \text{len} (\beta)$ and recall that $|AA_2|=|BB_2| = 10x$. Observe moreover that $\text{len} (\eta) \leq x$. Using the fomulae in \cref{thm:convex_triangles} (it is convenient to apply them to the right-angled triangle which is half of $AA_2D$) we get the identity
\begin{equation}\label{e:identita-tan}
\tan (5x) = \tfrac{\sqrt{3}}{2} \tan \tfrac{|AD|}{2}\, .
\end{equation}
In particular $5x < \frac{|AD}{2}$. Since $\frac{\tan x}{x}$ is monotone increasing, we then get
\[
\frac{\tan 5x}{5x} < \frac{\tan \tfrac{|AD|}{2}}{\tfrac{|AD|}{2}}
\]
and dividing the two sides of the latter inequality by the two sides of \eqref{e:identita-tan} we immediately conclude $y < \tfrac{10x}{\sqrt{3}}$. We can then bound
\begin{equation}\label{e:lineare-sotto}
|AA2| + |BB2| - (|A_2D|+|DA|+ \text{len} (\eta) + |CB2|) > (19 - 10\sqrt{3}) x\, .
\end{equation}
We next observe that the region $Q$ is foliated by a one parameter family of arcs $\{\eta_s\}_{s\in [0, |AD|]}$ which are at distance $s$ from $AB$. Since each of these arcs is shorter than $AB$, the coarea formula immediately gives
\begin{equation}\label{e:area-di-Q}
\mathcal{H}^2 (Q) \leq x |AD|
\leq \frac{10}{\sqrt{3}} x^2\, .
\end{equation}

We next wish to bound the area of the triangle $AA_2D$ (which is the same as the area of the triangle $BB_2C$). We split the triangle $AA_2D$ into two equal right angled triangles by drawing the altitude from $D$ to the edge $AB$. Take the one which has $A$ among its vertices and denote it by $T$. $T$ has internal angles $\frac{\pi}{2}$ (at the foot of the altitude) and $\frac{\pi}{6}$ at the vertex $A$. If we denote by $\theta$ the angle at the vertex $D$, then Gauss-Bonnet gives the formula
\[
\mathcal{H}^2 (T) = \theta - \frac{\pi}{3}\, .
\]
On the other hand
\[
\cos \theta = \cos \tfrac{|AA_2|}{2} \sin \tfrac{\pi}{6}\, .
\]
We now claim that
\begin{equation}\label{e:area-di-T}
\mathcal{H}^2 (T) \leq \tfrac{|AD|^2}{8} \leq \tfrac{25 x^2}{2}\, .
\end{equation}
This is implied equivalent to $\theta \leq \frac{\pi}{3} + \frac{100 x^2}{8}$, which by the monotonicity of the cosine, is equivalent to
\[
\cos \left(\tfrac{\pi}{3} + \tfrac{|AD|^2}{ 8}\right) \leq \cos \tfrac{|AA_2|}{2} \cos \tfrac{\pi}{6}\, .
\]
From the cosine duplication formula, the latter is in fact equivalent to
\[
\tfrac{1}{2} \cos \tfrac{|AD|^2}{8} \leq \tfrac{1}{2} \cos \tfrac{|AD|}{2} + \tfrac{\sqrt{3}{2}} \sin \tfrac{|AD|^2}{8}
\]
Using the known inequalities $1-\frac{t^2}{2} \leq \cos t \leq 1$ and $\sin t \geq t - \frac{t^3}{3!}$, we immediately see that, since $|AD| = 10 \text{len} (\beta) \leq 0.1$, the inequality above holds.

Using \eqref{e:lineare-sotto}, \eqref{e:area-di-Q}, and \eqref{e:area-di-T}, we see that \eqref{e:disuguaglianza-area-lunghezza} is implied by
\[
(19-10 \sqrt{3}) x \geq \left(\tfrac{20}{\sqrt{3}} + 100\right) x^2\, ,
\]
namely by $x \leq (19-10\sqrt{3})(100 + \frac{20}{\sqrt{3}})$. Since $x\leq 0.01$, the latter is equivalent to $0.01 \leq (19-10\sqrt{3})(100 + \frac{20}{\sqrt{3}})$, which can be easily verified by a calculator.

\subsection{Computations for \texorpdfstring{\cref{prop:effective-short-edge-2}}{short-edge-2}}
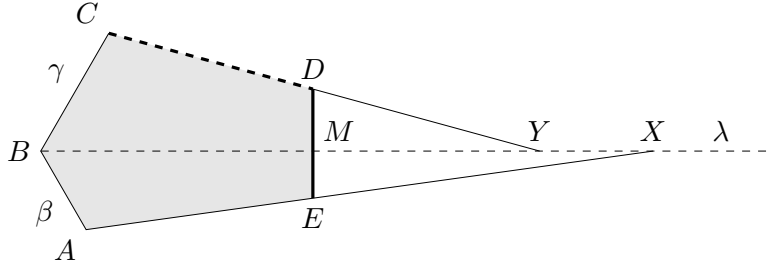
\begin{figure}[htbp]
        \begin{tikzpicture}[scale=3]
            \coordinate (A) at (0,0);
            \coordinate (B) at (-0.4/2, {0.4*sqrt(3)/2});
            \coordinate (C) at (0.1, {sqrt(3)/2});
            \coordinate (M) at (1, {0.4*sqrt(3)/2});
            \coordinate (D) at (1, {0.4*sqrt(3)/2 + 0.6*sqrt(3)/(2*1.9)});
            \coordinate (E) at (1, {0.4*sqrt(3)/2 * (1-1.5/2.5)});
            \node[below left] at (A) {$A$};
            \node[above left] at (C) {$C$};
            \node[left] at (B) {$B$};
            \node[above right] at (M) {$M$};
            \node[above] at (2.8, {0.4*sqrt(3)/2}) {$\lambda$};
            \fill[gray!20] (A) -- (B) -- (C) -- (D) -- (E) -- cycle;
            \draw[very thin] (2.5, {0.4*sqrt(3)/2}) -- (A) --node[below left] {$\beta$} (B) --node[above left] {$\gamma$} (C);
            \draw[very thin] (D) -- (2, {0.4*sqrt(3)/2});
            \draw[very thick, dashed] (C) -- (D);
            \draw[very thick] (D) -- (E);
            \node[below] at (E) {$E$};
            \node[above] at (D) {$D$};
            \node[above] at (2, {0.4*sqrt(3)/2}) {$Y$};
            \node[above] at (2.5, {0.4*sqrt(3)/2}) {$X$};
            \draw[very thin, dashed] (B) -- (3, {0.4*sqrt(3)/2});
        \end{tikzpicture}
        \caption{The competitor used to prove \cref{prop:effective-short-edge-2}. The dashed segment is removed from the net and replaced with the thick segment. The gray area is the region spanned by the homotopy from the original net to the competitor. The points $X$ and $Y$ are intersections of the equator $\lambda$ with the equators containing the segments $\alpha$ and $\delta$. }\label{f:competitor-2-bis}
    \end{figure}

To complete the proof of the proposition we give here the arguments for the claims (a), (b), and (c). It is convenient to mark the points $X$ and $Y$ closest to $B$ where the equator $\lambda$ intersects the two equators $\bar\alpha$ and $\bar\delta$ containing, respectively, the edges $\alpha$ and $\delta$, cf. \cref{f:competitor-2-bis}. Moreover, before proving (a), we can simply mark the points $D$ and $E$ as the closest ones to $M$ where $\bar\lambda$ intersects $\bar{\alpha}$ and $\bar{\delta}$. The claim (b) is then obvious, since by construction all internal angles of the pentagon $ABCDE$ are smaller than $\pi$.

In order to prove (a) we just need to check that the segments $|AE|$ and $|CD|$ are shorter, respectively, than $\text{len} (\alpha)$ and $\text{len} (\delta)$. Because the latter are strictly larger than $10 (|AB|+ |BC|) = 2 |BM|$, it suffices to show that $|AE| \leq 2 |BM|$ and $|CD|\leq 2 |BM|$. The argument is the same and we focus therefore on the first inequality.
Observe that, since sum of the angles $X\hat{A}B$ and $A \hat{B}X$ measure $\frac{2\pi}{3}$ and $\frac{\pi}{3}$ (and in particular they sum to $\pi$), elementary geometry implies that
\begin{align}
|AX|+|BX| & = \pi \label{e:furbo1}\\
|AB| & > B\hat{X} A \label{e:furbo2}\\
|BX| & > |AX|\, .\label{e:furbo3}
\end{align}
We moreover observe that
\[
|EM|= \dist\, (E, XB) \leq \dist\, (A, XB) = B\hat{X} A \leq |AB|
\]
(the second inequality follows because $|EX| < |AX| < \pi$ by \eqref{e:furbo1} and \eqref{e:furbo3}). In particular
\[
|AE|\leq |BM|+|BA|+|ME|\leq |BM|+2|AB|\leq 2 |BM|\, .
\]
We are thus left with proving (c). Observe that the argument just given above also implies
\begin{align}
|CD| - |DE| &\geq |BM|-|CB|-|DM| - |DE| \geq |BM| - 2 |CB| - |DE|\nonumber\\
&\geq |BM|- 3 |BC| - |AB| \geq \tfrac{2}{5} |BM|\, .\label{e:stima-lineare}
\end{align}
We will prove below that
\begin{equation}\label{e:area-ineq}
\mathcal{H}^2 (BMEA) \leq |BM| \cdot |BA|\, .
\end{equation}
Assuming the latter fact for a moment, the same inequality applied the quadrilateral $BCDM$ then gives
\begin{equation}\label{e:stima-quadratica-2}
\mathcal{H}^2 (ABCDE) \leq |BM| \cdot (|AB|+|BC|) \leq \tfrac{1}{5} |BM|^2\, .
\end{equation}
In particular, taking into account \eqref{e:stima-lineare} and \eqref{e:stima-quadratica-2}, (c) holds if $|BM|<1$, which is certainly the case because $|BM| = 5 (|AB|+|BC|) < \frac{1}{10}$.

\begin{figure}[htbp]
        \begin{tikzpicture}[scale=2]
          \draw (-3.5,0) -- (3.5,0);
          \node[above] at (3.3, 0) {$\lambda$};
          \coordinate (X) at (3,0);
          \coordinate (X') at (-3,0);
          \coordinate (C) at (0,4);
          \tkzDrawArc(C,X')(X);
          \node[above] at (X) {$X$};
          \node[above] at (X') {$X'$};
          \draw[dashed] (-2,-1) -- (2,-1);
          \node[right] at (2,-1) {$\partial D_r (N)$};
          \coordinate (M) at (1,0);
          \coordinate (B) at (-0.5,0);
          \coordinate (N) at (0,-2.5);
          \node[below] at (N) {$N$};
          \node[above] at (M) {$M$};
          \node[above] at (B) {$B$};
          \draw (M) to[out=270, in=50] (N);
          \node[above left] at (0.95,-0.95) {$E$};
          \draw (B) to [out=270, in=110] (N);
          \draw (B) -- (0.4,{4-sqrt(25-0.16)});
          \node[above] at  (0.5,{4-sqrt(25-0.16)}) {$A$};
        \end{tikzpicture}
        \caption{The argument for the inequality \eqref{e:area-ineq}.}\label{f:stima-area}
    \end{figure}
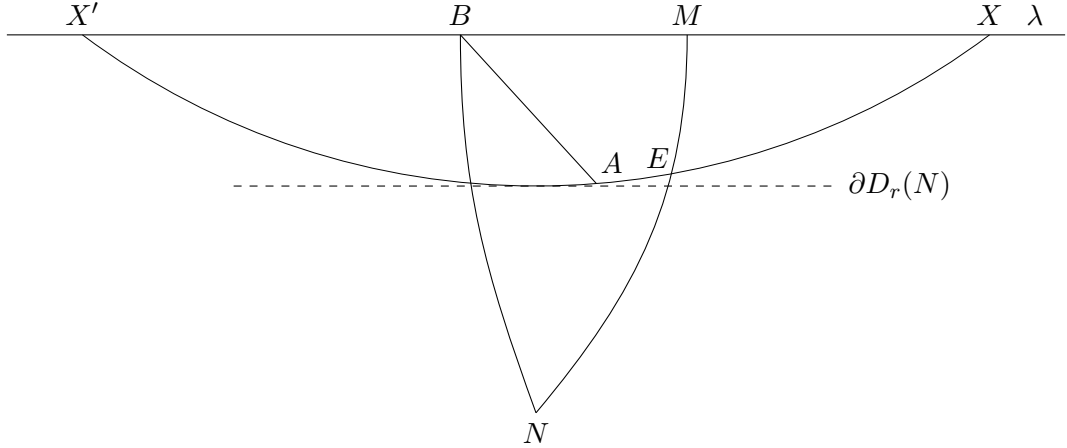

In order to show \eqref{e:area-ineq}, consider the equator $\lambda$ and the equator $\bar\alpha$ containing the edge $\alpha$ and let $X$ and $X'$ be their two intersections. The quadrilateral $BMEA$ is contained in the convex region delimited by two halves of these equators connecting $X$ and $X'$, cf. \cref{f:stima-area}, and we denote by $\eta$ the half-equator (of the two) which is contained in $\bar\alpha$. Denote further by $H$ the hemisphere bounded by $\lambda$ which contains $\eta$, by $N$ the center of $H$, and by $r< \frac{\pi}{2}$ the radius such that the disk $D_r (N)$ centered at $N$ of radius $r$ is tangent to $\eta$. By convexity $\eta$ touches the disk from the outside. In particular we conclude that
\[
BMEA \subset \{r < \dist (\cdot , N) < \tfrac{\pi}{2}\}\, .
\]
The point of tangency between $\gamma$ and $\partial D_r (N)$ is given by the midpoint of $\eta$ and it follows therefore that
\begin{equation}\label{e:compute-r}
r = \tfrac{\pi}{2} - A\hat{X}B\, .
\end{equation}
We now let $NB$ and $NM$ be the shortest geodesic segments connecting $N$ with, respectively, $B$ and $M$: these segments meet the equator $\lambda$ perpendicularly at $B$ and $M$. Because the internal angle of $ABME$ at $B$ is $\frac{\pi}{3} < \frac{\pi}{2}$ and the one at $M$ is $\frac{\pi}{2}$, if we denote by $T$ the convex triangle with vertices $N$, $B$, and $M$ we have
\[
BMEA \subset T \cap \{r < \dist (\cdot , N) < \tfrac{\pi}{2}\} =: \Lambda\, .
\]
Using that the function $p\mapsto \dist\, (p, N)$ is smooth in $\Lambda$ and its gradient has modulus $1$, we can use the coarea formula to compute
\[
\mathcal{H}^2 (BMEA) \leq \mathcal{H}^2 (\Lambda) = \int_r^{\frac{\pi}{2}} \mathcal{H}^1 (T \cap \{\dist (\cdot, N) =\rho\})\, d\rho \leq |BM| (\tfrac{\pi}{2} - r)\, ,
\]
where we have used that the function $\rho \mapsto  \mathcal{H}^1 (T \cap \{\dist (\cdot, N) = \rho\})$ reaches its maximum at $\rho=\frac{\pi}{2}$. Using \eqref{e:compute-r} we then get $\mathcal{H}^2 (BMEA)\leq |BM| \cdot A\hat{X}B$, which then implies \eqref{e:area-ineq} by \eqref{e:furbo2}.

\section{Proof of \texorpdfstring{\cref{l:quadrilateral-diagonal}}{l:quadri}}\label{a:quadrilateral-diagonal}

\begin{figure}
        \begin{tikzpicture}[scale=1]
        \draw[very thick] (2,1) -- (2,-1) -- (-2,-1) -- (-2,1) -- (2,1);
        \node[below right] at (2,1) {$A$};
        \node[above left] at (-2,1) {$D$};
        \node[below right] at (2,-1) {$B$};
        \node[below left] at (-2,-1) {$C$};
        \draw[very thin] (-2,-1) -- (2,1);
        \node[below right] at (0,0) {$\lambda$};
        \node[above right] at (-1.5,-1.08) {$\theta$};
        \draw[very thick, dashed] (2,1) -- (4,3);
        \node[above right] at (4,3) {$Q$};
        \draw[very thick, dashed] (-2,-1) to [out=45, in=210] (4,3);
        \node[below right] at (3,2) {$\gamma_A$};
        \node[above left] at (1.5,1.5) {$\gamma_C$};
        \end{tikzpicture}
        \caption{The quadrilateral $ABCD$ and the triangle $ACQ$.}\label{f:quadrilateri-e-diagonali}
    \end{figure}
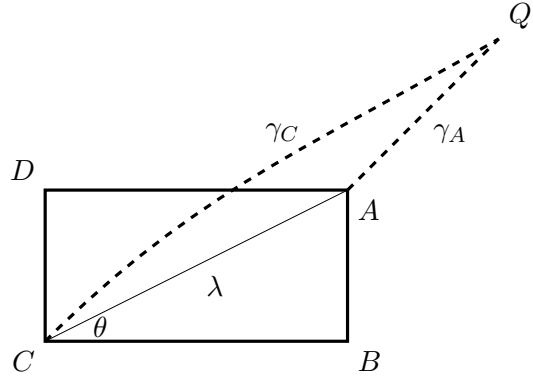

Consider a quadrilateral $ABCD$ as in the lemma, with vertices labeled clockwise. Assume that the equators $\gamma_A$ and $\gamma_C$ bisect the angles at $A$ and $C$. These equators coincide if and only if the diagonal $AC$ bisects both angles, which happens if and only if the quadrilateral is also equilateral. We can therefore draw the diagonal connecting $A$ and $C$ and assume, without loss of generality, that $\theta := A\hat{C}B$ is strictly smaller than $\frac{\pi}{3}$. Of course the angle $B\hat{A}C$ measures $\frac{2\pi}{3}-\theta$. We can explicitly compute the cosine of the length $\lambda$ of $AC$ as
\[
\cos \lambda = \frac{-\tfrac{1}{2} + \cos(\theta) \cos\left(\frac{2\pi}{3} - \theta\right)}{\sin(\theta) \sin\left(\frac{2\pi}{3} - \theta\right)}\,
\]
using the three angles of the convex triangle $ABC$ via the spherical law of cosines, cf. \cref{f:quadrilateri-e-diagonali}
We apply the product-to-sum identities to the numerator and denominator. For the cosines:
\begin{align*}
    \cos(\theta) \cos\left(\tfrac{2\pi}{3} - \theta\right) &= \tfrac{1}{2} \left[ \cos\left(2\theta - \tfrac{2\pi}{3}\right) + \cos\left(\tfrac{2\pi}{3}\right) \right]
    = \tfrac{1}{2} \left[ \cos\left(2\theta - \tfrac{2\pi}{3}\right) - \tfrac{1}{2} \right]
\end{align*}
For the sines:
\begin{align*}
    \sin(\theta) \sin\left(\tfrac{2\pi}{3} - \theta\right) &= \tfrac{1}{2} \left[ \cos\left(2\theta - \tfrac{2\pi}{3}\right) - \cos\left(\tfrac{2\pi}{3}\right) \right]
    = \tfrac{1}{2} \left[ \cos\left(2\theta - \tfrac{2\pi}{3}\right) + \tfrac{1}{2} \right]
\end{align*}
Substituting these expansions back into the equation for $\cos \lambda$ and
multiplying the numerator and denominator by 4 yields the simplified closed-form expression:
\begin{equation*}
    \cos \lambda = \frac{2\cos\left(2\theta - \frac{2\pi}{3}\right) - 3}{2\cos\left(2\theta - \frac{2\pi}{3}\right) + 1}
\end{equation*}

Consider now one of the two points $Q$ of intersections of $\gamma_A$ and $\gamma_C$, more precisely the first which we encounter on the equator $\gamma_A$ if we start at $A$ and exit the quadrilateral $ABCD$. $Q$, $A$, and $C$ form a convex triangle in which the angles adjacent to $AC$ measure $\frac{2\pi}{3}$ and $\frac{\pi}{3}-\theta$, cf. \cref{f:quadrilateri-e-diagonali}. The cosine of the angle $\alpha= C\hat{Q}A$ can thus be computed
applying again the spherical law of cosines to the triangle $CQA$:
\begin{equation*}
    \cos \alpha = -\cos\left(\tfrac{2\pi}{3} + \theta\right) \cos\left(\tfrac{\pi}{3} - \theta\right) + \sin\left(\tfrac{2\pi}{3} + \theta\right) \sin\left(\tfrac{\pi}{3} - \theta\right) \cos \lambda\, .
\end{equation*}
As above we evaluate the products using the sum and difference of the angles. The sum is $\pi$ and the difference is $2\theta + \frac{\pi}{3}$:
\begin{align*}
    \cos\left(\tfrac{2\pi}{3} + \theta\right) \cos\left(\tfrac{\pi}{3} - \theta\right) &= \tfrac{1}{2}\left[\cos\left(2\theta + \tfrac{\pi}{3}\right) + \cos(\pi)\right]
    = \tfrac{1}{2}\left[\cos\left(2\theta + \tfrac{\pi}{3}\right) - 1\right]
\end{align*}
\begin{align*}
    \sin\left(\tfrac{2\pi}{3} + \theta\right) \sin\left(\tfrac{\pi}{3} - \theta\right) &= \tfrac{1}{2}\left[\cos\left(2\theta + \tfrac{\pi}{3}\right) - \cos(\pi)\right]
    = \tfrac{1}{2}\left[\cos\left(2\theta + \tfrac{\pi}{3}\right) + 1\right]
\end{align*}
Substituting these into the equation for $\cos \alpha$ and using $\cos\left(2\theta + \frac{\pi}{3}\right) = - \cos (\frac{2\pi}{3}-2\theta)$, we get
\begin{equation*}
    \cos \alpha = \frac{1 + \cos\left(\frac{2\pi}{3}-2\theta\right)}{2} + \frac{1 - \cos\left( \frac{2\pi}{3}-2\theta\right)}{2} \left( \frac{2\cos\left(\frac{2\pi}{3}-2\theta\right) - 3}{2\cos\left(\frac{2\pi}{3}-2\theta\right) + 1} \right)\, .
\end{equation*}
After a few simple algebraic manipulations we find
\begin{equation*}
    \cos \alpha = \frac{8\cos\left(\frac{2\pi}{3}-2\theta\right) - 2}{2\left(2\cos\left( \frac{2\pi}{3}-2\theta\right) + 1\right)} = \frac{4\cos\left(\frac{2\pi}{3}-2\theta\right) - 1}{2\cos\left(\frac{2\pi}{3}-2\theta\right) + 1}\, .
\end{equation*}
The claim of the lemma is then that $\alpha$ uniquely determines $\theta$. Set then $x = \cos (\frac{2\pi}{3}-2\theta) \in (-\frac{1}{2},1)$. Given that $\cos$ on $[0, \pi]$ is an homeomorphism onto $[-1,1]$ and the fact that $0<\alpha < \pi$, our claim is equivalent to the injectivity of the function
\[
x\mapsto \tfrac{4x-1}{2x+1} = 2 - \tfrac{3}{2x+1}
\]
on the interval $(-\frac{1}{2}, 1)$. This is however an obvious fact.

\bibliographystyle{plainurl}
\bibliography{references.bib}

\end{document}